\documentclass[hidelinks,pdftex,phd,final]{pittetd} 
\usepackage{graphicx}
\usepackage[utf8]{inputenc}
\usepackage{indentfirst}
\usepackage[mathlines]{lineno}
\usepackage{xcolor}
\usepackage{listings, lstcoq}
\usewithpatch{amsmath, amssymb, amsfonts, amscd, amsthm, mathrsfs}
\usepackage{pdflscape}
\usepackage{cleveref}
\usepackage{subfig}
\usepackage{caption} 
\usepackage{pifont}
\usepackage{mathpartir}
\usepackage{eso-pic}
\usepackage{microtype}
\usepackage{changepage}
\usepackage{etoolbox} 
\usepackage{enumitem}
\hypersetup{hidelinks}
\InputIfFileExists{amsmath.pit}{%
            \ClassInfo{pittetd}{Patch for amsmath loaded}}{%
            \ClassInfo{pittetd}{No patch found for amsmath}}
\InputIfFileExists{amsthm.pit}{%
            \ClassInfo{pittetd}{Patch for amsthm loaded}}{%
            \ClassInfo{pittetd}{No patch found for amsthm}}
\makeatletter
\renewcommand{\@biblabel}[1]{[#1]\hfill}
\makeatother
\makeatletter
\newcommand\fake@math{}
\def\fake@math#1\){[math]}
\makeatother

\makeatletter
\def\@hangfrom#1{\setbox\@tempboxa\hbox{{#1}}%
      \hangindent 0pt
      \noindent\box\@tempboxa}
\makeatother
\newtheorem{theorem}{Theorem}

\newtheorem{lemma}[theorem]{Lemma}
\newtheorem{proposition}[theorem]{Proposition}

\newtheorem{corollary}[theorem]{Corollary}
\newtheorem{definition}[theorem]{Definition}
\newtheorem*{example}{Example}
\newtheorem*{remark}{Remark}
\newtheorem{problem}{Problem}
\newtheorem*{conjecture}{Conjecture}
\numberwithin{theorem}{section}

\def\R{\mathbb{R}}
\def\Z{\mathbb{Z}}
\def \h{\mathfrak{h}}
\def\H{\mathcal{H}}

\def\star{\mathfrak{h}^{\ast}}
\def\truncstar{\mathfrak{h}^{\ast \ast}} 

\def\SL{\mathrm{SL}_2}
\def\ac{\mathrm{ac}}
\def\val{\mathrm{val}}
\def\sl{\mathfrak{sl}_2}
\def \SOtwo{\mathrm{SO}(2,1)}
\def \sotwo{\mathfrak{so}(2,1)}

\def \SU{\mathrm{SU}(1,1)}
\def\su{\mathfrak{su}(1,1)}
\def\SO{\mathrm{SO}_2(\R)}
\def\Ad{\mathrm{Ad}}
\def\Jsu{J_{\mathfrak{su}}}
\def\ad{\mathrm{ad}}
\def\packing{\mathcal{P}}
\def\LL{\mathbf{L}}
\def\Kbal{\mathfrak{K}_{bal}}
\def\tr{\mathrm{trace}}
\def\O{\mathcal{O}}
\def\OX{\mathcal{O}_{X}}
\def \delx{\frac{\partial}{\partial x}}
\def \dely{\frac{\partial}{\partial y}}
\newcommand{\mattwo}[4]{\left(
\begin{array}{cc}
#1 & #2 \\
#3 & #4 \\  
\end{array}
\right)
}
\newcommand{\dettwo}[4]{\left|
\begin{array}{cc}
#1 & #2 \\
#3 & #4 \\  
\end{array}
\right|
}
\newcommand{\matthree}[9]{
\begin{array}{ccc}
#1 & #2 & #3 \\
#4 & #5 & #6 \\
#7 & #8 & #9 \\
\end{array}
}
\newcommand{\bracks}[2]{
\left\langle #1 , #2 \right\rangle
}

\newcommand{\Kccs}{\mathfrak{K}_{ccs}}
\newcommand{\RR}{\mathfrak{R}}

\newcommand{\states}{S}
\newcommand{\E}{\mathbb{E}}
\newcommand{\actions}{A}
\newcommand{\bind}{\mathsf{bind}}
\newcommand{\ret}{\mathsf{ret}}

\newcommand{\giry}{P}
\newcommand{\policy}{\pi}
\newcommand{\transition}{T}
\newcommand{\reward}{r}
\newcommand{\reals}{\mathbb{R}}
\newcommand{\integers}{\mathbb{Z}}

\newcommand{\bellman}{\mathbf{B}}

\newcommand{\mdplen}{n}
\newcommand{\step}{k}
\newcommand{\libname}{CertRL }

\makeatletter
\patchcmd{\BR@backref}{\newblock}{\newblock(on page~}{}{}
\patchcmd{\BR@backref}{\par}{)\par}{}{}
\makeatother

\newcommand*\linenomathpatch[1]{%
  \cspreto{#1}{\linenomath}%
  \cspreto{#1*}{\linenomath}%
  \csappto{end#1}{\endlinenomath}%
  \csappto{end#1*}{\endlinenomath}%
}

\linenomathpatch{equation}
\linenomathpatch{gather}
\linenomathpatch{multline}
\linenomathpatch{align}
\linenomathpatch{align*}
\linenomathpatch{alignat}
\linenomathpatch{flalign}

\chapterfloats 

\title[On the Reinhardt Conjecture and Formal Foundations of Optimal Control]{On the Reinhardt Conjecture and Formal Foundations of Optimal Control}

\author{Koundinya Vajjha}
\degree{M.Sc.,~University of Western Ontario, 2018 \\ M.Math.,~Indian Statistical Institute, 2016 \\
B.Math.,~Indian Statistical Institute, 2014}

\school{Department of Mathematics}

\date{July 25th 2022}


\subject{Dissertation}

\begin{document}
\maketitle



\committeemember{Professor Thomas C. Hales, Department of Mathematics, University of Pittsburgh}

\committeemember{Professor Kiumars Kaveh, Department of Mathematics, University of Pittsburgh}
\committeemember{Professor Anna Vainchtein, Department of Mathematics, University of Pittsburgh}
\committeemember{Professor Greg Kuperberg, Department of Mathematics, University of California, Davis}

\school{Department of Mathematics}
\makecommittee
\copyrightpage                     

\begin{abstract}
We describe a reformulation (following Hales (2017)) of a 1934 conjecture of Reinhardt on pessimal packings of convex domains in the plane as a problem in optimal control theory. Several structural results of this problem including its Hamiltonian structure and Lax pair formalism are presented. 

General solutions of this problem for constant control are presented and are used to prove that the Pontryagin extremals of the control problem are constrained to lie in a compact domain of the state space. 

We further describe the structure of the control problem near its singular locus, and prove that we recover the Pontryagin system of the multi-dimensional Fuller optimal control problem (with two dimensional control) in this case. We show how this system admits logarithmic spiral trajectories when the control set is the circumscribing disk of the 2-simplex with the associated control performing an infinite number of rotations on the boundary of the disk in finite time. 

We also describe formalization projects in foundational optimal control \textit{viz.,} model-based and model-free Reinforcement Learning theory. Key ingredients which make these formalization novel \textit{viz.,} the {Giry monad} and {contraction coinduction} are considered and some applications are discussed. 
\end{abstract}




\tableofcontents

\listoftables                      

\listoffigures                

\section{List of Notation}

\begin{description}[align=left,labelwidth=3cm]
\item[$\packing$] Packing of bodies in $\R^n$.
\item[$\delta(K,\packing)$] Packing density of a packing consisting of congruent copies of the body $K$. 
\item[$C$] The Cayley transform.
\item[$\delta(K)$] Best packing density of $K$.
\item[$\mu$] Lebesgue measure in $\R^2$.
\item[$\mathbf{L}$] Lattice in $\R^2$.
\item[$\delta(K,\mathbf{L})$] Density of the lattice packing of $K$.
\item[$\delta_L(K)$] Best lattice packing density of $K$.
\item[$\Kccs$] Convex, centrally symmetric discs in $\R^2$.
\item[$\Kbal$] A subset of discs in $\Kccs$. See Definition \ref{def:Kbal}.
\item[$\det(\mathbf{L})$] Determinant of the lattice $\mathbf{L}$. 
\item[$\delta_{\min}$] Worst best-packing density among all domains in $\Kccs$.
\item[$K_{\min}$] Global minimizer of the Reinhardt problem.
\item[$E$] Multi-points.
\item[$\sigma$] Multi-curves.
\item[$e_j^*$] The sixth roots of unity. 
\item[$\Ad$] The adjoint representation. 
\item[$\ad$] Infinitesimal generator of the adjoint action. 
\item[$\Ad^*$] The coadjoint representation.
\item[$\ad^*$] Infinitesimal generator of the coadjoint action. 
\item[$(g,X)$] State variables of the Reinhardt control problem. 
\item[$(\Lambda_1,\Lambda_2)$] Costate variables of the Reinhardt control problem. 
\item[$\Lambda_R$] Reduced costate variable of $\Lambda_2$.
\item[$J$] Infinitesimal generator of rotations in $\sl(\R)$. 
\item[$\Jsu$] Infinitesimal generator of rotations in $\su$. 
\item[$H_K$] The critical hexagon of the domain $K$.
\item[$h_K$] Inscribed centrally symmetric hexagon of the domain $K$.
\item[$\Delta(K)$] Minimal determinant of the domain $K$. 
\item[$\h$] The Poincar\'{e} upper half-plane. 
\item[$\star$] Star domain. 
\item[$\h^{**}$] Compactification of the star domain.
\item[$v_j$] State-dependent curvatures.
\item[$u_j$] Control variables.
\item[$U_T$] Simplex control set. 
\item[$U_C, U_I, U_r$] Various circular control sets (see Section~\ref{sec:control-sets}).
\item[$\alpha,\beta$] Control set parameters (see Table~\ref{tab:control-sets}).
\item[$Z_u$] Affine transformation mapping control set into $\sl(\R)$. 
\item[$R$] Rotation matrix $\exp(J\pi/3)$. 
\item[$\O_J,~\OX$] Adjoint orbit through $J$ and $X$.
\item[$T_X\OX,~T_X^*\OX$] Tangent and cotangent space to $\OX$ at $X$. 
\item[$TG,~T^*G$] Tangent and cotangent bundles of a Lie group $G$.
\item[$\lambda_{cost}$] Multiplier of the Pontryagin system. 
\item[$u^*(t)$] Optimal control. 
\item[$\H$] Full Hamiltonian of the Reinhardt control problem.
\item[$\H^*$] Maximized Hamiltonian.
\item[$\overrightarrow{\mathcal{G}}$] Hamiltonian vector field corresponding to $\mathcal{G}$. 
\item[$\Lambda_{sing}$] Singular locus.
\item[$\mathbf{J}^{\tau}$] Momentum map corresponding to an infinitesimal group of symmetries $\tau$.
\item[$\frac{\delta}{\delta X}$] Functional derivative with respect to $X$.
\item[$\RR$] Sesquilinear extension of the real part to $\mathbb{C}\times\mathbb{C}$: $\RR(u,v) := \mathrm{Re}(u\bar{v})$. 
\item[$\mu(w,z)$] Stands for $[w] - (\beta/\alpha) \RR(w,z)$ where $[w] = \sqrt{1 + |z|^2}$.
\item[$w,b,c$] Hyperboloid coordinates corresponding to $X,\Lambda_1,\Lambda_R$ respectively. 
\item[$\{\cdot,\cdot\}_{c}$] Fuller system Poisson bracket (See Section \ref{sec:fuller-system}).
\item[$\{\cdot,\cdot\}_{ex}$] Extended space Poisson bracket (See Appendix \ref{sec:poisson-bracket}).

\end{description}



\newpage
\begin{center}
    \tiny{R.S.}
\end{center}
\topskip0pt
\vspace*{\fill}
\begin{figure}[h]
    \center
    \includegraphics[scale=1.3]{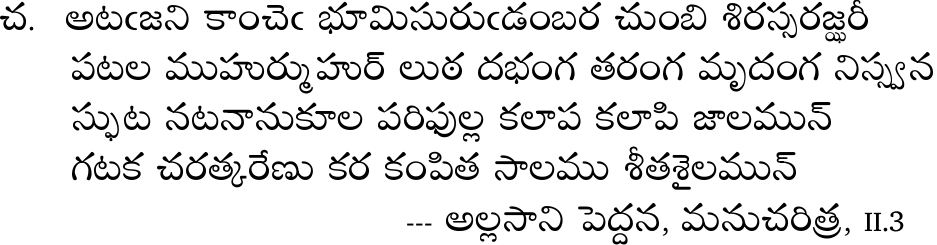}
\end{figure}
\smallskip
\medskip
\begin{quote}
    \centering 
    \footnotesize
        \emph{He witnessed the snow-capped mountain, its tall peak kissing the expanse of sky,  the gorge of which  gave birth to a waterfall, turning on itself, rumbling like the beating of a drum as it crashed onto the rocks below; 
        to which peacocks danced in time, feathers unfurled;
        while elephants roamed the mountain slopes shaking the sal trees with their trunks...\\
        \flushright
        --- Allasani Peddana, The Story of Manu, II.3}
\end{quote}
\vspace*{\fill}
\newpage
\phantomsection
\preface
This thesis would not have been possible without the vision, guidance and extraordinary mentorship of Thomas Hales. It has been the privilege of my life to have personally witnessed his approach and philosophy to mathematics. Not only did he teach me an huge amount of mathematics, but also how to develop an idea from a vague hunch into a crystallized, polished proof; when to bury an idea and move on; to accept mathematics for what it is and, above all, to appreciate its beauty. 

On June 17th, 2019 in Hanoi, he asked me if I ``want[ed] to prove the Reinhardt conjecture'', a question by which I was mildly taken aback. I replied that I would definitely want to try. In the months that followed, the problem, initially daunting, would always display remarkable structure\footnote{A structure which was not unlike that of the mythical mountain in the \textit{Story of Manu}...} and would always find a way to reel us back in. Some of my fondest memories over the last four years included time spent with this problem and for that, I am grateful. 

Thanks are to Avi Shinnar, Barry Trager, Vasily Pestun and Nathan Fulton at IBM Research for the internship and collaboration. With their help, not only did I formalize a lot of mathematics in the Coq proof assistant, but also found in them extremely knowledgeable mentors who I could rely on for support and discussion.

I also thank Greg Kuperberg, Anna Vainchtein and Kiumars Kaveh for agreeing to be on my thesis committee. I thank Velemir Jurdjevic for many helpful discussions and Larry Bates for patient answers to my many questions over email.

Thanks are also to friends at the Pitt math department: Luis, Sushmita, Arshia, Farjana, Grishma, Yujie, Mark and Anthony for having made my time spent in Thackeray Hall all the more enjoyable.

Lastly, I would not have imagined any of this without the love, support and patience of Amma, Dad, Monanna, Soni and Vihasi. Special thanks are to Amma: I wrote you this thesis in exchange for all the free food over the years. Finally, I dedicate this thesis to my grandmother Saraswati Mokkapati who --- like her namesake --- is refulgent, ineffable, and an Angel.


\phantomsection


\chapter{Introduction}%
Optimal Control problems solve for a policy driving an agent in an environment over a period of time such that a cost function is optimized. While this definition is broad, it encompasses several problems in both pure and applied mathematics which have been studied intensely for centuries. 

Optimal control theory takes its roots from the Calculus of Variations which was developed in the 18th century by Euler and Lagrange for solving optimization problems arising in mechanics. A very fruitful next century followed, with important contributions by Legendre, Hamilton, Jacobi and Weierstrass. Control theory was then solidified as an independent subject following the works of Bolza, McShane, Pontryagin and Bellman in the 20th century. The work of the latter two authors, Pontryagin and Bellman, has had a remarkable influence on control theory: Pontryagin was responsible (along with his school) for codifying the important Pontryagin Maximum Principle (PMP), which is considered a major milestone in control theory. Bellman's work brought forth the paradigm of dynamic programming, and his principle of optimality heralded the advent of what we today call Reinforcement Learning. 
To elaborate: a typical solution of an optimal control problem proceeds, after the problem is set up, by appealing to certain \textit{necessary conditions} which give information about the structure of the minimizer. In most cases, these necessary conditions give enough information about the minimizer to determine it completely. The PMP and the Hamilton-Jacobi-Bellman equation are examples of such necessary conditions.

Decades of research has resulted in the study of several variations of the base optimal control problem we mentioned above:
\begin{itemize}
    \item The \textit{policy} or \textit{control law} which we solve for may be state-dependent (closed-loop control) or just time-dependent (open-loop control).
    \item The \textit{environment} or \textit{state-space} may just be a finite set or may be an infinite manifold.
    \item The \textit{agent} or \textit{dynamical system} may transition in the state space stochastically or deterministically. 
    \item The \textit{period of time} we consider may be finite or infinite, and the time instances in which transition happens may be at discrete intervals or may be continuous.
    \item The \textit{cost function} may be time itself (time optimal control problem) or may be some other quantity. 
\end{itemize}

We now consider two examples which we hope shall illustrate the rich variety of problems that are possible to consider in this framework.  

\begin{example}
\begin{figure}
    \centering
    \includegraphics[scale=0.4]{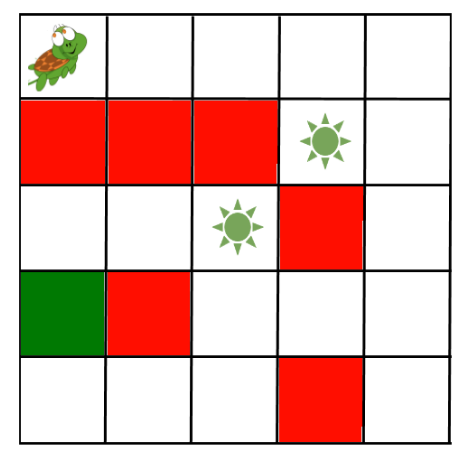}
    \caption{Turtle in a gridworld.}
    \label{fig:turtle}
\end{figure}

\normalfont The first example is the situation in Figure \ref{fig:turtle} showing a turtle in a $5 \times 5$ grid. Each square contains an associated reward for visiting it --- the white squares give zero rewards, red squares give negative rewards, the green stars have a positive reward, and time stops when the turtle reaches the green square. Our objective is to prescribe a list of actions (a \textit{policy} telling the turtle to move up, down, left or right at each time step) which the turtle should take so as to optimize its \textit{long-term reward} (the weighted average of all the rewards it accumulates following a policy). An important note is that the transitions of the turtle are \textit{stochastic}: at a particular state, given an action, the turtle's next state is only determined upto a probability. 
The study of problems of this nature forms the bulk of Reinforcement Learning (RL) theory. 
\end{example}

\begin{example}
\begin{figure}
    \centering
    \includegraphics[scale=1.5]{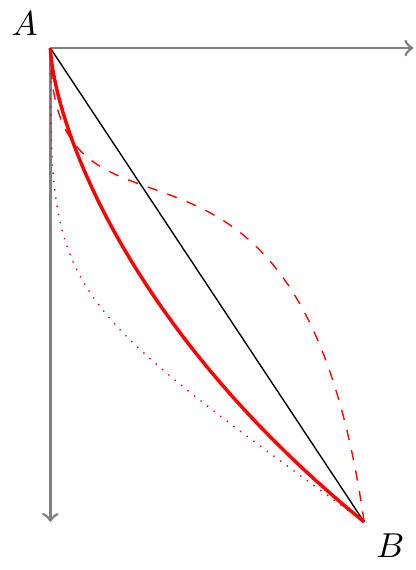}
    \caption{Brachistochrone problem.}
    \label{fig:brach}
\end{figure}
  
\normalfont The next example is the classical Brachistochrone problem shown in figure \ref{fig:brach}. Given two points on a wall, our objective is to join them by a curve so that a point mass, purely under the influence of gravity, starts at the first point and reaches the second point in the least amount of time by sliding on the curve without friction. We can see that in this problem is deterministic with the control law (policy) being a curve joining the two points in question, and the cost is the time taken by the mass to travel from the first point to the second point under gravity. Additionally, we see that the time is measured continuously.  
\end{example}

The Brachistochrone problem differs from the problems in Reinforcement Learning in that the state space of the problem is allowed to be a \textit{manifold} (which here is the plane $\mathbb{R}^2$). The field of Geometric Optimal Control theory --- which concerns itself with the study of control problems on manifolds --- provides powerful necessary conditions for optimality which have proven to be useful in solving different optimization problems, including the Brachistochrone. For more examples, see \cite{jurdjevic_1996,jurdjevic2016optimal,sussmann2002brachistochrone}.

These two examples are illustrative of the problems we consider in this thesis, which fall broadly into the areas of Reinforcement Learning and Geometric Optimal Control theory. 

\section{Formal Verification of Optimal Control}\label{sec:rein-learning}
Recent high-profile applications of Reinforcement Learning include beating the world's best players at Go \cite{alphagonature}, competing against top professionals in Dota \cite{openaidota}, improving protein structure prediction \cite{deepmindprotein}, and automatically controlling complex robots \cite{DBLP:journals/corr/GuHLL16}. These successes motivate the use of RL algorithms in safety-critical settings where it is important to ensure that the implementations of these algorithms are bug-free, to whatever extent possible.

Indeed, there is now a whole sub-field of Reinforcement Learning called ``safe RL'' which is concerned with designing RL algorithms in which the agents provably avoid certain regions of the state space~\cite{Junges16}.

In order to be certain that these algorithms do what they advertise, one should trust (among other things) that the proofs of correctness have no gaps and that the implementation of the algorithm in code conforms to the specification on paper. While actual numbers may vary, it is reported that commercial software has on average 1 to 25 bugs per 1000 lines of code~\cite{mcconnell2004code}. Even if not all of these bugs are serious, there is still a small chance of a serious one going unnoticed and potentially leading to damage or even loss of life. 

The most promising method of mitigating such disasters is via the technique of \textit{formal verification}, which refers to mathematically proving that a piece of code conforms to it's specification. At a high level, deductive formal verification proceeds by first embedding the piece of code as a mathematical model into (a computer implementation of) a logic and then proves that this model conforms to a specification (also expressed in the logic). These computer implementations of the logic (called \textit{proof assistants}) may include a reasoning layer to generate, manipulate and discharge \textit{proof obligations}, whose resolution imply the conformance of the code with its specification. Some logics (such as the Calculus of Inductive Constructions) have computational semantics, which means that their implementations also double up as programming languages. This makes it possible to program, specify and prove properties about computer code entirely within proof assistants. Notable successes within computer science which employ formal verification are the CompCert project which produced a formally verified optimizing C compiler and the seL4 project which produced a formally verified microkernel.

From a mathematician's point-of-view, a formal proof is a proof written down so that all logical inferences employed can be traced down all the way to the fundamental axioms of mathematics, without any appeals to intuition. The base logic which some proof assistants (such as Lean, Coq and Agda) implement are powerful and expressible enough to encode large swathes of mathematics. Despite this, formal proofs received little attention from professional mathematicians until the last decade. But today, encouragingly, we see mathematicians adopting and using proof assistants to formally verify their own work and also in building and maintaining libraries of formal mathematics \footnote{This is exemplified by the Lean Theorem Prover's math library $\mathsf{mathlib}$ which in the period between 2017-2022 has amassed $676,000$ lines of code.}.

Part of the appeal of formal verification, setting aside the correctness guarantees, is its ability to force the user to think through arguments down to first principles. This sometimes results in proofs which are simpler than their pen-and-paper equivalents. Frequently, implementation issues force maintainers of libraries of formal mathematics to state results which hold in maximal generality. Special cases of theorems then follow by applying the general theorem in their appropriate contexts. This bypasses the need to formalize the same argument twice. 

The reasoning within a proof assistant is ``close to the logic'' and so it behooves us to design proofs which conform to the constructs of the base logic. For example, certain primitives (such as inductive types) have good out-of-the-box support: every time an inductive type is declared in Lean/Coq, the associated recursors and properties they satisfy are automatically created and added to the local context. 
Consequently, a proof which is closer to the logic of the proof assistant is easier to formalize than one which is based on traditional set theory. 

Two examples of this (which we elaborate on further in Chapter III) are the Giry Monad and Metric Coinduction. These constructs expose the essential logical (coinductive) structure of long-term values associated to infinite horizon Markov Decision Processes. This is a great example of formal proofs providing a clarity to a well-established theory by viewing it from a different point-of-view.  

After describing the use of these techniques in formal verification projects, we shall also describe our work towards a formal proof of convergence of model-free Reinforcement Learning methods. All such convergence proofs utilize stochastic approximation in one form or another and so we formalize a very general stochastic approximation algorithm.
All results in Chapter III are joint work with Barry Trager, Avi Shinnar, Vasily Pestun and Nathan Fulton. 

\section{From the Reinhardt Conjecture to Optimal Control}\label{sec:opt-control}

A class of problems which discrete geometers are interested in concern minimizing or maximizing the \textit{packing density} $\delta(K,\packing)$ of a \emph{body} $K \subset \R^n$, that is a compact, connected set\footnote{Careful definitions of these terms are given at the beginning of Chapter 2.}. A body in $\R^2$ shall be called a \textit{disc}. The packing density roughly corresponds to the proportion of space taken up by congruent copies of a body $K$ when they are arranged according to the packing $\mathcal{P}$ in Euclidean space $\R^n$. Since $\delta(K,\packing)$ is a two variable function, different flavours of this question may be posed: we may restrict the classes $K \in \mathfrak{K}$ under consideration, or we may restrict the type of packings $\mathcal{P} \in \mathfrak{P}$. 

For example, the \textit{sphere packing problem} fixes the body $K$ to be $B^n$ (the unit ball in $\R^n$) and asks us to determine $\delta(K) := \sup_{\mathcal{P}} \delta(B^n,\packing)$, where the supremum ranges over all possible packings $\mathcal{P}$. Finding $\delta(K)$ for an arbitrary body $K$ is an extremely hard problem in general, even in low dimensions. In fact, it is no exaggeration to say that these problems are some of the most difficult in all of mathematics\footnote{For $n=3$, the sphere packing problem is the Kepler conjecture which was open for nearly 400 years.}.
\begin{figure}[ht]
    \centering
    \includegraphics[scale=1.1]{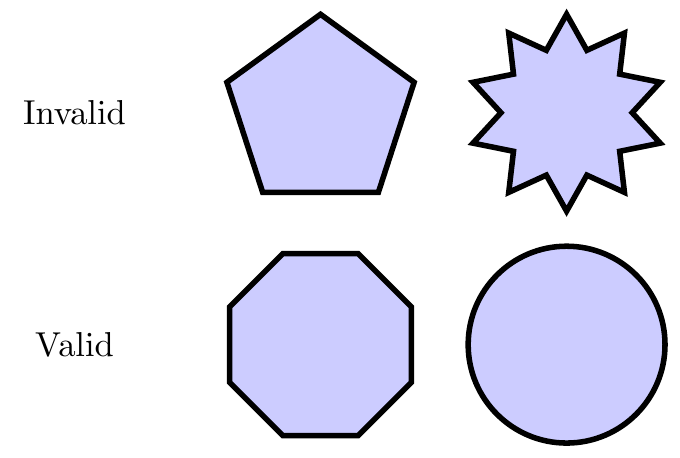}
    \caption{Valid and Invalid Convex Centrally Symmetric Discs.}
    \label{fig:class-ccs}
\end{figure}

\subsection{The Reinhardt problem}
In our situation, we fix the class of bodies $\mathfrak{K}$ to be the class $\mathfrak{K}_{ccs}$ of convex, centrally symmetric discs in the plane $\R^2$. Examples of discs belonging to $\mathfrak{K}_{ccs}$ are shown in \cref{fig:class-ccs}. We are now asked to determine the quantity \[\inf_{K \in \mathfrak{K}_{ccs}} \delta(K) = \inf_{K \in \mathfrak{K}_{ccs}} \sup_{\mathcal{P}} \delta(K,\packing)\]
and also that convex, centrally symmetric disc which whose best packing density achieves this minimum. Since affine transformations of a disc do not change its packing density, the answer to this problem is only unique upto an affine transformation. 
So, we need to find \textit{(upto an affine transformation) that unique convex, centrally symmetric disc in the plane whose best packing density is the worst among all such bodies}. 

\begin{figure}[ht]
    \centering
    \includegraphics[scale=0.43]{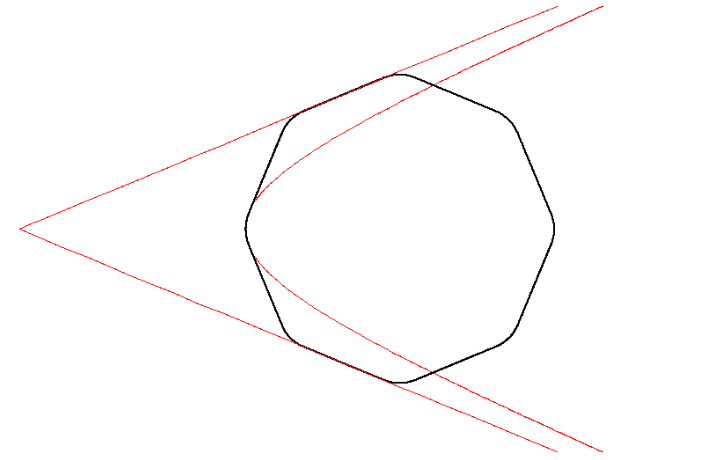}
    \caption{Construction of the Smoothed Octagon.}
    \label{fig:smoothed}
\end{figure}

While a plausible first guess for the minimizer is the circular disk in the plane, it turns out that there is a candidate which is slightly worse: In 1934, Karl Reinhardt \cite{reinhardt1934dichteste} conjectured that the minimum is achieved by the so-called \textit{smoothed octagon} pictured in \cref{fig:smoothed}. 

\begin{conjecture} [Reinhardt \cite{reinhardt1934dichteste}, Mahler \cite{mahler1947minimum}]
The smoothed octagon achieves the least best packing density among all other convex, centrally symmetric discs in the plane. It's density is given by 
\[ \inf_{K \in \mathfrak{K}_{ccs}} \delta(K) = \frac{8 - \sqrt{32} - \ln 2}{\sqrt{8} - 1} \approx 0.902414 \]
\end{conjecture}

As the picture shows, it is constructed by clipping the vertices of a regular octagon by hyperbolic arcs which are tangent to the sides making up the vertex, and which are asymptotic to the sides adjacent to those. 

\subsection{History of the Reinhardt Problem}

The earliest mention of the Reinhardt problem is as Problem 17 in \S27 of a 1923 book of Wilhelm Blaschke~\cite{blaschke1945differentialgeometrie}, where it is called \textit{Courant's conjecture}. This conjecture states that the best packing of all convex centrally symmetric discs has the packing density of the disk (whose best packing density in the plane is $\frac{\pi}{\sqrt{12}}$) as a greatest lower bound. 

Less than a decade later, Courant's conjecture was shown to be false by Reinhardt in his 1934 paper, by his construction of the smoothed octagon. After this paper, Kurt Mahler was also led to the smoothed octagon in a series of papers in 1946-47. The first paper \cite{mahler1947minimum} used methods of optimization theory (specifically, the calculus of variations) to refute Courant's conjecture by proving the existence of a disc whose packing density was worse than the disk. In the same paper, he states:

\begin{quote}
    It seems highly probable from the convexity condition, that the boundary of an extreme convex domain consists of line segments and arcs of hyperbolae. So far, however, I have not succeeded in proving this assertion.
\end{quote}

In a follow up paper, Mahler gives an explicit construction of the smoothed octagon~\cite{mahler1947area}. The term ``smoothed octagon'' appears explicitly in a follow up paper by Mahler and Ledermann in 1949~\cite{ledermann1949lattice}.

In our analysis of the Reinhardt conjecture, we largely follow Kurt Mahler's exposition in the articles cited above. Our reasons for doing so are because Reinhardt's original paper has no available English translation, and Mahler's presentation in terms of lattices is much easier to follow. 

Papers by Ennola~\cite{ennola1961lattice} in 1961 and Tammela~\cite{tammela1970estimate} in 1970 show that $\inf_{K \in \mathfrak{K}_{ccs}} \delta(K) > 0.8925$.  Nazarov \cite{nazarov1988reinhardt} proved that the smoothed octagon is a local minimum in the space of convex discs equipped with the Hausdorff metric. Hales \cite{hales2011} recast the Reinhardt problem as a problem in the calculus of variations. However, as of 2022, the full Reinhardt conjecture is still out of reach, having remained open since 1934. 

\subsection{Reinhardt Optimal Control Problem}
\begin{figure}[!htbp]
\centering

\caption{Multi-curves generating the smoothed octagon.}\label{fig:multi-curve-oct}
\subfloat[]{\includegraphics[width=0.33\textwidth]{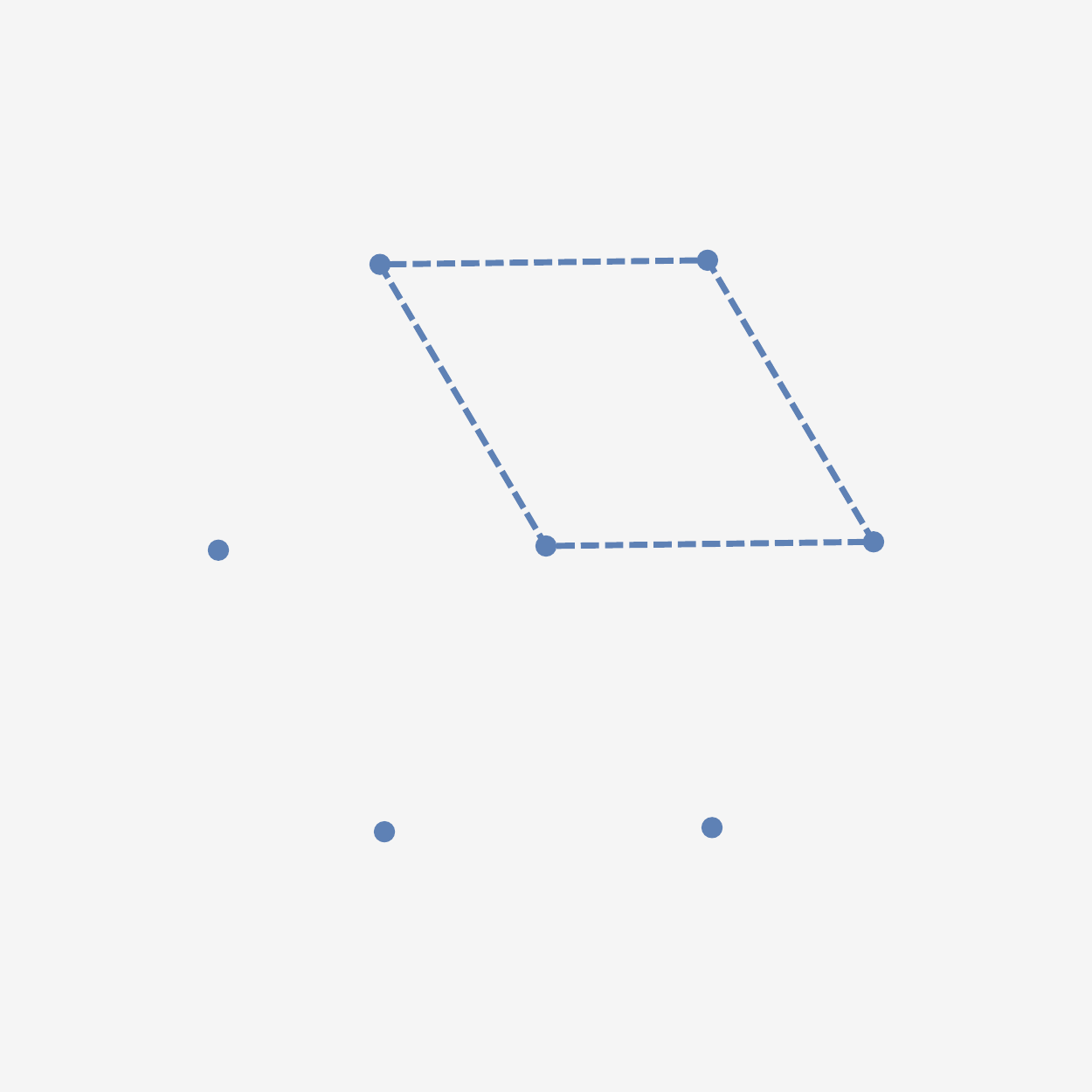}}\hfill
\subfloat[]{\includegraphics[width=0.33\textwidth]{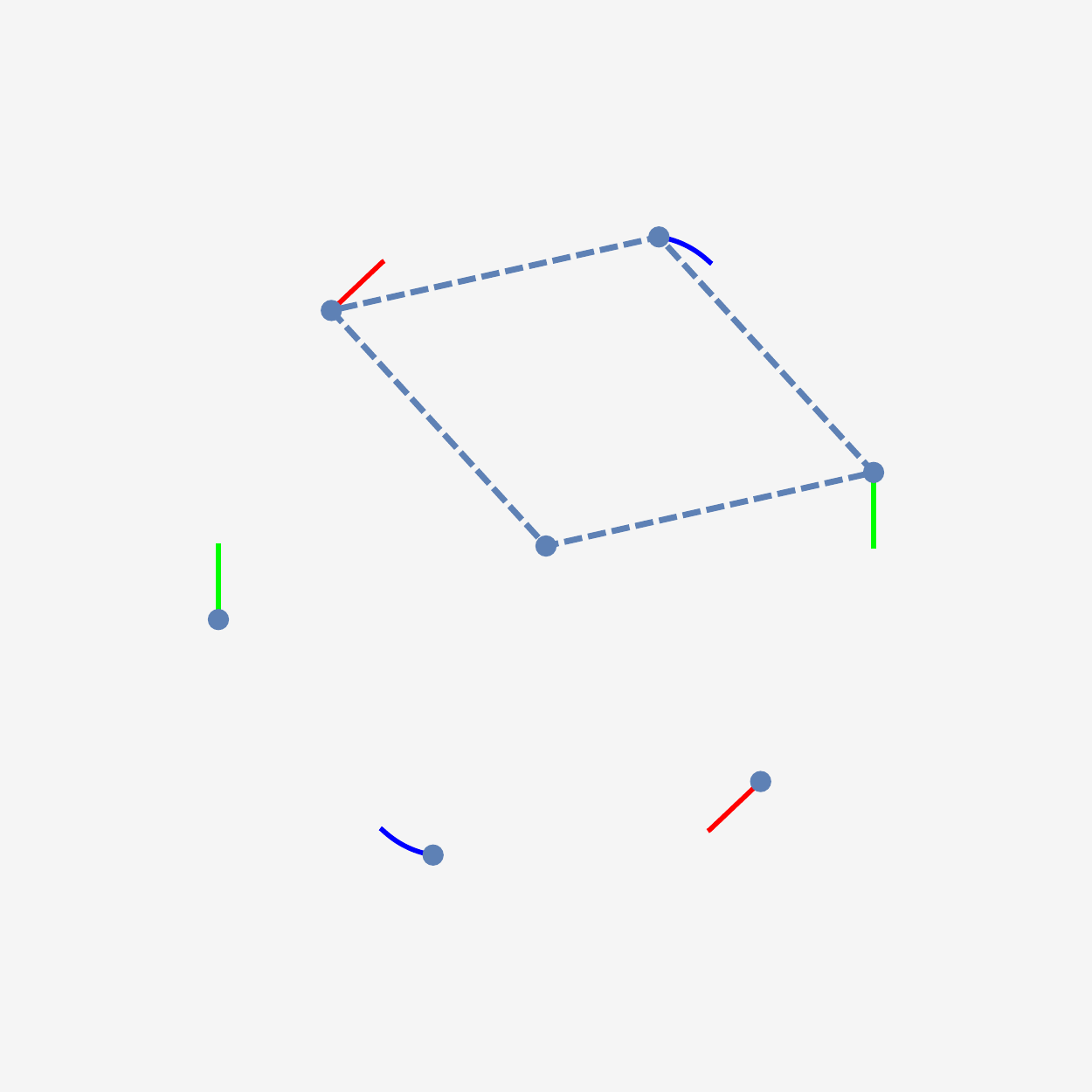}}\hfill
\subfloat[]{\includegraphics[width=0.33\textwidth]{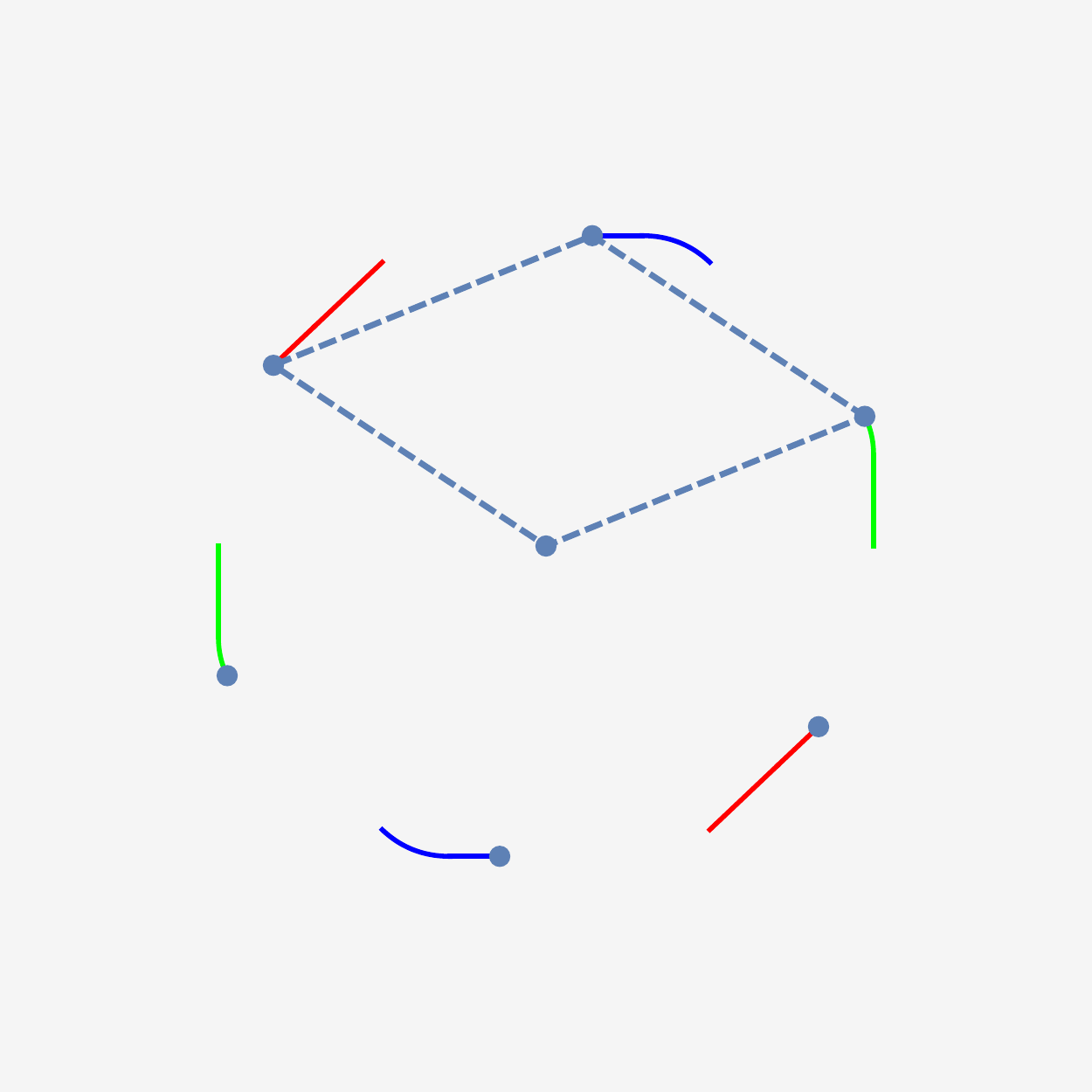}}

\subfloat[]{\includegraphics[width=0.33\textwidth]{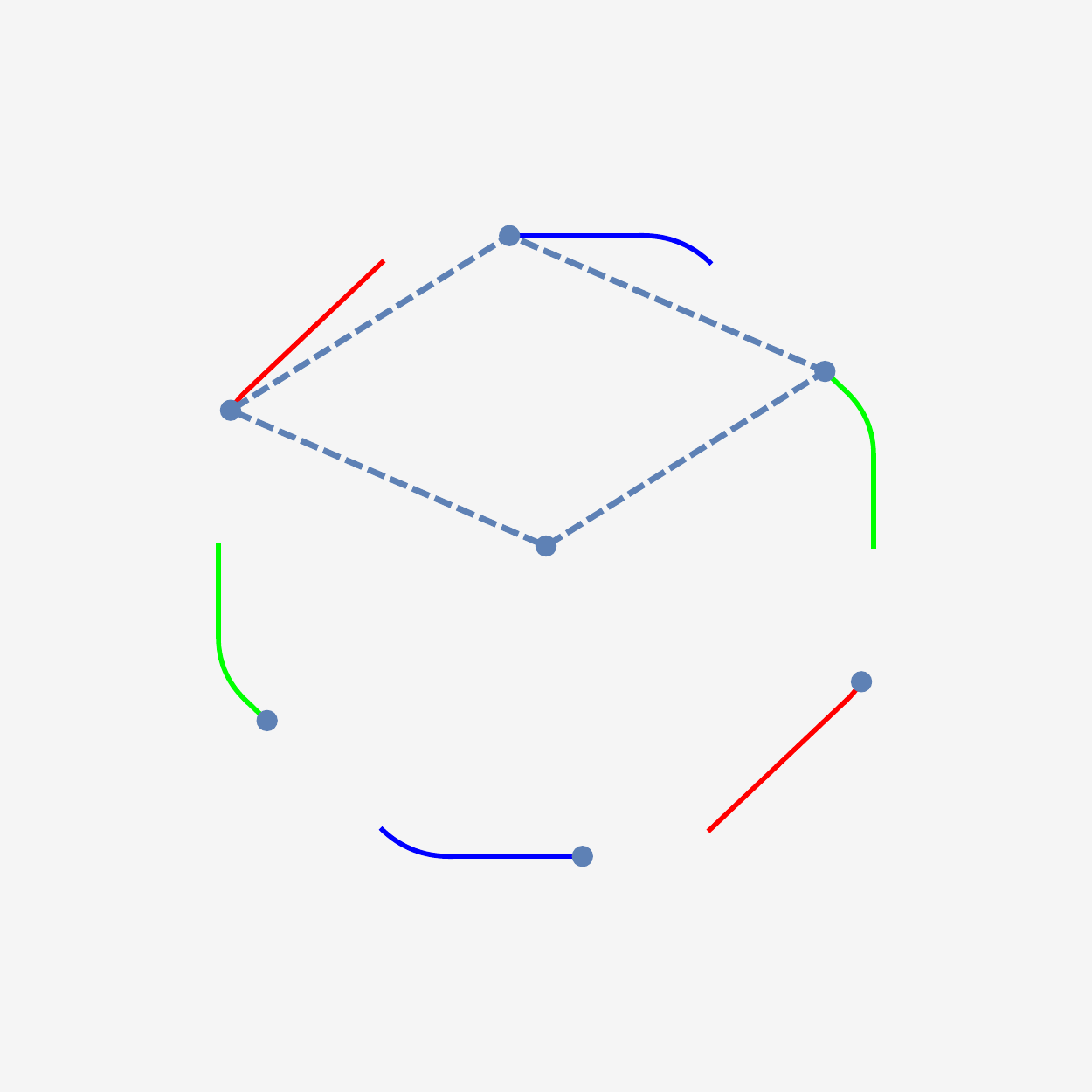}}%
\hspace*{0.005\textwidth}%
\subfloat[]{\includegraphics[width=0.33\textwidth]{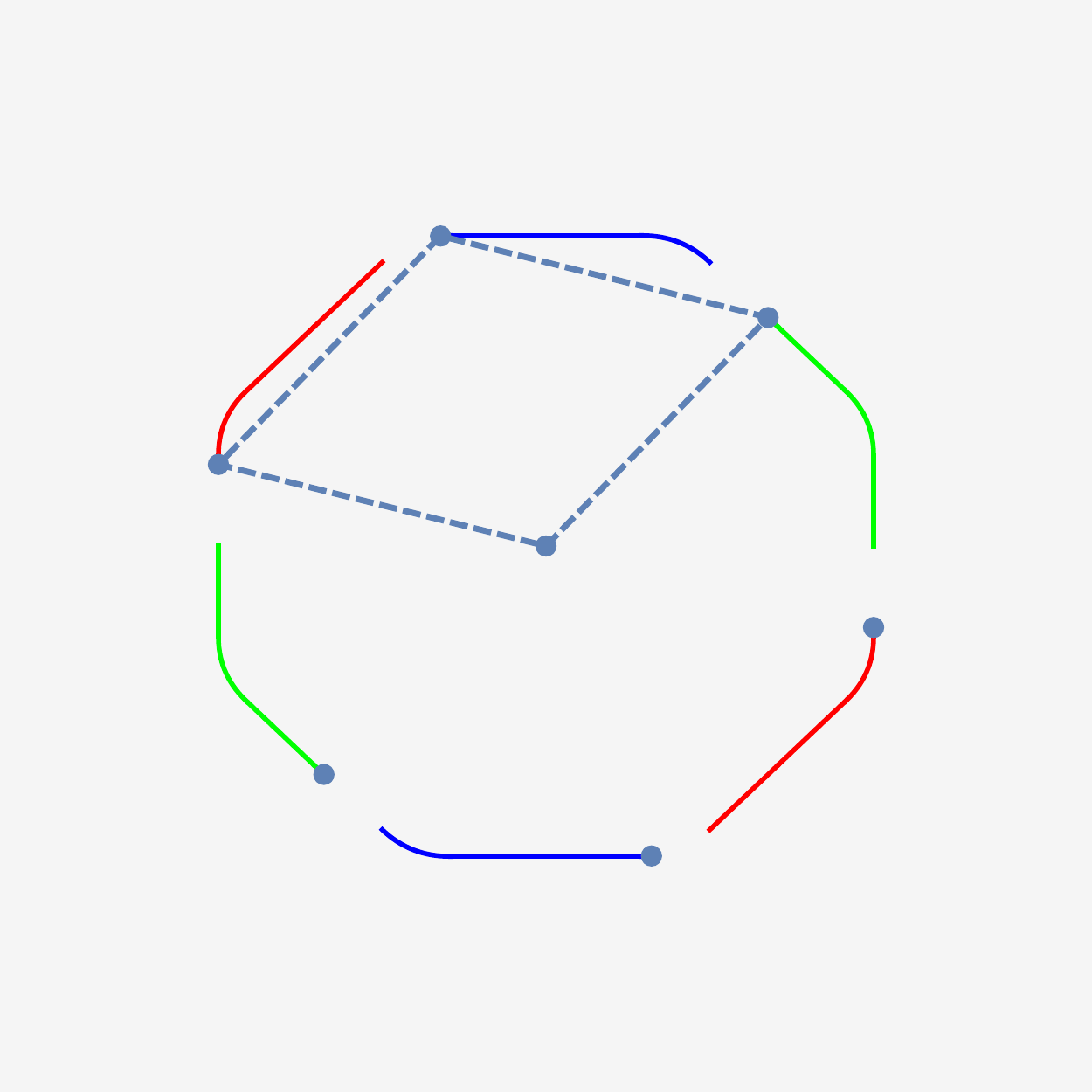}}
\subfloat[]{\includegraphics[width=0.33\textwidth]{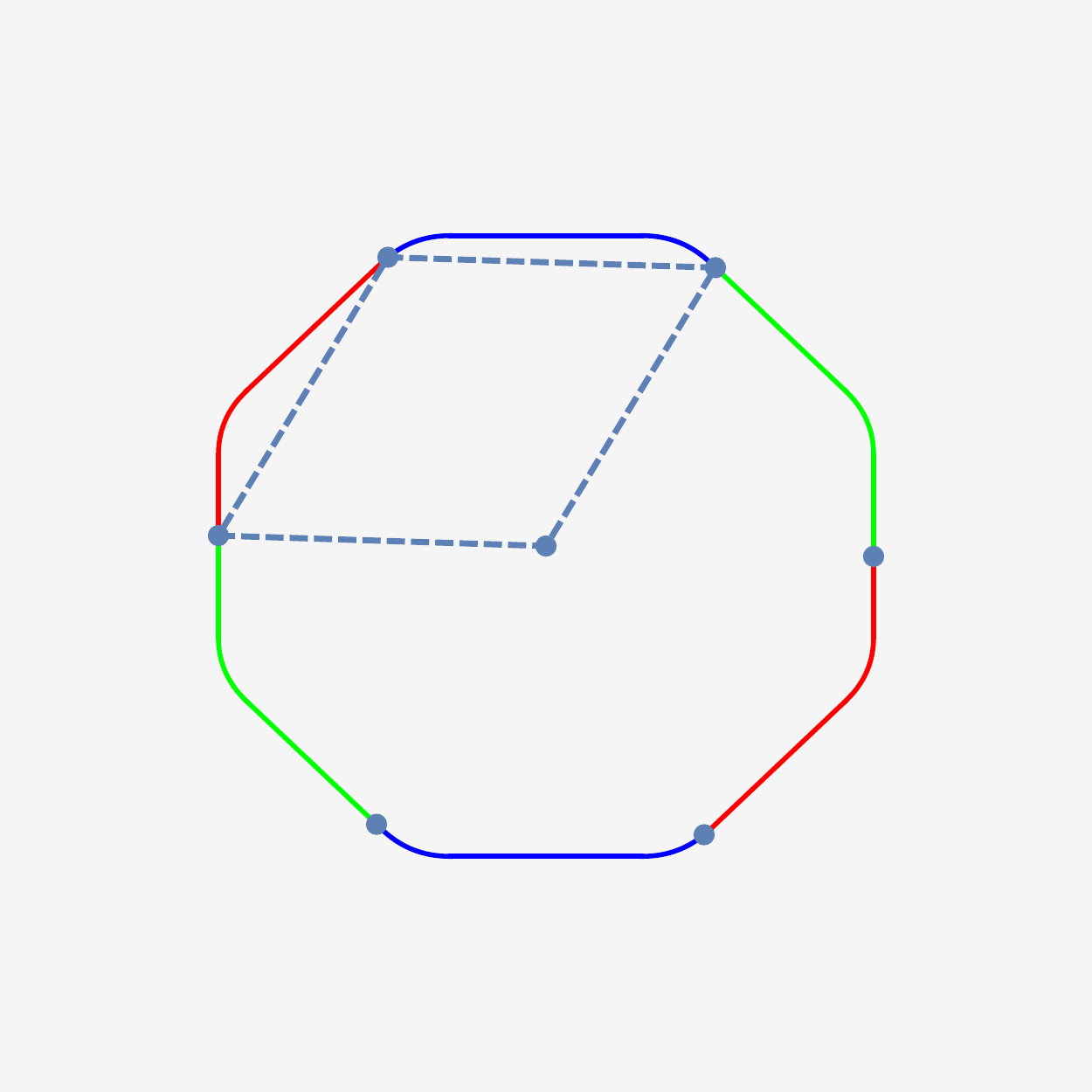}}

\end{figure}

In our opinion, the most promising line of attack to prove the Reinhardt conjecture begins with a paper Hales~\cite{hales2017reinhardt} in 2017 in which the Reinhardt problem is reduced to an optimal control problem on the tangent bundle of the Lie group $\SL(\R)$. The control set for this problem is the 2-simplex in $\R^2$. The properties proved by Reinhardt himself in 1934 play an essential role in this reformulation. As an example, Reinhardt proves that the boundary of the minimizer is generated by three curves which close up seamlessly, the points which generate them moving so that the area of the parallelogram with these three vertices and the origin remains fixed. This is shown in Figure \ref{fig:multi-curve-oct}.

The points are also required to move so that convexity of the final disc is not violated. Convexity is required \textit{locally} via a local curvature positivity condition and \textit{globally}, imposing conditions on the tangents to the six curves. The curvatures of these curves play a role in determining the packing density of the resultant disc in the plane. The control problem reformulation takes all of these into account.

The smoothed octagon was shown to be a critical point of the optimization problem (a Pontryagin extremal) given by a bang-bang control.

\begin{theorem}[Hales~\cite{hales2017reinhardt}]
The smoothed octagon is part of a family of Pontryagin extremals (smoothed $(6k+2)$-gons for every positive integer $k$) of the Reinhardt Optimal Control problem. The associated control is a bang-bang control.
\end{theorem}


Bang-bang controls arise frequently in optimal control problems and are characterized by the abrupt switching of the control between the extremes of a control set. The smoothed octagon also displays such extreme behaviour: the curvature of the curves which make up the boundary alternate between not curved at all (straight lines) to being as curved as possible (hyperbolic arcs). This insight explains its shape. 

Chapter II of this thesis will focus on exploring the the Reinhardt Optimal Control Problem (ROC), and will highlight its remarkable structure through its connections with hyperbolic geometry, Hamiltonian mechanics, conservation laws, and the theory of chattering control.

First order necessary conditions for optimality of an optimal control problem are given by the Pontryagin Maximum Principle (PMP). The PMP states that the optimal trajectory is given by a projection of the lifted controlled trajectory (living in the cotangent bundle). The lifted controlled trajectory is the Hamiltonian flow of the so-called \textit{maximized Hamiltonian}, which is the pointwise maximum of a control-dependent Hamiltonian function on the control set.

Our main insight is to change the shape of the control set (and hence also the control problem) from the 2-simplex to its circumscribing disk. This modification allows us to manufacture a conserved quantity by exploiting the rotational symmetry of this disk. This we do by appealing to a control theoretic version of the classical Noether theorem. 

\begin{theorem}
For the control problem with the circumscribing disk control set, the optimal trajectory satisfies a conservation law, which we call the angular momentum. 
\end{theorem}

This conservation law gives us valuable information about the optimal control. We use this to analyze the behaviour of the problem near the so-called \textit{singular locus}. Singular extremals in optimal control problems arise when the maximization condition in the Pontryagin maximum principle fails to give a unique minimizer for an interval of time. For our scenario, we have the following characterization of singular extremals.

\begin{theorem}
The circular disc in the plane is a singular extremal of the Reinhardt Optimal Control problem. 
\end{theorem}

Initial conditions in the extended state space giving rise to the circular disc in the plane is termed the singular locus. In the 2017 article, Hales proves a characterization of Pontryagin extremals. 

\begin{theorem}[Hales~\cite{hales2017reinhardt}]
All Pontryagin extremals of the Reinhardt control problem which do not pass through the singular locus are given by bang-bang controls with finitely many switches.
\end{theorem}

This means that the only to approach the singular locus is through \textit{chattering}. Chattering control happens when the control function performs abrupt and increasingly quicker transitions in the extremes of the control set in order to approach the singularity. They were first studied in a problem of Fuller and was considered a pathology for a time, but were proven to be \textit{ubiquitous} in a very precise sense by Kupka~\cite{kupka2017ubiquity} in the 1990s. The main result if this thesis is the recovery of the Fuller optimal control system in a neighbourhood of the singular locus. 

\begin{theorem}
Near the singular locus, the orders of growth of the state and costate variables (discarding the Lie group term) are 1,2 and 3 respectively.
\end{theorem}

Using this result, we perform a truncation of the PMP system by expanding in terms of increasing orders of growth and then throwing away all higher-order terms. We then recover the PMP system of the Fuller optimal control problem (of odd chain length).  

All results in Chapter II are joint work with Thomas Hales.

\chapter{The Reinhardt Conjecture}
We shall call a compact set in $\R^n$ a \textit{body}, and a body in $\R^2$ shall be called a \textit{disc}. 
By a \textit{centrally symmetric} disc in the Euclidean plane, we mean a compact subset  $K$ of $\mathbb{R}^2$ such that if $x \in K$ then $-x \in K$. Here, and throughout this chapter, we assume the center of symmetry is the origin $O = (0,0)$. A disc $K$ is a \textit{convex} disc, if for all $x, y \in K$, the line segment through $x,y$ is also in $K$. That is, $\{ \lambda x + (1 - \lambda) y ~|~ 0 \leq \lambda \leq 1 \} \subset K$. We denote by $\Kccs$ the class of all convex and centrally symmetric discs in the plane $\R^2$, which have the origin as the center of symmetry. 

A family of convex discs in $\R^2$ is called a \textit{packing} if any two distinct bodies in the family have a disjoint interior. We can now define the \textit{packing density} and \textit{best packing density} of a packing in $\R^2$. Intuitively the packing density corresponds to the proportion of the plane taken up by the packing.

\begin{definition}[Best packing density]
Let $B_t$ be a ball of radius $t$ in $\R^2$ centered at the origin and let $\mu$ be the Lebesgue measure on $\R^2$. The \textit{upper} and \textit{lower density} of a packing $\packing$ are defined to be
    \[ 
        \limsup_{t \rightarrow \infty} \frac{1}{\mu(B_t)}\sum_{C \in \packing} \mu(C \cap B_t) \qquad \mathrm{and} \qquad
        \liminf_{t \rightarrow \infty} \frac{1}{\mu(B_t)}\sum_{C \in \packing} V(C \cap B_t)
    \]
respectively. If they both exist and coincide, the common number is called the \emph{density of the packing} $\packing$ and is denoted $\delta(K,\packing)$. Given a convex body $K$ we define the \emph{best packing density} as the packing density formed with congruent copies of $K$:
\[
\delta(K) := \sup \{ \delta(K,\packing)\ | \ \mathcal{\packing} \text{ is a packing with congruent copies of }K \}
\]
\end{definition}

A \textit{lattice} is simply a discrete additive subgroup of $\R^2$ of full rank.  An important class of packings are \textit{lattice packings}, which consist of lattice translates of of a convex disc $K$: If $\LL$ is a lattice in $\R^2$ and $K$ is a fixed convex disc, then we consider the \textit{packings of translates of} $K$ under $\LL$ (called the \textit{lattice packing} of $K$). So $K + l \in K + \LL$ is a lattice translate of the convex disc $K$. We can now similarly define the lattice packing density and best lattice packing density. 

\begin{definition}[Best lattice packing density]
We define the \textit{upper} and \textit{lower densities} of a lattice packing of congruent copies of a convex disc $K$ to be 
\[ \limsup_{t \rightarrow +\infty} \frac{\sum_{l \in \LL}\mu(B_t \cap (K + l))}{\mu(B_t)} \qquad \mathrm{and} \qquad  \liminf_{t \rightarrow +\infty} \frac{\sum_{l \in \LL}\mu(B_t \cap (K + l))}{\mu(B_t)}\]
respectively. If they both exist and coincide, then the common number is called the \emph{density of the lattice packing}, denoted by $\delta(K,\LL)$ and the \emph{best lattice packing density} is defined as:
\[
\delta_L(K) := \sup \{ \delta(K,\LL) \ | \ \LL \text{ a lattice in } \R^2 \}
\]
\end{definition}

\begin{remark}\normalfont\leavevmode
\begin{itemize}
    \item It can be proved (see Exercise 3.2 of Pach \& Agarwal~\cite{pach2011combinatorial}) that for a convex disc $K$, one can always find a packing $\mathcal{P}$ such that $\delta(K,\mathcal{P})$ exists and is equal to $\delta(K)$. 
    \item It can also be proved (see Corollary 30.1 of Gruber~\cite{gruber2007convex}) that, given a convex disc $K$ and a lattice $\LL$, the upper and lower lattice packing densities of the packing given by lattice translates of $K$ by $\LL$ (provided that the $\LL$-translates of $K$ have disjoint interiors) coincide and are both equal to: \[\delta(K,\LL) = \frac{\mu(K)}{\det(\LL)}\]
    where $\det(\LL)$ is the determinant of the lattice $\LL$ (see Definition \ref{def:lattice-determinant}).
\end{itemize}
\end{remark}

Now consider the following quantity:
\[ \delta_{\min} := \inf_{K \in \Kccs} \delta(K) \]
So $\delta_{\min}$ is the worst best packing density among all convex centrally symmetric discs in $\R^2$.

Reinhardt's problem now is to explicitly describe a $K_{\min} \in \Kccs$ for which $\delta(K_{\min}) = \delta_{\min}$, and also determine this least packing density. Reinhardt suggested a specific candidate, called the \emph{smoothed octagon} to have the worst packing. The smoothed octagon is a regular octagon whose vertices have been clipped by hyperbolic arcs (shown in Figure \ref{fig:smoothed}).  

\begin{conjecture} [Reinhardt \cite{reinhardt1934dichteste}, Mahler \cite{mahler1947minimum}]
The smoothed octagon achieves the least best packing density among all other convex, centrally symmetric discs in the plane. Its density is given by 
\[ \frac{8 - \sqrt{32} - \ln 2}{\sqrt{8} - 1} \approx 0.902414 \]
\end{conjecture}

We aim to resolve the Reinhardt conjecture by restating it as a problem in control theory. To do this, we rely on a few geometric facts characterizing our minimizer $D_{min}$ which we collect in the next section.

\section{Background Results}
Even in the early 1900s the connection between convex, centrally symmetric discs and lattices was realized and explored, most prominently by Minkowski, who proved his famous theorem on lattice points. Apart from this classic theorem, Minkowski also proved other results which are central to our approach to the Reinhardt conjecture. Before stating these results, we set up a few definitions.

\begin{definition}[Support line]
For a convex disc $K$, a support line is a line containing at least one point of $K$ but does not separate any points of $K$. 
\end{definition}

\begin{definition}[Admissible lattice]
For a $K \in \Kccs$ centered at the origin, a lattice $\LL$ is called $K$-admissible if no point of $\LL$ other than $O=(0,0)$ lies in the interior of $K$. 
\end{definition}

\begin{definition}[Determinant of a Lattice]\label{def:lattice-determinant}
 For any lattice $\LL = \{ rE_0 + sE_1~|~ r,s \in \mathbb{Z} \}$ with basis $E_0,E_1 \in \R^2$, the \emph{determinant} $\det(\LL)$ is equal to $|E_0 \wedge E_1|$, which is simply the absolute value of the determinant of the $2\times2$ matrix having $E_0$ and $E_1$ as columns. This is also sometimes called the \emph{covolume} of the lattice $\LL$.
\end{definition}

\begin{definition}[Minimal determinant]
For $K\in\Kccs$, the quantity \[ \Delta(K) := \inf_{K-\mathrm{admissible}} \det(\LL) \] where the infimum runs over all $K$-admissible lattices is called the \emph{minimal determinant} of the convex disc $K$.
\end{definition}

\begin{definition}[Critical lattice]
 A lattice is called \emph{critical} for a disc $K$ if its determinant is equal to the minimal determinant of the disc $K$. 
\end{definition}

\begin{definition}[Irreducible disc]\label{def:irreducible}
 A disc $K \in \Kccs$ is called \emph{irreducible} if every boundary point of $K$ lies on a critical lattice of $K$.
\end{definition}

Note that this is not the original definition of irreducibility of a convex disc. We choose our definition based on Lemma 3 of \cite{mahler1947minimum}. A convex disc $K \in \Kccs$ which is irreducible is what Hales~\cite{hales2011} calls a ``hexameral domain'', a terminology which is later dropped in \cite{hales2017reinhardt}.

\subsection{The hexagons \texorpdfstring{$h_K$}{hK} and \texorpdfstring{$H_K$}{HK}.}

Minkowski proved the following theorem which gives conditions under which points on $K$ give rise to critical lattices. We follow the presentation by Mahler~\cite{mahler1947minimum}.

\begin{theorem}[Minkowski~\cite{minkowski1907diophantische}, Mahler~\cite{mahler1946lattice}]\label{thm:minkowski-inscribe}
Let $\LL$ be a critical lattice of a disc $K \in \Kccs$. Then $\LL$ contains three points $E_0,E_1,E_2$ on the boundary of $K$ such that (i) $E_0,E_1$ is a basis of the lattice $\LL$, and (ii) $OE_0E_1E_2$ is a parallelogram of area $\det(\LL) = |E_0 \wedge E_1| = \Delta(K)$, the minimal determinant of $K$. Conversely, if $E_0,E_1,E_2$ are three points on the boundary of $K$ such that $OE_0E_1E_2$ is a parallelogram, then the area of this parallelogram is not less than $\Delta(K)$ and is equal to $\Delta(K)$ if and only if the lattice with basis $E_0,E_1$ is critical. 
\end{theorem}


Since centrally symmetric hexagons can be decomposed into three parallelograms, the above result shows that a critical lattice of a disc $K \in \Kccs$ gives rise to an \textit{inscribed centrally symmetric hexagon} $h_K$ within our disc $K$ so that $\Delta(K) = \frac{1}{3}\mu(h_K)$ which is \textit{minimal} in the sense that 
\[ 
\mu(h_K) = \inf_{h} \mu(h)
\] where the infimum is taken over all hexagons with $h$ with vertices $\pm E_0,\pm E_1,\pm E_2$ on the boundary of $K$ and with $E_0 - E_1 + E_2 = 0$. In 1947, Mahler~\cite{mahler1947minimum} proved an analogous result which talks about \textit{circumscribed hexagons} of $K$:

\begin{theorem}[Mahler~\cite{mahler1947minimum}]\label{thm:mahler-circumscribe}
Let $K \in \Kccs$ be a disc which is not a parallelogram; let $\LL$ be a critical lattice of $K$ and let $\pm E_0,\pm E_1,\pm E_2$ with $E_0 - E_1 + E_2= 0$ be lattice points of $\LL$ on the boundary of $K$. Then there are unique support lines $\pm L_0,\pm L_1, \pm L_2$ of $K$ at these points such that
\begin{enumerate}
    \item no two of these lines coincide
    \item the hexagon $H_K$ bounded by the support lines is of area $\mu(H_K) = 4\Delta(K)$
    \item each side of $H_K$ is bisected at the lattice point $\pm E_i$ where it touches the boundary of $K$.
    \item the hexagon $H_K$ is \emph{minimal} in the sense that 
\[ 
\mu(H_K) = \inf_{H} \mu(H)
\]
where the infimum is taken over the set of all hexagons $H$ bounded by supporting-lines $\pm L_1',\pm L_2',\pm L_3'$ of the disc $K$.
\end{enumerate}
\end{theorem}

Thus, the critical lattice of any $K \in \Kccs$ gives rise to hexagons $h_K$ and $H_K$ whose areas are related as
\begin{equation}\label{eq:outer-inner-area}
    \Delta(K) =  \frac{1}{3}\mu(h_K) = \frac{1}{4}\mu(H_K).
\end{equation}

\begin{definition}[Balanced/Critical Hexagon]
  For a disc $K \in \Kccs$, the hexagon $H_K$ is called its \emph{balanced hexagon} or \emph{critical hexagon}.
\end{definition}

\begin{proposition}
If $K \in \Kccs$ then $h_K,H_K \in \Kccs$. 
\end{proposition}
\begin{proof}
By construction.
\end{proof}

\begin{figure}
    \centering
    \includegraphics{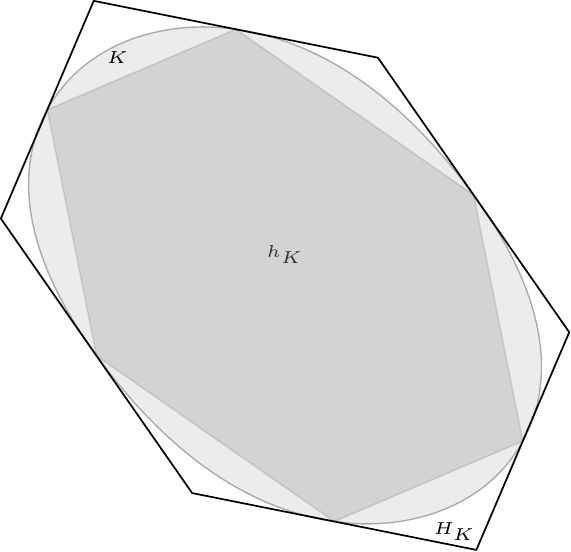}
    \caption{Critical hexagons for an ellipse.}
    \label{fig:hexagons}
\end{figure}

The hexagons $H_K$ and $h_K$ for an ellipse are shown in Figure \ref{fig:hexagons}. Our interest in the critical hexagon of a convex disc is primarily due to this theorem:

\begin{theorem}[Pach \& Agarwal~\cite{pach2011combinatorial}]\label{thm:packing-construction}
For a disc $K \in \Kccs$:
\begin{enumerate}
    \item The critical hexagon $H_K$ is the hexagon of smallest area circumscribing the disc $K$.
    \item a packing of density $\delta(K)$ of congruent copies of $K$ is constructed by tiling the plane with copies of $H_K$ and inscribing $K$ inside each copy. 
\end{enumerate}
\end{theorem}
\begin{proof}
Any hexagon of smallest area circumscribing $K$ can be chosen to be centrally symmetric by Theorem 2.5 of Pach \& Agarwal~\cite{pach2011combinatorial}. We get the first statement by Theorem \ref{thm:mahler-circumscribe},(4). Since $H_K$ is centrally symmetric, it tiles the plane. The second statement is Corollary 3.6 of Pach \& Agarwal~\cite{pach2011combinatorial}.
\end{proof}

\subsection{Properties of \texorpdfstring{$K_{\min}$}{Dmin}.}

Our starting point is the following theorem of L. Fejes T\'{o}th which states that the best packing density and best lattice packing densities are actually equal for the class $\Kccs$: 
\begin{theorem}[Fejes T\'{o}th]\label{thm:fejes-toth-packing}
    If $K \subset \R^2$ is a \textit{centrally symmetric} convex disc, then  \begin{equation}\label{eq:fejes-toth}
    \delta(K) = \delta_L(K) = \frac{\mu(K)}{\mu(H_K)}
    \end{equation}
    where $\mu(K)$ is the area of $K$ and $H_K$ is the critical hexagon of $K$.
\end{theorem} 
The main import of this theorem is the first equality which states that the best packing density of a disc $K \in \Kccs$ is given by a lattice packing. Theorem \ref{thm:packing-construction},(2) shows that this packing arises by first tiling the plane with copies the centrally symmetric critical hexagon $H_K$, and then inscribing $K$ inside each copy. The second equality in \eqref{eq:fejes-toth} was known for at least a few decades before Fejes T\'{o}th. In fact, it was proven by Reinhardt himself (see the introduction in Mahler \& Ledermann~\cite{ledermann1949lattice}) in 1934.

So, using Theorem \ref{thm:fejes-toth-packing}, Reinhardt's conjecture is now reduced to the following problem:

\begin{problem}[Reinhardt conjecture]\label{pbm:reinhardt-conjecture}
Describe those $K_{\min} \in \Kccs$ for which
\[
\delta(K_{\min}) = \inf_{K \in \Kccs}\frac{\mu(K)}{\mu(H_K)}.
\]
\end{problem}

In this form, this problem was studied by multiple authors: In 1904, Minkowski~\cite{minkowski1907diophantische} established a lower bound. In 1912, Blaschke~\cite{blaschke1945differentialgeometrie} called \textit{Courant's conjecture} the statement that the elliptical disk (as a convex, centrally symmetric disc) is $K_{\min}$.
Reinhardt~\cite{reinhardt1934dichteste} proved that minimizer of this quantity exists, and proved several properties about it, including the fact that the ellipse is not the minimizer, refuting Courant's conjecture:

\begin{theorem}[Reinhardt~\cite{reinhardt1934dichteste}, Mahler~\cite{mahler1947area}]\label{thm:reinhardt-existence}
There exists a convex, centrally symmetric disc $K_{\min}$ for which $\delta(K_{\min}) = \delta_{\min}$ and it has the following properties:
\begin{enumerate}
    \item $K_{\min}$ is not the ellipse. 
    \item $K_{\min}$ is an irreducible disc. 
    \item The boundary of $K_{\min}$ has no corners \emph{i.e.,} the boundary of $K_{\min}$ has at all points a unique support line. 
    \item Every point on the boundary of $K_{\min}$ lies on a \emph{unique} critical lattice. (Equivalently, every point on the boundary of $K_{\min}$ is the midpoint of a side of a unique critical hexagon.)
\end{enumerate}
\end{theorem}

We remark that the quantity $\mu(K)/\mu(H_K)$ in Problem~\ref{pbm:reinhardt-conjecture} is affine invariant and so there is no loss of generality in fixing the area of $\mu(H_K)$ and then considering the minimization problem to be over all $K \in \Kccs$ which have $H_K$ with that fixed area. Mahler~\cite{mahler1947area} chooses $\mu(H_K) = 4$, while Hales~\cite{hales2011} chooses $\mu(H_K) = \sqrt{12}$. We choose the latter, for reasons which will become apparent below. 
Using this and Theorem \ref{thm:reinhardt-existence}, we can reduce our problem to the following:


\begin{definition}\label{def:Kbal}
Define $\Kbal \subset \Kccs$ to be the subset of all irreducible discs whose boundary is subject to the conditions imposed in Theorem \ref{thm:reinhardt-existence}. 
\end{definition}
The set $\Kbal$ is nonempty since Reinhardt proved that the smoothed octagon and the circular disc both belong to this set. 

\begin{problem}[Reinhardt conjecture]\label{pbm:reinhardt-reduction}
Describe that $K_{\min} \in \Kbal$ for which
\[
\delta(K_{\min}) = \inf_{\{K \in \Kbal~|~\mu(H_K)=\sqrt{12}\}} \mu(K).
\]
\end{problem}
It would be prudent to get a better description of the subset $\Kbal$ to get a more well-defined minimization problem. This is what we proceed to do in the next subsection.

\subsection{Boundary Parameterization}
Let $K \in\Kbal$ be an arbitrary convex, centrally symmetric disc. By Theorem \ref{thm:minkowski-inscribe}, each point on the boundary of $K$ has two \textit{companion points} which, together with the origin, give rise to a parallelogram of area $\Delta(K)$. These three points give rise to three other points by central symmetry. Thus, every point $E_0$ on the boundary of $K$ gives rise to five other points $E_1,E_2,-E_0,-E_1,-E_2$ (vertices of $h_{K}$) also on the boundary. Let us order these points $E_0,E_1,\dots,E_5$. Then we have that $E_0 + E_2 + E_4 = 0$, $E_4 = -E_1$ and $|E_0 \wedge E_2| = \Delta(K)$, which is a fixed quantity. Furthermore, the area of the critical hexagon at each point is the same. By equation \eqref{eq:outer-inner-area} and Theorem \ref{thm:minkowski-inscribe}, the area of the critical hexagon is equal to a fixed fraction of the area of the parallelogram formed by the points $O,E_0,E_1$ and $E_2$. 

Since the boundary of $K$ does not contain any corners, we can parameterize the boundary by a regular $C^1$ curve $t \mapsto \sigma_0(t)$ with the convention that the boundary is parameterized in the counter-clockwise direction. We shall call this the \emph{positive orientation}. Then, at each time instant $t$, by the above discussion, the point $\sigma_0(t)$ gives rise to other points $\sigma_1(t),\sigma_2(t),\sigma_3(t),\sigma_4(t),\sigma_5(t)$ which are subject to the following conditions at each instant $t$:
\begin{align*}
    \sigma_0(t) + \sigma_2(t) + \sigma_4(t) &= 0\\
    \sigma_4(t) &= -\sigma_1(t) \\
    \sigma_0(t) \wedge \sigma_2(t) &= \Delta(K_{\min}) = \frac{\sqrt{3}}{2}
\end{align*}

This inspires the following definition, following Hales~\cite{hales2011}:
\begin{definition}[Multi-point and multi-curve]\label{def:multi-pt-multi-curve}
  A function $E : \Z/6\Z \to \R^2$ such that 
  \begin{equation}\label{eq:multi-point-properties}
      E_j + E_{j+2} + E_{j+4} = 0, \quad E_j = -E_{j+3}, \quad E_j \wedge E_{j+2} = \frac{\sqrt{3}}{2}
  \end{equation}
  is called a \emph{multi-point}. An indexed set of $C^1$ curves $\sigma :\Z/6\Z  \times [0,t_f] \to \R$ is a \emph{multi-curve} if for all $t \in [0,t_f]$, $\sigma(t)$ is a multi-point.
\end{definition}

\begin{example}
\normalfont\leavevmode
\begin{itemize}
    \item The sixth roots of unity $E_j = e_j^* = \exp\left(\frac{2\pi ij}{6}\right) \in \mathbb{C}$, viewed as points in $\R^2$, is an example of a multi-point. 
    \item For an example of a curve in $\R^2$ which satisfies the multi-curve properties in equation \eqref{eq:multi-point-properties}, see Section \ref{sec:hypotrochoids} in the Appendix. 
\end{itemize}

\end{example}

\subsubsection{Regularity Properties of Multi-Curves}
We would like to characterize the discs in $\Kbal$. To this end, first we establish the following Lemmas:


\begin{lemma}[Reinhardt~\cite{reinhardt1934dichteste}, Mahler~\cite{mahler1947area}, Hales~\cite{hales2011}]
    If $t \mapsto \sigma_0(t)$ is a positively oriented $C^1$ regular curve parameterizing the boundary of $K \in\Kbal$, then so is $\sigma_j(t)$ for $j=1,\dots,5$.
\end{lemma}
\begin{proof}
Given that $\sigma_0(t)$ is continuous, that $\sigma_1(t)$ and $\sigma_2(t)$ are continuous is proven in Lemma 9 of Mahler~\cite{mahler1947irreducible}. Given $\sigma_0'(t)$ is continuous, $\sigma_2'(t)$ is proven to be continuous in Lemma 11 of Hales~\cite{hales2011}. Similar statements for other curves follow by symmetry. 
\end{proof}

\begin{lemma}[Hales~\cite{hales2011}]\label{lem:tangents-lipschitz}
    Let $\sigma(t)$ denote a multi-curve parameterization of the boundary of $K \in\Kbal$.  Assume that the curve $\sigma_0$ is parameterized according to arc-length. Then, the tangents $\sigma_j'(t)$ are Lipschitz continuous for all $j$.
\end{lemma}
\begin{proof}
This is Lemma 18 of Hales~\cite{hales2011}.
\end{proof}

\begin{corollary}\label{lem:multicurve-c2}
The functions $\sigma_j'(t)$ are differentiable almost-everywhere.
\end{corollary}
\begin{proof}
Follows by Rademacher's theorem and Lemma \ref{lem:tangents-lipschitz}.
\end{proof}
Henceforth, we shall assume that the curve $t \mapsto \sigma_0(t)$ is parameterized according to arc-length. 
\subsubsection{Convexity of Multi-Curves}
We need to impose conditions on the curves $\sigma$ which make up the boundary of a $K \in\Kbal$ so that they enclose a convex disc, since this is not guaranteed \textit{a priori}.

\begin{figure}
    \centering
    \includegraphics{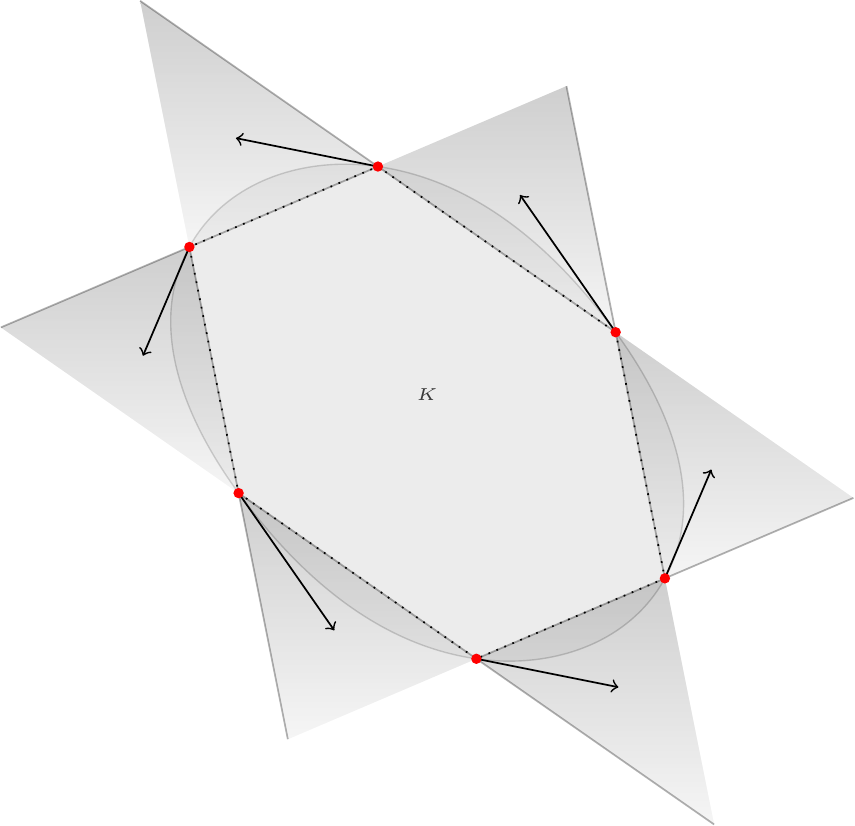}
    \caption{Global convexity condition for the ellipse.}
    \label{fig:ellipse-global-convex}
\end{figure}

\begin{lemma}[Star conditions]\label{lem:star-conditions}
    Let $K \in\Kbal$ be a disc with boundary parameterized by the multi-curve $\sigma$. Let $\sigma(t) = E$ be a multi-point on the boundary of the disc $K$. The convexity condition on $K$ forces the tangent vectors  $\sigma_j'(t)$ at time $t$ to point into the open cone with sides determined by the triangle with vertices $E_j,E_{j+1}$ and $E_j + E_{j+1}$ for each $j \in \Z/6\Z$.
    The tangent vector $\sigma_j'$ is also never parallel to the secant lines through $E_j$ and $E_{j+1}$. 
\end{lemma}
\begin{proof}
This situation is depicted in Figure~\ref{fig:ellipse-global-convex}.
This is asserted in Hales~\cite{hales2011,hales2017reinhardt} and is called the ``star condition''.
By convexity of the disc $K$, at any time $t$ the hexagon $h_K \subset K$ (as $h_K$ is the convex hull of the points $\{\sigma_j(t)\}$). Now, the vector $\sigma_j'(t)$ cannot point into the hexagon, because it would then create a non-convex piece of the curve $\sigma_j$ by continuity. Dually, it cannot point beyond the line through $\sigma_j(t) + \sigma_{j+1}(t)$ and $\sigma_j(t)$, as that would force $\sigma_{j+2}'(t)$ to point inward. 

If the vector $\sigma_j'(t)$ points along the edges of the triangle, then it would have to remain pointing in that same direction, as it cannot point inward (by the above argument) or outward (as then it would not be convex). This would then create a corner and $K$ would not belong to $\Kbal$. 

\end{proof}

Apart from the star conditions, there is another condition on the curvature of the boundary curve which needs to be imposed:

\begin{lemma}[Curvature constraint]\label{lem:local-curvature-constrt}
For a multi-curve $\sigma$ parameterizing the boundary of a disc $K \in\Kbal$, we have the following condition almost everywhere:
    \begin{equation}\label{eq:local-convexity}
    \sigma_j'(t) \wedge \sigma_{j}''(t) \ge 0 \quad j \in \Z/6\Z
\end{equation}
\end{lemma}
\begin{proof}
This is well-defined since we have seen in Lemma~\ref{lem:multicurve-c2} that the functions $\sigma_j'(t)$ are differentiable almost everywhere. 
A well-known theorem (see Proposition 3.8 of  Shifrin~\cite{shifrin2015differential}) states that a simple closed regular plane curve is convex if and only if its orientation can be chosen in such a way so that its curvature $\kappa \ge 0$ everywhere. The assertion follows.
\end{proof}

\subsection{A characterization of \texorpdfstring{$\Kbal$}{B}}
Summarizing, we have the following conditions characterizing the boundary of an arbitrary disc $K \in\Kbal$:
\begin{enumerate}
    \item The boundary of $K$ is parameterized by six regular $C^1$ curves $\sigma_j(t)$. 
    \item The functions $\sigma_j'$ are differentiable almost-everywhere. 
    \item For each $t$, $\{\sigma_j(t)\}$ is a multi-point. 
    \item For almost all $t$, we have the curvature constraint $\sigma_j'(t) \wedge \sigma_{j}''(t) \ge 0$.
    \item For each $t$, the vector $\sigma_j'(t)$ points into the open cone determined by the triangle with vertices $\sigma_j(t),\sigma_{j+1}(t)$ and $\sigma_j(t) + \sigma_{j+1}(t)$.
    \item The six curves $\{\sigma_j\}$ close up seamlessly so that the boundary of $K$ is a simple closed curve. 
\end{enumerate}
In our restatement of the Reinhardt conjecture (as Problem \ref{pbm:reinhardt-reduction}) we indicated that we are searching over all discs in $K \in \Kbal$ whose critical hexagon area $\mu(H_K) = \sqrt{12}$. By affine invariance, there is no loss of generality in assuming that the sixth roots of unity $E = \{e_j^*\}$ lie on the boundary of the disc $K$. We make this multi-point the \emph{initial} multi-point that the six curves $\sigma$ start out at. This explains the reason why we choose $\sqrt{12}$, since this is the area of the regular hexagon which has the sixth roots of unity as the midpoints of each side. 
Following Hales~\cite{hales2011}, we call this representation the \emph{circle representation of the disc} $K$. 

\begin{problem}[Reinhardt conjecture]\label{pbm:reinhardt-reduction-2}
Describe that $K_{\min} \in \Kbal$ which is in circle representation for which
\[
\delta(K_{\min}) = \inf_{\{K \in \Kbal~|~ K\text{ in circle representation}\}} \mu(K)
\]
where $\Kbal \subset \Kccs$ consists of all irreducible discs whose boundary is parameterized by multi-curves $\sigma$ subject to conditions listed above.
\end{problem}

\section{Derivation of the control problem}
Now that we have a better description of the set $\Kbal$, we can use this to restate the Reinhardt conjecture as an optimal control problem. 
Recalling from the introduction the characteristic features of optimal control problems, we now aim to 
\begin{enumerate}
    \item recast the system above as a dynamical system on a \emph{state space} (which, in our case, will be a manifold),
    \item define a well-defined \emph{cost functional},
    \item determine a well-defined \textit{control parameter}.
\end{enumerate}This will be our focus in this section.    
\subsection{State Dynamics in the Lie Group}\label{sec:lie-group-dynamics}
The multi-curve conditions give rise to a curve in $\SL(\R)$ in the following sense: 

\begin{theorem}[Mahler~\cite{mahler1947area}, Hales~\cite{hales2011}]\label{thm:defn-g}
  Let $\sigma$ be a $C^1$ multi-curve parameterizing the boundary of a disc $K \in\Kbal$. Then $\sigma$ determines a $C^1$ curve $g$ in $\SL(\R)$. 
\end{theorem}
\begin{proof}
Given a multi-curve $\sigma$, let $t_0 < t_1$ be two different time instants so that $\sigma(t_0)$ and $\sigma(t_1)$ are multi-points on the boundary of $K$. If we can exhibit a $2 \times 2$ real matrix $g = g(t_0,t_1)$ so that 
\begin{align*}
\sigma_0(t_0) &=  g\sigma_0(t_1) \\
\sigma_2(t_0) &=  g\sigma_2(t_1)
\end{align*}
then we are done since this implies $\sigma_j(t_0) = g\sigma_0(t_1)$ for all $j \in \Z/6\Z$ by the other multi-point conditions. However, a unique such matrix can always be found whenever $\sigma_0(t_1) \wedge \sigma_2(t_1) \ne 0$, which is evidently the case here. The identity $g\sigma_0(t_1)\wedge g\sigma_2(t_1) = \det(g)\sigma_0(t_1) \wedge \sigma_2(t_1)$ and the multi-point condition forces $\det(g) = 1$, which means that $g \in \SL(\R)$. Thus, if we have a continuous curve of multi-points $\sigma(t)$, we get a unique induced continuous curve $g(t) \in \SL(\R)$ for $t \in [0,t_f]$.
\end{proof}

\begin{remark}
\normalfont \leavevmode
\begin{itemize}
    \item The above proof also implies that if $t_0,t_1,t_2$ are three time instants with $t_2 = t_0+t_1$ then $g(t_2,t_0) = g(t_2,t_1)g(t_1,t_0)$.
    \item The key fact which we use later is that if $K$ is in circle representation, then the multi-curves $\sigma_j$ are given by $\sigma_j(t) = g(t)e_{2j}^*$ where $g(t) \in \SL(\R)$.
\end{itemize}
\end{remark}

Similarly, the tangents $\{\sigma_j'(t)\}$ give rise to corresponding elements in the Lie algebra $\sl(\R)$.
\begin{definition}\label{def:X-defn}
For $g$ as above, and for every $t$, define $X(t) \in \mathfrak{gl}_2(\R)$ to be such that $g'(t) = g(t) X(t)$.
\end{definition}

\begin{theorem}[Hales~\cite{hales2011}]\label{thm:defn-X}
Assume that we have a disc $K$ in circle representation with boundary parameterized by a multi-curve $\sigma$. Let $g$ be the induced curve in $\SL(\R)$. Then
\begin{enumerate}
    \item $\sigma_j'(t) = \Ad_{g(t)}(X(t))\sigma_j(t)=g(t)X(t)g(t)^{-1}\sigma_j(t)$ for all $j \in \Z/6\Z$.
    \item $X$ is Lipschitz continuous. 
    \item $X(t) \in \sl(\R)$.
\end{enumerate}
\end{theorem}
\begin{proof}
First of all, the matrix $X\in\sl(\R)$ since if $g : [0,t_f] \to \SL(\R)$ is any differentiable curve, then $X := g^{-1}g' \in \sl(\R)$. We then have:
\[
\sigma_j'(t) = g'(t)e_{2j}^* = g(t)X(t)e_{2j}^*=g(t)X(t)g(t)^{-1}g(t)e_{2j}^*=\Ad_{g(t)}(X(t))\sigma_j(t),
\]
where $\Ad_{g(t)}(X(t))=g(t)X(t)g(t)^{-1}$ is the adjoint representation.
Since $\sigma_j'(t) = g(t) X(t)e_{2j}^*$, $\sigma_j'$ is Lipschitz by Lemma \ref{lem:tangents-lipschitz} and $g(t)$ being a $C^1$ curve on a compact interval is bounded, this shows that $X(t)$ is Lipschitz. 


\end{proof}
\begin{corollary}[Hales~\cite{hales2011,hales2017reinhardt}]\label{cor:X-props}
For a disc $K$ in circle representation, we have the following properties for $X$ defined in Definition \ref{def:X-defn}:
\begin{itemize}
    \item Setting $X = \mattwo{a}{b}{c}{-a} \in \sl(\R)$, we then get $\sqrt{3}|a| < c$ and $3b+c<0$. In particular, $c>0$.
    \item $\det(X(t))>0$ for all $t$.
\end{itemize}
\end{corollary}
\begin{proof}

Assume that we start at time $t=0$, at the multipoint given by the sixth roots of unity $\{e_j^*\}$.
The star conditions in Lemma \ref{lem:star-conditions} imply that the tangent vector $\sigma_j'(0) = X(0)e_{2j}^*$ lies in between the secant lines joining $e_j^*$ and $e_{j\pm1}^*$. This means that, taking $j=0,5$
\begin{align}
    re_0^* + s(e_0^* -e_5^*) &= Xe_5^*\\
    p(e_1^*-e_0^*) + qe_1^* &= Xe_0^*
\end{align}
and we must have $p,q,r,s > 0$. Solving these systems of linear equations for $p,q,r,s$, we get \[\sqrt{3}|a| < c \quad 3b+c<0.\] This proves that $c>0$. Using these relations, we have that $\det(X(0)) = -bc-a^2 \ge -bc -\frac{c^2}{3} = \frac{-c(3b+c)}{3} > 0$. 
The same results hold for all time, using Theorem \ref{thm:defn-g} by shifting time by multiplying with a curve in $\SL(\R)$.  
\end{proof}


We have one equation for our state space dynamics viz., equation $g'=gX$. Since we are deriving dynamics in the Lie group and Lie algebra, we shift to a more convenient choice of parameterization.

\begin{proposition}\label{prop:reparam-lipschitz}
Let $s$ denote the arc-length parameter of $\sigma_0$ and let $\tilde{X}(t)$ denote the matrix $X$ with respect to the reparameterizion of $g$ to a time variable $t$ making $\det(\tilde{X}(t))=1$. Then we have that $\tilde{X}(t)$ is a Lipschitz continuous function.
\end{proposition}
\begin{proof}
Reparameterize $\tilde{g}(t) := g(s(t))$. In this new parameterization we must have $d\tilde{g}(t)/dt= \tilde{g}(t)\tilde{X}(t)$.
We have by the chain rule:
\begin{align}
  \frac{d\tilde{g}(t)}{dt} = \frac{d}{dt}g(s(t)) &= 
    \frac{d}{ds}g(s(t))\frac{ds}{dt} \\
    &= g(s)X(s)\frac{ds}{dt}
\end{align}
So that $\tilde{X}(t) = X(s(t))\frac{ds}{dt}$. Now $\det(\tilde{X}(t))=1$. So we have:
\[
\frac{ds}{dt} = \frac{1}{\det(X(s))}
\]
which gives us the required reparameterization equation. By Corollary \ref{cor:X-props} we have $\det(X(s)) > 0$ and so this equation is well-defined. Recall that we have that $X(s)$ is Lipschitz by Theorem $\ref{thm:defn-X}$ and $\det(X(s))$ is bounded since it a continuous function on a compact interval. This proves that $\tilde{X}$ is Lipschitz.
\end{proof}

We shall abuse notation and also denote by $t$ the new parameterization, with the understanding that $t$ runs so that $\det(X(t)) = 1$.

\begin{corollary}\label{cor:X-diff-ae}
With respect to the parameterization making $\det(X(t)) = 1$, the curve $X(t)$ is differentiable almost everywhere. 
\end{corollary}
\begin{proof}
This follows from Rademacher's theorem and Proposition \ref{prop:reparam-lipschitz}.
\end{proof}

\begin{corollary}[Hales~\cite{hales2017reinhardt}]\label{cor:param-equiv}
The parameterization which makes $\det(X(t))$ a constant also makes $X + X^{-1}X' \in \sl(\R)$ and vice versa.
\end{corollary}
\begin{proof}
This is immediate from the identity 
\[
\frac{d}{dt}\det(X(t)) =\det(X) ~\tr\left(X + X^{-1}X'\right).
\]

\end{proof}
\subsection{The Cost Functional}
We now compute the cost functional in terms of the quantity $X$. From Problem \ref{pbm:reinhardt-reduction-2}, we see that the quantity to be minimized is the area of a disc in $\Kbal$. Our strategy is to compute this area using Green's theorem, which requires us to compute the pullback of the 2-form $xdy-ydx$ in $\R^2$. 

\begin{lemma}\label{lem:area-helper-1}
Let $g : [0,t_f] \to \SL(\R)$ be a path so that $g' = gX$ as above and let $v \in \R^2$. Define $\gamma : [0,t_f] \to \R^2$ by $\gamma(t) := g(t)v$. Let $\omega = xdy - ydx$ be a 2-form in $\R^2$, then we have the following formula for the pull-back:
\[
\gamma^*\omega = -v^t J X v dt
\]
where $J = \mattwo{0}{-1}{1}{0}$.
\end{lemma}
\begin{proof}
If we write $\gamma(t) = (\gamma_1(t),\gamma_2(t))$, then 
\begin{align*}
\gamma^*\omega &= \omega(\gamma(t))\\
&= (\gamma_1(t)\gamma_2'(t) -\gamma_2(t)\gamma_1'(t))dt\\
&= (\gamma \wedge \gamma')dt \\
&= (gv \wedge gXv) dt \\
&= (v \wedge Xv) dt \\
&=-v^t JX v dt
\end{align*}
since $\det(g) = 1$ and $u \wedge Xv = -u^tJXv$ for all $X \in \sl(\R)$ and $u,v \in \R^2$. 
\end{proof}
The lemma above enables us to compute pullbacks of $\omega$ by the multicurves $\sigma_j$. Indeed, if we fix an initial $v \in \R^2$ on the boundary $\partial K$ of an arbitrary disc $K \in\Kbal$ is a simple closed curve parameterized by the curves $\sigma$, given by $\sigma_j(t) := g(t)(R^{2j}v)$. Here $R$ is the rotation matrix giving a counterclockwise rotation by an angle of $\pi/3$.

\begin{lemma}\label{lem:area-helper-2}
Let $Y \in \sl(\R)$ be an arbitrary matrix. Then we have 
\[
JY + (R^2)^tJYR^2 + (R^4)^tJYR^4 = \frac{3\tr(JY)}{2}I_2
\]
where $J = \mattwo{0}{-1}{1}{0}$ is the infinitesimal generator of rotations.
\end{lemma}
\begin{proof}
This is a simple computation. 
\end{proof}

We now derive a formula for the area of $K$. 
By Green's theorem and the Lemmas proven above, we have:
\begin{align*}
    \mu(K) = \iint_K d\mu &= \frac{1}{2}\oint_{\partial K}\omega \\
    &= \frac{1}{2}\int_{0}^{t_f}\gamma^*\omega dt \\
    &= \int_{0}^{t_f}\sigma_0^*\omega +\sigma_2^*\omega +\sigma_4^*\omega dt \quad \text{(by Lemma \ref{lem:area-helper-1})}\\
    &= \int_{0}^{t_f}-v^t JX v -(R^2v)^t JX (R^2v)
    - (R^4v)^t JX (R^4 v) dt \\
    &= -\int_{0}^{t_f} v^t \left(JX + (R^2)^tJXR^2 + (R^4)^t JXR^4\right)v dt \\ 
    &=-\int_{0}^{t_f} |v|^2 \frac{3\tr(JX)}{2} dt \quad \text{(by Lemma \ref{lem:area-helper-2})}
\end{align*}
This disc $K$ is in circle representation, which means that the point $e_0^* = (1,0)$ is on the boundary (along with the other sixth roots of unity). This means that we can take $v = e_0^* = (1,0)$ and so:
\begin{equation}\label{eq:cost}
\mu(D) = -\frac{3}{2}\int_{0}^{t_f} \tr(JX) dt
\end{equation}
which is the quantity to be minimized.  

\subsection{Control Sets}\label{sec:control-sets}
We now investigate the \emph{control} parameters which affect this cost. It is intutively obvious that the curvature of the curves making up the boundary of the disc $K \in\Kbal$ affect its area. So, it makes sense to consider the curvatures to play the role of the controls. Problems in which curvature plays the role of a control are well-studied in the literature, the Dubins-Delauney problem~\cite{jurdjevic2014delauney} being one a prominent such example. Other examples, such as Kirchoff's problem and the elastic problem are discussed in \cite{jurdjevic2016optimal}.

To begin, recall that Lemma \ref{lem:local-curvature-constrt} says that $\sigma_j'(t) \wedge \sigma_j''(t) \ge 0$ almost everywhere in $t$.  Now since our disc is in circle representation, we have $\sigma_j(t) = g(t)e_{2j}^*$. Then we have, for $j \in \Z/6\Z$:
\begin{align}
v_j(X) := \sigma_j'(t)\wedge\sigma_j''(t) &= g'(t)e_{2j}^*\wedge g''(t) e_{2j}^* \nonumber \\
&= gXe_{2j}^*\wedge (gX^2+gX') e_{2j}^* \nonumber \\
&= Xe_{2j}^*\wedge (X^2+X') e_{2j}^* \nonumber \\
&=e_{2j}^*\wedge \left(X+X^{-1}X'\right) e_{2j}^* \ge 0 \label{eq:state-dep-curvature}
\end{align}
almost everywhere in $t$. Here we used Proposition~\ref{prop:reparam-lipschitz} and Corollaries \ref{cor:X-props} and \ref{cor:X-diff-ae}. We call $v_j(X)$ the \emph{state-dependent curvature} as it depends on where we are in the state space.
\begin{lemma}[Hales~\cite{hales2011}]\label{lem:curvature-pos}
Let $g : [0,t_f] \to \SL(\R)$ be a twice-differentiable curve. Then there exists a $j$ so that $v_j(X) > 0$.
\end{lemma}
\begin{proof}
Taking $X = \mattwo{a}{b}{c}{-a}$, assuming $v_0(X),v_2(X)=0$ a short calculation shows that:
\[v_1(X) = e_{4}^*\wedge \left(X+X^{-1}X'\right) e_{4}^* = \frac{-3(a^2+bc)}{\sqrt{3}a + c} = \frac{3}{\sqrt{3}a + c}\]
almost everywhere. This is strictly positive due to the star conditions in Corollary \ref{cor:X-props}.
\end{proof}

Since the state-dependent curvatures $v_j(X)$ depend on $X$ and $X'$ in general, they are not suitable as control variables for our control problem. To this end, we introduce normalizations of the state dependent curvatures as follows: 

\begin{definition}[Control variables]\label{def:control-variables}
For each $j=0,1,2$, define a normalization of the state-dependent curvatures as
\[
u_j(t) := \frac {v_j}{v_1+v_2+v_3}
\] 
We will denote $v = v(X) := v_1 + v_2 + v_3$ so that $v_j(X) = v(X) u_j$. 
\end{definition}

Note that $v > 0$ by Lemma \ref{lem:curvature-pos}. The control variables $u_i$ evidently satisfy $0\le u_i \le 1$ and $u_0 + u_1 + u_2 = 1$. Note that the control variables are functions of time.

\begin{definition}[Triangular control set]\label{def:triang-control-set}
The \emph{triangular} or \emph{simplex} control set is the set \[U_T := \left\{(u_0,u_1,u_2)~|~ 0\le u_i \le 1, ~\sum_i u_i = 1\right\}\]
which is just the 2-simplex in $\R^3$.
\end{definition}

Later in the thesis, we will also be interested in an enlargement of the triangular control set to its circumscribing disk. 

\begin{definition}[Circular control sets]\label{def:circ-control-set}\leavevmode
\begin{itemize}
    \item The \emph{circumscribed} or \emph{disk} control set is the set $U_C$ which is the circumscribing disk of the 2-simplex in $\R^3$:
    \[
    U_C := \left\{(u_0,u_1,u_2)~|~ 0\le u_i \le 1, ~\sum_i u_i = 1,~\sum_i u_i^2 \le 1\right\}
    \]
    \item The \emph{inscribed} control set is the set $U_I$ which is the inscribed disk of the 2-simplex in $\R^3$:
     \[
    U_I := \left\{(u_0,u_1,u_2)~|~ 0\le u_i \le 1, ~\sum_i u_i = 1,~\sum_i u_i^2 \le \frac{1}{2}\right\}
    \]
    \item Later (in Section \ref{sec:su11-control-sets}), we shall also be interested in control sets which interpolate between $U_I$ and $U_C$. For $1/3 \le r^2$ we define:
    \[
    U_r := \left\{(u_0,u_1,u_2)~|~u_0 + u_1 + u_2 = 1,~u_0^2+ u_1^2 + u_2^2 \le r^2 \right\}
    \]
    We shall denote the boundary of these sets by $\partial U_C$, $\partial U_I$ and $\partial U_r$ respectively\footnote{Here, we do not mean the boundary as a subset of $\R^3$, but rather as a subset of the plane $\{(u_0,u_1,u_2)~|~\sum_{i=0}^3 u_i = 1\}$}. 
\end{itemize}

\end{definition}
\begin{remark}\normalfont
Controls $(u_0,u_1,u_2)$ from $U_C \setminus U_T$ give rise to boundary curves $\sigma_j$ of discs which fail to be convex. This is a consequence of the star conditions. 
\end{remark}

Our primary motivation in considering these control sets is the following: we can observe that the triangular control set is invariant under $\Z/3\Z$-rotations while the disk control set is invariant under rotations by the circle group $S^1$. The latter is important for our investigations, as it will allow us to employ Noether's theorem to derive a \textit{first integral} or a conserved quantity of the dynamics. A conserved quantity is useful because it facilitates a reduction in dimension of the original problem --- if the group of symmetries is large enough, reduction by that group may even lead to a direct solution.

However, our control sets are in $\R^3$, which is not the state space of our control problem. As equation \eqref{eq:state-dep-curvature} shows, the dynamics for $X$ seems to be related to the curvatures. Thus, our first task is to map these control sets into the Lie algebra $\sl(\R)$ and to derive a precise relation between the \emph{control matrices} and the evolution of the matrix $X$. 

To this end, we map the control set $U_T$ into the Lie algebra $\sl(\R)$. The following transformation
\begin{equation}\label{eq:Z0}
Z_u(u_0,u_1,u_2) = \left(\begin{matrix}
    \frac{u_1 - u_2}{\sqrt{3}} & \frac{u_0 - 2u_1 - 2u_2}{3} \\
    u_0 & \frac{u_2 - u_1}{\sqrt{3}}
    \end{matrix}\right)    \in \sl(\R)
\end{equation}
maps $U_T$ into $\sl(\R)$. Note that the control matrix $Z_u$ has been chosen to make
\begin{equation}\label{eq:def-Z0}
u_j = e_{2j}^* \wedge Z_u e_{2j}^*  \qquad j=0,1,2
\end{equation}

Thus, our control problem consists of finding a time-dependent control function constrained to the set $U_T$, which assigns a \textit{proportion} of curvature to the six curves (which, up to central symmetry are actually three curves) making up the boundary of a disc $K \in\Kbal$ at each time instant so that its area is minimized. 

\subsection{Lie Algebra Dynamics}\label{sec:lie-algebra-dynamics}

Let us now derive the control-dependent dynamics for $X$. Henceforth, we shall adopt the notation $\bracks{A}{B} := \tr(AB)$ for any two matrices $A,B \in \sl(\R)$. This form is a nondegenerate invariant bilinear form on $\sl(\R)$. Let us first collect a number of results about matrices in $\sl(\R)$ which we shall need. All of these are elementary and so we admit them without proofs.

\begin{lemma}\label{lem:sl2-lemmas}
We have the following results about matrices in $\sl(\R)$:
\begin{enumerate}
\item If $X \in \sl(\R)$, $\bracks{X}{X} = -2\det(X)$.
\item If $A,B$ are any two matrices in $\sl(\R)$, and $E$ is a multi-point
\[
E_{j} \wedge A E_{j} = E_{j} \wedge B E_{j},\quad j=0,1,2
\]
then $A=B$.
\item If $A,B \in \sl(\R)$ then $AB + BA = \bracks{A}{B} I_2$.
\item For any matrices $A,B,C \in \sl(\R)$ we have $\bracks{A}{[B,C]} = \bracks{[A,B]}{C}$.
\end{enumerate}
\end{lemma}

\begin{theorem}[Dynamics for $X$]\label{thm:X-dynamics-sl2}
The dynamics for $X$ (which is control-dependent) is given by:
\[
X' = \frac{\left[Z_u,X\right]}{\bracks{Z_u}{X}}
\]
where $[A,B] = AB - BA$ is the matrix commutator of any two matrices $A,B \in \sl(\R)$. 
\end{theorem}
\begin{proof}
By Corollary \ref{cor:param-equiv}, we get that $X + X^{-1}X' \in \sl(\R)$. 
From equations \eqref{eq:state-dep-curvature} and \eqref{eq:def-Z0} we get that 
\begin{align*}
u_j &=  e_{2j}^* \wedge Z_u e_{2j}^* \\
v_j &=  e_{2j}^* \wedge \left(X+X^{-1}X'\right)e_{2j}^* 
\end{align*}
for each $j$. Since, by Definition \ref{def:control-variables}, we have $v u_j = v_j$, by Lemma \ref{lem:sl2-lemmas},(2) we get that:
\[
X+X^{-1}X' = v Z_u
\]
from which we get 
\begin{equation}\label{eq:X-orig}
X' = X(v Z_u - X)
\end{equation}
Since $\tr(X') = 0$, we get $v = \bracks{X}{X}/\bracks{X}{Z_u}$.
Now we have:
 \begin{align}
        X' &= v X Z_u - X^2 \nonumber  \\
        &= \frac{\langle X,X \rangle}{\langle X,Z_u \rangle} X Z_u - \frac{\langle X,X\rangle}{2} I \qquad \mathrm{from~Lemma ~\ref{lem:sl2-lemmas},(3)} \nonumber\\ 
        &=\frac{\langle X,X \rangle}{\langle X,Z_u \rangle} X Z_u - \frac{\langle X,X\rangle}{2} \left( \frac{X Z_u + Z_u X }{\langle X,Z_u \rangle} \right) \qquad \mathrm{from~Lemma~\ref{lem:sl2-lemmas},(3)}\nonumber \\
        &= \frac{\langle X,X \rangle}{2\langle X,Z_u \rangle} \left( XZ_u - Z_u X\right) \nonumber \\ 
        &= \frac{-\langle X,X \rangle}{2\langle X,Z_u \rangle}\left[Z_u, X\right] = \frac{\left[Z_u,X\right]}{\bracks{Z_u}{X}}\quad\mathrm{by~ Lemma~\ref{lem:sl2-lemmas},(1)} \nonumber 
    \end{align}
    which proves the required. Note that we also have $v = -2/\bracks{X}{Z_u}$ as an application of Lemma \ref{lem:sl2-lemmas},(1).
\end{proof}

\begin{corollary}\label{cor:X-det-const}
The quantity $\bracks{X}{Z_u}$ for a fixed control matrix $Z_u$ is a constant of motion along $X$. 
\end{corollary}
\begin{proof}
The quantity in question is constant since:
\begin{equation}\label{eq:X-Z0-constant}
\langle X, Z_u \rangle' = \langle X', Z_u \rangle = \frac{\langle [Z_u, X], Z_u \rangle}{\langle X,Z_u \rangle} = 0
\end{equation}
where we have used the fact that $\langle [A,B],C \rangle = \langle A,[B,C] \rangle$.
\end{proof}

\begin{remark}\label{rmk:X-det-one}
\normalfont \leavevmode
\begin{enumerate}
    \item The equation $A'(t) = [B(t),A(t)]$ where $A,B$ are time-dependent matrices is called the \emph{Lax equation} and $A,B$ so related are called a \emph{Lax pair}. Lax pairs are well-studied in the theory of integrable systems (see Perelomov~\cite{perelomov1990integrable}, Jurdjevic~\cite{jurdjevic2016optimal}, Babelon et al.~\cite{babelon2003introduction}). Lax representations of integrable systems are quite desirable since the evolution of a Lax pair is \emph{isospectral}, meaning that the spectrum of the matrix $A(t)$ is an invariant of motion.
    \item The dynamics for $X$ is Hamiltonian for a particular Hamiltonian defined on $\sl(\R)$, with respect to a Poisson structure on $\sl(\R)$ called the \emph{Lie-Poisson structure}. See Appendix \ref{sec:X-lie-poisson} for more details. 
    \item As explained in Perelomov~\cite[p.~52]{perelomov1990integrable}, the spectral invariants guaranteed by the dynamics for $X$ are ``trivial'' integrals, and so it is more accurate to consider the dynamics for $X$ as giving a control-dependent infinitesimal generator for the (co)adjoint action of $\SL(\R)$ on $\sl(\R)$, rather than to regard it as describing the dynamics of an integrable system. 

    \item The equation \eqref{eq:X-orig} appears in Hales~\cite{hales2017reinhardt}, where its Lax pair reformulation was not explicitly recognized.
\end{enumerate}
\end{remark}



\subsection{Initial and Terminal Conditions}\label{sec:X-init-term-conds}
We now have dynamics for $g$ and $X$ in the Lie group and Lie algebra respectively. We also have an associated cost objective. The only thing remaining is to specify initial and terminal conditions. 
Since our disc is in circle representation, this means that $\sigma_j(0) = e_{2j}^*$ so that we start out at the sixth roots of unity, and so we set $g(0) = I_2$. 
The initial condition for $X(0) =  X_0$ may be an arbitrary matrix in $\sl(\R)$ provided it satisfies the star conditions in Corollary \ref{cor:X-props}. 

The terminal conditions $g(t_f)$ should be such that the curves $\sigma_j$ close up seamlessly to form a simple closed curve:
\begin{equation}
g(t_f)e_{j}^* = g(0)e_{j+1}^* \Leftrightarrow g(t_f) = R
\end{equation}
where is the rotation matrix $R = \exp(J\pi/3)$. For terminal conditions on $X$, note that we have the following conditions on $g$ which we get by the remark following Theorem \ref{thm:defn-g} (and setting $t_0=0$ there):
\[
g(t+t_f)e_j^* = g(t)e_{j+1}^* = g(t)Re_j^*, \qquad g(t+t_f)=g(t)R
\]
Differentiating this, we get:
\[
X(t+t_f) = R^{-1}X(t)R
\]
which gives us $X(t_f) = R^{-1}X_0R$. 

\section{Reinhardt Optimal Control Problem}\label{sec:ROC}
Summarizing the discussion so far, we are finally ready to state the Reinhardt conjecture as an optimal control problem. Let us begin with a proposition:

\begin{proposition}[Jurdjevic~\cite{jurdjevic_1996}]\label{prop:trivialization}
Let $G$ be any real Lie group and let $\mathfrak{g}$ be its Lie algebra. Then we have $T^*G \cong G \times \mathfrak{g}^*$ and $TG \cong G \times \mathfrak{g}$.
\end{proposition}
\begin{proof}
For each $g \in G$ let $L_g(h) := gh$ be the left-multiplication map. We then have that the tangent map $TL_g : T_x G \to T_{xh}G$ for each $x \in G$. Let $TL_g^*$ denote the dual mapping of $TL_g$. It follows that $T^*L_g(p)$ evaluated at a vector $X_x \in T_xG$ is equal to $p\circ TL_g(X_x)$. So, in particular, at $x=g$, $T^*L_{g^{-1}}$ maps $T_e^*G$ onto $T_g^*G$. The correspondence $ T^*L_{g^{-1}}(p) \leftrightarrow (g,p)$ realizes $T^*G \cong G \times \mathfrak{g}^*$. In a similar way, we also have a correspondence $TL_g(q) \leftrightarrow (g,q)$ giving the trivialization $TG \cong G \times \mathfrak{g}$.  
\end{proof}

The above proposition and remark apply to $\SL(\R)$. Using this, we can state our control problem as:

\begin{problem}[Control Problem (ROC)]\label{pbm:reinhardt-optimal-control}
The discs in $\Kbal \subset \Kccs$ in circle representation arise via the following optimal control problem:
On the manifold $\SL(\R) \times \sl(\R) \cong T\SL(\R)$, consider the following free-terminal time optimal control problem: 
\begin{align}
g' &= gX \qquad \qquad \qquad \qquad g:[0,t_f] \to \SL(\R) \label{eq:trisystem-g}\\
X' &=  \frac{1}{\langle X,Z_u \rangle}\left[Z_u, X\right]  \qquad~ X:[0,t_f] \to \sl(\R) \label{eq:trisystem-X} \\
-\frac{3}{2}\int_0^{t_f}&\langle J,X\rangle dt \to \mathrm{min} \qquad \quad \quad J = \left(\begin{matrix} 
 0 & -1 \\ 1 & 0 \end{matrix}\right) \label{eq:trisystem-cost}
\end{align}
where the set of controls for this problem is the 2-simplex in $\R^3$ mapped inside the Lie algebra $\sl(\R)$ via the affine map $Z_u$:
\begin{align}
    Z_u : U_T &= \left\{(u_0,u_1,u_2) \ | \ \sum_i u_i = 1,  \ u_i \ge 0 \right\} \to \sl(\R) \label{eq:trisystem-U} \\
    Z_u(u_0,u_1,u_2) &= \left(\begin{matrix}
    \frac{u_1 - u_2}{\sqrt{3}} & \frac{u_0 - 2u_1 - 2u_2}{3} \\
    u_0 & \frac{u_2 - u_1}{\sqrt{3}}
    \end{matrix}\right)   \label{eq:trisystem-Z0}
\end{align}
With initial conditions: $g(0) = I_2 \in \SL(\R)$ and $X(0) = X_0 \in \sl(\R)$ satisfying the star conditions. Also, with terminal conditions: $g(t_f) = I_2$ and $X(t_f)=RX_0R^{-1}$ where $R = \exp(J\pi/3)$. 
\end{problem}

\begin{definition}[Jurdjevic~\cite{jurdjevic_1996}]\label{def:left-invariant-control-system}
An arbitrary optimal problem with control system $dg/dt = F(g,u)$ defined on a real Lie group $G$ with control functions $u(t) \in U \subseteq \R^m$ is said to be \emph{left invariant} if $TL_g F(h,u) = F(L_g h, u)$ for each $g,h \in G$. Here $L_g(h) = gh$ is the left-multiplcation map and $TL_g:T_hG \to T_{gh}G$ is its tangent map.

Also, we require the associated cost function $f(g,u)$ to be left-invariant: $f(L_g h,u) = f(g,u)$ for all $g,h$ and $u$. 
\end{definition}

The dynamics of the Reinhardt optimal problem at the group level is clearly left-invariant. The reason for making this observation is that left-invariance of a dynamical system on a Lie group means that we can reduce its dynamics to co-adjoint orbits of the associated Lie algebra.

Since the cost and dynamics for $X$ are independent of $g$, we note that the control problem only depends on the group element $g$ for endpoint (transversality) conditions. This means that we can focus on the control problem exclusively at the Lie algebra level, and postpone the endpoint condition on $g$ until the very last step. In later sections, we shall drop the group dynamics and exclusively focus on state/costate dynamics in the Lie algebra. 

\subsection{The Upper Half-Plane.}
Now that we have the optimal control problem fully stated, a natural next step would be to write down the necessary conditions for optimality of trajectories. But before we do that, we shall first try to cut down the state space of the problem.

Recall that the star conditions (Corollary \ref{cor:X-props}) on the matrix $X$ imply that $c_{21}(X)$ (the entry on the second row and first column) is positive. We also just imposed the condition $\det(X)=1$. We begin with a characterization of such matrices.

\begin{lemma}[Hales~\cite{hales2017reinhardt}]\label{lem:X-hyperbolic}
The set of matrices $X \in \sl(\R)$ with $\det(X(t))=1$ and $c_{21}(X) > 0$ is the adjoint orbit $\O_J := \{\Ad_g J ~|~ g\in \SL(\R) \}$ in $\sl(\R)$ of the element $J = \mattwo{0}{-1}{1}{0}$. 
\end{lemma}
\begin{proof}
Let $\mathcal{X}$ be the set of all such matrices. 
The adjoint orbit of $J$, $\O_J$ in $\sl(\R)$ consists of elements that look like $gJg^{-1}$ for $g \in \SL(\R)$. The Iwasawa decomposition of $\SL(\R)$ states that $g = nak$ where 
\begin{align*}
     k \in K &= \left\{\mattwo {\cos \theta} {-\sin \theta} {\sin \theta} {\cos \theta} \ \Bigr| \ \theta \in [0,2\pi)\right\} \cong \SO \\
     a \in A &= \left\{ \mattwo r 0 0 {1/r} \Bigr| \ r > 0 \right\} \\
     n \in N &= \left\{ \mattwo 1 p 0 1 \Bigr| \ p \in \R \right\}.
\end{align*}
Since we have $kJk^{-1}=J$, the element looks like \[naJa^{-1}n^{-1}=\left(
     \begin{array}{cc}
      -p r^2 & -\frac{p^2 r^4+1}{r^2} \\
      r^2 & p r^2 \\
     \end{array}
     \right).\]
Clearly, this matrix satisfies the hypotheses and so $\O_J \subseteq \mathcal{X}$.
A generic $X$ satisfying the hypotheses looks like \begin{equation}\label{X-upper-defn}
X = X_z := \mattwo{x/y}{-x^2/y-y}{1/y}{-x/y} \in \O_J
\end{equation}
for some $x$ and $y>0$. This matrix is in the form described above for $p = -1/x$ and $r=1/\sqrt{y}$. So $X \subseteq \O_J$.
\end{proof}
\begin{remark}\normalfont
Since $\sl(\R)$ admits a nondegenerate invariant inner product, its coadjoint orbits and adjoint orbits coincide. See Appendix \ref{sec:kirillov} and also Chapter 5 of Jurdjevic~\cite{jurdjevic_1996}.
\end{remark}

\subsubsection{Transfer of Dynamics to the Upper Half-Plane}
By the Orbit-Stabilizer theorem, we have $\O_J \cong \SL(\R)/\SO$, since the stabilizer of $J$ in $\SL(\R)$ under the conjugation action (the centralizer) is $\SO$. Viewing this in a different way, the group $\SL(\R)$ acts on the upper half-plane $\h = \left\{x+iy~|~y>0\right\}$ by fractional linear transformations, with stabilizer of $i \in \h$ being given by $\SO$. Thus, the quotient is isomorphic to the the Poincar\'{e} upper half-plane $\h$. Putting all of this together, we have

\begin{lemma}\label{lem:def-phi}
The following map $\Phi$ is a isomorphism:
\begin{align*}
     \Phi : \mathfrak{h} &\to \mathcal{O}_J \\
     z = x+iy &\mapsto \mattwo {x/y} {-(x^2 + y^2)/y} {1/y} {-x/y} =: X_z
\end{align*} 
\end{lemma}
\begin{remark}\leavevmode\normalfont
\begin{itemize}
    \item Note that $\mathcal{O}_J = \mathcal{O}_{X_z}$ as $X_z \in \mathcal{O}_J$.
    \item  We shall write $X$ in place of $X_z$ for simplicity, bearing in mind that $\Phi$ is surjective.
    \item  Note that $X$ is a regular semisimple element of the Lie algebra $\sl(\R)$ as the element $J$ is. 
\end{itemize}
\end{remark}

This map $\Phi$ allows us to move back and forth between the upper half-planes and the adjoint orbit in the Lie algebra $\sl(\R)$. Also, the map $\Phi$ is more than just a bijection --- we show later that this map is actually an  anti-symplectomorphism and use this to transfer the state and costate dynamics from the Lie algebra to the upper half-plane. But first, we compute the tangent map $T \Phi$ at a $z \in \h$.

\begin{lemma}\label{lem:half-plane-lie-algebra-iso}
For any $X \in \sl(\R)$, we have 
 \[ T_X \OX = \{[Y,X] \ | \ Y \in \sl(\R) \}\cong \sl(\R)/\sl(\R)_X\]
where $\sl(\R)_X$ is the coadjoint isotropy algebra of the element $X$ (the centralizer of $X$) and $T_X \OX$ denotes the tangent space to $\OX$ at $X$. 
\end{lemma}
\begin{proof}
Note that $\OX = \{\Ad_g X \ | \ g \in \SL(\R) \}$. We have to describe tangent vectors to $\OX$. For any $Y \in \sl(\R)$, $\Ad_{\exp(tY)} X$ is a curve in $\OX$. Thus, the tangent vectors to this curve are computed as:
\[ 
\frac{d}{dt} \Ad_{\exp(tY)}X\Bigr|_{t=0} = \ad_Y X := [Y,X] \in T_X\OX
\]
This calculation is actually finding infinitesimal generator of the adjoint action. Next, note that the subspace $\{[Y,X] \ | \ Y \in \sl(\R) \}$ has dimension either 1, 2 or 3. It cannot be 1 dimensional since the vectors $[J, X]$ and $[J_1,X]$ are linearly independent (where $J_1 = \mattwo{0}{1}{1}{0}$). It cannot be three dimensional since that would mean that the map $\ad_X$ is onto and so, in particular, there would be a matrix $Y$ such that $\ad_X(Y) = X$, but such a $Y$ cannot exist, as a calculation shows. Thus, the subspace is 2-dimensional. 

The centralizer of $X$ lies in the kernel of $\ad_X$ and has dimension 1 since $X$ is a regular element and so the quotient space $\sl(\R)/\sl(\R)_X$ is 2-dimensional. This gives the required isomorphism. 
\end{proof}

\begin{lemma}\label{lem:tangent-maps}
We have the following expression for the tangent map $T \Phi$ and its inverse $(T \Phi)^{-1}$:
\begin{align}
     T \Phi : T_z \h &\to T_{X}\OX \cong \sl(\R)/\sl(\R)_{X} \\
     (k_1 \delx + k_2 \dely) &\mapsto \mattwo {k_2 /2y}{(yk_1 - k_2 x)/y}{0}{-k_2/2y} \mod \sl(\R)_{X}
\end{align}
\begin{align}
     (T \Phi)^{-1} : T_{X}\OX &\to T_z\h \\
     \left[\mattwo{a}{b}{c}{-a}\right] &\mapsto 
     \left(b + 2a x - c x^2 + cy^2 \right)\delx + (2 y(a - c x))\dely
\end{align}
\end{lemma}
\begin{proof}
We have, at $z = x+iy$:
\begin{align}
     T \Phi (k_1,k_2) &= \frac{d}{dt}{X}(x + t k_1, y + t k_2)\Bigr|_{t=0} \\
     &= k_1 \frac{\partial {X}}{\partial x} + k_2 \frac{\partial {X}}{\partial y} \in T_{X}\OX
\end{align}
We know by the previous lemma that there exists a matrix $Y_z$ such that \[T \Phi(k_1,k_2) = k_1 \frac{\partial {X}}{\partial x} + k_2 \frac{\partial {X}}{\partial y} = [Y_z(k_1,k_2),{X}]\]  Using this equation and solving for this matrix $Y_z$ gives us the following:
\begin{align}
Y_z &= \underbrace{Y_{12}\mattwo {x} {-x^2 -y^2} {1} {-x}}_{\in \sl(\R)_{X}}+ \mattwo {k_2 /2y}{(yk_1 - k_2 x)/y}{0}{-k_2/2y} \\
&\equiv \mattwo {k_2 /2y}{(yk_1 - k_2 x)/y}{0}{-k_2/2y} \mod \sl(\R)_{X}
\end{align}
So, for any arbitrary vector $(k_1,k_2) \in T_z \h$ we get a canonical representative inside the quotient space $\sl(\R)/\sl(\R)_{X}$. 

The map $T \Phi$ is a linear isomorphism. The inverse is the following:
\begin{equation}\label{eq:inv-tangent-map}
(T \Phi)^{-1} \left[\mattwo{a}{b}{c}{-a}\right]= \left(
    \left(b + 2a x - c x^2 + cy^2 \right)\delx + (2 y(a - c x))\dely
          \right)
          \in T_z\h
\end{equation}
where $ \left[\mattwo{a}{b}{c}{d}\right]\in \sl(\R)/\sl(\R)_{X} $. This can be verified directly. 
\end{proof}

We can now use these lemmas to compute the dynamics of the control problem in the upper half-plane\footnote{Although there is a more direct way of computing the dynamics of $x$ and $y$, this approach will be later extended to account for the costate dynamics of the control problem.}.
\begin{theorem}\label{thm:X-dynamics}
The ODE for $X$ viz.,
\[
     X' = \frac{[Z_u,X]}{\langle Z_u,X\rangle}
\]
transforms as
\begin{align}
x' &= f_1(x,y;u) := \frac{y \left(2 a x+b-c x^2+c y^2\right)}{2 a x+b-c x^2-c y^2} \\
y' &=  f_1(x,y;u):= \frac{2 y^2 (a-c x)}{2 a x+b-c x^2-c y^2}     
\end{align}
in upper half-plane coordinates. Here $Z_u = \mattwo{a}{b}{c}{-a}$.
\end{theorem}
\begin{proof}
We need to compute the tangent map $(T \Phi)^{-1}$ at the element $\left[\frac{Z_u}{\langle Z_u,X \rangle}\right] \in \sl(\R)/\sl(\R)_X$. 
To do this, we use equation \eqref{eq:inv-tangent-map}:
\[
(T \Phi)^{-1}\left[ \frac{Z_u}{\langle Z_u,X\rangle}\right] = \frac{1}{2 a x+b-c x^2-c y^2}\left(
     \begin{array}{c}
     b + 2a x - c x^2 + cy^2 \\
     2 y(a - c x)\\  
     \end{array}
     \right)\in T_z\h
\]
which gives us the required. 

This also shows that
\begin{equation}\label{eq:tangent-map-f1f2}
     T\Phi(f_1,f_2) = \left[\frac{Z_u}{\langle Z_u,X \rangle}\right] \in T_X\OX.
\end{equation}
\end{proof}

Thus, we have transferred the Lie algebra dynamics to the upper half-plane. More specifically, as we have mentioned before, this result is indicative of the fact that the coadjoint orbit through $J$ in $\sl(\R)^*$ and the upper half-plane are symplectomorphic. Thus, it is entirely equivalent to study the control problem in the Lie algebra picture or the half-plane picture. We may also transfer the dynamics to other models of hyperbolic geometry: for example, the Poincar\'{e} disk model or the hyperboloid model. Each picture has its advantages, with some simplifying equations while others are better since the symmetries are more apparent. 

For now, let us turn to deriving the cost functional in half-plane coordinates.

\subsubsection{The Cost Functional in Half-Plane Coordinates}
The cost of the control problem in equation \eqref{eq:cost} has the following expression in upper half-plane coordinates:

\begin{align}
    \frac{-3}{2}\int_0^{t_f}\bracks{J}{X}dt &=    \frac{-3}{2}\int_0^{t_f}\tr\left(\mattwo{0}{-1}{1}{0}\mattwo {x/y} {-(x^2 + y^2)/y} {1/y} {-x/y}\right)dt \nonumber \\
    &= \frac{3}{2}\int_0^{t_f}\frac{x^2+y^2+1}{y}dt \rightarrow \text{min} \label{eq:cost-upper-half-plane}
\end{align}

This reinterpretation makes the circular symmetry of the cost function apparent: the level sets of $(x^2+y^2+1)/y$ are concentric circles centered at the point $i$ in the upper half-plane. Thus, the cost is $\SO$-invariant.
\begin{remark}\normalfont
The Poincar\'{e} upper half-plane is conformally equivalent to other models of hyperbolic geometry such as the Poincar\'{e} disk and the hyperboloid model. The cost functional derived above can also be derived in these models. In the disk model, $\mathbb{D}=\{w \in \mathbb{C}~|~|w| < 1\}$, for example, the cost of a path $w : [0,t_f] \to \mathbb{D}$ becomes:
\[
3\int_0^{t_f}\frac{1+|w|^2}{1-|w|^2}dt \rightarrow \text{min}
\]
In the hyperboloid model, the cost functional becomes the integral of the height function. See Hales~\cite{hales2017reinhardt}.
\end{remark}
\subsubsection{The Star Domain in the Upper Half-Plane}
\begin{figure}[htbp]
    \centering
    \includegraphics{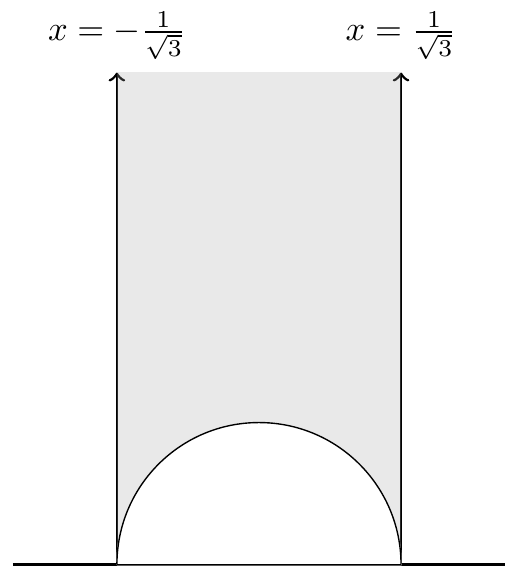}
    \caption{The star domain in the upper half-plane.}
    \label{fig:star-dom}
\end{figure}

We can now prove our first state space reduction result:

\begin{theorem}[\cite{hales2017reinhardt}]\label{prop:star-domain}
The dynamics of the Reinhardt control problem is constrained to an ideal triangle $\star$ in the upper half-plane. Thus, the new state-space of the control problem is $\SL(\R)\times \star$. 
\end{theorem} 
\begin{proof}
The star conditions on $X$ in Corollary \ref{cor:X-props} applied to $X = X_z$ give us the conditions on $x$ and $y$ which imply that the admissible trajectories should be constrained to be in:
\[
\star := \left\{x+iy \in \h~|~-\frac{1}{\sqrt{3}}<x<\frac{1}{\sqrt{3}},~\frac{1}{3}<x^2 + y^2\right\}
\]
which defines the interior of an ideal triangle in the upper half-plane. The vertices of this triangle are the lines $z = \pm\frac{1}{\sqrt{3}}$ and $z = \infty$. A picture of the star domain is show in Figure \ref{fig:star-dom}.
\end{proof}

In fact, we can now prove an equivalent star condition.
\begin{corollary}\label{cor:sl2-star-condition}
The star conditions on $X$ hold if and only if $\bracks{Z_u(u)}{X} < 0$ for all controls $u \in U_T$.
\end{corollary}
\begin{proof}
Note that $Z_u$ is an affine map and $U_T$ is a convex set. Thus, $Z_u(U_T)$ is also convex. Thus, it is enough to check that $\bracks{Z_u(u)}{X}<0$ for $u$ at the vertices of $U_T$. But these are precisely the star conditions. The converse is true simply by specializing $u$ to be the vertices of $U_T$. 
\end{proof}

Summarizing the results so far, we have parameterized the boundary of discs in $\Kbal \subset \Kccs$ as $U_T$-controlled paths $(g(t),z(t)) \in \SL(\R)\times \star$ subject to the terminal conditions. Our task is to find a \emph{control law} $u(t) \in U_T$ which minimizes the area enclosed by the resulting curve, given in hyperbolic coordinates by equation \eqref{eq:cost-upper-half-plane}.

\subsubsection{Compactification of the Star Domain}\label{sec:compactification}
While the star domain is a reduction of the state space, it is an ideal triangle with one vertex at infinity. Thus, it is open and unbounded in the upper half-plane. Our task in this section will be to explore a further reduction of this admissible region.

Empirical observations show that if $z \in \star$ is close to the boundary lines, then the corresponding critical hexagon (constructed in Lemma \ref{lem:pt-to-hexagon}) is close to a parallelogram. Classical results of Mahler and Reinhardt state that the only convex discs in $\Kccs$ with a parallelogram for a critical hexagon are parallelograms themselves. This suggests that there is a neighbourhood of the boundary of the star domain which gives rise to convex discs in $\Kccs$ whose packing density is close to being $1$ and so these discs can be excluded from consideration, since they are never optimal for our control problem. We make this intuition precise presently.

Our hope is that if can cut down enough of the state space, we can constrain ourselves to a compact region inside it so that eventually a computer simulation of the dynamics would become feasible. 

\begin{lemma}\label{lem:pt-to-hexagon}
At any time instant $t \in [0,t_f]$, consider $z(t) = z \in \star$ which is induced by the curve parameterizing the boundary of a disc $K$ in $\Kbal$. Then $z$ gives rise to a critical hexagon for $K$. We denote this as usual by $H_K(z)$.
\end{lemma}
\begin{proof}
The element $z$ in the upper half-plane determines a matrix $X$ in the (coadjoint orbit of $J$ in the) Lie algebra $\mathfrak{sl}_2$ via the isomorphism $\Phi$. The element
$X$ determines the six tangent directions $Xe_{2j}^*$ at the six points $\{e_j^*\}$ on the boundary of the disc $K$. From this, we take the six tangent lines, written as $e_j^* + s Xe_j^*$, as $s$ runs over the real numbers. By construction, the points of intersection of these six lines (which we denote $P_0,P_1,P_2,P_0',P_1',P_2'$) are the vertices of a critical hexagon for $K$. Here, we write $P_j'$ for $-P_j$ for $j=0,1,2$.
\end{proof}

\begin{figure}[htbp]
    \centering
    \includegraphics{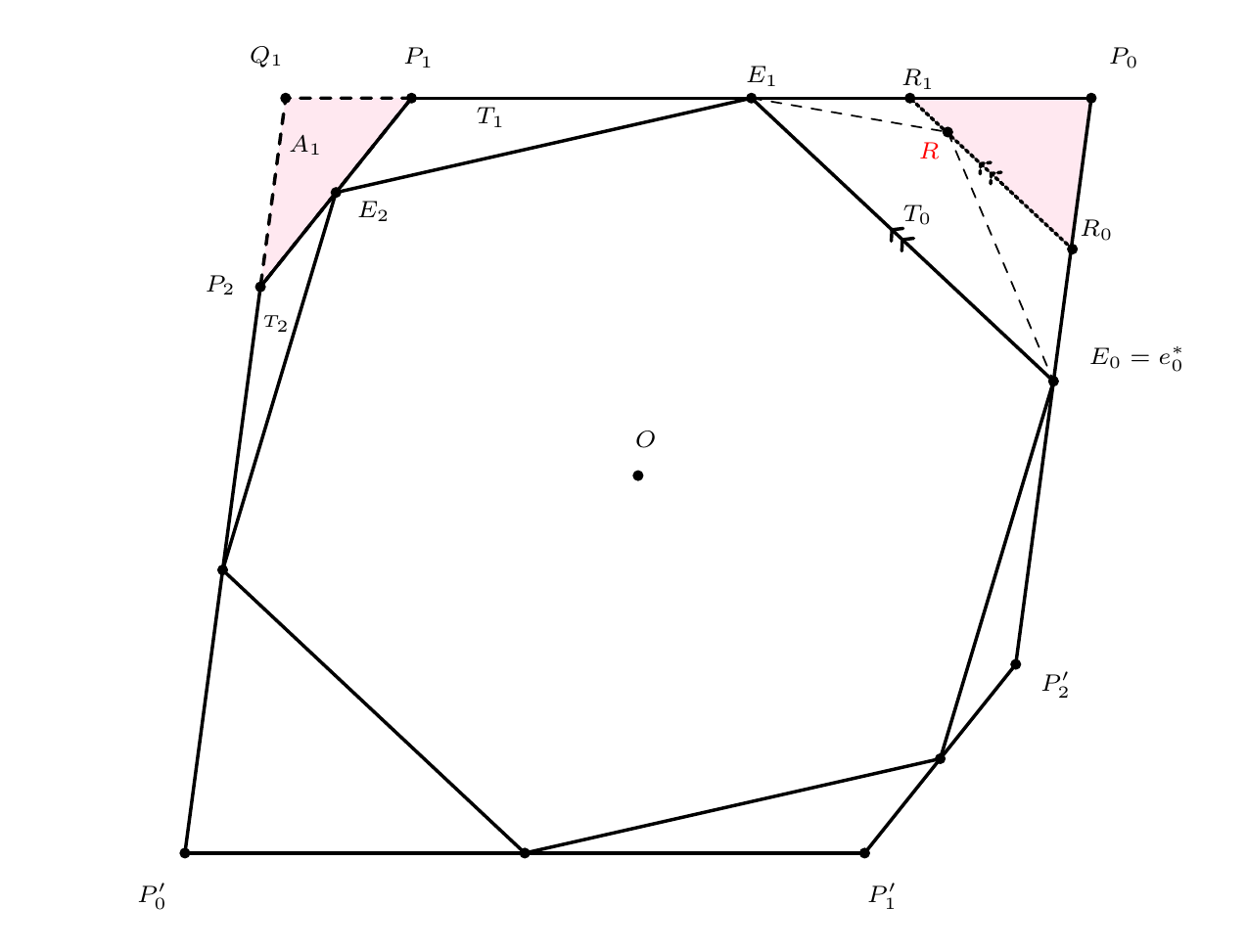}
    \caption{Critical Hexagon at a particular time instant.}
    \label{fig:compactification-result}
\end{figure}

The situation is depicted in Figure \ref{fig:compactification-result}. Recall that our disc is in circle representation, having the sixth roots of unity on its boundary. For notational convenience, we relabel $e_j^*=E_j$ for $j=0,1,2$ and add a prime subscript for the image of a point in the origin $O$, instead of the minus sign. 

We label triangles $T_0 = \triangle E_0 P_0 E_1,~ T_1 = \triangle E_1 P_1 E_2,~ T_2 = \triangle E_2 P_2 E_0'$. Let $Q_1$ denote the point of intersection of the lines ${P_0' P_2}$ and ${P_0 P_1}$. Similarly for $Q_0$ and $Q_2$. This now determines the triangles $A_1 = \triangle P_1 Q_1 P_2$, $A_2 = \triangle P_2 Q_2 P_0'$, and $A_0 = \triangle P_0 Q_0 P_1$ (the latter two triangles are not depicted in the figure). For our compactification result, we will need the areas of these triangles in terms of $z = x + iy$. 

\begin{lemma}\label{lem:triangle-areas}
Fix any point $z = x+iy \in \mathfrak{h}^{\ast}$. Let $\alpha(z) := \frac{1 + \sqrt{3}x}{y}$, $\beta(z) := \frac{1 - \sqrt{3}x}{y}$ and $\gamma(z) := \frac{3x^2 + 3y^2 - 1}{2y}$. Then we have: 

\[\mu(T_0) = \frac{\alpha \gamma}{4\sqrt{3}},
\qquad 
\mu(T_1) = \frac{\alpha \beta}{4\sqrt{3}},
\qquad
\mu(T_2) =  \frac{\beta \gamma}{4\sqrt{3}}
\]
and 
\[\mu(A_0) = \frac{\alpha^2}{\sqrt{3}},
\qquad 
\mu(A_1) = \frac{\beta^2}{\sqrt{3}},
\qquad
\mu(A_2) =  \frac{\gamma^2}{\sqrt{3}}
\]
where $\mu(\cdot)$ is the Lebesgue measure in $\R^2$.  Furthermore, we have $\mu(T_0) + \mu(T_1) +\mu(T_2) = \frac{\sqrt{3}}{4}$.  
\end{lemma}
\begin{proof}
The recipe in Lemma \ref{lem:pt-to-hexagon} gives us a method to start with the points $E_i$ and end up with points $P_i$ and $Q_i$ for $i=0,1,2$. Once we have coordinates of all these points, finding the areas of the associated triangles is straightforward.
\end{proof}
For ease of notation, we shall drop the Lebesgue measure $\mu(\cdot)$ in statements involving the area of domains in $\R^2$, since it can be inferred by context.

\begin{lemma}\label{lem:geom-thales}
\begin{figure}[t]
    \centering
    \includegraphics{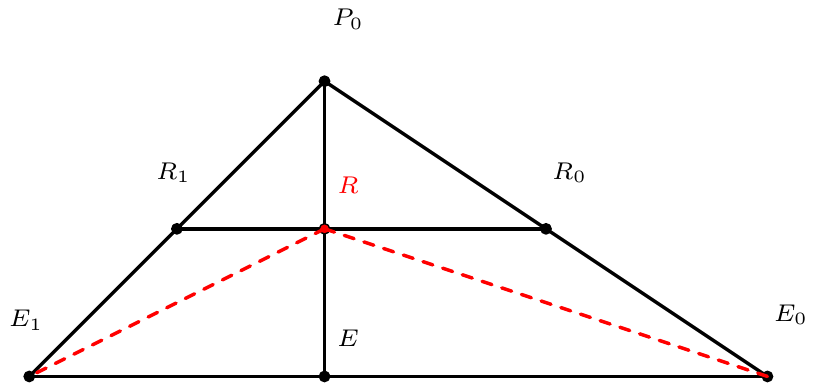}
    \caption{Area of cutoff triangles.}
    \label{fig:triangle-lemma}
\end{figure}
As shown in Figure~\ref{fig:triangle-lemma}, in triangle $\triangle E_0E_1P_0$, let $R_1$ and $R_0$ be points on $E_1P_0$ and $E_0P_0$ respectively such that $R_0R_1$ is parallel to $E_1E_0$. If $R$ is any point on $R_1R_0$, then \[\triangle E_0E_1R = \triangle E_0P_0E_1 - \sqrt{\triangle R_0P_0R_1 \cdot \triangle E_0P_0E_1}\]
\end{lemma}
\begin{proof}
By an affine transformation, we may assume the angle at $E_1$ is a right angle, $E_0E_1=1,~E_1P_0=1,~E_1R_1=r$ where $r \in (0,1)$. Then the identity to be proved is
\[
\frac{r}{2} =  \frac{1-\sqrt{(1-r)^2}}{2}
\]
which is immediate. 
\end{proof}

We now prove our main theorem. Let us denote by \[\h_i :=  \left\{z \in \star \ | \ T_i(z) \ge A_{i+1}(z)\right\},\quad i=0,1,2\]  Here and below, we consider the indices modulo 3. We shall derive a lower bound for all convex discs $K$ having $H_K(z)$ as a minimal midpoint hexagon, where $z$ belongs to the regions $\h_0, \h_1$ or $\h_2$. For $i=0,1,2$ define
\[ 
    \mathrm{area}_i(z) := \mathbf{I}_{\h_i}(z) \left(T_i(z) - \sqrt{A_{i+1}(z) \cdot T_i(z)}\right)  
\]

where $\mathbf{I}_{\h_i}(z)$ is the indicator function of the set $\h_i$. 

\begin{theorem}
If $K \in\Kbal$ is in circle representation and has $H_K(z)$ as a critical hexagon, where $z$ lies in the set $$ \h_0 \cup \h_1 \cup \h_2 = \{z \in \mathfrak{h}^{\ast} | \ T_0(z) \ge A_1(z) \ \mathrm{ or } \ T_1(z) \ge A_2(z) \ \mathrm{ or } \ T_2(z) \ge A_0(z) \}$$
then we have that 
\[\mu(D) \ge \frac{3\sqrt{3}}{2} + 2 \ \sum_{i} \mathrm{area}_i(z)\]
\end{theorem}
\begin{proof}
As in the Figure \ref{fig:compactification-result} above, let $P_0P_1P_2P_0'P_1'P_2'$ be the critical hexagon of an unseen disc $K$. The unseen disc $K$ is inscribed in this hexagon and passes through the points $\{e_j^*\}$ which are midpoints of its sides. 

First of all, we show that the set $\h_0 \cap \h_1 \cap \h_2$ is empty. That is, we cannot have $T_i(z) \ge A_{i+1}(z)$ for all $i=0,1,2$. Because if the inequalities hold, we have 
\[ 
\frac{\sqrt{3}}{4} = T_0 + T_1 + T_2 \ge A_0 + A_1 + A_2
\]
Substituting the expressions for $A_i$ from Lemma \ref{lem:triangle-areas}, this means that 
\[ 
\alpha^2 + \beta^2 + \gamma^2 \le \frac{3}{4}
\]
Again, from Lemma \ref{lem:triangle-areas} we have $T_0 + T_1 + T_2 = \frac{\sqrt{3}}{4}$ which means that $\alpha\beta+\beta\gamma+\gamma\alpha = 3$. So we have, by the Cauchy-Schwarz inequality, that 
\[ 
3 = \alpha\beta+\beta\gamma+\gamma\alpha \le \alpha^2 + \beta^2 + \gamma^2  \le \frac{3}{4}
\]
which is a contradiction. So we can only have the inequalities hold individually or pairwise. This gives us two cases:

\textbf{Case 1}: Without loss of generality assume that $z \in \h_0$ and $z \notin \h_1 \cup \h_2$. That means $T_0(z) \ge A_1(z)$. By the three-fold symmetry of the problem, the proof is identical (except for changing the indices) if we start with $z$ in the set $\h_1$ or $\h_2$. In Figure \ref{fig:compactification-result} above, this means that areas $\triangle E_0P_0E_1 \ge \triangle P_1Q_1P_2$. Construct a triangle $\triangle R_0P_0R_1$ such that $\mu(\triangle R_0P_0R_1) = \mu(\triangle P_1Q_1P_2)$ and such that the line segment $R_0R_1$ is parallel to $E_0E_1$. The triangles with equal area are shown in pink. 
We claim that there is at least one point $R$ on the line segment $R_0R_1$ which also lies on the boundary of the unseen disc $K$. This is because if there were no such point, then the disc $K$ would be strictly contained inside the octagon with vertices $R_0,R_1,P_1,P_2$ and their reflections. Thus, one can find points $S_1$ and $S_0$ on the line segments $E_0R_0$ and $E_1R_1$ respectively such that $K$ touches the segment $S_0S_1$. Thus, adding back the triangles $P_1Q_1P_2$ and $P_1'Q_1'P_2'$ one can construct a hexagon $S_0S_1Q_1S_0'S_1'Q_1'$ of strictly smaller area containing the unseen disc $K$. This contradicts the fact that $P_0P_1P_2P_0'P_1'P_2'$ was of minimal area. 

The above argument exhibits two points $R$ and $R'$ on $R_0R_1$ and $R_0'R_1'$ respectively which are also on the boundary of the disc $K$. Since $K$ was convex, it contains the convex hull $C$ of the points $E_0,E_1,E_2,R$ and their reflections.
Thus, we have 
\begin{align*}
\mu(D) \ge& \mu(C) \\
=&  \mu(E_0E_1E_2E_0'E_1'E_2'\}) + 2\triangle E_0RE_1\\
=& \frac{3\sqrt{3}}{2} + 2\triangle E_0RE_1 \\
=& \frac{3\sqrt{3}}{2} + 2\left(T_0 - \sqrt{A_1\cdot T_0}\right)\qquad \text{(using~Lemma~\ref{lem:geom-thales})}
\end{align*}

\textbf{Case 2}: Assume without loss of generality that $z \in \h_0 \cap \h_1$. We have $T_0(z) \ge A_1(z)$ and $T_1(z) \ge A_2(z)$. The above argument can also be adapted here again to exhibit four points (two new points along with their reflections) also on the boundary of the disc $K$, so that we have:
\[
\mu(D) \ge \frac{3\sqrt{3}}{2} + 2\left(T_0 - \sqrt{A_1\cdot T_0}\right) +  2\left(T_1 - \sqrt{A_2\cdot T_1}\right)
\]

Accounting for all cases, we have:
\[ \mu(D) \ge \frac{3\sqrt{3}}{2} + 2 \ \sum_{z \in \h_i} \mathrm{area}_i(z) \]
\end{proof}

\begin{figure}[htbp]
    \centering
    \includegraphics[scale=0.5]{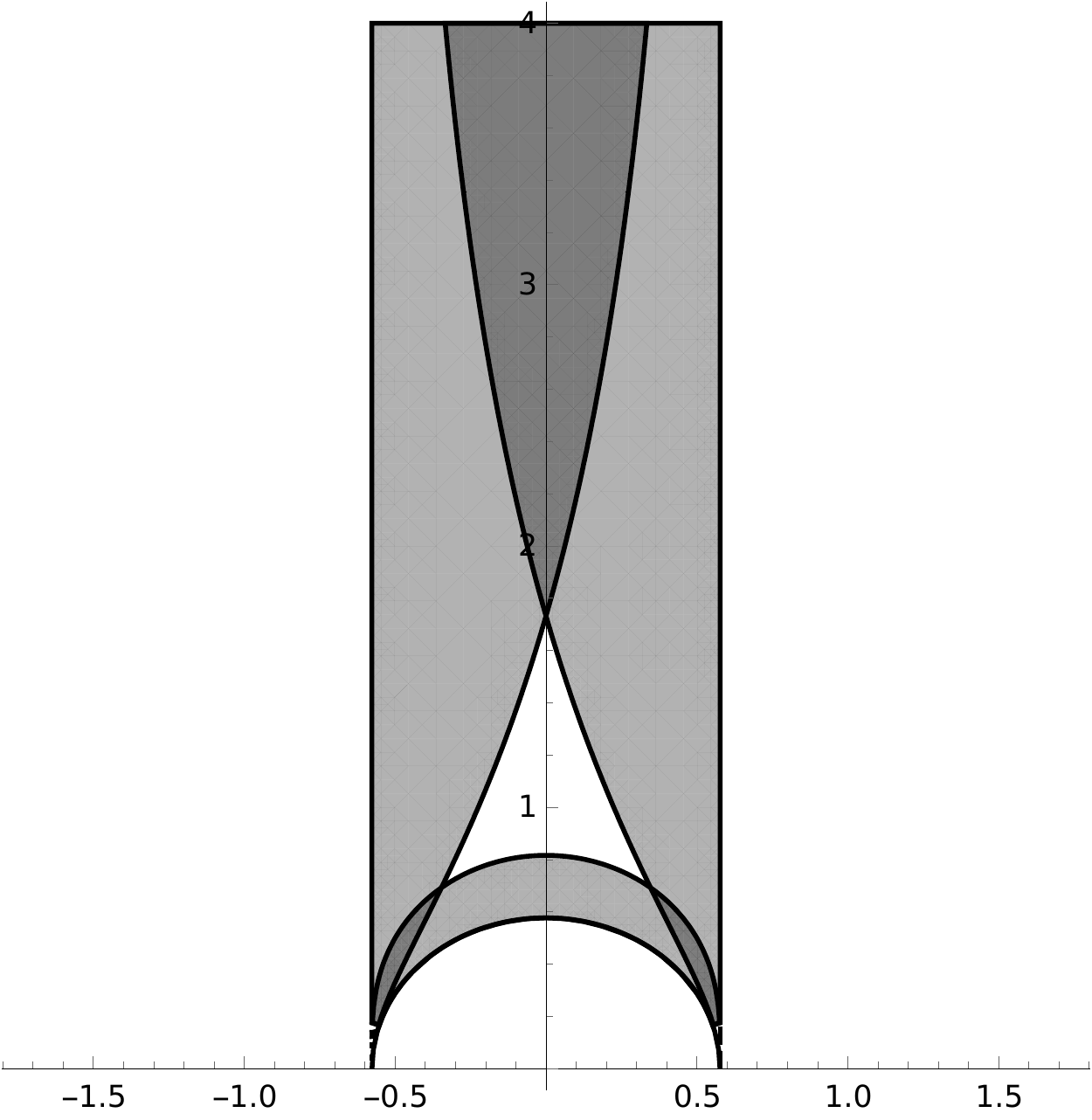}
    \caption{Covered boundary of the star domain.}
    \label{fig:union-contains-boundary}
\end{figure}

The sets $\h_0, \h_1, \h_2 \subset \star$ are given by:
\begin{align*}
\h_0 = \{ (x,y) \in \star \ | \ T_0 \ge A_1\} =&  \{(x,y) \in \star \ | \ 3 x^3 +3x y^2 - 7 \sqrt{3} x^2 + \sqrt{3}y^2  +15x-3\sqrt{3} \geq 0\} \\    
\h_1 = \{ (x,y) \in \star \ | \ T_1 \ge A_2\} =&\{(x,y) \in \star \ | \ -3x^4 -3y^4 -6x^2y^2 +x^2 + 2y^2\geq 0\} \\
\h_2 = \{ (x,y) \in \star \ | \ T_2 \ge A_0\} =&\{(x,y) \in \star \ | \ -3x^3 -3x y^2- 7\sqrt{3}x^2 +\sqrt{3}y^2 -15x -3\sqrt{3} \ge 0 \}
\end{align*}
A plot of $\h_0\cup\h_1\cup \h_2$ is shown in Figure \ref{fig:union-contains-boundary}.

Let us denote \[\mathcal{E}(z) := \frac{1}{\sqrt{12}}\left(\frac{3\sqrt{3}}{2} + 2 \ \sum_{z \in \h_i} \mathrm{area}_i(z) \right) = \frac{3}{4} +  \sum_{z \in \h_i} \frac{\mathrm{area}_i(z)}{\sqrt{3}} \] to be the lower bound of the packing density for discs $K \in \Kbal$ having $H_K(z)$ as a critical hexagon. Also, let $\delta_{oct}$ denote the packing density of the smoothed octagon. 
If $z \in \star$ is such that $\mathcal{E}(z) > \delta_{oct}$, then by the above theorem, the density of all convex discs having $H_K(z)$ as a critical hexagon is greater than $\delta_{oct}$ and so can be dropped from consideration. The next result shows that all but a compact subset of $\star$ can be excluded in this way. We start with a lemma.

\begin{lemma}\label{lem:compactification-lemma}\leavevmode
\begin{enumerate}
    \item For sufficiently large $y\gg 0$, the function $\mathrm{area}_0(x,y)$ is monotonically increasing in $x$ on $0<x<1/\sqrt{3}$.
    \item For sufficiently large $y\gg 0$, the function $\mathrm{area}_0(x,y) + \mathrm{area}_2(x,y)$ is monotonically increasing in $x$ on $0<x<1/\sqrt{3}$.
    \item For any $y$, the function $\mathrm{area}_0(x,y) + \mathrm{area}_2(x,y)$ takes its smallest value at $x=0$.
    \item If $\{z_n\}$ is a sequence going to infinity on the boundary $\partial\h_2$ given by 
    \[
    -3x^3 -3x y^2- 7\sqrt{3}x^2 +\sqrt{3}y^2 -15x -3\sqrt{3} = 0,
    \]
    the function $\mathcal{E}(z_n)$ goes to 1.  
\end{enumerate}
\end{lemma}
\begin{proof}
\begin{enumerate}
    \item We compute the partial derivative with respect to $x$ of $\mathrm{area}_0(x,y)$ to get:
    \begin{align*}
    24y^2\frac{\partial \mathrm{area}_0(x,y)}{\partial x} &= [6 x + 9 \sqrt{3} x^2] + \sqrt{3} [-1 + 3 y^2] \\ &+ \frac{[1 -\sqrt{3} x]\left[\left(15 x+\sqrt{3}\right) \left(3 x^2-1\right)+3 \left(9 x+\sqrt{3}\right) y^2\right]}{\sqrt{\left(\sqrt{3} x-1\right) \left(3 x^2-1\right) \left(3 x^2+3 y^2-1\right)}}
    \end{align*}
    For all $y>2$ and all $0<x<1/\sqrt{3}$, the terms in square parentheses above are positive, which means that the entire right hand side is positive. 
    \item Similarly, we compute the partial derivative \[M = \frac{\partial\left(\mathrm{area}_0(x,y) + \mathrm{area}_2(x,y)\right)}{\partial x}\] to get the following:
    \begin{align*}
12\sqrt{2}y^2M&= 
12x + \\& \frac{3(\sqrt{3}x+1)(\left(-5 x+1/\sqrt{3}\right) \left(3 x^2-1\right)+\left(-9 x+\sqrt{3}\right) y^2)
}{\sqrt{-(\sqrt{3} x+1) (3 x^2-1) (3 x^2+3 y^2-1)}} +
\\& \frac{3(-\sqrt{3}x+1)(\left(5 x+1/\sqrt{3}\right) \left(3 x^2-1\right)+\left(9 x+\sqrt{3}\right) y^2)
}{\sqrt{(\sqrt{3} x-1) (3 x^2-1) (3 x^2+3 y^2-1)}}.
    \end{align*}
Since the function $\mathrm{area}_0(x,y) + \mathrm{area}_2(x,y)$ is even in $x$, its derivative (the above expression) is odd in $x$. Also, this function goes to $\pm \infty$ as $x \to \pm 1/\sqrt{3}$ respectively and vanishes at $x = 0$.  It has no other roots in the interval $(0,1/\sqrt{3})$. These comments collectively show that $M$ is positive for $x$ in $(0,1/\sqrt{3})$. 
\item From the expression for the derivative $M$ derived above, we can see that $M$ vanishes when $x=0$. Since the function $\mathrm{area}_0(x,y) + \mathrm{area}_2(x,y)$ is even in $x$, and since we proved that it is increasing on $0<x<1/\sqrt{3}$, we get that it is decreasing on $-1/\sqrt{3}<x<0$. Thus, the critical point at $x=0$ is a minimum. 
\item Let $z=x+iy$. Along $\partial\h_2$, we get that $\mathrm{area}_2$ vanishes, and so $\mathcal{E}(z) = \mathrm{area}_0(z) + \mathrm{area}_2(z) = \mathrm{area}_0(z)$. We also have that the curve $\partial\h_0$ is asymptotic to $x = 1/\sqrt{3}$. Thus, if $y$ goes to infinity along this curve, we must have $x \to 1/\sqrt{3}$.
Let $y(x)$ denote the solution for $y$ in terms of $x$ of the equation determining the boundary of $\h_2$:
\[
-3x^3 -3x y^2- 7\sqrt{3}x^2 +\sqrt{3}y^2 -15x -3\sqrt{3} = 0,
\]
We choose the root for $y$ which is positive, since we are in the upper half-plane. We now can easily verify that $\lim_{x \to 1/\sqrt{3}}\mathrm{area}_0(x + i y(x)) = \frac{1}{4}$. The conclusion follows. 
\end{enumerate}
\end{proof}
\begin{theorem}
There exists a compact subset $\mathfrak{h}^{**} \subset \star$ such that all points of $\star \setminus \mathfrak{h}^{**}$ give rise to convex discs whose packing density is greater than that of the smoothed octagon. 
\end{theorem}
\begin{proof}
Note first that we have $T_0(z) = A_1(z) = 0$ if and only if $x = -\frac{1}{\sqrt{3}}$ and similarly for the other two cases. This shows that the boundaries of the star domain $\partial \star$ are contained in the set $\h_0 \cup \h_1 \cup \h_2$. This is shown in Figure~\ref{fig:union-contains-boundary}. 

We have from the previous theorem that for any $z \in \h_0 \cup \h_1 \cup \h_2$, the density $\delta(z)$ of any convex disc $K$ having $H_K(z)$ as a critical hexagon has the following bounds:
\begin{equation}\label{eq:delta-E}
1 \ge \delta(z) \ge \mathcal{E}(z) = \frac{3}{4} +  \sum_{z \in \h_i} \frac{\mathrm{area}_i(z)}{\sqrt{3}}
\end{equation}
We prove that as $z \rightarrow \partial\star$ we have $\mathcal{E}(z) \rightarrow 1$ uniformly.  
The star domain is an ideal triangle with one vertex at infinity. We prove our statement on the sides $x = \pm 1/\sqrt{3}$ and the vertex at infinity and conclude the same statement for the other vertices and side by appealing to symmetry. 
As a further reduction, we limit our analysis to the ``right half'' of the star domain \textit{viz.,} the region where $0 \le x < 1/\sqrt{3}$, since the other half is handled by symmetry. This region (non-exhaustively) consists of the regions $\h_0 \cap \h_2$ and $\h_0 \setminus (\h_0 \cap \h_2)$.

Now, let $\{z_n\}$ be any sequence in the right half of the star domain tending to $\infty$.
\begin{itemize}
    \item Assume that $z_n = x_n + iy_n$ goes to infinity while remaining inside $\h_0 \setminus (\h_0 \cap \h_2)$.  
    By Lemma \ref{lem:compactification-lemma},(1) we have $\mathcal{E}(z) = \mathcal{E}(x+iy)$ in this region is monotonically increasing in $x$ for all fixed and sufficiently large $y$.
    This means that:
    \[
    m(y_n) \le \sup_{0<x_n<1/\sqrt{3}}|\mathcal{E}(x_n+iy_n)| \le 1 
    \]
    for sufficiently large $n$, where $m(y_{n})$ is the value of $\mathcal{E}$ on the curve $\partial\h_0$. We also have that $\mathcal{E}(z) = 1$ when $x = 1/\sqrt{3}$. By Lemma \ref{lem:compactification-lemma},(4) we get  $m(y_{n}) \to 1$ as $n\to \infty$, and so we are done.   
    \item Assume that $z_n = x_n + iy_n$ goes to infinity while remaining inside $\h_0 \cap \h_2$. Then we have, by Lemma \ref{lem:compactification-lemma},(3) that, for all $n$:
    \[
    \sup_{-1/\sqrt{3}< x_n < 1/\sqrt{3}}|1-\mathcal{E}(x_n+iy_n)| = |1 -\mathcal{E}(iy_n)|
    \]
    We now show that this quantity is small as $n$ becomes larger. Within this region, $\mathcal{E}(z)= 3/4 + (\mathrm{area}_0(z) + \mathrm{area}_2(z))/\sqrt{3}$ and so the following limit proves the required:
    \[
       \lim_{y \to +\infty} \frac{\mathrm{area}_0(iy)+\mathrm{area}_2(iy)}{\sqrt{3}} = \lim_{y \to +\infty} \frac{y^2 \left(3y^2-2 \sqrt{6 y^2-2}\right)-y^2}{12 y^4} = \frac{1}{4}
    \]
\end{itemize}

The above results collectively prove that on each region we have that: for all $\epsilon > 0$ there exists a $y_c \gg 0$ such that for all $y > y_c$ and all $-1/\sqrt{3} <x<1/\sqrt{3}$, we have $|\mathcal{E}(x+iy)-1| < \epsilon$. So, choosing the largest $y_c$ among those that we get for each region, we may truncate the star domain by throwing away everything above that $y_c$, since all initial conditions in that region give a packing density arbitrarily close to 1. By symmetry, we may similarly truncate the cusps near the vertices on the $x$-axis. Thus, we have truncated the vertices of the star domain. This means that the $y$ can now only vary on a compact interval $[y_0,y_1]$ with $0<y_0<y_1<\infty$. 

 Now we turn to the sides. Note that we have:
\begin{align*}
    \mathcal{E}(z) &= \frac{3}{4} + \sum_{z \in \h_i}\frac{\mathrm{area}_i(z)}{\sqrt{3}} \\
    &\ge  \frac{3}{4} + \sum_{i = 0}^2\frac{\mathrm{area}_i(z)}{\sqrt{3}} \quad \text{(since } \h_0\cap\h_1\cap\h_2 = \emptyset) \\
    &= 1 - \left(\sqrt{\mu(A_0) \mu(T_2)} + \sqrt{\mu(A_2) \mu(T_1)} + \sqrt{\mu(A_1) \mu(T_0)}\right)\quad  (\text{by Lemma~\ref{lem:triangle-areas}}).
\end{align*}
This gives us \[|1 - \mathcal{E}(z)| \le |\sqrt{\mu(A_0) \mu(T_2)} + \sqrt{\mu(A_2) \mu(T_1)} + \sqrt{\mu(A_1) \mu(T_0)}|\]
The last expression is a continuous function in $x$ and $y$ since we have truncated $y$ to be away from zero. Also, we can check that this function goes to 0 as $x$ goes to the boundary for all $y$ uniformly, since $y$ varies in a compact set. This proves that: for all $\epsilon > 0$ there exists an $x_c$ such that for all $x > x_c$ and all $y \in [y_0,y_1]$, we have $|1 - \mathcal{E}(x+iy)| < \epsilon$. 
Equation \eqref{eq:delta-E}, we get $\delta(z) \rightarrow 1$ as $z$ goes to the boundary, by the sandwich theorem. Classical results of Mahler and Reinhardt show that we have $\delta(z) = 1$ on the boundary. So, by continuity, there is an open neighbourhood $S$ of the boundary contained within the set $\h_0 \cup \h_1 \cup \h_2$ such that $\delta(z) > \delta_{oct}$ for all $z \in S$. This shows that $\h^{**} := \star \setminus S$ is a closed and bounded set and so we are done. \end{proof}

\subsubsection{Control Problem in the Half-Plane}
  Evolution of $X(t)=h(t)Jh(t)^{-1}$ by adjoint action in $\mathcal{O}(J)$ corresponds to evolution by M\"{o}bius transformations of the corresponding element $z_0$ in the upper half-plane picture. 
    \[
    h(t)X_0 h(t)^{-1} \Longleftrightarrow h(t)\cdot z_0 = \frac{h_{11}(t)z_0 + h_{12}(t)}{h_{21}(t)z_0 + h_{22}(t)} 
    \]
    where
    \[
    h(t) = \left(\begin{matrix} h_{11}(t) & h_{12}(t) \\ h_{21}(t) & h_{22}(t) \end{matrix}\right) \in \SL(\R)
    \]
This means that the initial and terminal conditions derived in Section \ref{sec:X-init-term-conds} are transformed as:
\begin{align*}
X(0) = X_0 &\Longleftrightarrow z(0) = z_0, \\    
X(t_f) = R^{-1}X_0 R &\Longleftrightarrow z(0) = R^{-1}\cdot z_0
\end{align*}
where $R=\exp(J\pi/3)$ as usual.

Summarizing the results so far, we get the following reformulation of the Reinhardt conjecture from the coadjoint orbit of the Lie algebra to the upper half-plane picture.

\begin{problem}[Half-Plane Control Problem]\label{pbm:plane-optimal-control-problem}
On the set $\SL(\R)\times\star \subset T\SL(\R)$, consider the following free-terminal time optimal control problem: 
\begin{align*}
g' &= gX, \quad X = \mattwo {x/y} {-(x^2 + y^2)/y} {1/y} {-x/y}  \\[5mm]
x' &= f_1(x,y;u) := \frac{y \left(2 a x+b-c x^2+c y^2\right)}{2 a x+b-c x^2-c y^2} \\[5mm]
y' &=  f_2(x,y;u):= \frac{2 y^2 (a-c x)}{2 a x+b-c x^2-c y^2}         \\[5mm]
\frac{3}{2}&\int_0^{t_f}\frac{x^2+y^2+1}{y}dt \rightarrow \text{\normalfont min}  \\[5mm]
 & g:[0,t_f] \to \SL(\R),\quad x,y:[0,t_f] \to \star
\end{align*}
where the ODEs for $x,y$ are control-dependent:
\begin{equation*}
    a = \frac{u_2 - u_1}{\sqrt{3}} \quad
    b = \frac{u_0 - 2u_1 - 2u_2}{3} \quad
    c = u_0 
\end{equation*}
with  $(u_0,u_1,u_2) \in U_T$ which is the 2-simplex in $\R^3$. This problem has intial conditions $g(0) = I_2 \in \SL(\R)$ and $z(0) = z_0 \in \star$ and terminal conditions $g(t_f) = I_2$ and $z(t_f)=R^{-1}\cdot z_0$ where $R = \exp(J\pi/3)$.
Furthermore, the optimal trajectory is constrained to lie in a compact subset of the star region $\truncstar \subset \star$.
\end{problem}

As mentioned before, this problem is exactly equivalent to Problem \ref{pbm:reinhardt-optimal-control} because $O_X$ and $\h$ are anti-symplectomorphic (See Appendix \ref{sec:kirillov}).

We have finally reached the end of the reduction chain and have transformed a problem in discrete geometry to an optimal problem on $T\SL(\R)$. Already, we see that this problem is remarkably rich, with connections to Hamiltonian mechanics and hyperbolic geometry. 

\subsubsection{The Existence of Optimal Control}\label{sec:existence-opt-control}
The existence of optimal solutions to problems such as Problem \ref{pbm:plane-optimal-control-problem} is based on Filippov's theorem which gives conditions under which the corresponding \emph{attainable set} of the control system in question is compact. In informal terms, the attainable set corresponds to all the points in the manifold which are reachable provided one is allowed to move according to the control. Compactness of attainable sets imply the existence of optimal control. 

\begin{theorem}[Filippov]
    On a smooth manifold $M$, let $q' = f(q,u)$ be an optimal control system with an associated cost objective $\int_0^t \phi(q,u) dt \rightarrow \text{\emph{min}}$. Here $u \in U \subset \R^m$ which is a compact set. 
    Assume that the set $f(q,U) = \{ f(q,u)~|~u\in U\}$ is convex for each $q \in M$ and that there exists a compact set $K \subset M$ such that $f(q,u) = 0$ for $q \notin M, u\in U$. Then, the attainable sets are compact.
\end{theorem}

For the Reinhardt control problem, all assumptions of the Filippov theorem are shown to hold as follows:
\begin{itemize}
    \item the control set is $U_T$ defined in Definition \ref{def:triang-control-set} is obviously compact.
    \item the velocity sets are the convex hull of the velocities at the boundary of the control set, as shown in Hales~\cite{hales2017reinhardt} for the simplex $U_T$. 
    \item Recall that we have shown that the optimal trajectory is confined to a compact subset within the star domain \emph{viz.,} $\truncstar \subset \star$. We can multiply the right hand side of the vector fields in the star domain in Problem \ref{pbm:plane-optimal-control-problem} by the smooth cutoff function $\psi : \star \to \R$ with $\psi|_{\truncstar} = 1$ and zero everywhere else, to obtain a system with vector fields having compact support.  We can now apply Filippov's theorem to this new system and since there is a compact set in which both systems coincide, we obtain our desired conclusion on the common domain for the initial system. 
\end{itemize}
Thus, the optimal control exists for the plane optimal control system.  


\subsection{The Pontryagin Maximum Principle}\label{sec:PMP}
The Pontryagin Maximum Principle (PMP) is a powerful first-order necessary condition for optimality of solutions to an optimal control problem on a smooth manifold $M$ with control set $U \subseteq \R^m$ and free-terminal time. Given a control system 
\begin{align}
q' &= f(q,u) \in T_qM \\ 
q(0) &= q_0 \in M
\end{align}
on a manifold $M$ with an associated cost objective $$\min_{u \in U}\int_0^{t_f} \phi(q,u) dt$$
we assume that the control-dependent vector field $f(q,u)$ satisfies: 
\begin{enumerate}
    \item $q \mapsto f(q,u)$ is a smooth vector field on $M$ for any fixed $u\in U$.
    \item $(q,u) \mapsto f(q,u)$ is a continuous mapping for $q \in M, u \in \bar{U}$. 
\end{enumerate}
This optimal control system is denoted by the tuple $(M,U,f,\phi)$. We shall sometimes denote the vector field $f(u,q)$ for $u\in U$ by $f_u(q)$. 

The PMP first associates a \textit{control-dependent Hamiltonian} which is \emph{cost extended} on $T^*M$:
\[
H(q,p,u) = \langle p,f(q,u) \rangle_{*} + \lambda_{cost} \phi(q,u) \qquad p \in T^*_qM \quad \lambda_{cost} \in \R_{\le 0}
\]
where $\bracks{\cdot}{\cdot}_{*}$ is the natural pairing between a vector space and its dual and $\lambda_{cost}$ is a constant nonpositive scalar. Note that $H$ is linear in $p$. 

Let $u^*(t)$ denote the function defined implicitly as a function of $(x(t),p(t)) \in T^*M$ by\footnote{Note that it is quite possible that $u^*(t)$ may not be determined this way. This is when the extremal is called \textit{singular}. See Definition~\ref{def:normal-abonormal-singular-extremals}.}

\begin{equation}\label{eq:maximized-hamiltonian}
     H(q(t),p(t),u^*(t)) = \max_{u\in U}H(q(t),p(t),u) =: H^*(q(t),p(t),t)
\end{equation}
where $H^*$ is called the \textit{maximized Hamiltonian}. Then the PMP says that the extremals of the optimal control problem are projections (from $T^*M$ to $M$) of the flow-trajectory of $\overrightarrow{H^*}$, which is the Hamiltonian vector field corresponding to $H^*$ with respect to the canonical symplectic structure on the cotangent bundle $T^*M$. Integral curves of the vector field $\overrightarrow{H^*}$ satisfy:
\begin{equation}\label{eq:hamiltons-equations}
q' = \frac{\partial H^*}{\partial p} = f(x,u) \qquad p' = -\frac{\partial H^*}{\partial q}
\end{equation}
This trajectory $(u(t),q(t),p(t))$ in $T^*M$ is called the \textit{lifted controlled trajectory}.    

The PMP (for free terminal time periodic problems) also guarantees the following of lifted trajectories:
\begin{enumerate}
    \item Transversality conditions (endpoint conditions for the co-state trajectories) hold.
    \item The Hamiltonian $H(q,p,u^*)$ vanishes identically along the lifted controlled trajectory. 
    \item The covector $(\lambda_{cost},p(t)) \in \R_{\le 0} \times T^*_{q(t)}M$ is nonzero for all $t \in [0,t_f]$.  
    \item The quantity $\lambda_{cost}$ is constant, and when it is non-zero, may be taken to be $\lambda_{cost} = -1$ by rescaling the Hamiltonian. 
\end{enumerate} 

The lifted curves which satisfy the conditions of the PMP are called \emph{Pontryagin extremals} or simply \emph{extremals}.

\begin{definition}[Normal, Abnormal and Singular extremals]\label{def:normal-abonormal-singular-extremals}\leavevmode
\begin{enumerate}
    \item The extremal for which $\lambda_{cost} = 0$ is called an \emph{abnormal extremal}.
    \item The extremal for which $\lambda_{cost} \ne 0$ is called a \emph{normal extremal}. In the normal case a renormalization allows us to take $\lambda_{cost} = -1$. 
    \item If there is an open time interval in which equation \eqref{eq:maximized-hamiltonian} fails to uniquely determine the function $u^*(t)$, the extremal is called a \emph{singular extremal} or \emph{singular curve}.
\end{enumerate}
\end{definition}

Our strategy is to apply the maximum principle to our problem with the hope that these necessary conditions will provide us more information about the structure of the extremals. Recall the framework of our control problem in Problem \ref{pbm:reinhardt-optimal-control}, where we have dynamics occurring in the Lie group $\SL(\R)$ and the Lie algebra $\sl(\R)$. If we apply the PMP to this problem, the lifted trajectories live in \[T^*(\SL(\R) \times \sl(\R)) \cong \left(\SL(\R) \times \sl(\R)\right)\times \left (\sl(\R) \times \sl(\R)\right),\]
where we have used Proposition \ref{prop:trivialization} making $T^*\SL(\R) \cong \sl(\R)^* \times \SL(\R)$ and the identification $\sl(\R)^* \cong \sl(\R)$ via the trace inner product as in Appendix \ref{sec:X-lie-poisson}. So, for the state variables $p = (g,X) \in \SL(\R) \times \sl(\R)$ the PMP gives corresponding costate variables $q= (\Lambda_1,\Lambda_2) \in \sl(\R) \times \sl(\R)$. 

We now derive expressions for the Hamiltonian and the costate equations in both the Lie algebra coordinates and and the upper half-plane coordinates via the isomorphism described in Lemma \ref{lem:def-phi}. 

\subsubsection{Hamiltonian in the Lie Algebra and Half-Plane}\label{sec:max-ham}
Following our reformulation of the control problem as Problem \ref{pbm:reinhardt-optimal-control}, the Hamiltonian for this problem is composed of the Hamiltonians for the Lie group and the Lie algebra. 

The costate variable corresponding to the Lie group element $g \in \SL(\R)$ is denoted $\Lambda_1 \in \sl(\R)$. Ignoring the Lie algebra dynamics for a moment, we derive the Hamiltonian corresponding to the group element $g$. As pointed out earlier, the control problem is left-invariant (see Definition \ref{def:left-invariant-control-system}) and Hamiltonians of left-invariant systems are functions of $\sl(\R)^*$ only:

\begin{proposition}[Jurdjevic~\cite{jurdjevic_1996}]\label{prop:left-invariant-hams}
Consider an arbitrary left-invariant control system $dg/dt = F(g,u)$ on a real Lie group $G$, with $u \in U \subseteq \R^m$, with associated cost objective $\int_0^T p(g,u) dt \rightarrow \text{min}$. Let us also assume that the Lie algebra $\mathfrak{g}$ of $G$ is equipped with an invariant inner-product $\bracks{\cdot}{\cdot}$. Then the Hamiltonian function corresponding to this system looks like:
\[
H_F(g,p) = \bracks{p}{F(e,u)}, \quad p \in \mathfrak{g}
\]
where $e \in G$ is the group identity.
\end{proposition}
\begin{proof}
Any vector field $F(g,u)$ on a Lie group $G$ determines a linear function on the cotangent bundle $T^*G$ via the function $H_F(\xi) = \langle \xi,F(g,u)\rangle_*$ for $\xi \in T_g^* G$ (See Definition 6.9 of Jurdjevic~\cite{jurdjevic2016optimal}). Here $\bracks{\cdot}{\cdot}_*$ is the canonical pairing between a vector space and its dual. Now, since $T^*G \cong G \times \mathfrak{g}^*$ by Proposition \ref{prop:trivialization}, we have that $\xi = T^*L_{g^{-1}}(p)$ for some $p \in \mathfrak{g}^*$. Then we have
\[
H_F(g,p) = \bracks{p}{TL_{g^{-1}}\left(F(g,u)\right)}_* = \bracks{p}{TL_{g^{-1}}\left(T_gF(e,u)\right)}_* = \bracks{p}{F(e,u)}_*
\]
since the control system is left-invariant. Since $\mathfrak{g}$ is equipped with an invariant inner product $\bracks{\cdot}{\cdot}$, we have that $\mathfrak{g} \cong \mathfrak{g}^*$. Using this identification, we have $H_F(g,p) = \bracks{p}{F(e,u)}$.
\end{proof}
In our case, this means that the (cost-extended) Hamiltonian corresponding to the group is 
\[
\H_1(\Lambda_1,X) := \bracks{\Lambda_1}{X} - \frac{3}{2}\lambda_{cost}\bracks{J}{X}
\]

The Lie algebra part of the Hamiltonian, corresponding to the costate variable $\Lambda_2$, is
\[
\H_2(\Lambda_2,X;Z_u) := \frac{\bracks{\Lambda_2}{[Z_u,X]}}{\bracks{X}{Z_u}} = -\frac{\bracks{[\Lambda_2,X]}{Z_u}}{\bracks{X}{Z_u}}
\]
which suggests introducing a new variable $\Lambda_R := [\Lambda_2,X]$. The full Hamiltonian of the problem is now:
\begin{align}
\H(\Lambda_1, \Lambda_R, X ; Z_u) &:= \H_1(\Lambda_1,X) + \H_2(\Lambda_2,X;Z_u) \nonumber \\ 
&= \left\langle \Lambda_1 - \frac{3}{2}\lambda_{cost} J, X \right\rangle - \frac{\bracks{[\Lambda_2,X]}{Z_u}}{\bracks{X}{Z_u}} \nonumber \\ 
&= \left\langle \Lambda_1 - \frac{3}{2}\lambda_{cost} J, X \right\rangle  - \frac{\langle \Lambda_R,Z_u\rangle }{\langle X, Z_u \rangle} \label{eq:full-hamiltonian}
\end{align}

We now wish to transport the Hamiltonian from Lie algebra coordinates to upper half-plane coordinates. Recall that we have an isomorphism $\Phi : \mathfrak{h} \to \O_X \subset \sl(\R)$ defined in Lemma \ref{lem:def-phi}. This induces the tangent map $T\Phi : T_z\h \to T_X\O_X$ as described in Lemma \ref{lem:tangent-maps}. This also induces the dual (cotangent) map: $T^*\Phi:T^*_X\O_X \to T^*_z\h$.

\begin{lemma}\label{lem:cotangent-space-OX}
\[ 
     T^*_X \OX \cong X^{\perp} = \{ [W,X] \ | \ W \in \sl(\R)\}
\]
where $X^{\perp} = \{ Y \in \sl(\R) \ | \ \langle Y, X \rangle = 0 \}$.
\end{lemma}
\begin{proof}
By a general fact in elementary linear algebra we have:
\[ 
     T^*_X \OX \cong \left(\sl(\R)/\sl(\R)_X\right)^* \cong \sl(\R)_X^\circ
\]
where $\sl(\R)_X^\circ = \{ Y \in \sl(\R) \ | \ \langle Y, W \rangle = 0 \text{ for all } W \in \sl(\R)_X \}$ stands for the annihilator\footnote{Note that the annihilator is actually a subspace of the linear dual of a vector space, but here we identify the dual with the Lie algebra itself via the nondegenerate trace form.} of the centralizer of $X$. Since $X$ is a regular element in $\sl(\R)$, it's centralizer is one-dimensional and so $\sl(\R)_X^\circ = X^{\perp}$. 

It is clear that any Lie algebra element of the form $[W,X]$ is orthogonal to $X$ since $\langle [W,X],X \rangle = \langle W,[X,X]\rangle = 0$. 
Dimension counting again gives us that   $T^*_X \OX = \{ [W,X] \ | \ W \in \sl(\R)\} $.
\end{proof}

For the next theorem, we define $\nu := T^*\Phi(-\Lambda_R) = \nu_1 dx + \nu_2 dy \in T_z^*\h$.
\begin{theorem}[Hamiltonian in Upper Half-Plane]\label{thm:plane-hamiltonian}
    Let $\Phi$ be as in Lemma \ref{lem:def-phi} and let $\Lambda_R$ be as above.
    Then we have that the Hamiltonian $\mathcal{H}(\Lambda_1,\Lambda_R,X;Z_u)$ in equation \eqref{eq:full-hamiltonian}
     in upper half-plane coordinates becomes:
     \[
     H(\Lambda_1,x,y,\nu_1,\nu_2;u) =  \left\langle \Lambda_1 - \frac{3}{2}\lambda_{cost} J, X \right\rangle + \nu_1 f_1 + \nu_2 f_2 
     \]
\end{theorem}
\begin{proof}
The quantity $\nu \in T^*_z\h$ is well-defined since $\Lambda_R \in X^{\perp}$ and by Lemma \ref{lem:cotangent-space-OX}.
The only part of the Hamiltonian in equation \eqref{eq:full-hamiltonian} which will change is $\H_2(\Lambda_2,X;Z_u)$.
     \begin{align*}
    \H_2(\Lambda_R,X;Z_u) =      \frac{\langle -\Lambda_R,Z_u \rangle}{\langle Z_u,X \rangle} &=\left\langle -\Lambda_R, \frac{Z_u}{\langle Z_u,X \rangle}\right\rangle  \\
          &= \left\langle  -\Lambda_R, \frac{Z_u}{\langle Z_u,X \rangle} + \sl(\R)_X \right\rangle \\
          &=   \langle -\Lambda_R, T \Phi(f_1,f_2)\rangle \qquad \text{by \eqref{eq:tangent-map-f1f2}} \\
          &=\langle T^*\Phi(-\Lambda_R), (f_1,f_2)\rangle_* \\
          &= \langle \nu, (f_1,f_2) \rangle_*\\
          &= \nu_1 f_1 + \nu_2 f_2
    \end{align*}
    where we have used the definition of the cotangent map and the fact that $\Lambda_R \in X^{\perp}$, the annihilator subspace of $X$ in $\sl(\R)$.
\end{proof}

The maximum principle states that the extremals of the control problem are integral curves of the so-called maximized Hamiltonian. This is the pointwise maximum of the control-dependent Hamiltonian over the control set. For our immediate application, we take the control set to be the simplex $U_T$ (see Definition \ref{def:triang-control-set}):
\begin{align}
\H^*(\Lambda_1,\Lambda_2, X) &:= \max_{U_T}\H(\Lambda_1, \Lambda_R, X ; Z_u) \nonumber \\ 
&=\left\langle \Lambda_1 - \frac{3}{2}\lambda J, X \right\rangle + \max_{U_T} \frac{ \langle -\Lambda_R , Z_u \rangle}{\langle X,Z_u \rangle} \label{eq:max-ham}
\end{align}

\subsubsection{Costate Variables in Lie Algebra and the Half-Plane}\label{sec:costate-variables}
\begin{proposition}[Lie algebra costate variables]
The costate variables evolve as:
\begin{align}
\Lambda_1' &= [\Lambda_1, X] \\
\Lambda_R' &= \left(\left[P,\Lambda_R\right] - \frac{\langle \Lambda_R,Z_u \rangle}{\langle X, Z_u\rangle} \left[P,X\right]\right) + \left[ -\Lambda_1 + \frac{3}{2}\lambda_{cost} J ,X\right] \label{eq:lamR}
\end{align}
where $P = Z_u/\bracks X {Z_u}$ and $\Lambda_R := [\Lambda_2,X]$.
\end{proposition}
\begin{proof}
Let $\H^*$ denote the maximized Hamiltonian in equation \eqref{eq:max-ham}. Let \[(g,X,\Lambda_1,\Lambda_2) \in T^*(\SL(\R)\times \sl(\R)) \cong \SL(\R) \times \sl(\R) \times \sl(\R) \times \sl(\R)\] denote the state and costate variables.\footnote{Here we use the left-trivialization to write $T^*(\SL(\R)) \cong \SL(\R) \times \sl(\R)^*$ and use the invariant trace product to identify $\sl^*(\R)$ and $\sl(\R)$ as usual. } Since the Hamiltonian is independent of $\SL(\R)$, it is left-invariant. 

By the maximum principle, we have that the state and costate equations are Hamilton's equations in the cotangent bundle with respect to an appropriate symplectic form. For the $(g,\Lambda_1)$ pair, the dynamics is Hamiltonian with respect to the pullback of the canonical symplectic form on $T^*(\SL(\R))$ to $\SL(\R) \times \sl^*(\R)$. These equations are called the Euler-Arnold equations for a left-invariant Hamiltonian (See \cite[pp.~285]{cushman1997global}) and are given by:
\begin{align*}
    g'&=g\frac{\delta \H^*}{\delta \Lambda_1} = gX \\
    \Lambda_1' &= \ad^*_{\delta \H^*/\delta \Lambda_1} \Lambda_1 = -\ad_{X}\Lambda_1 = [\Lambda_1,X]
\end{align*}
where we identify $\sl(\R)^*$ with $\sl(\R)$ as usual, sending the $\ad^*$-operator to $\ad$-operator (note that the star on $\H^*$ has a different meaning than the star on $\ad^*$). The expression $\delta \H^*/\delta \Lambda_1$ denotes the functional derivative of $\H^*$ with respect to $\Lambda_1$ and is defined in Appendix \ref{sec:X-lie-poisson} in equation \eqref{eq:functional-derivative}. 
For the pair $(X,\Lambda_2) \in \sl(\R) \times \sl(\R)^*$, the dynamics is Hamiltonian with respect to the canonical symplectic structure on the trivial cotangent bundle $T^*(\sl(\R))$ which gives us Hamilton's equations in the usual form: 
\begin{align*}
    X' &= \frac{\delta \H^*}{\delta \Lambda_2} =  [P,X] \\
    \Lambda_2' &= -\frac{\delta \H^*}{\delta X} = -\Lambda_1 + \frac{3}{2}\lambda_{cost}J - [\Lambda_2, P] + \bracks {[\Lambda_2,P]} {X} P
\end{align*}
Using this, we can derive:
\begin{align*}
    \Lambda_R' &= [\Lambda_2',X] + [\Lambda_2,X']\\
    &=\left[-\Lambda_1 + \frac{3}{2}\lambda_{cost}J ,X\right] +  \bracks {[\Lambda_2,P]} {X} [P,X] - [[\Lambda_2, P],X] + [\Lambda_2,[P,X]] \\
    &= \left[-\Lambda_1 + \frac{3}{2}\lambda_{cost}J ,X\right] + \bracks {[\Lambda_2,P]} {X} [P,X] +[P,[\Lambda_2,X]] \\
    &= \left[-\Lambda_1 + \frac{3}{2}\lambda_{cost}J ,X\right] - \bracks {\Lambda_R} {P} [P,X] +[P,\Lambda_R] \\
    &= \left(\left[P,\Lambda_R\right] - \frac{\langle \Lambda_R,Z_u \rangle}{\langle X, Z_u\rangle} \left[P,X\right]\right) + \left[ -\Lambda_1 + \frac{3}{2}\lambda_{cost} J ,X\right] 
\end{align*}
\end{proof}

\begin{remark}\normalfont \leavevmode
Note that the variable $\Lambda_R$ is constrained to lie in the 2-dimensional subspace $\{A\in\sl(\R)~|~\bracks{A}{X} = 0\}$. 
Thus, we can consider $\Lambda_R$ to be the \emph{reduced} costate variable. 
\end{remark} 

\begin{corollary}[Jurdjevic~\cite{jurdjevic_1996,jurdjevic2016optimal}]\label{cor:lam1-gen-sol}
The variable $\Lambda_1(t)$ whose dynamics is given by a Lax pair, evolves in a coadjoint orbit of $\sl(\R)$ through the initial value $\Lambda_1(0)$ and its general solution is given by 
\[
\Lambda_1(t) = \Ad^{*}_{g(t)^{-1}}(\Lambda_1(0)) = g(t)^{-1}\Lambda_1(0)g(t)
\]
this also implies that $\det(\Lambda_1(t))$ is a constant of motion.
\end{corollary}
\begin{proof}
This can be verified by differentiating. We immediately get that its determinant is constant. 
\end{proof}


\begin{corollary} \label{cor:control-dep-ad-alternate}

In the ODE for $\Lambda_R$ in equation \eqref{eq:lamR}, the control dependent term has the following expression:
\begin{equation}
  \left(\left[P,\Lambda_R\right] - \frac{\langle \Lambda_R,Z_u \rangle}{\langle X, Z_u\rangle} \left[P,X\right]\right) =  \ad_{Z_u} \frac{\delta}{\delta Z_u} \frac{\langle \Lambda_R, Z_u\rangle}{\langle X, Z_u\rangle} = \ad_{Z_u} \frac{\delta \mathcal{H}_2}{\delta Z_u} 
\end{equation}
where $\frac{\delta}{\delta Z_u}$ denotes the functional derivative\footnote{See Section~\ref{sec:X-lie-poisson} for the definition of the functional derivative.}.
\end{corollary}
\begin{proof}
Look at the following computation:
\begin{align*}
     \ad_{Z_u} \frac{\delta}{\delta Z_u}\left( \frac{\langle \Lambda_R, Z_u\rangle}{\langle X, Z_u\rangle}\right) &= 
 \left[Z_u,\frac{\delta}{\delta Z_u}\left( \frac{\langle \Lambda_R, Z_u\rangle}{\langle X, Z_u\rangle}\right) \right] \\
     &=   \left[Z_u, \frac{\Lambda_R \langle X, Z_u \rangle - \langle \Lambda_R, Z_u\rangle X}{ \langle X, Z_u \rangle^2} \right] \\
     &=  \frac{\left[ Z_u, \Lambda_R \right]}{\langle X, Z_u\rangle}  -\frac{ \langle \Lambda_R, Z_u \rangle}{\langle X, Z_u\rangle^2} \left[ Z_u, X \right] \\
     &= [P,\Lambda_R]- \frac{\langle \Lambda_R, Z_u \rangle}{\langle X,Z_u\rangle}[P,X]
\end{align*}
\end{proof}

We now aim to transfer the costate ODE of $\Lambda_R$ to the upper-half plane. Let us begin with:

\begin{lemma}\label{lem:cotangent-maps}
     As 1-forms, we have the following expression for $T^*\Phi: T^*_X\O_X \to T^*_z\h$ and its inverse:
     \begin{align}
          T^*\Phi[W,X] &= -\langle W, dX\rangle = -\left\langle W, \frac{\partial X}{\partial x}dx +  \frac{\partial X}{\partial y} dy\right\rangle \\ 
          (T^*\Phi)^{-1}(v_1,v_2) &=v_1 dh_1 + v_2 dh_2
     \end{align}
     where $W\in \sl(\R)$, $(v_1,v_2)\in T_z\h$, $h_1 = b + 2a x - c x^2 + cy^2$, $h_2 = 2y(a-cx)$ and $dh_1,dh_2$ have the following expressions:
     \begin{align*}
          dh_i &= \left(\frac{\partial h_i}{\partial a} da + \frac{\partial h_i}{\partial b} db +\frac{\partial h_i}{\partial c} dc \right), \quad i=1,2
     \end{align*}
\end{lemma}
\begin{proof}
By the Lemma \ref{lem:cotangent-space-OX}, a typical element of $T_X^*\OX$ looks like $[W,X]$ where $W\in \sl$. 
By definition, the cotangent map is given by the following:
\begin{align*}
\langle \underbrace{T^*\Phi(\underbrace{[W,X]}_{\in T_X^*\OX})}_{\in T_z^*\h}, \underbrace{v}_{\in T_z\h} \rangle_{*}
&=
\langle \underbrace{[W,X]}_{\in T_X^*\OX}, \underbrace{T\Phi(v)}_{\in T_X\OX} \rangle_{*}
\end{align*}
where $\bracks{\cdot}{\cdot}_{*}$ is the natural pairing between the cotangent space and the tangent space. 
So, for any vector $v = (v_1,v_2) \in T_z\h$:
\begin{align*}
    \langle T^*\Phi[W,X],v \rangle_{*} &= \langle [W,X],T\Phi(v) \rangle\\
    &= \langle W,[X, T\Phi(v)] \rangle \\
&= \left\langle W,\left[X,  \mattwo {v_2 /2y}{(yv_1 - v_2 x)/y}{0}{-v_2/2y}\right] \right\rangle \\
&= - \left\langle W, v_1 \frac{\partial X}{\partial x} + v_2 \frac{\partial X}{\partial y} \right\rangle
\end{align*}
which gives us the required. 
For the inverse, if $[W] = \left[ \mattwo a b c {-a} \right] \in T_X\OX$:
\begin{align*}
     \langle T^*\Phi(v), [W]\rangle &= \langle v, T\Phi[W]\rangle \\
     &= \langle (v_1,v_2), (h_1,h_2)\rangle \\
     &= v_1 h_1 + v_2 h_2 \\
     &= v_1(b + 2a x - c x^2 + cy^2) + v_2(2y(a-cx)) \\
     &= 2(v_1x + v_2y)a+ v_1 b -(v_1(x^2 - y^2) + 2 v_2 x y)c \\
     &= v\cdot\left(\frac{\partial h_1}{\partial a},\frac{\partial h_2}{\partial a}\right) a +v\cdot\left(\frac{\partial h_1}{\partial b},\frac{\partial h_2}{\partial b}\right) b+v\cdot\left(\frac{\partial h_1}{\partial c},\frac{\partial h_2}{\partial c}\right) c \\
     &= v_1 \left(\frac{\partial h_1}{\partial a} a + \frac{\partial h_1}{\partial b} b +\frac{\partial h_1}{\partial c} c \right) +  v_2 \left(\frac{\partial h_2}{\partial a} a + \frac{\partial h_2}{\partial b} b +\frac{\partial h_2}{\partial c} c \right) \\
     &= v_1 dh_1(a,b,c) + v_2 dh_2(a,b,c)
\end{align*}
\end{proof}

Using this, we will transfer the ODE for $\Lambda_R \in T_X^*\OX$ into the ODE for $\nu \in T_z^*\h$. Recall that we have the following expression for the Hamiltonian.

\begin{lemma}\label{lem:helper-cotangent}
     Let $X$ be as before and denote \[Y_1 = \frac{y^2}{2}\frac{\partial X}{\partial x},\quad Y_2 = \frac{y^2}{2}\frac{\partial X}{\partial y},\quad X_1 = \frac{\partial X}{\partial x},\quad X_2 = \frac{\partial X}{\partial y}.\]
     We have the following facts:
     \begin{enumerate}
     \item The set $\left\{ X, Y_1, Y_2\right\}$ forms an orthogonal basis for $\sl(\R)$. 
     \item For all $i,j =1,2$, $\left\langle Y_i,X_j\right\rangle = \delta_{ij}$.
     \item For all $i=1,2$, we have:
     \[
          \left\langle \frac{\partial Y_i}{\partial x},X_2\right\rangle =           \left\langle \frac{\partial Y_i}{\partial y},X_1\right\rangle
     \]
     \item Along any trajectory $x' = f_1$ and $y'=f_2$ in $\mathfrak{h}$, we have:
     \[
     \langle Y_i', dX\rangle = \langle dY_i,X'\rangle    
     \]
     \end{enumerate}
Here, for a function $f$, we write $df$ to indicate the total deriviative with respect to $x,y$ variables and holding all other variables constant.  
\end{lemma}
\begin{proof}
    The first three are elementary calculations. For the fourth, for $i=1,2$:
    \begin{align*}
          \langle Y_i', dX\rangle &= \left\langle \frac{\partial Y_i}{\partial x} \frac{dx}{dt} +\frac{\partial Y_i}{\partial y} \frac{dy}{dt},
          \frac{\partial X}{\partial x} dx + \frac{\partial X}{\partial y} dy
          \right\rangle \\
          &= \left\langle \frac{\partial Y_i}{\partial x} \frac{dx}{dt},
          \frac{\partial X}{\partial x} dx 
          \right\rangle +  \left\langle \frac{\partial Y_i}{\partial y} \frac{dy}{dt},
          \frac{\partial X}{\partial y} dy 
          \right\rangle \\
          &= \langle dY_i, X'\rangle
    \end{align*}
    since the cross-terms vanish. 
\end{proof}

\begin{theorem}\label{thm:lamR-dynamics}
    The ODE for $\Lambda_R$ gets transformed to the ODE
    \[
    \nu' = -d\H = -\frac{\partial \H}{\partial x} dx - \frac{\partial \H}{\partial y} dy
    \]
    where $\H$ is the Hamiltonian derived in equation \eqref{eq:full-hamiltonian}. 
\end{theorem}
\begin{proof}
Throughout this proof, we make the identification of $X = X_z$ in upper half-plane coordinates via the isomorphism $\Phi$ defined in Lemma $\ref{lem:def-phi}$. This allows us to write $df$ (for $f$ a function of coordinates in the upper half-plane) in the same convention we adopt in Lemma \ref{lem:helper-cotangent}.

Since $X,Y_1,Y_2$ is a basis of $\sl(\R)$, the subspace $T_X^*\OX$ has $[Y_1,X],[Y_2,X]$ as a basis. So we have 
\[
\Lambda_R = [Y_1,X]\nu_1 + [Y_2,X]\nu_2
\]
We can now compute $T^*\Phi\Lambda_R =  T^*\Phi[Y_1,X]\nu_1+T^*\Phi[Y_2,X]\nu_2 = -\nu$ using the expression for $T^*\Phi$ in Lemma \ref{lem:cotangent-maps} (which is consistent with the definition of $\nu$ in Theorem \ref{thm:plane-hamiltonian}).

Now recall the expression for the Lie algebra Hamiltonian $\H$:
\begin{align*}
    \H &= \H_1 + \H_2 \\
    &=  \left\langle L, X \right\rangle - \bracks{\Lambda_R}{P}
\end{align*}
where $L = \Lambda_1 - \frac{3}{2}\lambda_{cost} J$ and $P = Z_u/\bracks{Z_u}{X}$. 

We now have:
\begin{align*}
    -d\H &= -d \left\langle L, X \right\rangle + d\bracks{\Lambda_R}{P} \\
     &= -\left\langle L, dX \right\rangle + \bracks{d\Lambda_R}{P} + \bracks{\Lambda_R}{dP} 
\end{align*}
We have $-\bracks{L}{dX} = T^*\Phi[L,X]$ and 
\[
\bracks{\Lambda_R}{dP} = \left\langle \Lambda_R, d\frac{Z_u}{\bracks{Z_u}{X}} \right\rangle = T^*\Phi\left(\bracks{P}{\Lambda_R}[P,X]\right)
\]
Using Lemma \ref{lem:helper-cotangent}, we now see:
\begin{align*}
    \bracks{d\Lambda_R}{P} &= \sum_{i=1,2} \bracks{d[Y_i,X]}{P}\nu_i \\
    &= -\sum_{i=1,2}\bracks{Y_i' + [Y_i,P]}{dX}\nu_i \\
    &= \sum_{i=1,2} T^*\Phi\left([Y_i' + [Y_i,P],X]\nu_i\right)
\end{align*}
So 
\[
-d\H = T^*\Phi\left( [L,X] + \bracks{P}{\Lambda_R}[P,X] + [Y_1' + [Y_1,P],X]\nu_1 + [Y_2' + [Y_2,P],X]\nu_2 \right)
\]

Now we have:
\begin{adjustwidth}{-1cm}{0pt}
\begin{align*}
    &T^*\Phi\left(\text{ODE for } \Lambda_R\right) \\ &= T^*\Phi(\left( \Lambda_R' +[\Lambda_R,P] + [L,X] + \bracks{P}{\Lambda_R}[P,X] \right) \\
    &= T^*\Phi\left( \sum_{i=1,2}([Y_i',X]\nu_i +  [Y_i,X']\nu_i + [Y_i,X]\nu_i'+  [[Y_i,X],P]\nu_i) + [L,X] + \bracks{P}{\Lambda_R}[P,X]\right) \\
    &= T^*\Phi\left(\sum_{i=1,2}[Y_i'+[Y_i,P],X]\nu_i + [Y_i,X]\nu_i'\right) + T^*\Phi\left([L,X] + \bracks{P}{\Lambda_R}[P,X] \right) \\
    &= -d\H -\nu'
\end{align*}
\end{adjustwidth}
Where we used the Jacobi identity. Thus, the ODE for $\Lambda_R$ is satisfied if and only if the ODE for $\nu$ is satisfied.
\end{proof}
\subsubsection{Transversality conditions}\label{sec:transversality}
For free terminal time optimal control problems, the PMP specifies \emph{transversality} conditions which are endpoint conditions which the extremals need to satisfy. 

On a manifold $M$, if $(q(t),u(t))$ is the projection of the lifted extremal trajectory $\xi(t) \in T^*_{q(t)}M$ in the cotangent bundle, then transversality requires that
\[ 
\bracks{\xi(t_f)}{v}_{*} = 0 \quad v \in T_{q(t_f)}M_f 
\]
where $T_{q(t_f)}M_f$ is the tangent space at $q(t_f)$ of the \emph{final submanifold} $M_f$ (which is the submanifold in which the terminal point $q(t_f)$ is allowed to vary in). 

More generally, if the initial point $q(0)$ is also allowed to vary in an \emph{initial submanifold} $M_0$, then transversality requires that the lifted extremal trajectory $\xi(t)$ in the cotangent bundle annihilates the vectors in the tangent spaces of the initial and terminal manifolds at the initial and terminal times.
\begin{align*}
\bracks{\xi(0)}{v}_{*} &= 0 \quad v \in T_{q(0)}M_0 \\
\bracks{\xi(t_f)}{v}_{*} &= 0 \quad v \in T_{q(t_f)}M_f 
\end{align*} 

Letting $R = \exp(J\pi/3)$, we have that our system is periodic up to a rotation by $R$. In this case, the initial and terminal submanifolds coincide after rotation by $R$ and so, transversality simply means that the lifted extremal trajectories are periodic functions modulo rotation by $R$. Within the Lie algebra $\sl(\R)$, rotation by $R$ amounts to $Y \mapsto \Ad_{R} Y = RYR^{-1}$.

For our system, we have $\xi(t) = \left(g(t),X(t),\Lambda_1(t),\Lambda_R(t)\right)$. We have already seen endpoint conditions for $g$ and $X$ in Section \ref{sec:X-init-term-conds}.  Collecting everything we get:
\begin{align*}
g(t_f) &= R \\
X(t_f) &= R^{-1} X(0) R \\
\Lambda_1(t_f) &= R^{-1}\Lambda_1(0)R \\
\Lambda_R(t_f) &= R^{-1}\Lambda_R(0)R
\end{align*}

\begin{remark}\normalfont\leavevmode
\begin{itemize}
    \item Note that the transversality conditions in the Reinhardt problem require the terminal time to satisfy 
    \[
    X(t_f) = R^{-1}X_0R
    \]
    where $R = \exp(J\pi/3)$. Thus, if we extend time, every optimal solution can be made into a periodic one:
    \[
    X(3t_f) = X_0
    \]
    The same requirement also holds of $g, \Lambda_1$ and $\Lambda_R$. Thus, the transversality requirement follows from a classification of periodic solutions of the lifted trajectories in the cotangent bundle.
    \item It is important to note that the transversality condition for $\Lambda_1$ can be dropped since we know its general solution and it follows from the transversality condition for $g$. See Corollary \ref{cor:lam1-gen-sol}. 
\end{itemize}
\end{remark}
\subsubsection{Summary}
At this point, we have the state and costate equations fully stated in both the coadjoint orbit picture and the upper half-plane picture. We also have the transversality conditions stated. We collect them as: 

\begin{problem}[State-Costate Equations]\label{pbm:state-costate-reinhardt}
The Reinhardt control problem which was described in Problem \ref{pbm:reinhardt-optimal-control} is an optimal control problem on the manifold $M = \SL(\R) \times \sl(\R)$. This problem has the the Hamiltonian:
\[
\H(\Lambda_1,\Lambda_R,X;Z_u) = \left\langle \Lambda_1 - \frac{3}{2}\lambda_{cost} J, X \right\rangle  - \frac{\langle \Lambda_R,Z_u\rangle }{\langle X, Z_u \rangle}
\]
and we lift the state trajectories described by $(g,X)$ in Problem \ref{pbm:plane-optimal-control-problem} to the following Hamiltonian system on the cotangent bundle of $M$. 

\begin{align*}
g'&=gX \\
X' &= [P,X] \\
\Lambda_1' &= [\Lambda_1, X] \\
\Lambda_R' &= \left(\left[P,\Lambda_R\right] - \frac{\langle \Lambda_R,Z_u \rangle}{\langle X, Z_u\rangle} \left[P,X\right]\right) + \left[ -\Lambda_1 + \frac{3}{2}\lambda_{cost} J ,X\right] 
\end{align*}
where $P = \frac{Z_u}{\langle Z_u,X\rangle}$. The transversality conditions described in Section \ref{sec:transversality} hold and the pair $(X,\Lambda_R)$ have equivalent dynamics in the upper half-plane described in Theorems \ref{thm:X-dynamics} and \ref{thm:lamR-dynamics}. The control dependent part of the Hamiltonian also has an expression in the upper half-plane picture as Theorem \ref{thm:plane-hamiltonian}.
\end{problem}

\subsection{Particular Solutions}
We have a well-defined control problem on the cotangent bundle and we now turn to establishing particular solutions of this system. We start with the easiest case, where the control is constant. 

\subsubsection{Solutions for constant control}

In this section, we keep the control matrix $Z_u$ constant and derive general solutions to the state and costate equations. Recall that this means that $\bracks{Z_u}{X}$ is a constant of motion (by equation \eqref{eq:X-Z0-constant}), and hence $P = \frac{Z_u}{\bracks{Z_u}{X}}$ is also constant. So, for $g(0) = I_2$ and $X(0)=X_0$ (or, equivalently $z(0)=z_0$), denote $P_0 := \frac{Z_u}{\bracks{Z_u}{X_0}}$. The general solutions for $(g,X)$ are:
\begin{align} 
g(t) &= \exp(t(X_0 + P_0))\exp(-t P_0) \label{eq:const-control-g} \\
z(t) &= \exp\left(tP_0 \right) \cdot z_0 \iff X(t) = \exp(tP_0)X_0\exp(-tP_0) = \Ad_{\exp(tP)}X_0 \label{eq:const-control-z-X}
\end{align}

As previously noted in Corollary \ref{cor:lam1-gen-sol}, the general solution for $\Lambda_1$ is:
\[
\Lambda_1(t) = \Ad_{g(t)^{-1}}(\Lambda_1(0)) = g(t)^{-1}\Lambda_1(0)g(t)
\]

We also have a rather complicated expression for the general solution for $\Lambda_R$. 
\[
\Lambda_R(t) = \Ad_{\exp(tP)}S(t)
\]
where 
\begin{align*}
    S(t) &:= \Lambda_R(0) - [L_P(t) + \ell(t)P,X_0] \\
    \ell(t) &:= \int_0^t \bracks{P}{\Lambda_R(0) - [L_P(s),X_0]}ds \\
    L_P(t) &:= \int_0^t \Ad_{\exp(-(X_0+P)s)}\Lambda_1(0) - \frac{3}{2}\lambda_{cost}\Ad_{\exp(-Ps)}J ds
\end{align*}
Note that if $S_0$ is the particular solution for $\Lambda_R(0) = 0$, then the general solution adds an affine term:
\[
S(t) = S_0(t) + \Lambda_R(0) + t\bracks{P}{\Lambda_R(0)}[P,X_0]
\]
This can be verified by differentiating. 

\subsubsection{Singular solutions and Singular locus}\label{sec:singular-locus}
We have stated earlier that the optimal control matrix $Z_u^*$ is implicitly determined by the Hamiltonian maximization condition in equation \eqref{eq:max-ham}. Singular solutions arise when this maximization condition fails to produce a unique candidate for the matrix $Z_u^*$ over an entire time interval. Throughout this section, we denote $J = \mattwo{0}{-1}{1}{0}$ as usual. The following theorem appears in the upper half-plane picture in Hales~\cite{hales2017reinhardt}.

\begin{theorem}
If a lifted extremal has $\Lambda_R$ vanishing identically on an interval $[t_1,t_2]$, then this induces the control to switch to the center of the control set for $t \in (t_1,t_2)$. This means that $Z_u^*(t) = \frac{1}{3}J$ in this interval and this determines an arc of the circle in $\Kccs$ as a singular extremal.
\end{theorem}
\begin{proof}
Assume that, along an extremal curve, for all $t \in [t_1,t_2]$ we have $\Lambda_R(t) \equiv 0$. Note that this means that equation \eqref{eq:maximized-hamiltonian} does not involve the control matrix $Z_u^*$ and so the maximization fails to uniquely determine the control matrix in this interval. Thus, the lifted extremal in this interval is singular (according to Definition \ref{def:normal-abonormal-singular-extremals}). 

In this interval we also have $\Lambda_R'(t) = 0$ on $(t_1,t_2)$. 
The maximum principle states that the Hamiltonian vanishes identically along the lifted extremal (See Section \ref{sec:PMP}). Thus, we have:
\[
\H(\Lambda_1,\Lambda_R,X; Z_u) = \langle \Lambda_1 - \frac{3}{2}\lambda J, X \rangle - \frac{\langle \Lambda_R,Z_u\rangle }{\langle X, Z_u \rangle} =  \langle \Lambda_1 - \frac{3}{2}\lambda J, X \rangle = 0 \quad \text{ for } t\in [t_1,t_2]
\]
where we abbreviate $\lambda_{cost} = \lambda$. By Lemma \ref{lem:sl2-lemmas},(3) this implies:
\begin{align}
\left(\Lambda_1 - \frac{3}{2}\lambda J\right)X +  X\left(\Lambda_1 - \frac{3}{2}\lambda J\right) = 0 
\end{align} 
Now, since $ \Lambda_R'(t) = 0$, we have:
\[ 
\Lambda_R'(t) =  \left[ -\Lambda_1 + \frac{3}{2}\lambda J ,X\right] = 0 \quad \text{ for } t \in (t_1,t_2).
\]
Putting these together this means that $\left(\Lambda_1 - \frac{3}{2}\lambda J\right)X = 0$. Since $X^2$ is a non-zero multiple of the identity, this proves that $\Lambda_1 \equiv \frac{3}{2}\lambda J$ on $(t_1,t_2)$. This proves that $\Lambda_1' = [\Lambda_1,X] = \frac{3}{2}\lambda\left[ J,X \right]= 0$ on $(t_1,t_2)$.

At this point, note that we must have $\lambda \ne 0$ for otherwise we will have $\lambda = 0$, $\Lambda_R \equiv 0$ and $\Lambda_1 \equiv 0$ on this interval, contradicting the non-vanishing of the costate variables in the maximum principle (see Section \ref{sec:PMP}). 
Thus, we have $[J,X] = 0$. Projecting down to the upper half-plane (\textit{i.e.,} making $X = X_z$ via the isomorphism $\Phi$ in Lemma \ref{lem:def-phi}) and solving this equation we get $x = 0, y = 1$ and so $X \equiv J$ on $t \in (t_1,t_2)$. The unique control which achieves this is $(u_0,u_1,u_2) = (1/3,1/3,1/3)$ which is the common centroid of the control sets $U_I,U_T,U_C$.

Now, the curve $g(t) = \exp(Jt)$ satisfies $g'=gX=gJ$ and this is the curve in $\SL(\R)$ which gives rise to the circle in the packing plane as a convex centrally symmetric disc, assuming $g(0) = I_2$. 
\end{proof}

Thus, in the \textit{singular locus} of the cotangent space, we must necessarily have:
\begin{align*}
g(t) &= g_0\exp(Jt) \qquad \quad X(t) = J \leftrightarrow z(t) = i \\
\Lambda_1(t) &= \frac{3}{2}\lambda_{cost}J \qquad \quad \Lambda_R(t) = 0 \\
\end{align*}
where $g(0) = g_0$. 

\begin{definition}[Singular Locus]
The region of the extended state space $T^*(\SL(\R) \times \sl(\R))$ given by 
\[
\Lambda_{sing} = \left\{\left(g_0,\frac{3}{2}\lambda_{cost}J,J,0\right)~|~g_0\in\SL(\R)\right\} \subset T^*(\SL(\R) \times \sl(\R))
\]
is called the \emph{singular locus}. In the star domain of the upper half-plane, the singular locus corresponds to the point $z=i \in \star$. 
\end{definition}

\begin{remark}\normalfont
Note that $\Lambda_{sing}$ gives the initial conditions corresponding to the circle in $\Kccs$ (which has $g_0 = I_2 \in \SL(\R)$) up to a transformation in $\SL(\R)$. 
\end{remark}

When we are using the simplex control set $U_T$, we might have a case where $\Lambda_R \ne 0$, but where the Hamiltonian is maximized along an edge of the simplex. A long argument in Hales~\cite{hales2017reinhardt} proves:
\begin{theorem}[Hales~\cite{hales2017reinhardt}]\label{thm:no-singular-arcs}
The global optimizer of the Reinhardt Optimal Control problem does not contain any singular arcs \emph{i.e.,} it does not have the control switching to $(1/3,1/3,1/3)$ or an edge of the control simplex for any interval of time.  
\end{theorem}

In other words, the globally optimal solution cannot stay at the singular locus for any finite interval of time. 

\subsubsection{Smoothed \texorpdfstring{$(6k+2)$}{6k+2}-gons}\label{sec:6k+2-gons}
In the 2017 article, Hales constructs a family of Pontryagin extremals of the control problem, all given by a bang-bang control. In this section we shall briefly report the statements of the main results of that article. 

\begin{definition}[Bang-bang control]
A control function is said to be \emph{bang-bang} if it is characterized by discontinuous switching between the extremes of the control set. 
\end{definition}

\begin{lemma}\label{lem:exp-of-sl2}
If $X \in \sl(\R)$ and $\det(X)\ne 0$, then $\exp(X) \in \SL(\R)$ and also (for nonzero $X$),
\begin{equation}\label{eq:matrix-exp-sl2}
\exp(t X) = (\cosh{(t\sqrt{-r})}) I_2 + \frac{\sinh{(t\sqrt{-r})}}{\sqrt{-r}}X
\end{equation}
where $r = \det(X)$.
\end{lemma}
\begin{proof}
If $X\in \sl(\R)$, we have 
\[
\det(\exp (X)) = e^{\tr (X)} = 1
\]
since $\tr(X) = 0$. This shows that $\exp(X) \in \SL(\R)$. 
By the Cayley-Hamilton theorem, we have:
\[ 
X^2 + r I = 0
\]
By induction using the above, we get an appropriate formula for $X^{2n}$ and $X^{2n + 1}$. Substituting these into the series expansion $\exp (X) := \sum_{n=0}^\infty \frac{X^n}{n!}$ we get
\[ 
\exp(X) = (\cosh{\sqrt{-r}}) I_2 + \frac{\sinh{\sqrt{-r}}}{\sqrt{-r}}X 
\]
which gives the required. 
\end{proof}

Let us label the vertices of the control set $U_T \subset \R^3$ as $e_0 = (1,0,0)$, $e_1 = (0,1,0)$ and $e_2 = (0,0,1)$. 

\begin{lemma}[Hales~\cite{hales2017reinhardt}]\label{lem:const-control-trajectories}
If the control stays at one of the vertices $e_1,e_2 \in U_T$ of the simplex control set, the corresponding trajectories in the star domain are straight lines. If the control is $e_0$, the trajectory (in the upper half-plane) is an arc of a circle. 
\end{lemma}
\begin{proof}
By equation \eqref{eq:const-control-z-X}, we have that the trajectory in the upper half-plane $z(t)$ is given by fractional linear transformations. So, with the $e_1$ or $e_2$ control, and by equation \eqref{eq:Z0}, we get $c_{21}(Z_u) = 0$, which, along with Lemma \ref{lem:exp-of-sl2}, proves that
$c_{21}(\exp(tP_0)) = 0$. So, the linear fractional transformation by the matrix $\exp(tP_0)$ is actually an affine transformation. 

The $e_0$ trajectory is obtained from $e_1$ by rotation by the matrix $R = \exp(J\pi/3)$. The conclusion follows since we know that fractional linear transformations map straight lines to circles.
\end{proof}

\begin{lemma}[Hales~\cite{hales2017reinhardt}]\label{lem:control-face-lemma}
Assume that our control set is a compact, convex set in the affine plane $\{(u_0,u_1,u_2)~|~\sum u_i =1\}$. The controls maximizing the Hamiltonian in equation \eqref{eq:max-ham} are constrained to a face\footnote{A nonempty convex subset $F$ of a convex set $C$ is called a \emph{face} if $\alpha x + (1 - \alpha) y \in C$ for $x,y \in C$ and some $0 < \alpha < 1$ implies $x,y \in C$.} of the control set. 
\end{lemma}
\begin{proof}
We see that the control-dependent part --- denoted $\H_\h$ --- of the Hamiltonian in equation \eqref{eq:max-ham} is a ratio of two linear functions. 
If $su_1+(1-s)u_2$ for $0\le s \le 1$ is a segment within the control set, then the control-dependent part becomes a function of $s$ as: 
\[
\H_\h(s) = \frac{as+b}{cs+d}
\]
the derivative of this expression has fixed sign. Thus, the Hamiltonian is monotonic along every segment within the control set and so the set of maximizers is constrained to a face of the convex control set. 
\end{proof}

Thus, if we consider the control set $U_T$, the set of maximizers are either the entire set (which coincides with the singular case), or an edge or a vertex. 

\begin{figure}[htbp]
        \centering
        \includegraphics[scale=0.40]{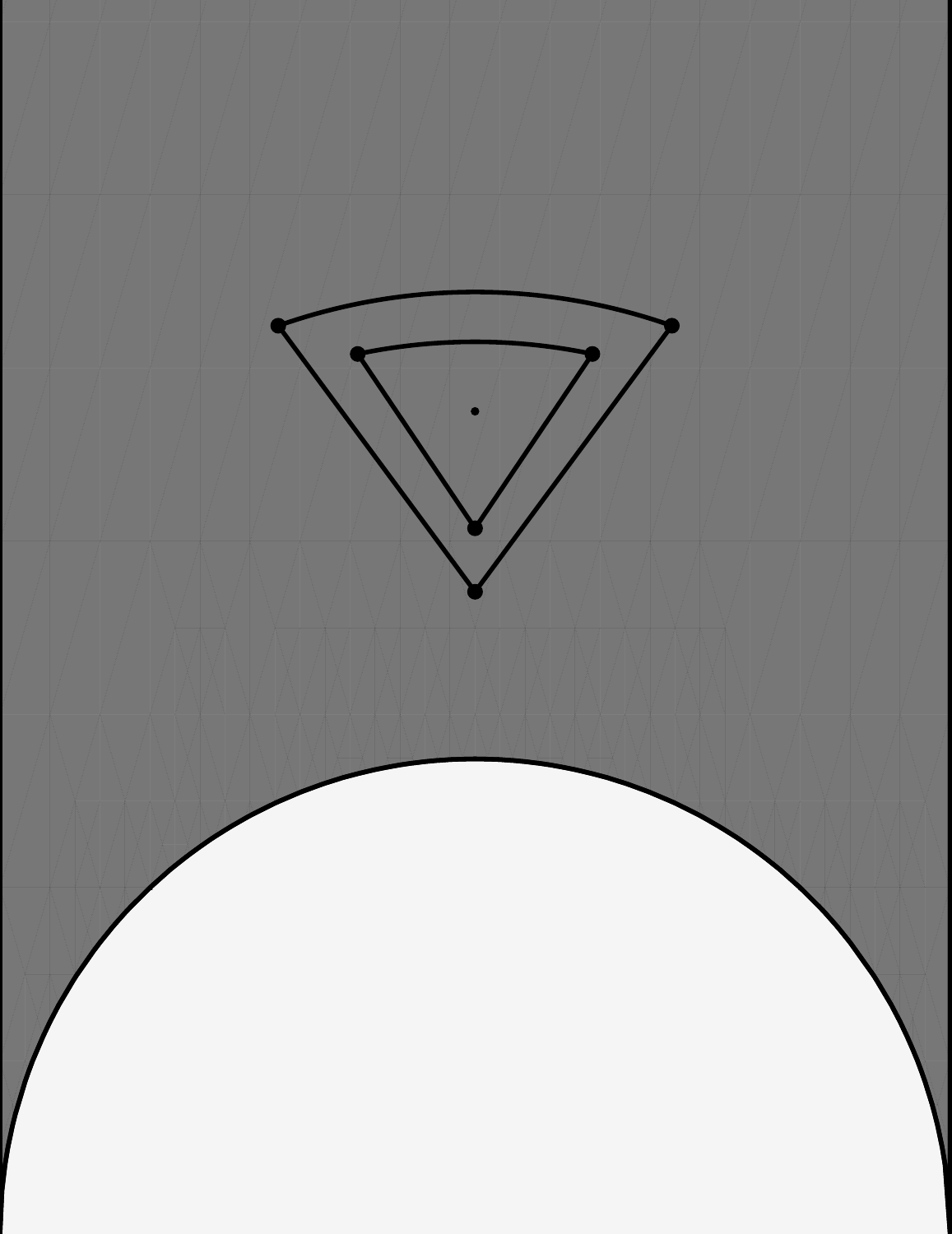}
        \caption{Trajectories of the smoothed 8-gon and 14-gon in \texorpdfstring{\(\star\)}{h*}}
        \label{fig:8-14-gons}
    \end{figure}

The smoothed octagon is part of a family of smoothed \((6k+2)\)-gons, which are all given by bang-bang controls. The trajectory of the smoothed 8-gon and the 14-gon are shown in Figure \ref{fig:8-14-gons}. These are equilateral triangles in $\h$ which shrink toward the singular locus $i \in \star$ as $k \to \infty$. 

\begin{theorem}[Hales~\cite{hales2017reinhardt}]\leavevmode\label{thm:bang-bang}
\begin{itemize}
    \item For each positive integer $k$, the smoothed $6k+2$-gon is a Pontryagin extremal given by a bang-bang control. 
    \item Pontryagin extremals which do not pass through the singular locus $\Lambda_{sing}$ are given by a bang-bang control with finitely many switches.
\end{itemize}
\end{theorem}

\section{Circular control set}\label{sec:circular-control}
 Up until now, we have been exclusively working with the simplex control set $U_T$. Now, we change our control set to be $U_C$, the circumscribing disk of $U_T$ as described in Definition \ref{def:circ-control-set}. This change is motivated by the following theorem. 
 
\begin{theorem}[Hales~\cite{hales2017reinhardt}]\label{thm:symmetries}
Let $\H$ be the Hamiltonian of the Reinhardt optimal control problem. Assuming that the control set is closed under rotations (i.e., if $\rho \in \SO$ is a rotation and $Z_u$ is in the control set, then so is $\Ad_\rho Z_u$), the Hamiltonian $\H$ is invariant under the action of the subgroup $\SO$ of $\SL(\R)$. 
\end{theorem}
\begin{proof}
The Hamiltonian depends on the quantities $X,\Lambda_1,\Lambda_R$ and the control matrix $Z_u$. Ignoring $Z_u$ for the time being, if $\rho \in \SO$ is an arbitrary rotation, then let us see how these quantities transform. $\rho$ acts on trajectories in $\h$ by fractional linear transformations and so, equivalently, we have $X \mapsto \Ad_\rho X= \rho X \rho^{-1} =: \tilde{X}$ in the adjoint orbit picture. Now, $\Lambda_1$ transforms as $\Lambda_1 \mapsto \rho \Lambda_1 \rho^{-1} =: \tilde{\Lambda_1}$ since these transformed quantities satisfy the same ODE:
\[\tilde{\Lambda_1}' = \Ad_\rho \Lambda_1' = \Ad_\rho [\Lambda_1,X] = [\Ad_\rho\Lambda_1,\Ad_\rho X] = [\tilde{\Lambda_1},\tilde{X}]\] 

Similarly, we can also see that $\Lambda_R$ transforms as $\Lambda_R \mapsto \rho \Lambda_R \rho^{-1}$.
Now, the control set transforms as $Z_u \mapsto \rho Z_u \rho^{-1}$, which may, in general, fall outside the control set given by $U_T$. So, if we modify the control set so that it does not, a simple computation now using the expression for the control-dependent Hamiltonian in equation \eqref{eq:full-hamiltonian} shows that it is unchanged under these transformations by $\rho$. 
\end{proof}

\begin{remark}\leavevmode\normalfont
Following the discussion at the end of Section \ref{sec:control-sets}, the control set $U_T$ is only symmetric up to discrete $\Z/3\Z$-rotations and not under general $\SO$ rotations. 
\end{remark}

\subsection{Conserved quantity for the circular control set.}\label{sec:angular-momentum}
Recalling the discussion at the end of Section \ref{sec:control-sets}, we note that we enlarge the control set $U_T$ to $U_C$ to have continuous $\SO$-symmetry in order to manufacture a conserved quantity for the dynamics. This is achieved by Noether's theorem. 

The version of Noether's theorem which we use for optimal control is the one described by Sussmann~\cite{sussman}. Before recalling the statement, we shall begin with a few definitions. 
We let $\mathcal{P} = (M,U,f,\phi)$ denote a generic optimal control system satisfying the regularity conditions as in Section \ref{sec:PMP}.  Also, we assume that the cost functional $\phi$ is independent of the control as in our case. In what follows, we shall denote the vector field $f(u,q) \in T_{q}M$ for $u\in U$ by $f_u(q)$. The definitions below are all from Sussmann~\cite{sussman}. 

\begin{definition}[Symmetry of a control system]
Let $G$ be a Lie group with Lie algebra $\mathfrak{g}$ and let $\mathcal{P} = (M,U,f,\phi)$ be an optimal control system. A \emph{symmetry} of this optimal control system is a diffeomorphism $\sigma : V_1 \to V_2$ where $V_1,V_2 \subset M$ are open such that for every $u \in U$ there exist $u_1,u_2 \in U$ for which $d\sigma(x)(f_u(x)) = f_{u_1}(\sigma(x))$ and $d\sigma(x)(f_{u_2}(x)) = f_{u}(\sigma(x))$ for all $x \in V_1$.
\end{definition}

\begin{definition}[Infinitesimal Group of Symmetries]
An \emph{infinitesimal group of symmetries} of a control system $\mathcal{P}$ is a smooth action $\tau : \mathfrak{g} \to \Gamma^{\infty}(TM)$, which assigns to every Lie algebra element $\xi \in \mathfrak{g}$ a smooth vector field on $M$, such that every diffeomorphism $\exp(t\tau(\xi))$ is a symmetry of $\mathcal{P}$. 
\end{definition}

\begin{definition}[Momentum Map]
To an infinitesimal group of symmetries $\tau : \mathfrak{g} \to \Gamma^{\infty}(TM)$ of an optimal control system $\mathcal{P}$, we associate the \emph{momentum map} $\mathbf{J}^\tau : T^*M \to \mathfrak{g}^*$ given by
\[
\mathbf{J}^\tau(x,z)(\xi) = \bracks{z}{\tau(\xi)(x)}_*
\]
where $x \in M,\ z\in T_x^*M$ and $\tau(\xi)(x) \in T_x M$. 
\end{definition}

See Abraham~\&~Marsden~\cite{abraham2008foundations} for the general theory of momentum maps in symplectic geometry.

\begin{theorem}[Sussmann-Noether Theorem]
Assume that $\mathcal{P} = (M,U,f,\phi)$ is an optimal control system as above, let $\mathfrak{g}$ be the Lie algebra of the Lie group $G$ and let $\tau : \mathfrak{g} \to \Gamma^{\infty}(TM)$ be an infinitesimal group of symmetries of $\mathcal{P}$. Let $\Xi$ be a lifted controlled trajectory (guaranteed by the PMP) in $T^*M$. Then the function $\mathbf{J}^\tau : T^*M \to \mathfrak{g}^*$ is constant along the trajectory $\Xi$.
\end{theorem}

\begin{problem}[Circular Control Problem]
We start with Problem \ref{pbm:state-costate-reinhardt} and enlarge the control set $U_T$ to its circumscribing disk $U_C$. This control problem is called the \emph{Circular Control Problem}.  
\end{problem}
Let us now apply the Sussmann-Noether theorem to the circular control problem.

\begin{theorem}[Angular Momentum]\label{thm:angular-momentum}
Because of the $\SO$ symmetry, we have that $\bracks{J}{\Lambda_1 + \Lambda_R}$ is conserved along the optimal trajectory.
\end{theorem}
\begin{proof}
Our optimal control system consists of $(g,X,\Lambda_1,\Lambda_2) \in T^*(T\SL(\R))$. 
Note that our Lie group here is $\SO$ and hence its Lie algebra is one-dimensional $\mathfrak{so}_2(\R) = J\R$, where $J = \mattwo{0}{-1}{1}{0}$ is the infinitesimal generator of rotations. This infinitesimal symmetry gives rise to rotations  $\exp(J\theta) \in \SO$, which in turn give rise to the following action on our manifold: 
\begin{align*}
\SO \times T^*(T\SL(\R)) &\to T^*(T\SL(\R)) \\
(\exp(J\theta), (g,X,\Lambda_1,\Lambda_2)) &\mapsto (\exp(J\theta)g\exp(-J\theta), \Ad_{\exp(J\theta)}X, \Ad_{\exp(J\theta)}\Lambda_1, \Ad_{\exp(J\theta)}\Lambda_2)
\end{align*}
where the action on $g$ is by inner automorphisms and the rest are given by the adjoint action on each copy of $\sl(\R)$. (Note that, throughout, we make the identification $\sl(\R)^* \cong \sl(\R)$ via the nondegenerate trace form.)
These are \emph{symmetries} by the proof of Theorem \ref{thm:symmetries} and also since the rotation action on the control matrix $Z_u$ is 
\[
Z_u \mapsto \Ad_{\exp(J\theta)} Z_u \in U_C
\]
Now, we have to calculate the momentum map corresponding to action. 
The momentum map is computed by the canonical pairing between the costate variables $(\Lambda_1,\Lambda_2) \in T_{(g,X)}^*(T\SL(\R))$ and tangent vectors
in $T_{(g,X)}(T\SL(\R))$ giving the infinitesimal rotation action. The first component of this tangent vector is given by:
\begin{align*}
\frac{d}{d\theta}\exp(J\theta)g\exp(-J\theta)\Bigr|_{\theta=0} &= \underbrace{Jg - gJ}_{\in T_g\SL(\R)} \\ 
&=\underbrace{g}_{\in \SL(\R)}\underbrace{(g^{-1}Jg - J)}_{\in \sl(\R)}
\end{align*}
thus, we can identify the first component with $g^{-1}Jg - J$ in the Lie algebra $\sl(\R)$. The second component is given by 
\[
\frac{d}{d\theta}\Ad_{\exp(J\theta)}X = \ad_J X = [J,X]
\]
which is already in the Lie algebra $\sl(\R)$. 
So, putting all this together, we get the momentum map:
\begin{align*}
\mathbf{J}\left((g,\Lambda_1,X,\Lambda_2) \right) &= \bracks{(\Lambda_1,\Lambda_2)}{(\Ad_{g^{-1}}J - J, [J,X])}_{*} \\
&= \bracks{\Lambda_1}{\Ad_{g^{-1}}J - J} + \bracks{\Lambda_2}{[J,X]} \\     
&= \bracks{\Lambda_1}{\Ad_{g^{-1}}J} -\bracks{\Lambda_1}{J} - \bracks{[\Lambda_2,X]}{J} \\
&= \bracks{\Lambda_1}{\Ad_{g^{-1}}J} - \bracks{J}{\Lambda_1 + \Lambda_R}
\end{align*}
Where, as usual, $\bracks{\cdot}{\cdot}_*$ denotes the natural pairing between a vector space and its dual. 
Note however, from Corollary \ref{cor:lam1-gen-sol}, that $\Lambda_1(t) =  g^{-1}\Lambda_1(0)g(t)$ and so $\bracks{\Lambda_1}{\Ad_{g^{-1}}J}$ is a constant. 

Thus, by the Sussmann-Noether theorem $\bracks{J}{\Lambda_1 + \Lambda_R}$ is a constant of motion along the optimal trajectory of the circular control problem. 
\end{proof}

\begin{remark}\normalfont \leavevmode
\begin{itemize}
    \item Since it is the conserved quantitiy arising from a rotational symmetry, $\bracks{J}{\Lambda_1 + \Lambda_R}$ will be called the \emph{angular momentum}.
    \item The spurious constant in the expression for the momentum map is a consequence of the action of $\SO$ on $\SL(\R)$ by inner automorphisms. This means that we could also modify the action to be $g \mapsto g \exp(-J\theta)$ to get a valid conserved quantity.  
    \item We also get the same conserved quantity for a control set $U_r$ which is a disk of arbitrary radius $r$ provided $r^2 > 1/3$. 
    (See Definition \ref{def:circ-control-set}.) 
\end{itemize}
\end{remark}

Conserved quantities are very useful since they introduce constraints and cut down the dimension of the problem. But an immediate application is that they give us a constraint on the optimal control matrix in terms of the state-costate variables, which is what we aim to derive. But before we do, we make a quick detour to understand the structure of the control sets when viewed in $\su$.

\subsection{Control sets in \texorpdfstring{$\su$}{su(1,1)}.}\label{sec:su11-control-sets}
From Section \ref{sec:sl2-exceptional}, we know that the Lie algebra $\su$ and $\sl(\R)$ are isomorphic. This isomorphism is given by the Cayley transform:
\[
C^{-1} \sl(\R) C = \su
\]

This means that the control matrix $Z_u$, given in equation \eqref{eq:Z0}, viewed in $\su$ becomes:
\begin{equation}\label{eq:Z0-su11}
\frac{1}{3}\mattwo{-i}{2(u_0 + \zeta u_1 + \zeta^2 u_2)}{2(u_0 + \zeta^2 u_1 + \zeta u_2)}{i}
\end{equation}

Now, let $\partial U_r$ denote the boundary of the disk control set $U_r$ of radius $r$ in $\R^3$:
\[
U_r := \left\{(u_0,u_1,u_2)~|~u_0 + u_1 + u_2 = 1,~u_0^2+ u_1^2 + u_2^2 \le r^2 \right\}
\]
We have the following lemma:

\begin{lemma}\label{lem:circle-plane-bijection}
The set $\partial U_r$ and the circle $\left\{z\in \mathbb{C}~|~|z|^2 = \left(\frac{3r^2-1}{2}\right)\right\}$ are in bijection.
\end{lemma}
\begin{proof}
Consider the affine plane $\Pi = \left\{(u_0,u_1,u_2)~|~u_0 + u_1 + u_2 = 1 \right\}$ and consider the map $L :  \Pi \to \mathbb{C} $, defined by $L(u_0,u_1,u_2) := (u_0 + \zeta u_1 + \zeta^2 u_2)$. This map is the restriction of the linear map 
\[
(u_0,u_1,u_2) \mapsto \left(\matthree{1}{1}{1}{1}{\zeta}{\zeta^2}{1}{\zeta^2}{\zeta}\right)\left(\begin{array}{ccc}
u_0 \\
u_1 \\
u_2
\end{array}\right) \in \{(1,z,\bar{z}~|~z\in \mathbb{C}\}
\]
which has non-zero determinant and so $L$ is an isomorphism of affine planes. This isomorphism restricts to a bijection between the circles $\partial U_r$ and  $\{z\in \mathbb{C}~|~|z|^2 = \left(\frac{3r^2-1}{2}\right)\}$ since 
\begin{align*}
|(u_0 + \zeta u_1 + \zeta^2 u_2)|^2 &= \left(\frac{3}{2}(u_0^2 + u_1^2 + u_2^2) - \frac{(u_0+u_1+u_2)^2}{2}\right) \\ 
&= \left(\frac{3r^2 - 1}{2} \right).
\end{align*}
\end{proof}

The above lemma shows that if the control set is $U_r$, then we can take $Z_u$ in general to be:
\begin{equation}\label{eq:general-Z0}
    Z_u =\mattwo{-i \alpha}{\beta z}{\beta \bar{z}}{i \alpha} \in \su
\end{equation}
where $\alpha\in \R$ and $\beta \in \R,~ z \in \mathbb{C}$ with $|z| = 1$ are the polar coordinates of the entry $c_{12}(Z_u)$. With this notation, $\alpha = \frac{1}{3}$ and $\beta = \frac{2}{3}|u_0 + \zeta u_1 + \zeta^2 u_2|$. Using all this, we get:
\begin{equation}\label{eq:det-Z0}
    \det(Z_u) = (\alpha^2 - \beta^2) = \frac{1 - 4\left(\frac{3r^2-1}{2}\right)}{9} = \frac{1 - 2r^2}{3}
\end{equation}

Table \ref{tab:control-sets} shows the square of the radii $r^2$ of the circumscribing disk, inscribed disk and center of the control simplex $U_T$ and their corresponding radii when viewed as disks in the complex plane (following Lemma \ref{lem:circle-plane-bijection}). \newline

\begin{table}[htbp]
    \centering
    \begin{tabular}{|c|c|c|c|c|}
        \hline
        $\det(Z_u)$ & Relation to simplex $U_T$ (in $\R^3$) & $r^2$ & Radius of disk in $\mathbb{C}$ = $(3r^2-1)/2$ & $(\alpha,\beta)$  \\
        \hline
        $- \frac{1}{3}$ & Circumscribing disk $U_C$ & 1 & 1 & (1/3,2/3)\\
        $0$ & Inscribed disk $U_I$ & $1/2$ & 1/4 & (1/3,1/3)\\
        $\frac{1}{9}$ &  Center of disk $(1/3,1/3,1/3)$& $1/3$ & 0 & (1/3,0)\\
        \hline
    \end{tabular}
    \caption{Various control sets and their radii.}
    \label{tab:control-sets}
\end{table}

\subsection{Quadratic equation for optimal control matrix.}\label{sec:quadratic-equation}
Henceforth, let us denote the optimal control matrix for the control problem in Problem \ref{pbm:state-costate-reinhardt} (with the control set $U_r$) as $Z_u^*$. We now hope to derive a constraint on $Z_u^*$ from angular momentum conservation. 

We begin with a few lemmas:

\begin{lemma}\label{lem:boundary-circle-control}
The set of maximizers of the Hamiltonian (considered as a function of the control) when the control set is the disk $U_r$ is either the entire disk or just a point on $\partial U_r$.
\end{lemma}
\begin{proof}
This is a direct corollary of Lemma \ref{lem:control-face-lemma}.
\end{proof}




\begin{theorem}
Let $Z_u^*$ be the optimal control matrix for the control problem in  Problem \ref{pbm:state-costate-reinhardt} (with the control set $U_r$). We then have\footnote{Note that the $\Jsu$ in this theorem is an element of $\su$ and so is given by \[\Jsu = C^{-1}\mattwo{0}{-1}{1}{0}C = \mattwo{-i}{0}{0}{i}.\]}:
\[
\left \langle \Jsu, \ad_{Z_u}\frac{\delta}{\delta Z_u}\frac{\langle \Lambda_R, Z_u\rangle}{\langle X, Z_u\rangle} \right\rangle = 0
\]
where $\frac{\delta}{\delta Z_u}$ is the functional derivative with respect to $Z_u$.
\end{theorem}
\begin{proof}
We give two proofs of this fact:

For the first proof, we differentiate the angular momentum expression (in Theorem \ref{thm:angular-momentum}) with respect to time, to get:
\begin{align}
    0 = \langle \Jsu, \Lambda_1' + \Lambda_R' \rangle &= \left\langle \Jsu, [\Lambda_1, X] + [P^*,\Lambda_R] + \left[-\Lambda_1 + \frac{3}{2}\lambda \Jsu, X \right] - \langle P^*,\Lambda_R \rangle [P^*,X] \right\rangle \nonumber \\ 
    &=\left\langle \Jsu, [P^*,\Lambda_R] - \langle P^*,\Lambda_R\rangle [P^*,X] \right\rangle  \label{eq:deriv-ang-momentum}\\
    &= \left \langle \Jsu, \ad_{Z_u^*}\frac{\delta}{\delta Z_u}\frac{\langle \Lambda_R, Z_u^*\rangle}{\langle X, Z_u^*\rangle} \right\rangle \qquad \text{(by Corollary \ref{cor:control-dep-ad-alternate})} \nonumber
\end{align}
where $P^* = Z_u^*/\bracks{Z_u^*}{X}$.

For the second proof, note that, from Lemma \ref{lem:boundary-circle-control}, the Hamiltonian is maximized at a point on $\partial U_r$. By the form of the control matrix $Z_u^*$, on the boundary of $U_r$ we get two constraints: 
\[ 
\bracks{\Jsu}{Z_u^*} = \mathrm{constant}, \qquad \bracks{Z_u}{Z_u^*} = \mathrm{constant}
\]
The Hamiltonian maximization of the PMP can be considered as a constrained maximization problem subject to the above two constraints. The functional derivatives of the above two constraints are $\Jsu$ and $2Z_u$ respectively. By Lagrange multipliers, we get that the derivative $\delta \H/\delta Z_u$ should lie in the plane containing the span of the derivatives of the constraints \emph{viz}., $\Jsu$ and $Z_u^*$. That is, $\ad_{Z_u^*}(\frac{\delta \H}{\delta Z_u})$ is in the span of $\ad_{Z_u^*}\Jsu$. Thus, 
\[
\bracks{\ad_{Z_u^*}\frac{\delta \H}{\delta Z_u}}{\Jsu} = 0
\]
which gives us the required. 
\end{proof}
We also have the following result, which shows that the angular momentum and the Hamiltonian are in involution with respect to the Poisson bracket $\{\cdot,\cdot\}_{ex}$ on the extended state space $T^*T\SL(\R)$. 
\begin{proposition}
The angular momentum $\bracks{\Jsu}{\Lambda_1 + \Lambda_R}$ Poisson commutes with the Hamiltonian on $T^*(T\SL(\R))$, provided the control set is rotationally invariant.
\end{proposition}
\begin{proof}
We use the Poisson bracket on the manifold $T^*(T\SL(\R))$ derived in Section \ref{sec:poisson-bracket}, which we recall here. If $F,G$ are two left-invariant smooth functions on $T^*T\SL(R)$, their extended space Poisson bracket is:
\[
\{F,G\}_{ex} := \left\langle \Lambda_1, \left[ \frac{\delta F}{\delta\Lambda_1}, \frac{\delta G}{\delta\Lambda_1} \right] \right\rangle
 + \bracks{\frac{\delta F}{\delta X}}{\frac{\delta G}{\delta \Lambda_2}} -  \bracks{\frac{\delta F}{\delta \Lambda_2}}{\frac{\delta G}{\delta X}}
\]
Now, let us set $\mathcal{A} = \bracks{\Jsu}{\Lambda_1 + \Lambda_R}$ to be the angular momentum. Let us now compute:
\begin{align*}
    \{\mathcal{A},\mathcal{H}\}_{ex} &= \frac{\bracks{\Jsu}{[\Lambda_2,[Z_u^*,X]]}}{\bracks{X}{Z_u^*}} + \frac{\bracks{\Jsu}{[X,[\Lambda_2,Z_u^*]]}}{\bracks{X}{Z_u^*}} + \frac{\bracks{\Jsu}{[X,Z_u^*]}\bracks{\Lambda_R}{Z_u^*}}{\bracks{X}{Z_u^*}^2} \\ 
    &= \frac{\bracks{\Jsu}{[[X,\Lambda_2],Z_u^*]}}{\bracks{X}{Z_u^*}} + \frac{\bracks{\Jsu}{[X,Z_u^*]}\bracks{\Lambda_R}{Z_u^*}}{\bracks{X}{Z_u^*}^2} \qquad \text{(by the Jacobi identity)} \\
    &= \left\langle \Jsu, [P^*,\Lambda_R] - \langle P^*,\Lambda_R\rangle [P^*,X] \right\rangle = 0 \quad \text{(by equation \eqref{eq:deriv-ang-momentum})}
\end{align*}
where $P^* = Z_u^*/\bracks{Z_u^*}{X}$.
\end{proof}

Simplifying equation \eqref{eq:deriv-ang-momentum} above, we get a more symmetric expression for the optimal control matrix, which is homogeneous in $Z_u^*$:
\begin{equation}\label{eq:symm}
    \langle \Jsu,[Z_u^*,X]\rangle\langle Z_u^*,\Lambda_R\rangle = \langle \Jsu,[Z_u^*,\Lambda_R]\rangle{\langle Z_u^*,X\rangle}
\end{equation}

This is the same as saying:
\begin{equation}
  \dettwo{ \langle [\Jsu,Z_u^*],X]\rangle}{\langle [\Jsu,Z_u^*],\Lambda_R\rangle}{\langle Z_u^*,X\rangle}{\langle Z_u^*,\Lambda_R\rangle} = 0
\end{equation}


\begin{proposition}
Angular momentum conservation gives the following constraint on the optimal control matrix:
\begin{equation}\label{eq:opt-control-constraint}
\langle Z_u^*,[\Lambda_R,X] \rangle = \frac{\langle Z_u^*, Z_u^* \rangle}{\langle Z_u^*,\Jsu \rangle}\langle \Jsu, [\Lambda_R,X]\rangle
\end{equation}
\end{proposition}
\begin{proof}
From \eqref{eq:symm} and Proposition \ref{prop:trace-quotient-sub} setting $A = [\Jsu,Z_u^*], B=Z_u^*,C =X , D =\Lambda_R$, we get:
\begin{equation}\label{eq:ang-mom-quotient}
\langle [\Jsu,Z_u^*],\Lambda_R \rangle\bracks {Z_u^*} {X} - \langle Z_u^*,\Lambda_R \rangle\bracks {[\Jsu,Z_u^*]} {X} = \frac{\left\langle [[\Jsu,Z_u^*],Z_u^*], [\Lambda_R, X] \right\rangle}{2} = 0
\end{equation}
Note that, we have:
\[
Z_u^* \Jsu + \Jsu Z_u^* = \langle Z_u^*, \Jsu \rangle I_2.
\]
Multiplying this by $Z_u^*$ on the left and rearranging gives:
\begin{equation}\label{eq:z0jz0}
 Z_u^* \Jsu Z_u^* = \det(Z_u^*) \Jsu + \langle Z_u^*, \Jsu \rangle Z_u^*.
\end{equation}
Now we have, by Proposition \ref{prop:brack-brack}:
\begin{align*}
[[\Jsu,Z_u^*],Z_u^*] &= -2 \det(Z_u^*) \Jsu - 2 Z_u^* \Jsu Z_u^* \\
&= -4 \det(Z_u^*)\Jsu - 2\langle Z_u^*, \Jsu \rangle Z_u^* \quad \text{ by \eqref{eq:z0jz0}}.
\end{align*}
Substituting this into \eqref{eq:ang-mom-quotient} we get
\[
2\det(Z_u^*) \langle \Jsu, [\Lambda_R, X] \rangle + \bracks {Z_u^*} {\Jsu} \langle Z_u^*, [\Lambda_R, X] \rangle =0
\]
which implies the required, once we use the fact that $-2\det(Z_u^*)= \langle Z_u^*, Z_u^* \rangle$.
\end{proof}

\begin{definition}[Augmented Determinant]
For $L = \mattwo{l_{11}}{l_{12}}{l_{21}}{-l_{11}} \in \sl$ and $\alpha, \beta \in \R$, define the augmented determinant $D(L,\alpha,\beta) := \alpha^2 l_{11}^2 + \beta^2 l_{12} l_{21}$. 
\end{definition}

\begin{proposition}\label{prop:control-quadratic}
The optimal control in the circular control case $U_r$ is determined by the roots of a quadratic equation. 
\end{proposition}
\begin{proof}
We know by Lemma \ref{lem:boundary-circle-control} that the optimal control matrix $Z_u^*$ should take values in the boundary of the disk. By equation \eqref{eq:general-Z0}, we can take 
\begin{equation}\label{eq:general-boundary-Z0}
Z_u^* = \left(\begin{matrix}
-i\alpha & \beta z \\ 
\frac{\beta}{z} & i \alpha
\end{matrix}\right)
\end{equation}
With this notation, \eqref{eq:opt-control-constraint} becomes:
\begin{equation}\label{eq:opt-control-expl}
    \langle Z_u^*,[\Lambda_R,X] \rangle = \frac{(\alpha^2 - \beta^2)}{\alpha}\langle \Jsu, [\Lambda_R,X]\rangle
\end{equation}
where $\alpha, \beta \in \R$ and $\alpha > 0$. 


Simplifying this, we get the following quadratic equation in $z$
\begin{equation}\label{eq:quadratic-equation}
\alpha[\Lambda_R,X]_{21} z^2 - 2i\beta [\Lambda_R,X]_{11}z + \alpha [\Lambda_R,X]_{12} = 0
\end{equation}
We then solve for $z$ to get two roots
\begin{equation}\label{eq:z-quad-roots}
z_{1,2} = i\left(\frac{\beta[\Lambda_R,X]_{11} \pm  \sqrt{D([\Lambda_R,X],\alpha,\beta)}}{\alpha[\Lambda_R,X]_{21}}\right)  
\end{equation}
where $D([\Lambda_R,X],\alpha,\beta) := \beta^2[\Lambda_R,X]_{11} + \alpha^2[\Lambda_R,X]_{21}[\Lambda_R,X]_{12}$ is the augmented determinant of $[\Lambda_R,X]$.
The root $z_1$ (with the positive sign in front of the square root) maximizes the Hamiltonian.
\end{proof}

This determines the optimal control explicitly as a function of the state-costate variables. We see how modifying the control set to be more symmetrical has resulted in a conserved quantity which has in turn given us valuable information about the optimal control.

\section{Inscribed Disk Control Set}\label{sec:inscribed-control}

In this section, we specialize to the case where the control set is the inscribed disk $U_I$ of the 2-simplex in $\R^3$. The parameters for this set are $\alpha = \beta = 1/3$. Our main result in this case is a classification of periodic abnormal extremals in Section \ref{sec:no-abnormal-extremals}.

\subsection{Simplifying the Costate}
Our first result in this case is a simplification of the ODE for $\Lambda_R$:

\begin{proposition}
In the case where the control set is the inscribed disk of the 2-simplex, the control dependent term in the ODE for $\Lambda_R$ vanishes:
\begin{equation}\label{eq:control}
[P^*,\Lambda_R] - \langle P^*,\Lambda_R\rangle [P^*,X] = 0
\end{equation}
\begin{proof}
The numerator of \eqref{eq:control} has the following expression:
\begin{equation}\label{eq:num-ctrl-dep-term}
\langle Z_u^*,X \rangle[Z_u^*, \Lambda_R]- \langle Z_u^*, \Lambda_R \rangle[Z_u^*,X] = 2 Z_u^* \Lambda_R X Z_u^* + \det(Z_u^*)[\Lambda_R,X] 
\end{equation}
which can be verified by Lemma \ref{lem:sl2-lemmas}. 
Now, as already mentioned, in the case of the inscribed control set we have $\alpha=\beta = 1/3$ and so we have $\det(Z_u^*) = 0$ by equation \eqref{eq:general-Z0}. Also, equation \eqref{eq:opt-control-expl} implies:
\begin{align}
\langle Z_u^*,  [\Lambda_R,X] \rangle = 0    
\end{align}
Now since $\langle \Lambda_R, X \rangle = 0$, we get $\bracks {Z_u^*} {\Lambda_R X} = 0$ (as $\Lambda_R X = - X \Lambda_R$ by Lemma \ref{lem:sl2-lemmas}). Thus, $Z_u^* \Lambda_R X + \Lambda_R X Z_u^*= 0$. Multiply by $Z_u^*$ on the right and use the fact that $\det(Z_u^*) = 0$ to conclude that $Z_u^* X \Lambda_R Z_u^* = 0$. Thus, both terms on the RHS of \eqref{eq:num-ctrl-dep-term} are zero. 
\end{proof}
\end{proposition}

Because of (\ref{eq:control}), the costate equations (\ref{eq:lamR}) simplify to:
\begin{align}
\Lambda_1' &= [\Lambda_1, X] \\
\Lambda_R' &= \left[ -\Lambda_1 + \frac{3}{2}\lambda \Jsu ,X\right]
\end{align}

\subsection{Simplifying the Hamiltonian}

\begin{lemma}\label{lem:gen-disc-formula}
We have:
\[
\Delta([\Lambda_R, X],\alpha,\beta) = -\alpha^2 \bracks{X}{X} \bracks{\Lambda_R}{\Lambda_R} - (\beta^2 - \alpha^2)\bracks{\Lambda_R}{\Jsu X}^2
\]
\end{lemma}
\begin{proof}
Compute both sides.
\end{proof}
\begin{lemma}
In the case where the control set is $U_r$, the control dependent term in the Hamiltonian has the following expression:
\begin{equation}
\frac{\bracks{\Lambda_R}{Z_u^*}}{\bracks{X}{Z_u^*}} = 
\frac{(\beta^2 - \alpha^2)\bracks{\Jsu}{\Lambda_R} \bracks{\Jsu}{X} \pm 2 \beta \sqrt{D([\Lambda_R,X],\beta,\alpha)}}{4 D(X,\alpha,\beta)}
\end{equation}
In particular, if $\alpha = \beta = 1/3$ (that is, the control set is $U_I$), we get:
\[
\frac{\bracks{\Lambda_R}{Z_u^*}}{\bracks{X}{Z_u^*}} = \pm \sqrt{\frac{\bracks{\Lambda_R}{\Lambda_R}}{2}}
\]
\end{lemma}
\begin{proof}
The proof is given by using Lemma \ref{lem:gen-disc-formula}, the general expression for $Z_u^*$ as given in equation \eqref{eq:general-Z0} and the expression for the roots in equation \eqref{eq:z-quad-roots} and computing both sides. 
\end{proof}

We choose the root with the positive sign of the square-root, in turn giving the maximum Hamiltonian. 
\begin{equation}\label{eq:simp-inradius-maxham}
\H^*(\Lambda_1,\Lambda_R,X) = \langle \Lambda_1 - \frac{3}{2}\lambda \Jsu, X \rangle + \sqrt{\frac{\langle \Lambda_R, \Lambda_R\rangle}{2}}
\end{equation}

Let us now turn to deriving the dynamics of the Lie algebra element $X$. 
\subsection{Simplifying the State Equations}
\begin{theorem}
With the inscribed control disk $U_I$, we find that the state equation for $X$ simplifies as follows:
\begin{align}
    X' = -\frac{\left[\Lambda_R,X\right]}{\sqrt{2 \langle \Lambda_R, \Lambda_R\rangle}} \label{eq:Xbeforeparam}
\end{align}
\end{theorem}
\begin{proof}
Recall that $X$ is an element of the Lie algebra, and we already have an expression for the Hamiltonian in equation \eqref{eq:simp-inradius-maxham}. Recall also that $\Lambda_R = [\Lambda_2, X]$ and the dynamics for $X$ is given by $X' = \frac{\delta \H^*}{\delta \Lambda_2}$ (See Section \ref{sec:costate-variables}).

First, we note that:
\begin{equation}\label{eq:helper}
    [\Lambda_R,X] = [[\Lambda_2,X],X] = [\Lambda_2 X - X \Lambda_2,X] = -4\Lambda_2-2\bracks{\Lambda_2}{X}X 
\end{equation}

Now we compute:
\begin{align*}
    X' &=\frac{\delta}{\delta \Lambda_2}\sqrt{\frac{\bracks{[\Lambda_2,X]}{[\Lambda_2,X]}}{2}} \\ 
    &= \frac{1}{\sqrt{2\bracks{\Lambda_R}{\Lambda_R}}} \frac{\delta}{\delta \Lambda_2}\frac{\bracks{[\Lambda_2,X]}{[\Lambda_2,X]}}{2} \\ 
    &=  \frac{1}{\sqrt{2\bracks{\Lambda_R}{\Lambda_R}}} \frac{\delta}{\delta \Lambda_2}\left( 2 \bracks{\Lambda_2}{\Lambda_2} + \bracks{\Lambda_2}{X}^2\right)\quad \mathrm{(by~ Proposition~\ref{prop:trace-quotient-sub})} \\
    &=  \frac{1}{\sqrt{2\bracks{\Lambda_R}{\Lambda_R}}} \left(4\Lambda_2 + 2\bracks{\Lambda_2}{X}X\right) \\ 
    &=  -\frac{[\Lambda_R,X]}{\sqrt{2\bracks{\Lambda_R}{\Lambda_R}}}\quad \text{(using equation \eqref{eq:helper})}
\end{align*}
\end{proof}

This gives a remarkable simplification of both the state and costate equations in the inscribed control case. In summary we have the following:

\begin{align*}
    g'&=gX \qquad \qquad \qquad \qquad \quad \Lambda_1' = [\Lambda_1,X] \\
    X'&=- \frac{[\Lambda_R,X]}{\sqrt{2\bracks{\Lambda_R}{\Lambda_R}}} \qquad \qquad \Lambda_R' = \left[-\Lambda_1 + \frac{3}{2}\lambda \Jsu, X \right]
\end{align*}

We shall now use this result to classify the abnormal extremals. 
\subsection{No abnormal extremals.}\label{sec:no-abnormal-extremals}
Recall that an extremal is \emph{abnormal} if $\lambda = 0$. In this case we get a further simplification:
\begin{align}
    g'&=gX \qquad \qquad \qquad \qquad \quad \Lambda_1' = [\Lambda_1,X] \nonumber \\
    X'&=- \frac{[\Lambda_R,X]}{\sqrt{2\bracks{\Lambda_R}{\Lambda_R}}} \qquad \qquad \Lambda_R' = \left[-\Lambda_1, X \right] \label{eq:lamR-abnormal-incircle}
\end{align}
Let us ignore the $\SL(\R)$ component $g$ for the moment and let us do a crude dimension count. Recall that we started out with a 9-dimensional system in $(X,\Lambda_1, \Lambda_2) \in \su^3$.
The following observations reduce the dimension of the system:
\begin{itemize}
    \item We reduced the dimension by 1 by our change of variable $\Lambda_R$ along with the realization that $\langle \Lambda_R, X \rangle = 0$ (See Section \ref{sec:max-ham}). 
    \item If $\lambda = 0$, then from equation \eqref{eq:lamR-abnormal-incircle} we get that $\Lambda_1' + \Lambda_R' = 0$ so that:
    \begin{equation}\label{eq:ang-mom}
    \Lambda_1 + \Lambda_R  = K
    \end{equation} where $K$ is a constant.  Note that this relation means also that $\langle K, X \rangle = \langle \Lambda_1, X \rangle$. This reduces the dimension of the system by 3. 
    \item We have, by the Lax pairs for $X$ and $\Lambda_1$ that $\det(X)$ and $\det(\Lambda_1)$ are constant. Equivalently, we have that $\langle X,X \rangle = -2$ and $\langle \Lambda_1, \Lambda_1 \rangle$ is a constant. This reduces the dimension of the system by 2. 
    \item The vanishing of the Hamiltonian means that \begin{equation}\label{eq:ham-vanish}
    \H^*(\Lambda_1,\Lambda_R,X) = \langle \Lambda_1, X\rangle + \sqrt{\frac{\langle \Lambda_R, \Lambda_R\rangle}{2}} = 0
    \end{equation}
which reduces the dimension of the system by one. 
\end{itemize}
At this point, we have a system with two degrees of freedom. 
Using equation \eqref{eq:ang-mom}, we get:
\begin{align}
    \Lambda_R' &= \left[\Lambda_R - K, X \right]    \label{eq:abnormal-system-lamR} \\
    X' &= -\frac{\left[\Lambda_R,X\right]}{\sqrt{2 \langle \Lambda_R, \Lambda_R\rangle}} \label{eq:abnormal-system-X}
\end{align}
where $K$ is a constant. 

\begin{lemma}\label{lem:su11-invariance}
If $K,X,\Lambda_R$ is a solution to equations \eqref{eq:abnormal-system-lamR} and \eqref{eq:abnormal-system-X}, then so is $\Ad_g K,\Ad_g X, \Ad_g \Lambda_R$ (where $g \in \SU$) and they also satisfy the relations listed at the beginning of the section.
\end{lemma}
\begin{proof}
This is simple to verify once we use the fact that the trace inner product is an invariant inner product and using properties of Lie brackets and inner products. 
\end{proof}

We now classify periodic solutions to this system of equations. In fact, we prove that periodic solutions do not exist (assuming that $\Lambda_R(t_0) \ne 0$ for some $t_0$). We begin with a lemma:

\begin{lemma}\label{lem:lamR-X-props}
We also have the following identities relating $X$ and $\Lambda_R$. 
\begin{enumerate}
    \item $[X, [X, \Lambda_R]] = -4 \Lambda_R $
    \item $[\Lambda_R,[\Lambda_R,X]] = 2 \langle \Lambda_R, \Lambda_R \rangle X$
    \item $c_K := \langle K, K \rangle$ and $c_R := \bracks {\Lambda_R} {\Lambda_R} - 2 \bracks {K} {\Lambda_R} $ are both constant.
    \item $\bracks {[[K,X],X]} K = -4 \bracks K K - 2 \bracks X {K} ^2$
    \item $K$ is not zero.
\end{enumerate}
\end{lemma}
\begin{proof}
If $K=0$ then we have $\bracks{\Lambda_R}{\Lambda_R} = 0$ using equation \eqref{eq:ham-vanish}. This implies that $\Lambda_R = 0$ since the trace inner-product is positive definite on the 2-dimensional subspace consisting of vectors orthogonal to $X$, contradicting the non-triviality condition of the maximum principle. 
The constancy of $c_K$ is trivial and $c_R$ is constant because its derivative vanishes after using \eqref{eq:abnormal-system-lamR}. Now, we have
\[ 
 \bracks {[[K,X],X]} K = 2\bracks X X \bracks K K - 2 \bracks X K^2 = -4 \bracks K K - 2 \bracks X K ^2
\]
by Propositions \ref{prop:brack-brack} and \ref{prop:brack-brack-inner-prod} and the fact that $\bracks X X = -2$ which shows (5). For the rest, since the identities are $\SU$-invariant, it is enough to test them at $X=\Jsu$ and $\Lambda_R = \mattwo 1 0 0 {-1}$, since $X$ belongs to the adjoint orbit of $J$ and $\Lambda_R$ lies in the subspace orthogonal to the subspace spanned by $X$. 
\end{proof}

Lemma \ref{lem:su11-invariance} suggests that we look at functions in the algebra of $\SU$-invariant functions. To this end, consider the following expressions:
\begin{align*}
    u &:= \langle X,K \rangle \\
    v &:= \bracks {[\Lambda_R,X]} K \\
    w &:= \sqrt{2 \bracks {\Lambda_R} {\Lambda_R}}
\end{align*}

The crude dimension count we did earlier suggests that we should have an extra identity relating $u,v,w$ apart from the Hamiltonian vanishing condition. 

\begin{lemma}\label{lem:deriv-u-v-w}
We have the following expressions for the derivatives of $u,v,w$:
\begin{align*}
    u' &= -\frac{v}{w} \\
    w' &= \frac{2v}{w} \\
    v' &= (2 c_R - w^2) + 4c_K + 2 u^2 - w u
\end{align*}
We also have the following polynomial relations between $u,v,w$:
\begin{align}
    2u + w &= 0 \label{eq:relHam} \\
    4v^2 - 4u^2w^2 - 8c_K w^2 + (w^2 - 2 c_R)^2 &= \mathrm{constant} \label{eq:relPhi}
\end{align}

\end{lemma}
\begin{proof}
The Hamiltonian vanishing condition \eqref{eq:ham-vanish} shows that $2u + w = 0$. 

For the rest, we have:

\begin{align*}
    u' &= \bracks X K ' = \bracks {X'} K = \frac{\bracks {[\Lambda_R, X]} K}{\sqrt{2 \bracks {\Lambda_R} {\Lambda_R}}} = -\frac{v}{w} \\
    w' &=(\sqrt{2 \bracks {\Lambda_R} {\Lambda_R}})' = \frac{2  \bracks {\Lambda_R} {\Lambda_R}'}{w} = \frac{2 v}{w} \\
    v' &= \bracks {[\Lambda_R, X]} K' \\ 
    &= \bracks {[[\Lambda_R,X],X]} K - \bracks {[[K,X],X]} K - \frac{\bracks {[\Lambda_R, [\Lambda_R,X]]} K }{\sqrt{2 \bracks {\Lambda_R} {\Lambda_R}}} \\
    &= \underbrace{-4 \bracks {\Lambda_R} K}_{\text{by Lemma \ref{lem:lamR-X-props}, (1)}} + \underbrace{ 4 \bracks K K + 2 \bracks X K ^2}_{\text{by Lemma \ref{lem:lamR-X-props}, (4)}} - \underbrace{\frac{ 2 \bracks {\Lambda_R}{\Lambda_R} \bracks K X }{\sqrt{2 \bracks {\Lambda_R} {\Lambda_R}}}}_{\text{by Lemma \ref{lem:lamR-X-props}, (2)}} \\
    &= \underbrace{(2 c_R - w^2)}_{\text{by Lemma \ref{lem:lamR-X-props}, (3)}} + 4c_K + 2 u^2 - w u
\end{align*}
Finally, to verify \eqref{eq:relPhi}, we differentiate its left hand side and use the relations of $u',v',w'$:
\begin{align*}
     (4v^2 - 4u^2w^2 - 8c_K w^2 + (w^2 - 2 c_R)^2)'
     &= 16vc_R - 8vw^2 + 32c_K v + 16u^2 v - 8 vwu \\
     &-4(-2 uvw + 4 u^2 v) - 32c_K v + 8(w^2 - 2 c_R)v  \\
     &= 0.
\end{align*}\end{proof}
We can now prove:
\begin{lemma}\label{thm:v-affine}
The function $v$ is an affine function.
\end{lemma}
\begin{proof}
Note that if we eliminate $w$ from the relations \eqref{eq:relPhi} and \eqref{eq:relHam} we get
\begin{equation}\label{eq:v2-expression}
v^2 = (8c_K + 4c_R)u^2 - c_R^2
\end{equation}
Differentiating this, we get:
\begin{align*}
    v' &= \frac{2(8c_K + 4c_R)}{2v}uu' \\
    &=\frac{2(8c_K + 4c_R)}{2v} \times \frac{-w}{2} \times \frac{-v}{w} \\
    &=(4c_K + 2c_R)
\end{align*}
which is a constant. (Here we used expressions for $u$ and $u'$ derived in Theorem \ref{lem:deriv-u-v-w}.) This proves that $v$ is affine.
\end{proof}

From the remark at the end of Section \ref{sec:transversality}, we recall that we are primarily interested in \emph{periodic} solutions to the state-costate ODEs. From our analysis so far, we can prove that periodic solutions in the abnormal system do not exist. 

\begin{theorem}
In the case where the control set is the inscribed disk $U_I$ of the 2-simplex, the state-costate equations do not have any periodic solutions. 
\end{theorem}
\begin{proof}
If the solutions are periodic, then $v$ shall be a periodic function of the same period. However, Theorem \ref{thm:v-affine} shows that this can only happen provided $(2c_K + c_R) = 0$ \emph{i.e.,} if $v$ is constant. From \eqref{eq:v2-expression}, we get that 
\[
v^2 + c_R^2 = 0
\]
which proves that $v = c_R = c_K = 0$. We also have that $u = \bracks{X}{K}$ and $w^2 = \bracks{\Lambda_R}{\Lambda_R}$ are both constant. By Lemma \ref{lem:lamR-X-props}, we have that $K \ne 0$ but $c_K = \bracks{K}{K} = 0$. So, up to conjugation by an element in $\SU$, $K$ looks like $\frac{\pm1}{2}\mattwo{i}{1}{1}{i}$. Now if we have (by Lemma \ref{lem:def-phi}) that $X = X_z$,  $\bracks{X}{K}$ constant gives us that $y = \mathrm{Im}(z)$ is constant.  Using the fact that $X$ is periodic, this implies that $X'(t_0) = [\Lambda_R(t_0),X(t_0)] = 0$ for some $t_0$. Also, $\bracks{\Lambda_R(t_0)}{X(t_0)} = 0$. These relations imply that $2 \Lambda_R(t_0) X(t_0) = 0$, which in turn implies that $\Lambda_R(t_0)=0$. Since $\bracks{\Lambda_R}{\Lambda_R}$ is constant, we have that $\Lambda_R$ is identically zero. This means that we are in the singular locus, which gives a contradiction. 
\end{proof}
\section{Singular locus}

The singular locus, as shown in Section \ref{sec:singular-locus} is the region of the cotangent space $T^*(T\SL(\R))$ given by:
\[
\Lambda_{sing} = \left\{\left(g_0,-\frac{3}{2}J,J,0\right)~|~g_0\in\SL(\R)\right\} \subset \SL(\R) \times \sl(\R) \times \sl(\R)\times \sl(\R)
\]


In Hales~\cite{hales2017reinhardt} it is shown that:
\begin{itemize}
    \item Pontryagin extremals which avoid the singular locus are given by bang-bang controls with finitely many switches (See Theorem \ref{thm:bang-bang}).
    \item the global optimal trajectory of the Reinhardt control problem with the control set $U_T$ cannot stay within the singular locus for any finite interval of time (See Theorem \ref{thm:no-singular-arcs}).
\end{itemize} 

This means that if we are considering the control set $U_T$, for a Pontryagin extremal to approach the singular locus, the control should switch infinitely many times around the boundary of $U_T$ in a finite interval of time. This is the \emph{chattering} phenomenon (see Fuller~\cite{fuller1963study} and Zelikin \& Borisov~\cite{zelikin2012theory}). 

If we use a control set $U_r$ which has a smooth boundary, we would expect the optimal control to perform infinitely many rotations along the boundary $\partial U_r$ to approach the singular locus in finite time. A system with exactly this behaviour is described in Manita \& Ronzhina~\cite{manita2022optimal} and the associated trajectory is spiral-like. The results of this section are motivated by that paper.

These results warrant a study of the behaviour of the system near the singular locus. To this end, we introduce convenient coordinates and re-express the state, costate and optimal control equations in these coordinates. Throughout this section, unless otherwise specified, we work with the circular control sets $U_r$ and we assume that $\alpha > 0$ (See Section \ref{sec:su11-control-sets}). 

\subsection{Hyperboloid Coordinates}

We begin with the following lemma:

\begin{lemma}\label{lem:hyperboloid-coords}
If $A = \mattwo{a}{b}{c}{-a}$ is an arbitrary matrix in $\sl(\R)$, then it is conjugate by the Cayley transform to the following matrix in $\su$.
\begin{equation}
     = \mattwo{i\sqrt{\det(A) + |p|^2}}{p}{\overline{p}}{-i\sqrt{\det(A) + |p|^2}} \in \su
\end{equation}
where $p = \frac{(b+c)}{2} + i a$.
\end{lemma}
\begin{proof}
If $A \in \sl$, then we have, by the Cayley transform:
\[
C^{-1}AC = \left(
\begin{array}{cc}
 i\frac{(b-c)}{2}  & \frac{b+c}{2}+ i a \\
 \frac{b+c}{2} -i a & -i\frac{(b-c)}{2}  \\
\end{array}
\right)
\]
Now, we have that $\mathrm{det}(A) = -a^2 - bc$ and $\sqrt{\mathrm{det}(A) + |p|^2} = \frac{b-c}{2}$ which gives the required.
\end{proof}

This shows that the matrix $A = \mattwo{a}{b}{c}{-a}$ is parameterized by a single complex number $p = \frac{(b+c)}{2} + i a$ and $\det(A)$.

\begin{remark} \normalfont \leavevmode
\begin{itemize}
    \item For a complex number $p$, we shall denote $[p] := \sqrt{1 + |p|^2}$. 
    \item  We shall write
    \begin{equation*}
        \RR(u,v) := \mathrm{Re}(\overline{u}v)
    \end{equation*}
    to denote the real part of a complex number extended to $\mathbb{C} \times \mathbb{C}$ as sesquilinear forms.
\end{itemize}
\end{remark}

\begin{definition}[Hyperboloid coordinates]
The \emph{hyperboloid coordinates} of the matrix $A = \mattwo{a}{b}{c}{-a} \in \sl(\R)$ are $p = \frac{b+c}{2} + i a$ and $\det(A)$. 
\end{definition}

Now, using the above, we shall transform the ODE for $X,\Lambda_1,\Lambda_R$ into ODEs given by hyperboloid coordinates. Recall that we have the constraints $\bracks{X}{X} = -2$, $\bracks{\Lambda_1}{\Lambda_1} = \mathrm{const.} $ and $\bracks{X}{\Lambda_R} = 0$. This implies that $\det(X) = 1$ and we assume that $\det(\Lambda_1) := d$. Using these and Lemma \ref{lem:hyperboloid-coords}, we can write\footnote{In general, the determinant of $\Lambda_1$ may be positive, negative or zero. But the trajectory moves leaving this determinant unchanged. So, we recover the $d = \det(\Lambda_1) < 0$ by negating $\Lambda_1$. The case $d=0$ needs to be handled to account for the full optimal control problem. Here, we work with $d>0$ since this is the case for the singular locus. }:

\begin{equation}\label{eq:st-costate-hyperboloid}
    X := \mattwo{-i[w]}{w}{\overline{w}}{i[w]} \quad
    \Lambda_1 := \sqrt{d} \mattwo{i[b]}{b}{\overline{b}}{-i[b]} \quad
    \Lambda_R := \mattwo{-i\frac{\RR(c,w)}{[w]}}{c}{\overline{c}}{i\frac{\RR(c,w)}{[w]}}
\end{equation}
for complex numbers $c,b,w \in \mathbb{C}$, which ensure that the constraints are satisfied. (Note that $(w,[w])$ lie on the upper sheet of the hyperboloid $\{(w,t) \in \mathbb{C} \times \R~|~t^2 - |w|^2 = 1 \}$, which justifies our nomenclature.)

For a control matrix given by equation \eqref{eq:general-boundary-Z0} (which is already in $\su$ by the results of Section \ref{sec:su11-control-sets}) and using Corollary \ref{cor:sl2-star-condition}, we may derive the star condition in these coordinates as: 
\begin{equation*}
    \mu(w,z) := [w]-\frac{\beta}{\alpha}\RR(w,z) =  \frac{\bracks{X}{Z_u}}{-2\alpha}  > 0
\end{equation*}
for $|z| < 1$. 

Another advantage of switching to hyperboloid coordinates is that, at the singular locus, we have $w=b=c=0$. Let us now turn to deriving the optimal control in hyperboloid coordinates.\footnote{As usual, we take $J_{\su} = C^{-1}JC = \mattwo{-i}{0}{0}{i}$.}

\begin{remark}\normalfont
In what follows, we write:
\begin{equation}
    f_1(t) = f_2(t) + o(e(t))
\end{equation}
to mean $e(t)$ is non-zero on a punctured neighbourhood of the singular locus and the limit of $(f_1 - f_2)/e$ is 0 at $t$ goes to zero at the singular set. 
\end{remark}

\begin{lemma}\label{lem:maximizing-root}
On a punctured neighbourhood of the singular locus, the root of the quadratic equation giving the optimal control is $z = c/|c| + o(1)$.
\end{lemma}
\begin{proof}
For the circular control set $U_r$, the Hamiltonian maximization occurs on the boundary of the control set, and if the optimal control matrix is given by:
\begin{equation*}
Z_u^* = \left(\begin{matrix}
-i\alpha & \beta z \\ 
\frac{\beta}{z} & i \alpha
\end{matrix}\right)
\end{equation*}

Then by Proposition \ref{prop:control-quadratic}, angular momentum conservation gives the following quadratic equation in $z$ (see equation \eqref{eq:quadratic-equation}):
\begin{equation*}
\alpha[\Lambda_R,X]_{21} z^2 - 2i\beta [\Lambda_R,X]_{11}z + \alpha [\Lambda_R,X]_{12} = 0
\end{equation*}

On a punctured neighbourhood of the singular locus, we may assume $c\ne0$. If we set $w_1 := \bar{c}w/|c|$ and $z_1 := \bar{c}z/|c|$, then the quadratic equation becomes
\[
a_1 + 2b_1 z_1 - \bar{a}_1z_1^2 := 
(2 + |w_1|^2 - w_1^2) + 2 \beta_1(w_1-\bar{w}_1)[w_1]z_1 - (2+ |w_1|^2 - \bar{w}_1^2)z_1^2 = 0
\]
A calculation shows that the Hamiltonian is not constant in $z$, which implies that the discriminant is non-zero and the two roots are distinct. Now we have that, near the singular locus, the coefficients are $|a_1| = 2 + o(1) > 0$ and $b_1 = o(1)$. Using this, we get that the discriminant of this quadratic in $z_1$ is a real number, and equals:
\[
\Delta := 4|a_1|^2+4b_1^2 = 16 + o(1) > 0
\]
From Proposition \ref{prop:control-quadratic}, the maximizing root of the quadratic is the one with the positive sign in front of the square root. So, the maximizing root $z$ of the quadratic in \eqref{eq:quadratic-equation} is
\begin{equation}
    z^* = \frac{c}{|c|}z_1^* = \frac{c}{|c|}\frac{2b_1 + \sqrt{\Delta}}{2\bar{a}_1} =\frac{c}{|c|} + o(1)
\end{equation}
which gives the required. 
\end{proof}

The ODEs in Lie algebra coordinates for $X,\Lambda_1$ and $\Lambda_R$ were derived in Sections \ref{sec:costate-variables} and \ref{sec:lie-algebra-dynamics}. We now write them out in the hyperboloid coordinates we have defined above. This will enable us to better understand the dynamics near the singular locus. 

\begin{theorem}[ODE in Hyperboloid Coordinates]
In hyperboloid coordinates, the dynamics for $X,\Lambda_1$ and $\Lambda_R$ look like:
\begin{align}
    b' &= 2i\left([b]w +b[w]\right) \label{eq:eqn-b} \\
    w' &= i\frac{w - \beta_1 [w] z^*}{\mu(w,z^*)} \label{eq:eqn-w} \\
    c' &=i \frac{\left(-[w] + \beta_1 \RR(w,z^*)\right)\left(-[w]c + \beta_1 \RR(c,w) z^* \right)}{[w] \mu(w,z^*)^2} \nonumber \\ 
    &- i\frac{\left(-\RR(c,w) + \beta_1 [w] \RR(c,z^*)\right)\left(-w + \beta_1 [w] z^*\right)}{[w] \mu(w,z^*)^2} \nonumber \\
    &-i\left((-3+2[b]\sqrt{d})w + 2 b \sqrt{d}[w]\right) \label{eq:eqn-c}
\end{align}
where $\beta_1 = \beta/\alpha$ and $z^*$ is the maximizing root of the quadratic equation for the optimal control. 
\end{theorem}
\begin{proof}
This is an elementary computation using the hyperboloid form of the state and costate equations as in equation \eqref{eq:st-costate-hyperboloid} and substituting them into the state and costate ODE derived in Sections \ref{sec:lie-algebra-dynamics} and \ref{sec:costate-variables}. 
\end{proof}

Similarly, we have expressions for the Hamiltonian and the angular momentum as:
\begin{align}
\mathcal{A} &= 2 \sqrt{d}[b] - 2 \frac{\RR(w,c)}{[w]} \label{eq:hyperboloid-ang-mom} \\ 
\H &= (2\sqrt{d}\RR(w,b) + (2\sqrt{d}[b]-3)[w])-\frac{\RR(w-\beta_1z^*,c)}{\mu(w,z^*)[w]} \label{eq:hyperboloid-ham}
\end{align}

\subsection{Valuation-Truncation}

It can be checked that naive linearization of the system $\eqref{eq:eqn-b},\eqref{eq:eqn-w},\eqref{eq:eqn-c}$ around the singular locus $w=b=c=0$ failed to give anything interesting, which suggests that the rate of growth in $t$ of these quantities near the singular locus is not all the same. Thus, we need a mechanism to track the rate of growth of each quantity. To this end, we introduce \emph{valuations} as a way to systematically track these orders of growth.

\begin{definition}[Valuations and Angular Component]
For $\epsilon>0$, if $f:(0,\epsilon) \to \mathbb{C}$ is any function, we say that $f$ has \emph{valuation $k \in \mathbb{Z}_{\ge 0}$ at $t = 0$} (denoted $\val(f)$) if there exists a smallest integer $k$ for which 
\[
\lim_{t \to 0}\frac{f(t)}{t^k} \ne 0
\] This non-zero limit is called the \emph{angular component} of $f$, which we shall denote $\mathrm{ac}(f)$. 
If $f \equiv 0$, then we set $\val(f) = +\infty$.
\end{definition}

We have the following properties of valuations:

\begin{lemma}[Properties of Valuations]\label{lem:valuation-props-1}
If $f,g : (0,\epsilon) \to \mathbb{C}$ are two complex-valued functions, then:
 \begin{enumerate}
     \item if $\val(f)$ and $\val(g)$ exist, then so does $\val(fg)$ and \[\val(fg) = \val(f) + \val(g).\]
     thus, if $\val(f)$ exists, then so does $\val(f^n)$ and is equal to $n\,\val(f)$ for $n\in \mathbb{N}$. 
     \item if $\val(f), \val(g)$ and $\val(f+g)$ exist, then \[\val(f+g) \ge \min(\val(f),\val(g)).\] 
      Thus, if  $f + g = 0$, then  $\val(f) = \val(g)$ if $\val(f)$ exists.
     \item if $\val(f), \val(g)$ exist and $\val(f) \ne \val(g)$ then $\val(f+g)$ exists and \[\val(f+g) = \min(\val(f),\val(g)).\]
    \item if $\val(f) \ne 0$ exists and if $f$ is differentiable then $\val(f')$ also exists and is equal to $\val(f)-1$. 
 \end{enumerate}
\end{lemma}
\begin{proof}
\begin{enumerate}
    \item 
Let $k_1 = \val(f)$ and $k_2 = \val(g)$ so that $\lim_{t \to 0}f(t)/t^{k_1} = \ac(f)$ and $\lim_{t \to 0}g(t)/t^{k_2} = \ac(g)$, both non-zero. Then, $\val(fg)$ exists since 
\[
\lim_{t\to 0}\frac{f(t)g(t)}{t^{k_1 + k_2}}=\ac(f)\ac(g) \ne 0
\]
So, $\val(fg) \le k_1 + k_2$. We can see that $k_1 + k_2$ is the smallest index for which this happens, for otherwise we could find a smaller index for at least one of $f$ or $g$ for which the defining property of valuations is true. Thus, $\val(fg) = \val(f) + \val(g)$. 
\item 
 In the case where $\val(f) = \val(g) = k$, the angular components of $f$ and $g$ may sum to zero, and so $\val(f+g)$ may not exist in that case. But if we assume that $k_3 = \val(f+g)$ exists, then we have 
\[
 \lim_{t\to 0}\frac{f(t)}{t^{k}}+\lim_{t\to 0}\frac{g(t)}{t^{k}} =  \lim_{t\to 0}\frac{f(t)+g(t)}{t^{k}} 
\]
which exists and may be zero. By definition this means that $k_3 \ge k$, which gives $\val(f+g) \ge \min(\val(f), \val(g))$.
\item If $k_1 = \val(f)$, $k_2 = \val(g)$ are distinct:
\[
\lim_{t\to 0}\frac{f(t)+g(t)}{t^{\min(k_1,k_2)}}= \begin{cases}
\ac(f) \ne 0 & \text{if } k_1 < k_2 \\
\ac(g) \ne 0 & \text{if } k_2 < k_1
\end{cases}
\]
thus, $\val(f + g)$ exists and is equal to $\min(\val(f),\val(g))$.
\item 
Now, if $f$ is differentiable and $k=\val(f) \ge 1$ exists, then we have:
\[
0 \ne \lim_{t \to 0}\frac{f(t)}{t^k} =\frac{1}{k} \lim_{t \to 0}\frac{f'(t)}{t^{k-1}}
\]
by L'H\^{o}pital's rule, which proves that $\val(f') = \val(f) - 1$.

\end{enumerate}

\end{proof}

We also need a lemma about how valuation interacts with the real and imaginary part bilinear forms:
\begin{lemma}\label{lem:valuation-props-2}
If $f,g : (0,\epsilon) \to \mathbb{C}$ are functions so that $|f|,|g|$, $\RR(f,g)$ and $\RR(f,ig)$ have valuations, then:
\begin{enumerate}
    \item $\val(\RR(f,ig)) = \val(\RR(if,g))$
    \item we have $\val(\RR(f,g)) \ge \val|f| + \val|g|$ and $\val(\RR(f,ig)) \ge \val|f| + \val|g|$.
    \item If $g$ is differentiable, then 
    \[ \val(\RR(g',g)) = 2 \val|g| - 1.
    \]
    \item if $|f|$ has positive valuation, then $[f]$ has valuation zero. 
\end{enumerate}
\end{lemma}
\begin{proof}
\begin{enumerate}
    \item The first statement is obvious. 
    \item Note that we have the identity:
    \[
    |f|^2|g|^2 = |\overline{f}|^2|g|^2 = \mathrm{Re}(\overline{f}g)^2 + \mathrm{Im}(\overline{f}g)^2 = \RR(f,g)^2 + \RR(if,g)^2
    \]
    Then we get:
    \[
    \val|f| + \val|g| \le \min(\val(\RR(f,g)),\val(\RR(f,ig)))
    \]
    which proves the required. 
    \item 
    To prove this statement, we differentiate the equation $|g|^2 = g \overline{g}$ to get:
    \[
    2|g|'|g| = \RR(g,g')
    \]
    from this we get $\val(\RR(g,g')) = 2\val(|g|)-1$ since $\val(g') = \val(g)-1$ by Lemma \ref{lem:valuation-props-1},4. 
    \item This is because we have the binomial expansion: 
    \[
    [f] = \sqrt{1 + |f|^2} = 1 + |f|^2/2 + \dots
    \] 
    If $|f|$ has positive valuation, then as $t \to 0$, $[f] \to 1$ and so $[f]$ has valuation zero.  
\end{enumerate}


\end{proof}

We now use this apparatus to derive a truncation of the system in equations $\eqref{eq:eqn-b},\eqref{eq:eqn-w},\eqref{eq:eqn-c}$ near the singular locus $w=b=c=0$. We assume that $|w|,|b|,|c|$ have non-zero and finite valuations. 

\begin{theorem}[Valuation-Truncation]\label{thm:valuation-truncation}
Assume that $w,b,c : (0,\epsilon) \to \mathbb{C}$ are optimal trajectories \emph{(}evolving as described by equations \eqref{eq:eqn-b},\eqref{eq:eqn-w},\eqref{eq:eqn-c}\emph{)} and tending toward the singular locus $w=b=c=0$. Assume that $|w|,|b|,|c| : (0,\epsilon) \to \R$ have finite, non-zero valuations and that the control set $U_r$ has parameters so that $\beta/\alpha \ne 0$. Then, $\RR(b,c),\RR(c,w),\RR(b,w),\RR(b,ic),\RR(c,iw)$ and $\RR(b,iw)$ have valuations and
\begin{align*}
    \val|w| &= 1,\quad \val|b|=2,\quad\val|c|=3 \\
    \val\RR(b,w)&=3,\quad\val\RR(c,w)=4,\quad\val\RR(b,c)=5 \\
    \val\RR(b,iw)&=3,\quad\val\RR(c,iw)=4,\quad\val\RR(b,ic)=5
\end{align*}
\end{theorem}
\begin{proof}
Our strategy of proof will be to use the conservation laws, and the ODEs for $c,b,w$ to systematically obtain information about their valuations as equalities and inequalities. 

Let us start with the angular momentum in equation \eqref{eq:hyperboloid-ang-mom}. Note that, at the singular locus, $d = \det(\Lambda_1) = 9/4$ and $\Lambda_R = 0$. So, the angular momentum is $\bracks{J}{\Lambda_1 + \Lambda_R} = 3$. The angular momentum equation reads:
\[
\mathcal{A} = 2 \sqrt{d}[b] - 2 \frac{\RR(w,c)}{[w]} = 3[b] - 2\frac{\RR(w,c)}{[w]}
\]
Now, since $\val(\mathcal{A} - 3) = +\infty$, we have that $\val(\RR(w,c))$ exists. So, we get by properties of valuations, that:
\[
\val\left(\frac{\RR(w,c)}{[w]}\right) = \val(\RR(w,c)) = \val([b] - 1)
\]
Where we used the fact that $\val(1/[w]) = 0$, since it can be expanded as a series in $|w|$ about $w=0$.  
For similar reasons, since $\mathrm{ac}([b]) = 1$, we get that $\val([b]-1) = \val|b|^2 = 2 \val|b|$. 
So we get:
\begin{equation}\label{eq:val-b-ge}
2\val|b| = \val|b|^2 = \val \RR(c,w) \ge \val|c| + \val|w|
\end{equation}
using Lemma \ref{lem:valuation-props-2},(2). 

Lemma \ref{lem:maximizing-root} shows us that the maximizing root of the quadratic equation in equation \eqref{eq:quadratic-equation} is $z^* = c/|c| + o(1)$, so that the control matrix becomes: 
\[Z_u^* = \left(\begin{matrix}
-i\alpha & \beta c/|c| \\ 
\beta\overline{c}/|c| & i \alpha
\end{matrix}\right) + o(1)
\]
This means that the Hamiltonian in equation \eqref{eq:hyperboloid-ham} becomes:
\begin{equation}\label{eq:val-ham}
 0 = \mathcal{H} = (3\RR(w,b) + (3[b]-3)[w])-\frac{\RR(w-\beta_1c/|c|,c)}{\mu(w,c/|c|)[w]} 
\end{equation}
Let us analyze each summand on the right hand side above. The last term is:
\begin{align*}
\left(\frac{\RR(w-\beta_1c/|c|,c)}{\mu(w,c/|c|)[w]} \right) - \beta_1 |c| = \frac{|c|\RR(w,c)(-1+\beta_1^2)}{|c|[w]-\beta_1 \RR(w,c)} = o(|c|)
\end{align*}
Since, by equation \eqref{eq:val-b-ge}, we have $\val(\RR(w,c)) > \val(|c|)$. For the second term of equation \eqref{eq:val-ham}, we have:
\[ 
(3[b] - 3)[w] = \frac{3 |b|^2}{2} + o(|b|^2)
\]
once we expand $[b]$ and $[w]$. Now, again by \eqref{eq:val-b-ge}, we have $\val|b|^2 = \val(\RR(w,c)) > \val|c|$ and so the previous two equations show that: 
\[
(3[b]-3)[w] + \frac{\RR(w-\beta_1c/|c|,c)}{\mu(w,c/|c|)[w]} = \beta_1 |c| + o(|c|)
\]
These remarks show that $\val(\RR(w,b))$ exists and by \eqref{eq:val-ham}, since $\H = 0$, we get:
\begin{equation}
\val(3\RR(w,b)) = \val\left( (3[b]-3)[w]) - \frac{\RR(w-\beta_1c/|c|,c)}{\mu(w,c/|c|)[w]} \right)    
\end{equation}
Thus, 
\begin{equation}\label{eq:val-c-ge}
\val|c| = \val(\RR(w,b)) \ge \val|w| + \val|b|
\end{equation}
We have used up the conservation laws to gain information about the valuations of $c,b,w$. Let us now turn to their ODEs. We shall expand the right hand sides of the ODEs and retain the least valuation term in each case. 

From \eqref{eq:val-c-ge} and \eqref{eq:val-b-ge}, we conclude that $\val|b| > \val |w|$. This, and the ODE for $b$ in equation \eqref{eq:eqn-b} gives:
\begin{equation*}
    b' = 2i([b]w + b[w]) = 2iw + o(|w|)
\end{equation*}
Now, from Lemma \ref{lem:valuation-props-1},4 we get:
\begin{equation}\label{eq:valuation-b'}
    \val|b|' = 2 \val|b|-1 = \val(\RR(b,iw))
\end{equation}

Similarly, the ODE for $w$, in equation \eqref{eq:eqn-w} gives:
\[
w' = i\left(\frac{w - (\beta_1 w c/|c|)}{[w] - \beta_1 \RR(w,c/|c|)}\right) = -i\beta_1\frac{c}{|c|} + o(1)
\]
and, using Lemma \ref{lem:valuation-props-1},4
\begin{equation}\label{eq:valuation-w}
    2\val|w| - 1 + \val|c| = \val(\RR(w,ic))
\end{equation}
Lastly, the ODE for $c$ in equation \eqref{eq:eqn-c} gives:
\begin{equation*}
    c' = -3ib + o(|b|)
\end{equation*}
which we derive by expanding terms out in order of increasing valuation and then apply \eqref{eq:val-c-ge} to retain the most significant term. This also gives:
\begin{equation}\label{eq:valuation-c}
    2\val|c| - 1 = \val(\RR(b,ic))
\end{equation}

Let us now turn to finding the valuations of $\RR(w,c)$. 
We have:
\begin{equation*}
    \RR(w,c)' = \RR(w',c) + \RR(w,c')
\end{equation*}
we now use the the ODEs for $w$ and $c$ and properties of valuations to get:
\[
\val(\RR(w,c))-1 = \min(\val|c|,\val(\RR(w,ib))) \le \val|c|
\]
Using the fact that $\val|w| + \val|c| \le \val(\RR(w,c))$, we get $\val|w| \le 1$. 

Now, from \eqref{eq:valuation-w} we have:
\[
2\val|w| + \val|c| -1 = \val(\RR(w,ic)) \ge \val|w| + \val|c|
\]
from which we get that $\val|w| \ge 1$. 
This proves that $\val|w| = 1$.

Now we turn to finding the valuation of $\RR(b,c)$:
\begin{equation*}
    \RR(c,b)' = \RR(c',b) + \RR(c,b')
\end{equation*}
we now use the the ODEs for $b$ and $c$ and properties of valuations to get:
\begin{equation}\label{eq:valuation-Rbc}
    \val(\RR(c,b))-1 = \val(\RR(w,ic)) = 2 \val|w| + \val|c| - 1
\end{equation}

And since we have $\val|c| + \val|b| \le \val(\RR(c,b))$, we get $\val|b| \le 2\val|w|$. Now, from equation \eqref{eq:valuation-b'}, we have $2 \val|b| - 1 = \val(\RR(b,iw))$, which implies $\val|w| \le \val|b| - 1$. This we get:
\[
\val|w| + 1 \le \val|b| \le 2 \val|w|
\]
which, since $\val|w| = 1$, gives $\val|b| = 2$ and so $\val(\RR(b,iw)) = 3$. 

From equation \eqref{eq:val-c-ge}, we have
\[
\val|c| = \val(\RR(b,w)) \ge \val|b| + \val|w| = 3
\]
Now, we also have, from equation \eqref{eq:val-b-ge} that $\val|c| \le 2\val|b| - \val|w| = 3$. Thus, we have proven that $\val|c| = \val(\RR(b,w)) = 3$. From equation \eqref{eq:val-b-ge} again, we get $\val(\RR(c,w)) = 4$. From equation \eqref{eq:valuation-Rbc}, we get $\val(\RR(b,c)) = 5$. The conclusion follows. 





\end{proof}

The previous theorem implies that near the singular locus, the $c,b,w$ system looks like:
\begin{align*}
    c' &= -3ib +o(t^2) \\
    b' &= 2iw + o(t)  \\
    w' &= -i\beta_1 \frac{c}{|c|} + o(1) 
\end{align*}
using this, we create a \textit{truncated system} of ODEs given by: 
\begin{align}
    c_F' &= -3ib_F \label{eq:eqn-cF} \\
    b_F' &= 2iw_F  \label{eq:eqn-bF} \\
    w_F' &= -i\beta_1 \frac{c_F}{|c_F|}  \label{eq:eqn-wF}
\end{align}
where the little-oh terms are discarded. This new system governs the dynamics of the Reinhardt system very close to the singular locus. We also have the following expressions for the truncated Hamiltonian and angular momentum:
\begin{align}
    \H_c &:=2\RR(w_F,b_F) + \beta_1|c_F| \\
    \mathcal{A}_c &:= \frac{2}{3}|b_F|^2 - 2 \RR(c_F,w_F)
\end{align}
\subsection{Hamiltonian nature of truncated system.}\label{sec:fuller-system}
We introduce new coordinates $z_3 := c_F/6\beta_1$, $z_2 := c_F'/6\beta_1 = -ib_F/3\beta_1$, $z_1 := c_F''/6\beta_1 = w_F/\beta_1$, so that the truncated system in equations \eqref{eq:eqn-cF},\eqref{eq:eqn-bF},\eqref{eq:eqn-wF} becomes:
\begin{equation}\label{eq:fuller-system} 
    z_3' = z_2, \quad 
    z_2' = z_1, \quad
    z_1' = -i\frac{z_3}{|z_3|} 
\end{equation}

The truncated Hamiltonian and angular momentum are:
\begin{align*}
\mathcal{H}_c &:= \frac{i}{2}\left( z_2 \bar{z_1} - \bar{z_2}z_1\right) + \sqrt{z_3 \bar{z_3}} \\
\mathcal{A}_c &:= {z_2} \bar{z_2} - \left(z_1 \bar{z_3} + \bar{z_1}z_3\right)
\end{align*}
In this section, we study this system. 
\begin{theorem}[Fuller Hamiltonian System]
The system \eqref{eq:fuller-system} is Hamiltonian with respect to a non-standard Poisson bracket. The angular momentum is in involution with the Hamiltonian with respect to this bracket.
\end{theorem}
\begin{proof}
For smooth functions $A,B : \mathbb{C}^3 \to \R$, we define their non-standard Poisson bracket as:
\begin{equation}\label{eq:fuller-poisson-brack}
\{A,B\}_c := \frac{2}{i} \left(\frac{\partial A}{\partial z_1} \frac{\partial B}{\partial \bar{z_3}} - \frac{\partial A}{\partial \bar{z_3}} \frac{\partial B}{\partial z_1} \right) + \frac{2}{i}\left( - \frac{\partial A}{\partial z_2} \frac{\partial B}{\partial \bar{z_2}} +  \frac{\partial A}{\partial \bar{z_2}} \frac{\partial B}{\partial z_2} \right) + \frac{2}{i}\left(\frac{\partial A}{\partial z_3} \frac{\partial B}{\partial \bar{z_1}} -  \frac{\partial A}{\partial \bar{z_1}} \frac{\partial B}{\partial z_3} \right)
\end{equation}
we can now verify directly that equations \eqref{eq:fuller-system} become:
\begin{align*}
    z_1' &= \{z_1,\mathcal{H}_c\}_c \\
    z_2' &= \{z_2,\mathcal{H}_c\}_c \\
    z_3' &= \{z_3,\mathcal{H}_c\}_c
\end{align*}
which are Hamilton's equations for this Poisson bracket. We can also verify that $\{\H_c,\mathcal{A}_c\}_c = 0$. 
\end{proof}

The system \eqref{eq:fuller-system} is acted upon by the two-dimensional scaling group $\SO \times \mathcal{G}$ where $\mathcal{G} = \R_{>0}$. The action is given by:
\begin{equation*}
    (e^{i\theta},r)\curvearrowright(z_1(t),z_2(t),z_3(t)) := (e^{i\theta}r z_1(t/r),e^{i\theta}r^2 z_2(t/r),e^{i\theta}r^3 z_3(t/r))
\end{equation*}
We also have an involution given by time reversal:
\begin{equation*}
    \tau \curvearrowright (z_1(t),z_2(t),z_3(t)) := (\overline{z}_1(-t),\overline{z}_2(-t),\overline{z}_3(-t))
\end{equation*}
\begin{remark}\leavevmode\normalfont
\begin{itemize}
    \item The group $\mathcal{G}$ acts also on the Kepler dynamical system. Cushman \& Bates~\cite{cushman1997global} call this action the \emph{virial action}.
    \item It is noteworthy that the rotation group $\SO$ is a symmetry of the system and as a result of the classical Noether theorem, we recover $\mathcal{A}_c$ as a conserved quantity. 
\end{itemize}
\end{remark}

Thus, the truncated angular momentum and the Hamiltonian are \emph{exactly conserved} for the truncated system. As an aside, we note that the Poisson bracket in \eqref{eq:fuller-poisson-brack} arises via a symplectic structure on $\mathbb{C}^3$:

\begin{proposition}
Let $A,B$ be smooth, real-valued functions on $\mathbb{C}^3$. Consider the following symplectic form on $\mathbb{C}^3$:
\begin{equation}\label{eq:fuller-symplectic-form}
    \omega_c := \frac{(-1)^{j}}{2i}  dz_j \wedge d\bar{z}_{4-j}.
\end{equation}
where we employ the Einstein summation convention to sum over $j=1,2,3$. 
Let $\vec{\mathcal{A}}$ and $\vec{\mathcal{B}}$ denote the Hamiltonian vector fields of smooth functions $\mathcal{A},\mathcal{B}$ with respect to this symplectic form. Then we have:
\[
\{\mathcal{A},\mathcal{B}\}_c = \omega_c(\vec{\mathcal{A}},\vec{\mathcal{B}})
\]
\end{proposition}
\begin{proof}
We derive an expression for the Hamiltonian vector field of a smooth function $\mathcal{G}$ associated to the non-standard symplectic form $\omega_c$:
\begin{align}\label{eq:ham-vec-field}
    d\mathcal{G}(Y) = \omega_c(\overrightarrow{\mathcal{G}},Y)
\end{align}
for an arbitrary vector field $Y$ on $\mathbb{C}^3$.

To see this, set $Y = b_j \frac{\partial}{\partial z_j} +  \bar{b}_j \frac{\partial}{\partial \bar{z}j}$ and $\overrightarrow{\mathcal{G}} = a_j \frac{\partial}{\partial z_j} +  \bar{a}_j \frac{\partial}{\partial \bar{z}j}$ for arbitrary complex numbers $b_j, a_j \in \mathbb{C}$, where we employ the Einstein summation convention.

Then \eqref{eq:ham-vec-field} becomes:
\begin{align*}
     b_j \frac{\partial \mathcal{G}}{\partial z_j} +  \bar{b}_j \frac{\partial  \mathcal{G}}{\partial \bar{z}_j} &= \frac{1}{2i}\left[(-1)^{j} dz_j \wedge d\bar{z}_{4-j}\right]( a_j \frac{\partial}{\partial z_j} +  \bar{a}_j \frac{\partial}{\partial \bar{z}j}, b_j \frac{\partial}{\partial z_j} +  \bar{b}_j \frac{\partial}{\partial \bar{z}j}) \\
     &= \frac{1}{2i}\left[ (-1)^{j} (a_j \bar{b}_{4-j} - \bar{a}_{4-j} b_{j})\right]
\end{align*}
The above equation holds for all $b_j$. First choose $b_j = i$ and $b_k = 0$ for all $k \ne j$. Next, choose $b_j = 1$ and $b_k = 0$ for all $k \ne j$. Then we get the following two equations:
\begin{align*}
i\left[ \frac{\partial \mathcal{G}}{\partial z_j} -  \frac{\partial  \mathcal{G}}{\partial \bar{z}_j} \right]    = \frac{i}{2}\left[ (-1)^{j+1} (-i a_{4-j} -i\bar{a}_{4-j})\right] = (-1)^{j+1}\mathrm{Re}[a_{4-j}] \\
\frac{\partial \mathcal{G}}{\partial z_j} + \frac{\partial  \mathcal{G}}{\partial \bar{z}_j}  = \frac{i}{2}\left[ (-1)^{j+1} (a_{4-j} -\bar{a}_{4-j})\right] = (-1)^{j}\mathrm{Im}[a_{4-j}]
\end{align*}
from which, we get, for $j=1,2,3$:
\begin{align*}
a_{4-j} &= \mathrm{Re}[a_{4-j}] + i \mathrm{Im}[a_{4-j}] \\
&= (-1)^{j+1}i\left[ \frac{\partial \mathcal{G}}{\partial z_j} -  \frac{\partial  \mathcal{G}}{\partial \bar{z}_j} \right] + (-1)^j i\left[\frac{\partial \mathcal{G}}{\partial z_j} +  \frac{\partial  \mathcal{G}}{\partial \bar{z}_j} \right] \\ 
&= (-1)^j 2i \frac{\partial \mathcal{G}}{\partial \bar{z}_j}
\end{align*}
and similarly,
\begin{align*}
\overline{a}_{4-j} &= \mathrm{Re}[a_{4-j}] - i \mathrm{Im}[a_{4-j}] \\
&= (-1)^{j+1}i\left[ \frac{\partial \mathcal{G}}{\partial z_j} -  \frac{\partial  \mathcal{G}}{\partial \bar{z}_j} \right] + (-1)^{j+1} i\left[\frac{\partial \mathcal{G}}{\partial z_j} +  \frac{\partial  \mathcal{G}}{\partial \bar{z}_j} \right] \\ 
&= (-1)^{j+1} 2i \frac{\partial \mathcal{G}}{\partial z_j}
\end{align*}

Thus, we have:
\begin{align*}
    \overrightarrow{\mathcal{G}} &= (-1)^j2i\left(\frac{\partial \mathcal{G}}{\partial \overline{z}_{4-j}}\frac{\partial}{\partial z_j} - \frac{\partial \mathcal{G}}{\partial z_{4-j}}\frac{\partial}{\partial \overline{z}_j} \right) \\ 
    &= 
    (-1)^j\frac{2}{i}\left(\frac{\partial \mathcal{G}}{\partial z_{4-j}}\frac{\partial}{\partial \overline{z}_j} - \frac{\partial \mathcal{G}}{\partial \overline{z}_{4-j}}\frac{\partial}{\partial z_j} \right)
\end{align*}

From this, we get that the Poisson bracket of two functions $\H$ and $\mathcal{G}$ is:

\[
\{\H,\mathcal{G}\}_c = \overrightarrow{\mathcal{G}}\mathcal{H} = \omega_c\left(\overrightarrow{\H},\overrightarrow{\mathcal{G}}\right)
\]

and thus, we get the required. 

\end{proof}

The above results also hold true for the Fuller system of chain length $n$: 

\begin{definition}[Length-$n$ Fuller system]
The system of ODEs given by
\begin{equation}\label{eq:fuller-system-len-n} 
    z_n' = z_{n-1}, \quad 
    z_{n-1}' = z_{n-2} \quad
    \dots \quad
    z_2' = z_1, \quad
    z_1' = \gamma \frac{z_n}{|z_n|} 
\end{equation}
where $z_i \in \mathbb{C}$ for $i=1,\dots n$ and $\gamma \in i^n\R^{\times}$, where $\R^{\times} = \R \setminus \{0\}$. This system shall be called the \emph{Fuller system of length $n$}. 
\end{definition} 

\begin{remark}\normalfont
By scaling each $z_j \mapsto z_j/\gamma$, the constant $\gamma$ in equation \eqref{eq:fuller-system-len-n} scales by $\gamma \mapsto \gamma/|\gamma|$. Thus, there is no loss of generality in assuming that $|\gamma|=1$.
\end{remark}

We can also generalize the symplectic form to:

\begin{definition}[Fuller symplectic form]\label{def:fuller-symplectic-form}
On the manifold $\mathbb{C}^n$, we have the following symplectic form:
\begin{equation*}
    \omega_c := \sum_{j=1}^n \frac{(-1)^{j}}{2i}  dz_j \wedge d\bar{z}_{n-j+1} 
\end{equation*}
\end{definition}

\begin{theorem}
The length $n$ Fuller system in equation \eqref{eq:fuller-system-len-n} is the Hamiltonian vector field (with respect to the Fuller symplectic form) of the Hamiltonian:
\begin{equation*}
\mathcal{H}_n := \sum_{j=1}^n (-1)^j \RR(z_j,i^n z_{n-j})
\end{equation*}
The angular momentum:
\begin{equation*}
    \mathcal{A}_n := \sum_{j=1}^n (-1)^j\RR(i^{n+1}z_j,z_{n-j+1})
\end{equation*}
is conserved along this system. 
\end{theorem}
\begin{proof}
If $\mathcal{G}$ is a smooth function, the above proof generalizes in a straightforward way to give the following expression for the Hamiltonian vector field of $\mathcal{G}$:
\begin{equation*}
     \overrightarrow{\mathcal{G}} = \frac{2}{i}\sum_{j=1}^n (-1)^j\left(\frac{\partial \mathcal{G}}{\partial z_{n-j+1}}\frac{\partial}{\partial \overline{z}_j} - \frac{\partial \mathcal{G}}{\partial \overline{z}_{n-j+1}}\frac{\partial}{\partial z_j} \right)
\end{equation*}
Using this expression, we can compute the Hamiltonian vector field of $\H_n$ and we recover exactly the system \eqref{eq:fuller-system-len-n}. Differentiating $\H_n$ and $\mathcal{A}_n$ along the length-$n$ Fuller system shows that they are conserved. 
\end{proof}

\subsection{Log-spiral solutions.}\label{sec:log-spiral-solutions}

The system \eqref{eq:fuller-system} admits the following outward-moving logarithmic spiral solution:
\begin{align*}
    z_3^*(t) &= A_3 t^{3-i}= \frac{1}{10} t^{3-i}\\
    z_2^*(t) &= A_2 t^{2-i} = \frac{(3-i)}{10}t^{2-i} \\
    z_1^*(t) &= A_1 t^{1-i} = \frac{(2-i)(3-i)}{10}t^{1-i} \\
    u^*(t) &= \beta_1\frac{z_3^*(t)}{|z_3^*(t)|} 
\end{align*}

The solution $z^*(t) := (\bar{z_3}^*(T-t),\bar{z_2}^*(T-t),\bar{z_1}^*(T-t))$ gives an inward-moving logarithmic spiral, where $T$ is the time taken for the trajectory to fall into the singular locus. This can be verified by differentiating. We can also verify that $\H_c(z_3^*,z_2^*,z_1^*)=\mathcal{A}_c(z_3^*,z_2^*,z_1^*)=0$. 

To approach the singular locus, the optimal control $u^*(t)$ (for the inward log-spiral) performs an infinite number of rotations along the circle $\partial U_r$ in finite time $T$. 
\subsection{Other Fuller-type systems}
Fuller systems (over $\R$) were first described in Fuller~\cite{fuller1963study} and arises as the Pontryagin system of what is now called the Fuller optimal control problem. This problem can be described as 
\[
x'=y \quad y' = u \quad \int_0^\infty x^2 dt \to \min
\]
with initial conditions $x(0) = x_0$, $y(0) = y_0$ and $u \in [-1,1]$ is a control variable taking values in an interval. The optimal trajectory for this problem consists of an arc whose control switches infinitely many times at the extremes of the control set in a finite amount of time. 

Generalizations of the Fuller phenomenon are studied in the book of Zelikin \& Borisov~\cite{zelikin2012theory}. Problem 5.1 is Chapter 5 of this book is exactly the length $n$ Fuller system in equation \eqref{eq:fuller-system-len-n} specialized to $\R$.
This system is called the \textit{multi-dimensional Fuller problem} with 1-dimensional control. 

Our system in equation \eqref{eq:fuller-system-len-n}, taking inspiration for this terminology, will be called the \emph{multi-dimensional Fuller system with 2-dimensional control}. Problem 7.2 of Zelikin \& Borisov studies a Fuller problem with multidimensional control. In particular, equation (7.11) on page 230 of this book is exactly our system \eqref{eq:fuller-system-len-n} for $n=4$. Just as we did, Zelikin \& Borisov construct log-spiral solutions to the Fuller problem for 2-dimensional control, but leaves the exploration of other solutions as a research problem in Chapter 7. 

An important point is that the Fuller systems considered in the literature have even length. Our system, because of left-invariance, has odd chain length. The same remark also applies in our derivation of the extended state space Poisson bracket (see Section \ref{sec:poisson-bracket}). Also, in our case, the extra dimension for the control and the circular symmetry of the control set gives us an additional symmetry and thus another conservation law. 



\chapter{Formal verification of Optimal Control}
As we explained in the introduction, \textit{formal verification} refers to the practice of checking the validity of a mathematical proof all the way down to the fundamental axioms of mathematics. While it is possible to do this by hand (an example being Russell and Whitehead's \textit{Principia Mathematica}), it is very tedious and so the contemporary usage of the term \textit{formal verification} is usually taken to mean \textit{computer-assisted} formal verification. 

To check the validity of a proof of a proposition down to the axioms, we first need a computer implementation of a foundation of mathematics which is expressive enough for this task. While pen-and-paper mathematics assumes Zermelo-Fraenkel set theory with Choice (ZFC) as its foundations, intense research in the latter part of the 20th century led to the development of alternative foundations of mathematics which are arguably more suited for formalization. These alternative foundations are expressive enough to capture large swathes of mathematics and (some of them) have the added advantage of having \textit{computational semantics}, which means that their implementations double up as programming languages (we shall briefly describe such a foundation in \cref{sec:coq-intro}). Because of this, it is possible to program in these languages, write specifications about program behaviour and prove that they conform to the specification all within the same environment. Formally verified code is guaranteed to behave no more and no less than what the specification allows, assuming that the implementation and the foundations are consistent. 

We shall be specifically interested in \textit{interactive theorem proving} (ITP) which is distinct from \textit{automated theorem proving} (ATP). ITP is so named because proofs are constructed \textit{interactively}, with feedback provided to the user by the system after each step in a proof (it is for this reason that these systems are also called \textit{proof assistants}). ATP, on the other hand, is concerned with the automatic discovery of proofs of propositions via push-button techniques, typically of statements in propositional logic (SAT solvers) or quantifier-free fragments of first order logic (SMT solvers). The table below shows a few proof assistants along with the foundations they implement.

\begin{center}
\begin{tabular}{ |c|c| } 
 \hline
 Proof Assistant & Foundations \\
 \hline 
 Mizar & Tarski-Grothendieck Set Theory \\ 
 HOL Light & Simple Type Theory \\
 Coq & Calculus of Inductive Constructions \\
 Lean & Calculus of Inductive Constructions\\
 Agda & Unified Theory of Dependent Types \\
 PVS & Church's Theory of Types with Dependent Types \\
 Arend & Homotopy Type Theory with Interval Types \\
 \hline
\end{tabular}
\end{center}

\section{The Coq Proof Assistant.} \label{sec:coq-intro}
Our work in formal verification uses the Coq proof assistant which was has been in development since the 1980s. Coq has been used in a wide variety of formalization projects in both Computer Science and Mathematics: it was used to produce a formally verified C compiler in the CompCert project~\cite{leroy2016compcert} and was also used to formalize the proof of the Feit-Thompson odd order theorem~\cite{gonthier2013machine}.

\subsection{Type Theory of Coq} \label{sec:coq-type-theory}
Coq implements a type theory called the Calculus of Inductive Constructions (CIC) which is a purely syntactic formal system consisting of \textit{terms} and \textit{types}. Every term in Coq has a type. Types (which can be thought of as generalizations of data types in programming languages such as C) come with type formation rules (also called \textit{rules of inference}) which specify valid \textit{judgements} in the theory. These rules of inference state how to create, use and compute with types and terms. One of the core features of a type theory such as Coq is a reserved base sort $\mathsf{Prop}$ which is the type of all propositions (statements which can be subject to proof).

These rules (\textit{syntax}) are then \textit{interpreted} in different ways (\textit{semantics}). For example, the syntactical expression $a : A$ (read as ``the term $a$ is of type $A$'') can variously be interpreted as ``$a$ belongs to the set $A$'' or (if $A$ is a proposition) ``$a$ is a proof of the proposition $A$''. For another example, the expression $f : A \to B$, can mean ``$f$ is a function from $A$ to $B$'' or (if $A$ and $B$ are both propositions) ``$f$ is a proof that proposition $A$ implies proposition $B$''. 

The type theory of Coq supports \textit{dependent types}, which are types that depend on a particular term. They are written as $\prod_{s:S} A(s)$ to stand for the \textit{dependent function type} and $\sum_{s:S}A(s)$ to stand for the \textit{dependent pair type}. The syntactical expression $p : \prod_{s:S} A(s)$ can stand to mean ``$p$ is a proof of the proposition $\forall s, A(s)$'' or ``$p$ is a function which takes a term $s : S$ and returns a term which has type $p(s) : A(s)$'' according as what $A$ is. Similarly, the syntactical expression $p : \sum_{s:S} A(s)$ can stand to mean ``$p$ is a proof of the proposition $\exists s, A(s)$'' or ``$p$ is a function which takes a term $s : S$ and returns a pair $(s,p(s))$ where $p(s) : A(s)$'' according as what $A$ is. 

This interpretation we are describing is called \textit{propositions-as-types} and is part of a larger identification called the \textit{Curry-Howard correspondence} between proof calculi and models of computation. 

Syntactic theories such as CIC were carefully designed to be expressible enough to interpret logic and logical constructs and also to preserve computational meaning. In fact, it was shown that the calculus of inductive constructions (CIC) can be interpreted in Zermelo-Fraenkel Set Theory (ZFC) and vice versa, thus proving their equiconsistency \cite{werner1997sets}. Recent advances have shown that some extensions of CIC can have more sophisticated (homotopy-theoretic) semantics \cite{kapulkin2021simplicial}. These results showed that CIC (and its sister theories and extensions) could serve perfectly well as a foundation of mathematics. Their computational character means that so long as we do not extend CIC with non-computable axioms (such as the axiom of choice) these syntactic theories can also serve as models of computation and their implementations double up as programming languages. 

In the sections that follow, we use this type theoretic syntax in our informal presentation.

\subsection{Proving Theorems Using Coq}
\begin{figure}
    \centering
    \includegraphics[scale=0.5]{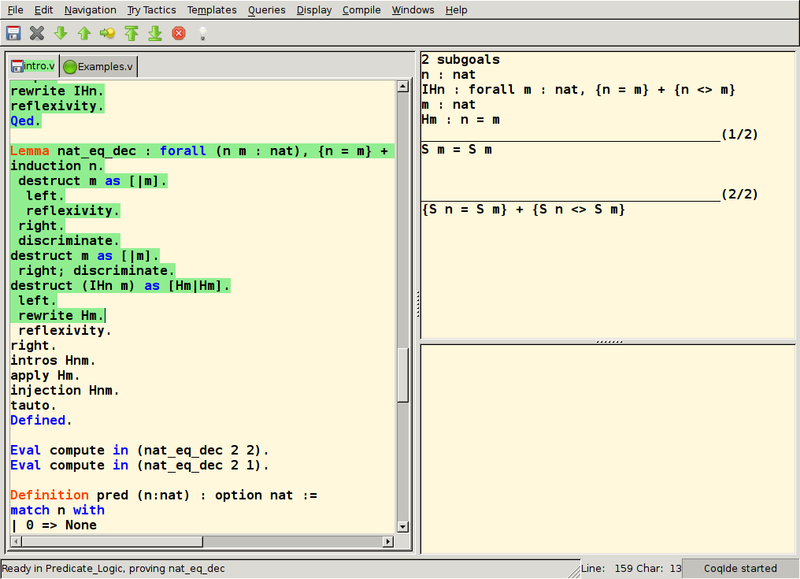}
    \caption{A Coq session.}
    \label{fig:coq-proof}
\end{figure}

As explained above, Coq is an interactive theorem prover: Propositions are stated in Coq syntax, type-checked and then loaded into the Coq session as \textit{goals} or \textit{proof obligations}, which are then proven using hypotheses which the system stores in the \textit{local context}. The proof may also invoke previously proven theorems which can be loaded into the environment (via an importing mechanism akin to Python or Haskell). Within a proof block, the user writes programs called \textit{tactics} which either make progress by transforming the goal into a new proof obligation or fail to make progress, upon which the system generates an error message informing the user of failure. This process continues until the goal is reduced to a triviality and a \textit{proof term} is constructed. Upon completion, the trusted kernel of the Coq system checks that this interactively constructed proof term does indeed prove the stated proposition. A screenshot of a Coq session is shown in Figure \ref{fig:coq-proof}. The highlighted text (in green) has been checked by the Coq system. 

This process is labour intensive and complicated proofs may well take months to formalize; but the payoff is that formal proofs come with the greatest possible guarantee of correctness and trust. Let us now turn to describing the part of mathematics our formalization is based on \textit{viz.,} Reinforcement Learning theory. 

\section{Optimal Control and Reinforcement learning} 

Reinforcement learning (RL) algorithms solve sequential decision making problems in which the goal is to choose actions that maximize a quantitative utility function \cite{bellman1954,howard1960dynamic,Puterman1994,sutton.barto:reinforcement}. Reinforcement learning algorithms produce, at a minimum, a \emph{policy} that specifies which action(s) should be taken in a given state. The primary correctness property for reinforcement learning algorithms is \emph{convergence}: in the limit, a reinforcement learning algorithm should converge to a policy that optimizes for the expected future-discounted value of the reward signal. Henceforth, whenever we say we ``proved an algorithm correct'', it shall mean that we formally proved in Coq that the algorithm in question has this property of converging to the optimal policy.

\begin{figure}
    \centering
    \includegraphics{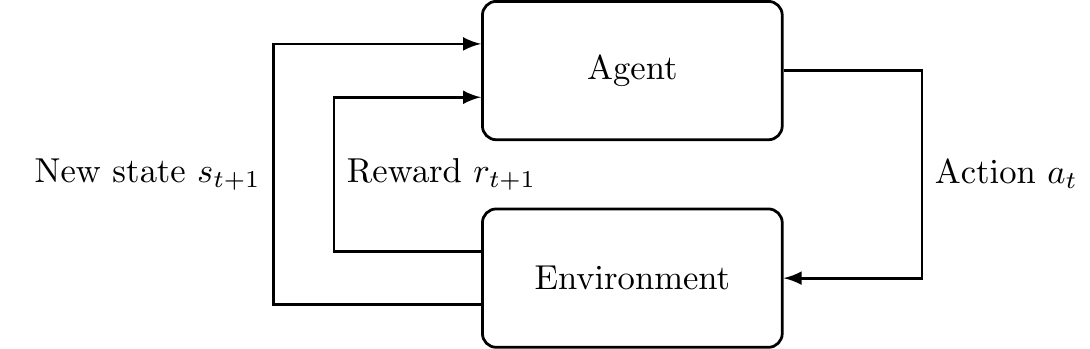}
    \caption{A typical RL loop.}
    \label{fig:rl-loop}
\end{figure}

Figure \ref{fig:rl-loop} shows a typical RL loop. At each time step $t$, the \textit{agent} (which is in a particular state $s_t$) interacts with the environment by performing an action $a_t$ available at that state. Time then progresses and the environment responds by rewarding the agent with a reward $r_{t+1}$ and the agent transitions to the new state $s_{t+1}$, upon which the loop continues.  These transitions of the agent are in general stochastic --- the next state of the agent may be predicted only up to a probability.

The class of RL algorithms may be divided into cases where the transition probability distribution of the agent is known or unknown. The former case is called \textit{Model-based RL}, while the latter is called \textit{Model-free RL}. We shall now present formalization projects which relate to both these cases. 

\subsection{Model-based RL --- CertRL}

\newcommand{\HTMLBase}{https://FormalML.github.io/CPP21/documentation/html}
\newcommand{\coqHTMLBase}{\HTMLBase}
\newcommand{\coqBaseModule}{FormalML.}
\newcommand{\fleur}{\ding{95}}
\newcommand{\coqtop}{\text{\href{https://github.com/IBM/FormalML}{\fleur}}}
\newcommand{\coqdef}[2]{\text{\href{\coqHTMLBase/\coqBaseModule#1.html\##2}{\fleur}}}

Our work in formalization of Model-based RL includes the CertRL library, which is a Coq library containing:
\begin{itemize}
  \item a formalization of Markov decision processes, long-term values,  optimal long-term value functions and the Bellman operator.
  \item a formal proof of convergence for value iteration and a formalization of the policy improvement theorem in the case of stationary policies, and
  \item a formal proof that the optimal value function for finitary sequences satisfies the  finite time analogue of the Bellman equation. 
  \item Novel use of the finitary Giry monad and a proof rule called contraction coinduction to streamline the formalization. 
\end{itemize}

All results announced is this section are joint work with Avi Shinnar, Barry Trager, Vasily Pestun and Nathan Fulton and are from the paper \cite{vajjha2021certrl}. Our presentation closely follows the work of Frank Feys, Helle Hvid Hansen and Larry Moss \cite{feys2018long}. Throughout the text which follows, hyperlinks to theorems, definitions and lemmas which have formal equivalents in the Coq development are indicated by a $\coqtop$.

The following sections elaborate on the above bulleted results more carefully. 

\subsubsection{Giry Monad}\label{sec:giry}
The Giry monad is a monad structure on the category of all measurable spaces, first described by Lawvere in \cite{lawvere1962category} and later explicitly defined by Giry in \cite{Giry1982ACA}. While the construction is very general (applying to arbitrary measures on a space), for our purposes it suffices to consider finitely supported probability measures. In this context, the Giry monad is called the \emph{finitary Giry monad}.

The Giry monad has been extensively studied and used by various authors because it has several attractive qualities that simplify (especially formal) proofs.
First, the Giry monad naturally admits a denotational monadic semantics for certain probabilistic programs \cite{DBLP:conf/popl/RamseyP02,DBLP:conf/lics/JonesP89,DBLP:conf/haskell/ScibiorGG15,audebaud2009proofs}. 
Second, it is useful for rigorously formalizing certain informal arguments in probability theory by providing a means to perform \textit{ad hoc} notation overloading \cite{tassarotti2021formal}. 
Third, it can simplify certain constructions such as that of the product measure \cite{eberl2015pdf}. 

\libname uses the Giry monad as a substitute for the stochastic matrix associated to a Markov decision process. 
This is possible because the \emph{Kleisli composition} of the Giry monad recovers the Chapman-Kolmogorov formula \cite{perrone2018categorical}. This formula is a classical result in the theory of Markovian processes that states the probability of transition from one state to another through two steps can be obtained by summing up the probability of visiting each intermediate state. Reasoning about probabilistic processes requires us to compose probabilities and so the Chapman-Kolmogorov formula plays a fundamental role, but requires formalizing and reasoning about matrix operations. Kleisli composition provides an alternative and more elegant mechanism for reasoning about compositions of probabilistic choices.

On a set $A$, let $P(A)$ denote the set of all finitely-supported probability measures on $A$  \coqdef{converge.pmf_monad}{Pmf}.
An element of $P(A)$ is a list of elements of $A$ together with probabilities. The probability assigned to an element $a : A$ is denoted by $p(a)$. 

In our development we state this as the record 

\begin{coq}
  Record Pmf (A : Type) := mkPmf {
  outcomes :> list (nonnegreal * A);
  sum1 : list_fst_sum outcomes = R1
 }.
\end{coq}

where \texttt{outcomes} is a list which stores all the entries of the type \texttt{A} along with their atomic probabilities. The field \texttt{sum1} ensures that the probabilities sum to 1. The Giry monad is defined in terms of two basic operations associated to this space: 
\begin{align*}
\mathsf{ret}: \ &A \rightarrow P(A) \ \coqdef{converge.pmf_monad}{Pmf_pure} \\   
& a \mapsto \lambda x:A, \ \delta_a(x)
\end{align*}
where $\delta_a(x) = 1$ if $a=x$ and $0$ otherwise. The other basic operation is
\begin{align*}
    \mathsf{bind} : \ &P(A) \rightarrow (A \rightarrow P(B)) \rightarrow P(B) \ \coqdef{converge.pmf_monad}{Pmf_bind} \\
    \mathsf{bind} \ &p \ f = \lambda b:B, \ \sum_{a \in A} f(a)(b) * p(a)
\end{align*}
In both cases the resulting output is a probability measure.
The above definition is well-defined because we only consider finitely-supported probability measures.
A more general case is obtained by replacing sums with integrals.

The definitions of $\bind$ and $\ret$ satisfy the following properties: 
\begin{align}
  \bind\ (\ret\ x)\ f & = \ret(f(x)) \ \coqdef{converge.pmf_monad}{Pmf_bind_of_ret} \\
  \bind\ p\ (\lambda x,\, \delta_x) & = p \ \coqdef{converge.pmf_monad}{Pmf_ret_of_bind}\\
  \bind\ (\bind\ p\ f)\ g & = \bind\ p\ (\lambda x, \bind\ (f x)\ g) \ \coqdef{converge.pmf_monad}{Pmf_bind_of_bind}
\end{align}

These \emph{monad laws} establish that the triple $(P,\bind,\ret)$ forms a monad. Think of $P(A)$ as the random elements of $A$ (\cite[page 15]{perrone2018categorical}).
In this paradigm, the set of maps $A \rightarrow P(B)$ are simply the set of maps with a random outcome. When $P$ is a monad, such maps are called Kleisli arrows of $P$.  

In terms of reinforcement learning, a map $f : A \rightarrow P(B)$ is a rule which takes a state $a:A$ and gives the probability of transitioning to state $b:B$. 

Suppose now that we have another such rule $g : B \rightarrow P(C)$. 
Kleisli composition puts $f$ and $g$ together to give a map $(f >=> g) : A \rightarrow P(C)$.
It is defined as:
\begin{align}
    f >=> g &:= \lambda x:A, \bind \ (f \ x)\ g \label{eq:firstguy} \\
             &= \lambda x:A, (\lambda c:C, \sum_{b : B} g(b)(c) *f(x)(b)) \\
             &= \lambda (x:A) \ (c:C), \sum_{b : B} f(x)(b)*g(b)(c)\label{eq:chap-kol}
\end{align}

The motivation for (\ref{eq:firstguy})--(\ref{eq:chap-kol}) is intuitive.
In order to start at $x:A$ and end up at $c:C$ by following the rules $f$ and $g$, one must first pass through an intermediate state $b:B$ in the codomain of $f$ and the domain of $g$. 
The probability of that point being any particular $b:B$ is $f(x)(b)*g(b)(c).$ So, to obtain the total probability of transitioning from $x$ to $c$, simply sum over all intermediate states $b:B$. This is exactly (\ref{eq:chap-kol}).
We thus recover the classical Chapman-Kolmogorov formula, but as a Kleisli composition of the Giry monad. 
This obviates the need for reasoning about operators on linear vector spaces, thereby substantially simplifying the formalization effort.

Indeed, if we did not use Kleisli composition, we would have to associate a stochastic transition matrix to our Markov process and manually prove various properties about stochastic matrices, which can quickly get tedious. With Kleisli composition however, our proofs become more natural and we reason closer to the metal instead of adapting to a particular representation. This is analogous to ``basis-free'' proofs in linear algebra where one reasons directly with linear operators instead of matrices associated to them. 

\subsubsection{Markov Decision Processes}
We refer to \cite{Puterman1994} for detailed presentation of the
theory of Markov decision processes.
Our formalization considers the theory of \emph{infinite-horizon discounted Markov decision processes with deterministic stationary policies}.



\begin{definition}[Markov Decision Process \coqdef{converge.mdp}{MDP}]\label{def:MDP}
  A Markov decision process consists of the following data:
  \begin{itemize}
    \item A nonempty finite type $\states$ called \emph{the set of states}.\footnote{There are various definitions of finite.  Our mechanization uses surjective finiteness (the existence of a surjection from a bounded set of natural numbers)\coqdef{utils.Finite}{Finite}, and assumes that there is a decidable equality on $\states$.  This pair of assumptions is equivalent to bijective finiteness.}
    
    \item For each state $s : \states$, a nonempty finite type $\actions(s)$ called \emph{the type of actions} available at state $s$. This is modelled as a dependent product type as described in Section \ref{sec:coq-type-theory}.
    \item A \emph{stochastic transition structure} $\transition: \prod_{s : S} (\actions(s) \to \giry(\states))$. Here $\giry(\states)$ stands for the set of all probability measures on $\states$, as described in Section \ref{sec:giry}.
    \item A \emph{reward function} $\reward: \prod_{s : S} (\actions(s) \to \states \to \reals)$ 
      where $\reward(s, a, s')$ is the reward
    obtained on transition from state $s$ to state $s'$ under action $a$. 
  \end{itemize}
\end{definition}
From these definitions it follows that the rewards are bounded in absolute value: since the state and action spaces are finite, there exists a constant $D$ such that 
\begin{equation}\label{eq:bdd_rewards}
  \forall (s \ s' : \states), (a : A(s)), |\reward (s, a, s')| \le D \ \coqdef{converge.mdp}{bdd}
\end{equation}

\begin{definition}[Decision Rule / Policy]\label{def:dec_rule}
  Given a Markov decision process with state space $\states$ and action space $\prod_{s : \states}\actions(s)$,
  \begin{itemize}
    \item   A function $ \policy: \prod_{s:\states} \actions(s)$ is called a \emph{decision rule}  \coqdef{converge.mdp}{dec_rule}. The decision rule is \emph{deterministic} \footnote{if the decision rule takes a state and returns a probability distribution on actions instead, it is called \emph{stochastic}.}. 
    \item   A \emph{stationary policy} is an infinite \emph{sequence} of decision rules: $(\pi,\pi,\pi,...)$ \coqdef{converge.mdp}{policy}. \emph{Stationary} implies that the same decision rule is applies at each step.
  \end{itemize}
\end{definition}

This policy $\policy$ induces a stochastic dynamic process on $\states$
evolving in discrete time steps $\step \in \integers_{\geq 0}$.
In this section we consider only stationary policies, and therefore
use the terms \emph{policy} and \emph{decision rule} interchangeably.

\subsubsection{Kleisli Composites in a Markov Decision Process}

Note that for a fixed decision rule $\policy$, we get a Kleisli arrow $\transition_{\policy}: \states \to \giry(\states)$ defined as $\transition_{\policy}(s) = \transition(s)(\pi(s))$.

Conventionally, $\transition_{\policy}$ is represented as a row-stochastic matrix
$(T_{\policy})^{s}{}_{s'}$ that acts on the probability co-vectors from the right, so that the row $s$
of $T_{\policy}$ corresponding to state $s$
encodes the probability distribution of states $s'$ after a transition from the state $s$. 

Let $p_{\step} \in P(\states)$ for $\step \in \integers_{\geq 0}$ denote a probability distribution on $S$ evolving
under the policy stochastic map $\transition_{\policy}$ after $\step$ transition steps, so
that $p_0$ is the initial probability distribution on $S$ (the initial distribution is usually taken to be $\ret \ s_0$ for a state $s_0$). These are related by 
\begin{equation}\label{eq:transition}
  p_{\step} = p_0 \transition_{\policy}^{\step}  
\end{equation}
In general (if $p_0 = \ret \ s_0$) the number $p_{\step}(s)$ gives the probability that starting out at $s_0$, one ends up at $s$ after $\step$ stages. So, for example, if $\step=1$, we recover the stochastic transition structure at the end of the first step \coqdef{converge.mdp}{bind_stoch_iter_1}. 

Instead of representing $T_{\policy}^{\step}$ as an iterated product of a stochastic matrix in our formalization, we recognize that \eqref{eq:transition} states that $p_{\step}$ is the $\step$-fold iterated Kleisli composite of $T_{\pi}$ applied to the initial distribution $p_0$ \coqdef{converge.mdp}{bind_stoch_iter}.
\begin{equation}
  p_{\step} = (p_0 >=> \underbrace{T_{\policy} >=> \dots >=> T_{\policy}}_{\step \ \text{times}} )
\end{equation}
Thus, we bypass the need to define matrices and matrix multiplication entirely in the formalization. 

\subsubsection{Long-Term Value of a Markov Decision Process}
 Since the transition from one state to another by an action is governed by a probability distribution $T$, there is a notion of expected reward with respect to that distribution.

\begin{definition}[Expected immediate reward]
  For a Markov decision process, 
  \begin{itemize}
  \item An \emph{expected immediate reward} to be obtained in the transition under action $a$
    from state $s$ to state $s'$
     is a function $\bar \reward: \states \to \actions \to \reals$ computed by averaging the
    reward function over the stochastic transition map to a new state $s'$
    \begin{equation}\label{eq:expt_im_rwd}
      \bar \reward(s, a) := \sum_{s' \in \states} \reward(s, a, s') \transition(s,a) (s')
    \end{equation}
    \item An \emph{expected immediate reward under a decision rule} $\policy$, denoted $\bar \reward_\policy: S \to \reals$ is defined to be:  
    \begin{equation} \label{eq:step_expt_reward}
      \bar \reward_{\policy} (s) := \bar \reward(s,\policy(s)) \  \coqdef{converge.mdp}{step_expt_reward} 
    \end{equation}
   That is, we replace the action argument in \eqref{eq:expt_im_rwd} by the action prescribed by the decision rule $\policy$.
    \item The \emph{expected reward} at time step $\step$ of a Markov decision process starting at initial state $s$, following policy $\pi$ is defined as the expected value of the reward with respect to the $\step$-th Kleisli iterate of $T_{\policy}$ starting at state $s$.
    \[  r^{\policy}_{\step}(s) := \mathbb{E}_{T_\policy^{\step}(s)}\left[\bar \reward_{\policy} \right] = \sum_{s' \in S}\left[ \bar \reward_{\policy}(s') T_{\policy}^{\step}(s)(s')\right] \ \coqdef{converge.mdp}{expt_reward}
    \] 
  \end{itemize}
  \end{definition}

The long-term value of a Markov decision process under a policy $\policy$ is defined as follows:

\begin{definition}[Long-Term Value]\label{def:ltv}
  Let $\gamma \in \reals, 0 \le \gamma < 1$ be a \emph{discount factor}, and $\policy = (\policy,\policy,\dots)$ be a stationary policy. Then $V_{\policy}: \states \to \reals$ is given by
  \begin{equation}
    V_{\policy}(s) = \sum_{\step=0}^{\infty} \gamma^{\step} r^{\policy}_{\step}(s) \ \coqdef{converge.mdp}{ltv}
  \label{eq:policyvalue}
  \end{equation}
\end{definition}

The rewards being bounded in absolute value implies that the long-term value function $V_{\policy}$ is well-defined for every initial state \coqdef{converge.mdp}{ex_series_ltv}.

It can be shown by manipulating the series in \eqref{eq:policyvalue} that the long-term value satisfies the Bellman equation:
\begin{align}
V_{\policy}(s) &= \bar v(s,\policy(s)) + \gamma \sum_{s' \in S} V_{\policy}(s') T_{\policy}(s)(s') \ \coqdef{converge.mdp}{ltv_corec}  \label{eq:ltv_corec}\\ 
&= \bar \reward_\policy(s) + \gamma \mathbb{E}_{T_\pi(s)} \left[V_{\policy}\right]
\end{align}
Abstracting the above, we can write the following definition:
\begin{definition}
Given a Markov decision process, we define the \emph{Bellman operator} as
\begin{align*}
  \bellman_{\policy}: \ (\states \to \R) &\to (\states \to \R)   \\
   W &\mapsto \bar \reward_\policy(s) + \gamma  \mathbb{E}_{T_\pi(s)} W
\end{align*}
\end{definition}

\begin{theorem}[Properties of the Bellman Operator \coqdef{converge.mdp}{bellman_op_monotone_le}]\label{thm:bellman_op_prop}
The Bellman operator satisfies the following properties:
\begin{itemize}
  \item As is evident from \eqref{eq:ltv_corec}, the long-term value $V_{\policy}$ is the fixed point of the operator $\bellman_\policy$ \coqdef{converge.mdp}{ltv_bellman_op_fixpt}.
  \item  The operator $\bellman_\policy$ is a contraction in the max norm on the space of functions $S \to \R$ \coqdef{converge.mdp}{is_contraction_bellman_op}. 
  \item The operator $\bellman_\policy$ is a monotone operator. That is, for any $W_1,W_2 : S \to \R$:
  \[\forall s, W_1(s) \le W_2(s) \Rightarrow \forall s, \bellman_\policy(W_1)(s) \le \bellman_\policy(W_2)(s) \]
\end{itemize}
\end{theorem}

The Banach fixed point theorem now says that $V_{\policy}$ is the unique fixed point of this operator and also gives us a way to approximate this fixed point. 

Let $V_{\policy, \mdplen}: \states \to \reals$ be the $\mdplen$-th iterate of the Bellman operator $B_\policy$. It can be computed by the recursion
relation 
\begin{equation}
  \begin{aligned}
   & V_{\policy, 0}(s_0) = 0 \\
   &  V_{\policy, \mdplen + 1}(s_0) = \bar \reward_{\policy}(s_0) + \gamma \mathbb{E}_{T_\pi(s_0)}V_{\policy, \mdplen} \qquad \mdplen \in \integers_{\ge 0} 
 \end{aligned}
 \label{eq:trecursion}
\end{equation}
where $s_0$ is an arbitrary initial state. The first term in the reward function $V_{\policy, \mdplen+1}$ for the process of length $\mdplen+1$ is the
sum of the reward collected in the first step (\emph{immediate reward}), and the remaining
total reward obtained in the subsequent process of length $\mdplen$ (\emph{discounted future reward}).
The $\mdplen$-th iterate is also seen to be equal to the $\mdplen$-th partial sum of the series \eqref{eq:policyvalue} \coqdef{converge.mdp}{ltv_part}.

The sequence of iterates $\{V_{\policy, \mdplen}\}|_{n = 0,1,2, \dots}$ is pointwise convergent and equals $V_{\policy}$, by the Banach fixed point theorem. 
\begin{equation}
  V_{\policy} = \lim_{\mdplen \to \infty} V_{\policy, \mdplen} \ \coqdef{converge.mdp}{bellman_op_iterate}
\end{equation} 

\subsubsection{Contraction Coinduction}
Just as we reformulated the transition matrix of a Markov process to be more naturally expressible in Coq (via the Giry monad), we can also similarly reformulate the Banach Fixed Point theorem.

Let $(X,d)$ denote a complete metric space with metric $d$. 
Subsets of $X$ are modeled by terms of the function type $\phi : X \rightarrow \mathsf{Prop}$.
Another interpretation is that $\phi$ denotes all those terms of $X$ which satisfy a particular property: whether or not an element $x : X$ belongs to $\phi$ is dictated by whether or not the proposition $\phi(x) : \mathsf{Prop}$ is true.

A \emph{Lipschitz map} \coqdef{converge.LM.fixed_point}{is_Lipschitz} is a mapping that is Lipschitz continuous; 
i.e., a mapping $F$ from $(X,d_X)$ into $(Y, d_Y)$ such that for all $x_1, x_2 \in X$ there is some $K \ge 0$ such that
\[
  d_Y(F(x_1), F(x_2)) \le K d_X(x_1, x_2).
\]
The constant $K$ is called a Lipschitz constant.

A map $F : X \rightarrow X$ is called a \textit{contractive map} \coqdef{converge.LM.fixed_point}{is_contraction}, or simply a \textit{contraction}, if there exists a constant $0 \leq \gamma < 1$ such that 
\[ 
d(F(u), F(v)) \leq \gamma d(u,v) \quad \forall u,v : X.
\]
Contractive maps are Lipschitz maps with Lipschitz constant $\gamma < 1$. 

The Banach fixed point theorem is a standard result of classical analysis which states that contractive maps on complete metric spaces have a unique fixed point.

\begin{theorem}[Banach fixed point theorem]
\label{thm:bfpt}
    If $(X,d)$ is a nonempty complete metric space and $F : X \rightarrow X$ is a contraction, then $F$ has a unique fixed point; 
    \textit{i.e.,} there exists a point $x^* \in X$ such that $F(x^*) = x^*$. This fixed point is $x^* = \lim_{\mdplen \rightarrow \infty} F^{(\mdplen)}(x_0)$ where $F^{(\mdplen)}$ stands for the $\mdplen$-th iterate of the function $F$ and $x_0$ is an arbitrary point in $X$.  
\end{theorem}

The Banach fixed point theorem generalizes to subsets of $X$.

\begin{theorem}[Banach fixed point theorem on subsets \coqdef{converge.LM.fixed_point}{FixedPoint}] 
\label{thm:bfpsubsets}
    Let $(X,d)$ be a complete metric space and $\phi$  a closed nonempty subset of $X$. Let $F : X \rightarrow X$ be a contraction and assume that $F$ preserves $\phi$. In other words, \[ 
     \phi (u) \rightarrow \phi (F(u))    
    \]
    Then $F$ has a unique fixed point in $\phi$; \textit{i.e.,} a point $x^* \in X$ such that $\phi (x^*)$ and $F(x^*) = x^*$. 
    The fixed point of $F$ is given by $x^* = \lim_{\mdplen \rightarrow \infty} F^{(\mdplen)}(x_0)$ where $F^{(\mdplen)}$ stands for the $\mdplen$-th iterate of the function $F$. 
\end{theorem}

We can now restate \cref{thm:bfpsubsets} as a proof rule:
\begin{equation}\label{eq:metricc}
    \inferrule{\phi \ \mathrm{ closed} \quad \exists x_0, \phi(x_0) \quad  \phi (u) \rightarrow \phi (F(u))}{\phi (\mathsf{fix} \ F \ x_0)} \ \coqdef{converge.mdp}{metric_coinduction}
\end{equation}

which one should read as ``if $\phi$ is a closed, non-empty subset of a complete metric space $X$ which is preserved by a contractive map $F$, then the fixed point of $F$ also belongs to the subset $\phi$''. 

This proof rule states that in order to prove that
some closed $\phi$ is a property of the fixed point of $F$,
it suffices to establish the standard inductive assumptions:
that $\phi$ holds for some initial $x_0$, 
and that if $\phi$ holds at $u$ then it also holds after a single application of $F$ to $u$. In this form, the Banach fixed point theorem is called \emph{Metric coinduction}\footnote{The rule \eqref{eq:metricc} is \emph{coinductive} because it is equivalent to the assertion that a certain coalgebra is final in a category of coalgebras. (See Section 2.3 of Kozen and Ruozzi \cite{KozenRuozzi})}.

\begin{definition}[Ordered Metric Space]\label{def:ordered-met-space}
    A metric space $X$ is called an \emph{ordered metric space} if the underlying set $X$ is partially ordered and the sets $\{z \in X | z \le y \}$ and $\{z \in X | y \le z\}$ are closed sets in the metric topology for every $y \in X$.
\end{definition}

Note that if $A$ is a finite set, the space of functions $A \to \R$ is a complete metric space with the metric induced by the  $L^{\infty}$ norm \coqdef{converge.mdp}{Rmax_norm}:
\begin{equation}\label{eq:rfct_norm}
    \|f\|_{\infty} = \max_{a \in A}|f(a)|    
\end{equation}
It is also ordered since it inherits the order from $\R$. 

For ordered metric spaces, metric coinduction specializes to \cite[Theorem 1]{feys2018long}, which we restate as  \cref{thm:contrc} below.

\begin{theorem}[Contraction coinduction]\label{thm:contrc}
    Let $X$ be a non-empty, complete ordered metric space. If $F : X \rightarrow X$ is a contraction and is order-preserving, then:
    \begin{itemize}
        \item $\forall x, F(x) \le x \Rightarrow x^* \le x$ \coqdef{converge.mdp}{contraction_coinduction_Rfct_le} and
        \item $\forall x, x \le F(x) \Rightarrow x \le x^*$ \coqdef{converge.mdp}{contraction_coinduction_Rfct_ge}
    \end{itemize}
    where $x^*$ is the fixed point of $F$.
\end{theorem}

We will use the above result to reason about Markov decision processes. Let us now consider the Value and Policy iteration algorithms. 

\subsubsection{The Convergence of Value Iteration}

In the previous subsection we defined the long-term value function
$V_{\policy}$ and showed that it is the fixed point of the Bellman operator. It is also the pointwise limit of the iterates $V_{\policy,\mdplen}$, which is the expected value of all length $\mdplen$ realizations of the Markov decision process following a fixed stationary policy $\policy$. 

The space of all decision rules is finite because the state and action spaces are finite \coqdef{converge.mdp}{dec_rule_finite}.

The above facts imply the existence of a decision rule (stationary policy) which maximizes the long-term reward. We call this stationary policy the \emph{optimal policy} and its long-term value the \emph{optimal value function}. 
\begin{equation}\label{eq:max_ltv}
    V_*(s) = \max_{\policy} \{V_\policy(s)\} \ \coqdef{converge.mdp}{max_ltv} 
\end{equation}

\textit{The aim of reinforcement learning is to have tractable algorithms to find the optimal policy and the optimal value function corresponding to the optimal policy}. 

Bellman's \emph{value iteration} algorithm is such an algorithm, which is known to converge asymptotically to the optimal value function. Value iteration is one of the most basic RL algorithms for calculating the optimal policy associated to an MDP. Despite its elementary nature, studying the convergence proof of Value iteration is instructive since it serves as a template for convergence proofs of other RL algorithms. In this section we describe this algorithm and formally prove this convergence property.

\begin{definition}\label{def:bellmanopt}
   Given a Markov decision process we define the \emph{Bellman optimality operator} as:
   \begin{align*}
    \hat \bellman:& (\states \to \reals) \to (\states \to \reals) \\
    &W \mapsto \lambda s, \max_{a \in A(s)}
    \left(\bar \reward(s,a)  + \gamma \mathbb{E}_{\transition(s,a)} [W]\right) \ \coqdef{converge.mdp}{bellman_max_op} 
   \end{align*}
\end{definition}
where the max is taken over all actions allowable at a state $s$. 

\begin{theorem}\label{thm:bellman_max_op_prop}
    The Bellman optimality operator $\hat \bellman$ satisfies the following properties:
    \begin{itemize}
      \item  The operator $\hat \bellman$  is a contraction with respect to the $L^{\infty}$ norm (\ref{eq:rfct_norm}) \coqdef{converge.mdp}{is_contraction_bellman_max_op}.
      \item The operator $\hat \bellman$ is a monotone operator. That is, 
      \[\forall s, W_1(s) \le W_2(s) \Rightarrow \forall s, \hat \bellman(W_1)(s) \le \hat \bellman(W_2)(s) \ \coqdef{converge.mdp}{bellman_max_op_monotone_le}\]
    \end{itemize}
\end{theorem}

Now we move on to proving the most important property of $\hat \bellman$:
the optimal value function $V_*$ is a fixed point of $\hat \bellman$.

By Theorem \ref{thm:bellman_max_op_prop} and the Banach fixed point theorem, we know that the fixed point of $\hat \bellman$ exists. Let us denote it $\hat V$. Then we have: 

\begin{theorem}[Lemma 1 of \cite{feys2018long} \coqdef{converge.mdp}{ltv_Rfct_le_fixpt}]\label{thm:ltv_Rfct_le_fixpt}
    For every decision rule $\sigma$, we have $V_\sigma\le \hat V$. 
\end{theorem}
\begin{proof}
    Fix a policy $\sigma$. 
    Note that for every $f : S \rightarrow \reals$, we have $\bellman_{\sigma}(f) \le \hat \bellman (f)$\coqdef{converge.mdp}{bellman_op_bellman_max_le}. In particular, applying this to $f = V_\sigma$ and using Theorem \ref{thm:bellman_op_prop}, we get that $V_\sigma = \bellman_\sigma(V_\sigma)\le \hat \bellman(V_\sigma)$. Now by contraction coinduction (Theorem \ref{thm:contrc} with $F = \hat \bellman$ along with Theorem \ref{thm:bellman_max_op_prop}) we get that $V_\sigma \le \hat V$.
\end{proof}

Theorem \ref{thm:ltv_Rfct_le_fixpt} immediately implies that $V_* \le \hat V$. To go the other way, we introduce the following policy, called the \emph{greedy} decision rule. 

\begin{equation}\label{eq:greedy} 
    \sigma_*(s) := \text{argmax}_{a \in A(s)}\left( \bar \reward (a,s) + \gamma \mathbb{E}_{T(s,a)}[\hat V]\right) \ \coqdef{converge.mdp}{greedy}
\end{equation}

We now have the following theorem:

\begin{theorem}[Proposition 1 of \cite{feys2018long} \coqdef{converge.mdp}{exists_fixpt_policy}] \label{thm:exists_fixpt_policy}
   The greedy policy is the policy whose long-term value is the fixed point of  $\hat \bellman$: 
   \[V_{\sigma_*} = \hat V\]
\end{theorem}
\begin{proof}
    We observe that $\bellman_{\sigma_*}(\hat V) = \hat V$ \coqdef{converge.mdp}{exists_fixpt_policy_aux}.    Thus, $\hat V \le \bellman_{\sigma_*}(\hat V) $.
    Note that we have $V_{\sigma_*}$ is the fixed point of $B_{\sigma_*}$ by Theorem \ref{thm:bellman_op_prop}. Now applying contraction coinduction with $F = \bellman_{\sigma_*}$, we get $\hat V \le V_{\sigma_*}$. 
    From Theorem \ref{thm:ltv_Rfct_le_fixpt} we get that $V_{\sigma_*} \le \hat V$.
\end{proof}

Theorem \ref{thm:exists_fixpt_policy} implies that $V_* \ge \hat V$ and so we conclude that $V_* = \hat V$ \coqdef{converge.mdp}{max_ltv_eq_fixpt}. Thus, the fixed point of the optimal Bellman operator $\hat \bellman$ exists and is equal to the optimal value function. 

Now with these results in mind, we can sketch the mains steps of the \textit{value iteration algorithm}:
\begin{enumerate}
  \item Initialize a value function $V_0 : \states \to \reals$. 
  \item Define $V_{\mdplen+1} = \hat \bellman V_\mdplen$ for $\mdplen\ge 0$. At each stage, the following policy is computed
  \[
    \policy_\mdplen(s) \in \text{argmax}_{a \in A(s)}\left( \bar \reward(s,a) + \gamma \mathbb{E}_{T(s,a)}[V_\mdplen]\right) 
  \]
  \item Check if $|V_{n+1} - V_{n}| < \theta$, a cutoff value.
  \item If no, go to II. If yes, return $V_n$ as our candidate optimal value function.  
\end{enumerate}

By the Banach Fixed Point Theorem, the sequence $\{V_\mdplen\}$ converges to the optimal value function $V_*$ \coqdef{converge.mdp}{bellman_iterate}. This means that the above algorithm is guaranteed to stop. It is important to note that the Value iteration algorithm computes the optimal policy implicitly (in step II of the algorithm sketch above) and so we also get the optimal policy as an output at the end of this algorithm. 

\subsubsection{The Convergence of Policy Iteration}
We have just seen that the convergence of value iteration is asymptotic: the iteration is continued until a fixed threshold is breached. Policy iteration is a similar iterative algorithm to compute the optimal policy that benefits from a more definite stopping condition. Define the \emph{$Q$ function} to be:
\[ 
    Q_\policy(s,a) := \bar \reward(s,a) + \gamma \mathbb{E}_{T(s,a)}[V_{\policy}].
\]
The \textit{policy iteration algorithm} proceeds in the following steps:
\begin{enumerate}
    \item Initialize the policy to $\pi_0$. 
    \item Policy evaluation: For $\mdplen \ge 0$, given $\pi_\mdplen$, compute $V_{\pi_\mdplen}$. 
    \item Policy improvement: From $V_{\pi_\mdplen}$, compute the greedy policy:
    \[
      \pi_{\mdplen+1}(s) \in \text{argmax}_{a \in A(s)}\left[ Q_{\pi_\mdplen}(s,a) \right]  
    \] 
    \item Check if $V_{\pi_\mdplen} = V_{\pi_{\mdplen+1}}$. If yes, stop.
    \item If not, repeat (2) and (3).
\end{enumerate}

This algorithm depends on the following results for correctness. We follow the presentation from \cite{feys2018long}.

\begin{definition}[Improved policy \coqdef{converge.mdp}{improved_tot}] \label{def:improved_policy}    
A policy $\tau$ is called an improvement of a policy $\sigma$ if for all $s \in S$ it holds that 
    \[ 
    \tau(s) = \mathrm{argmax}_{a \in A(s)}\left[ Q_{\sigma}(s,a)] \right]
    \]
\end{definition}
So, step (2) of the policy iteration algorithm simply constructs an improved policy from the previous policy at each stage. 

\begin{theorem}[Policy Improvement Theorem] \label{thm:policy_improvement}
    Let $\sigma$ and $\tau$ be two policies. 
    \begin{itemize}
        \item If $\bellman_\tau V_{\sigma} \ge \bellman_\sigma V_{\sigma}$ then $V_\tau \ge V_\sigma$ \coqdef{converge.mdp}{policy_improvement_1}.
        \item If $\bellman_\tau V_{\sigma} \le \bellman_\sigma V_{\sigma}$ then $V_\tau \le V_\sigma$ \coqdef{converge.mdp}{policy_improvement_2}.
    \end{itemize}
\end{theorem}

Using the above theorem, we have: 
\begin{theorem}[Policy Improvement Improves Values \coqdef{converge.mdp}{improved_has_better_value}] 
    If $\sigma$ and $\tau$ are two policies and if $\tau$ is an improvement of $\sigma$, then we have $V_{\tau} \ge V_{\sigma}$. 
\end{theorem}
\begin{proof}
    From Theorem \ref{thm:policy_improvement}, it is enough to show $\bellman_\tau V_{\sigma} \ge \bellman_\sigma V_{\sigma}$.
    We have that $\tau$ is an improvement of $\sigma$.
    \begin{align}
        \tau(s) &= \text{argmax}_{a \in A(s)}\left[ Q_{\sigma}(s,a) \right] \\
        &= \text{argmax}_{a \in A(s)}\left[ \bar \reward(s,a) + \gamma \mathbb{E}_{T(s,a)}[V_{\sigma}] \right]  \label{eq:tau_prop}
    \end{align}

    Note that 
    \begin{align*}
        \bellman_\tau V_{\sigma} &= \bar \reward(s,\tau(s)) + \gamma  \mathbb{E}_{T(s,\tau(s))} [V_{\sigma}] \\
        &= \max_{a \in A(s)}\left[\bar \reward(s,a) + \gamma \mathbb{E}_{T(s,a)}[V_{\sigma}]\right] \quad \mathrm{by} \ \mathrm{\eqref{eq:tau_prop}}\\
        &\ge \bar \reward(s,\sigma(s)) + \gamma \mathbb{E}_{T(s,\sigma(s))}[V_{\sigma}] \\
        &= \bellman_\sigma V_\sigma
    \end{align*}
\end{proof}

In other words, since $\policy_{\mdplen+1}$ is an improvement of $\policy_\mdplen$ by construction, the above theorem implies that  $V_{\policy_\mdplen} \le V_{\policy_{\mdplen+1}}$. This means that $\policy_\mdplen \le \policy_{\mdplen+1}$. 

Thus, the policy constructed in each stage in the policy iteration algorithm is an improvement of the policy in the previous stage. Since the set of policies is finite \coqdef{converge.mdp}{dec_rule_finite}, this policy list must at some point stabilize. 
Thus, the algorithm is guaranteed to terminate.

\subsubsection{Comments on Formalization}

Comparing Theorem \ref{thm:ltv_Rfct_le_fixpt} and Theorem \ref{thm:exists_fixpt_policy} with the the equivalent results from Puterman \cite[Theorem 6.2.2]{Puterman1994} demonstrates that \libname avoids reasoning about low-level $\epsilon - \delta$ details through strategic use of coinduction. The usefulness of contraction coinduction is also reflected in the formalization, sometimes resulting in Coq proofs whose length is almost the same as the English text.

Another important comparison is between a Coq coinductive proof and an English non-coinductive proof.
The Coq proof of the policy improvement theorem provides one such point of comparison. Recall that theorem states that a particular closed property (the set $\{x | x \le y\}$) holds of the fixed point of a particular contractive map (the Bellman operator).  
The most common argument --- presented in the most common textbook on reinforcement learning --- proves this theorem by expanding the infinite sum in multiple steps \cite[Section 4.2]{sutton.barto:reinforcement}. We reproduce this below:
\begin{theorem}[Policy Improvement Theorem \cite{sutton.barto:reinforcement}]
  Let $\pi,\pi'$ be a pair of deterministic policies such that, for all states $s$, 
  \begin{equation}\label{eq:policyimp}
  Q_\pi(s,\pi'(s)) \ge V_{\pi}(s)
  \end{equation} 
  then $V_{\pi'}(s) \ge V_{\pi}(s)$. 
\end{theorem}
\begin{proof}
  Starting with (\ref{eq:policyimp}) we keep expanding the $Q_{\pi}$ side and reapplying (\ref{eq:policyimp}) until we get $V_{\pi'}(s)$. 
  \begin{align*}
    V_{\pi}(s) &\le Q_{\pi}(s,\pi'(s)) \\
    &= \E_{\pi'}\left\{ r_{t+1} + \gamma V_{\pi}(s_{t+1})|s_t = s\right\} \\
    &\le \E_{\pi'}\left\{ r_{t+1} + \gamma Q_{\pi}(s_{t+1},\pi'(s_{t+1}))| s_t = s\right\} \\
    &= \E_{\pi'}\left\{ r_{t+1} + \gamma \E_{\pi'}\left\{ r_{t+2} + \gamma V_{\pi}(s_{t+2}) \right\} | s_t = s\right\}\\
    &=\E_{\pi'}\left\{ r_{t+1} + \gamma r_{t+2} + \gamma^2 V_{\pi}(s_{t+2}) \ |\ s_t = s \right\}\\
    &\le \E_{\pi'}\left\{ r_{t+1} + \gamma r_{t+2} + \gamma^2 r_{t+3} + \gamma^3 V_{\pi}(s_{t+3}) \ |\ s_t = s \right\}\\
    &\vdots \\
    &\le \E_{\pi'}\left\{ r_{t+1} + \gamma r_{t+2} + \gamma^2 r_{t+3} + \gamma^3 r_{t+4} \dots \ |\ s_t = s \right\}\\
    &= V_{\pi'}(s)
  \end{align*}
\end{proof}

At a high level, this proof proceeds by showing that the closed property $\{x | x \le y\}$ holds of each partial sum of the infinite series $V_{\pi}(s)$. By completeness and using the fact that partial sums converge to the full series $V_{\pi}$, this property is also shown to hold of the fixed point $V_{\pi}$. 

However, the coinductive version of this proof (Theorem \ref{thm:policy_improvement}) is simpler because it exploits the fact that this construction has already been done once in the proof of the fixed point theorem: the iterates of the contraction operator were already proven to converge to the fixed point and so there is no reason to repeat the construction again. Thus, the proof is reduced to simply establishing the ``base case'' of the (co)induction. 

\begin{coq}
Theorem policy_improvement_1 (σ τ : dec_rule M) :
(forall s, bellman_op γ ~ τ (ltv γ σ) s >= bellman_op γ ~ σ (ltv γ ~ σ) s)
  -> forall s, ltv γ ~ τ s >= ltv γ ~ σ s.
Proof.
  intros Hexpt. unfold bellman_op in Hexpt. unfold step_expt_reward in Hexpt.
  set (Hτ := ltv_bellman_op_fixpt γ ~ hγ τ ~ (ltv γ ~ τ)).
  rewrite Hτ.
  apply contraction_coinduction_Rfct_ge.
  - apply is_contraction_bellman_op ; auto.
  - apply bellman_op_monotone_ge ; auto.
  - unfold Rfct_ge. intros s.
    replace (ltv γ ~ σ s) with (bellman_op γ ~ σ (ltv γ ~ σ) s)
      by (now rewrite <-(ltv_bellman_eq_ltv γ ~ hγ σ)).
    apply Hexpt.
Qed.
\end{coq}

This shows how in some cases the act of formalization, since one is forced to think through the argument from first principles, the reasoning performed ``closer to the logic'' and this sometimes may result in simpler proofs as compared their traditional pen-and-paper equivalents. 

\subsection{Model-free RL --- Stochastic Approximation}
\renewcommand{\HTMLBase}{https://FormalML.github.io/ITP22/documentation/html}
\renewcommand{\coqHTMLBase}{\HTMLBase}
\renewcommand{\coqBaseModule}{FormalML.}

\renewcommand{\fleur}{\ding{95}}
\renewcommand{\coqtop}{\text{\href{https://github.com/IBM/FormalML}{\fleur}}}
\renewcommand{\coqdef}[2]{\text{\href{\coqHTMLBase/\coqBaseModule#1.html\##2}{\fleur}}}

We now turn to Model-free Reinforcement Learning. Recall that this means that we do not have any information on how the agent transitions or what rewards it receives (in other words, the stochastic transition structure $T$ and the reward function $r$ of a Markov Decision Process in Definition \ref{def:MDP} are not provided to us). We can only observe the transitions of the agent and are tasked with designing a Reinforcement Learning algorithm to get the optimal policy. 

This suggests a more general problem: \textit{how do we approximate a target if the target is give by an unknown function and direct observations of it are noise-corrupted}? 

More precisely, say we have an unknown function $M(\theta)$ which can be written as $M(\theta) = \mathbb{E}_\theta[Y]$ where $Y \sim p(y|\theta)$ is a $\theta$-dependent random variable. Assume further that we have access to an oracle which provides us samples $Y_1,Y_2,\dots$ of the random variable $Y$.
Given such a situation, how do we approximate/estimate the root $\theta^*$ of the equation $M(\theta) = b$?

To consider a special case of this problem, assume that $Y$ is independent of $\theta$. Then $M(\theta) = \E[Y]$ and our problem reduces to estimating $\E[Y]$ from samples $Y_1,Y_2,\dots$ of $Y$. In this case, our approximation would simply be the unbiased estimator $\frac{1}{n}\sum_k Y_k$ of the expected value of $Y$. So, in a sense, our problem is a generalization of the parameter estimation problem. 

This example of Kushner and Yin \cite{kushner2003stochastic} gives further motivation: assume that $M(\theta)$ is a function which measures the probability of success of a drug at dosage level $\theta$. We are interested, as an experimenter, to calculate the dosage level $\theta^*$ for which the probability of success $M(\theta^*) = v$. A naive approach would be to fix an arbitrary dosage level $\theta_1$, perform the experiment a number of times and calculate the success or failure of curing the disease and estimate the probability from that and compare that with the required probability $v$. If they are not comparable, then we adjust the dosage level accordingly to $\theta_2$ and repeat this process. 

One can imagine that if the current dosage level is not near the optimum, a lot of effort shall be expended in repeated wasteful experiments. A major insight of Robbins and Monro \cite{robbins1951stochastic} was that if $\theta_n$ is the dosage level at iteration $n$ and $Y_n$ is the probability of success at that level one, then one could use the following iterative procedure instead:
\[
    \theta_{n+1} = \theta_n + a_n(v - Y_n)
\]
They then proved that if 
 \begin{itemize}
        \item $\{a_n\}$ is chosen so that \begin{equation}\label{eq:learning-rate} a_n \to 0, \   \sum_{n=1}^{\infty} a_n = \infty, \   \sum_{n=1}^{\infty} a_n^2 < \infty,\end{equation}
        \item $Y$ is bounded with probability 1
        \item $M(\theta)$ is non-decreasing and
        \item $M'(\theta^*)$ exists and is positive
\end{itemize}
then $\theta_n \to 0$ in $L^2$.
In the Robbins-Monro scheduling assumption (\ref{eq:learning-rate}), it is
clear that the step-sizes have to converge to zero (otherwise the
model would fluctuate and never converge to the exact solution). The second assumption
$\sum_{n=1}^{\infty} a_n = \infty$ that says that the rates
should not converge to zero too fast is also sensible, as otherwise it
is easy to imagine a learning schedule with $a_n$ dropping to
zero so fast that the iterative process does not reach the root $\theta_{*}$ from an initial point
$\theta_0$ (for a concrete example, see \cite[pp. 5]{chen2006stochastic}). The third assumption, $\sum_{n=1}^{\infty} a_n^2 < \infty$, is
more subtle and technical, it primarily ensures that even in
situations when the noise-terms have self-correlation they would not
move the iterative process out of its track of converging with
probability 1 to the root.

One can test this iterative procedure in the case where $M(\theta) = \frac{1}{1 + \exp(-\theta)}$ and want to find the root $\theta^*$ of the equation $M(\theta) = \frac{1}{2}$ given that we can only observe $Y(\theta) = M(\theta) + U$ where $U \sim N(0,1)$ is a standard normal variate. 

\begin{figure}
    \centering
    \includegraphics[scale=0.34]{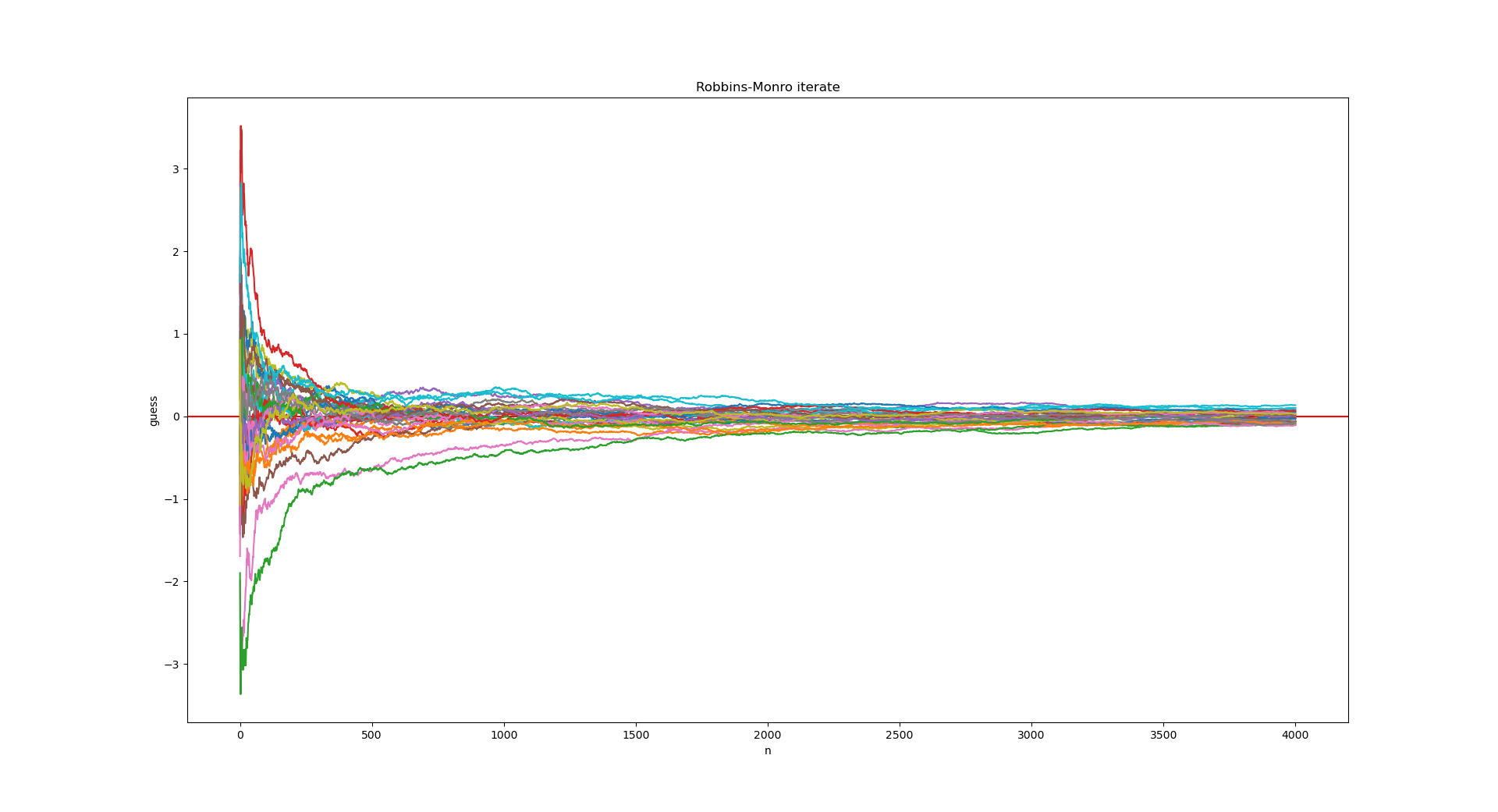}
    \caption{Testing the Robbins-Monro Procedure}
    \label{fig:robbins-monro}
\end{figure}

Figure \ref{fig:robbins-monro} shows the Robbins-Monro iterates to approximate the root of the above function for 30 random seed values in $[-2,2]$. In each case, with $a_n = n^{-0.8}$, 4000 iterates of the Robbins-Monro procedure are enough to produce a candidate root which seemingly asymptotically converges to the actual root $\theta^* = 0$. In certain situations, a slower decreasing of learning rate $\{a_n\}$ still leads to convergence, and is faster in practice.

This extensive digression on stochastic approximation serves to illustrate the fact that every model-free RL algorithm utilizes these methods in one form or another. Thus, to prove formal convergence of such algorithms, we must first formalize stochastic approximation theory. A recurring theme in the context of formal mathematics is to express things in \textit{maximal generality}. This is useful simply because it avoids repeated work since the special case can be derived from a generalization if we have the general theorem available in our library. 

Our motivation in tackling this problem was to pick a stochastic approximation result which was general enough to imply the results which are used in applications. In this process, we chose to formalize a stochastic approximation theorem of Aryeh Dvoretzky proven in 1956. 

We shall now describe this formalization effort. All results announced in this and the following sections are joint work with Avi Shinnar, Barry Trager and Vasily Pestun and appear in \cite{vajjha2022formalization}.

\subsubsection{History of Stochastic Approximation}
After Robbins and Monro's seminal 1951 paper, Kiefer and Wolfowitz \cite{10.1214/aoms/1177729392} took a similar
approach but considered the problem of estimating the parameter $x$ where the function
$M(x)$ has a maximum, and proved convergence in probability.

Wolfowitz \cite{10.2307/2236689} weakened the assumption of
Robbins-Monro about boundedness of $Y$: instead his version assumes only
that the variance of $Y$ is bounded uniformly over $\theta$, and $M(\theta)$ is
bounded, and with those assumptions Wolfowitz proves convergence in
probability.
Blum \cite{10.1214/aoms/1177728794} weakened further the assumptions
of Robbins-Monro and Wolfowitz and proved a substantially stronger
result, namely that the iterative sequence 
\begin{equation}
\label{eq:stoch-approx}
  \theta_{n+1} := \theta_{n} + a_n (b - Y_n)
\end{equation}
converges with probability
1 to the root of the equation $M(\theta) = b$. Blum requires the variance of $Y$ be uniformly bounded over $\theta$,
but he allows the expectation value $M(\theta) = \mathbb{E}_\theta[Y]$ to be bounded by a linear function of $\theta$ 
\begin{equation}
\label{eq:linear-bound}
  |M(\theta)| \leq A |\theta| +  B \quad A, B \geq 0
\end{equation}
instead of a constant.  Blum's proof is based on a version of
Kolmogorov's inequality adopted in a suitable way by Lo\`eve
\cite{zbMATH03068821} where instead of series of independent random
variables, a certain dependence was allowed but constrained by a
conditional expectation value. This extension of Kolomogorov's
inequality to the conditional situation was related to earlier works
of Borel, L\'evy and Doob about convergence with probability 1 of
certain stochastic processes.

Finally, the most general form of stochastic approximation was
formulated by Dvoretzky \cite{dvoretzky1956}. In the original Robbins-Monro
stochastic approximation (\ref{eq:stoch-approx}), the next value
$\theta_{n+1}$ is determined through the previous value $\theta_{n}$ and the sample $y_n$. Dvoretzky allowed himself more general hypotheses in
which $\theta_{n+1}$ is determined through a certain function that can take
as arguments complete history of all previous values
$\theta_{1}, \dots, \theta_{n}$ and not just the current sample $Y_{n}$. Concretely, let $T_{n}: \mathbb{R}^{n} \to \mathbb{R}$ be a measurable
real-valued function in $n$-variables. Consider the stochastic process
\begin{equation}
\label{eq:full_history}
  \theta_{n+1} := T_{n}(\theta_1, \dots, \theta_n) + W_n
\end{equation}
where $W_1, W_2, \dots$ are random variables, with $W_n$ dependent on
the previous history $\theta_1, \dots, \theta_n$ such that
\begin{equation}
\label{eq:noise}
\mathbb{E}(W_n|~\theta_1, \dots, \theta_n) = 0
\end{equation}

Another way to formulate Dvoretzky's  setup is to say that for any sequence of
random variables $X_1, X_2, \dots$ where we have conditional
probability distribution of $X_{n+1}$ dependent on the complete
history $X_1, \dots, X_n$, and then \emph{define}
\begin{equation}
  \begin{aligned}
    T(X_1, \dots, X_n) & := \mathbb{E}[X_{n+1} | X_1, \dots, X_n] \\
    W_n & := X_{n+1} - T(X_1, \dots, X_n)
  \end{aligned}
\end{equation}
in this way we automatically get the relation (\ref{eq:full_history}) (for $\theta_n = X_n$) with noise terms $W_n$ that satisfy (\ref{eq:noise}).

This setup specializes to the Robbins-Monro procedure (\ref{eq:stoch-approx}) upon setting
\begin{equation}
  \label{eq:RM-Dvoretzky}
  \begin{aligned}
& T(\theta_{1}, \dots, \theta_{n}) = \theta_{n} + a_n (b - M(\theta_n)) \\
& W_{n} = a_n(M(\theta_n) - Y_n)
  \end{aligned}
\end{equation}

For his hypotheses, Dvoretzky assumed that:
\begin{itemize}
\item there exists a point $\theta_{*}$ such that 
\begin{equation}
\label{eq:T-bound}
  |T_n(\theta_1, \dots, \theta_n) - \theta_{*}| \leq \max(\alpha_n, (1 + \beta_n) | \theta_n - \theta_{*}| - \gamma_n)
\end{equation}
where $\alpha_n, \beta_n, \gamma_n$ are sequences of non-negative real numbers with
\begin{align}
\label{eq:alpha-beta-gamma}
  \alpha_n \to 0 \\
  \sum_n \beta_n < \infty \\
  \sum_n \gamma_n = \infty 
\end{align}
\item The cumulative variance of the noise terms $W_n$ is bounded
  \begin{equation}
\label{eq:comul-noise}
    \sum_{n=1}^{\infty} \mathbb{E} [W_n^2] < \infty, \qquad \mathbb{E} [W_n | \theta_1, \dots, \theta_n] = 0
  \end{equation}
\end{itemize}
and proved that the iterative sequence (\ref{eq:full_history})
converges with probability 1 to the point $\theta_{*}$.

As the above discussion has shown, the Robbins-Monro paper spawned 
a huge literature on the analysis and applications of such stochastic algorithms. 
This is because the problem of estimating unknown parameters of a model from 
observed data is quite a fundamental one, with variants of this problem appearing in
one form or another in control theory, learning theory and other
fields of engineering.  Because of the pervasive reach of stochastic approximation methods,
any serious formalization effort of algorithms involving parameter 
estimation with an unknown underlying model will
eventually have to contend with formalizing tricky stochastic
convergence proofs.  We chose to formalize Dvoretzky's theorem as it
implies the convergence of both the Robbins-Monro and Kiefer-Wolfowitz
algorithms, various stochastic gradient descent agorithms and various
reinforcement learning algorithms such as Q-learning based on
Bellman's optimality operator.

\subsubsection{Dvoretzky's Theorem}

After Dvoretzky's original publication \cite{dvoretzky1956} of his
theorem and several very useful extensions, several shorter proofs
have been proposed.  A simplified proof was published by Wolfowitz
\cite{10.1214/aoms/1177728082} who like Blum relied on the conditional
version of Kolmogorov's law exposed by Lo\`eve
\cite{zbMATH03068821}. A third, more simplified proof was published by
Derman and Sacks \cite{derman1959dvoretzky}, who again relied on the
conditional version of Kolomogorov's law, streamlined the chain of
inequality manipulations with Dvoretzky's bounding series parameters
$(\alpha_n, \beta_n, \gamma_n)$ and used Chebyshev's inequality and the
Borel-Cantelli lemma to arrive at a very short proof. Robbins and
Siegmund generalized the theorem to the context where the variables
take value in generic Hilbert spaces using the methods of supermartingale
theory \cite{ROBBINS1971233}, as did Venter \cite{10.1214/aoms/1177699145}. For a survey see Lai \cite{10.2307/3448398}.  Dvoretzky himself published a revisited version in \cite{DVORETZKY1986220}.

We have chosen to formalise the proof following Derman and Sacks
\cite{derman1959dvoretzky} as this version appeared to us as being the
shortest and most suitable to formalize using constructions from
our library of formalized probability theory. We present complete formalization of the scalar version of Dvoretzky's theorem, with random variables taking value in $\R$.

Here is a full statement of Dvoretzky's theorem: 

\begin{theorem}[Regular Dvoretzky's Theorem~\coqdef{QLearn.Dvoretzky}{Dvoretzky_DS_simple_vec_theta}\label{th:Dvoretzky}] 
  Assuming the following:
  \begin{description}
 \item[$\mathsf{H}_1:$] Let $(\Omega, \mathcal{F}, P)$ be a probability space

 \item[$\mathsf{H}_2:$] For $n = 1,2, \dots$
  
 \item[$\mathsf{H}_3:$] Let $\mathcal{F}_{n}$ be an increasing sequence of sub $\sigma$-fields
  of $\mathcal{F}$

 \item[$\mathsf{H}_4:$] Let $X_{n}$ be $\mathcal{F}_n$-measurable random
  variables taking values in $\R$.

 \item[$\mathsf{H}_5:$] Let $T_{n}: \R^{n} \to \R$
  be a measurable function

 \item[$\mathsf{H}_6:$] Let $W_{n}$ be $\mathcal{F}_{n+1}$-measurable random variables
  taking values in $\R$ such that
      \[\quad X_{n+1} = T(x_1, \dots, x_n) + W_{n}\]
 \item[$\mathsf{H}_7:$] $ \mathbb{E}(W_n | \mathcal{F}_n) = 0$

 \item[$\mathsf{H}_8:$]  $  \sum_{n=1}^{\infty} \mathbb{E} W_n^2 < \infty$

 \item[$\mathsf{H}_9:$] Let $\alpha_n, \beta_n, \gamma_n$ be a series of real numbers such
  that

 \item[$\mathsf{H}_{10}:$]  $\alpha_n \geq 0$

 \item[$\mathsf{H}_{11}:$] $\beta_n \geq 0$

 \item[$\mathsf{H}_{12}:$] $\gamma_n \geq 0$

 \item[$\mathsf{H}_{13}:$] $\lim_{n = \infty} \alpha_n = 0$

 \item[$\mathsf{H}_{14}:$] $\lim_{n = \infty} \sum_{k=1}^{n} \beta_k < \infty$

 \item[$\mathsf{H}_{15}:$] $\lim_{n = \infty} \sum_{k=1}^{n} \gamma_k = \infty$

 \item[$\mathsf{H}_{16}:$] Let $x_{*}$ be a point in $\R$ such that for all $n = 1,2,\dots$
  and for all $x_1, \dots, x_n \in \R$,
  \[| T_n(x_1, \dots, x_n) - x_{*} | \leq \max(\alpha_n, (1 + \beta_n)|x_n - x_{*}| - \gamma_n)\]
\end{description}
Then the sequence of random variables $X_1, X_2, \dots $ converges with probability 1
  to $x_{*}$: 
 \[ \quad \quad  P\{\lim_{n=\infty} X_{n} = x_{*}\} = 1 \]
\end{theorem}

An increasing sequence $\mathcal{F}_n$ of sub-$\sigma$-fields of $\mathcal{F}$ (a filtration) formalizes a notion of
a discrete stochastic process moving forward in time steps $n$, where
$\mathcal{F}_n$ formalizes the history of the process up to the time
step $n$. Assuming an $\mathcal{F}_n$-measurable random variable $X_n$
means assuming a stochastic variable $X_n$ that is included into the history
up to the time step $n$.

  We have also formalized the extended version of Dvoretzky's theorem
  in which $\alpha_n, \beta_n, \gamma_n$ are promoted to real valued functions 
  and $T_n$ is promoted to be an $\mathcal{F}_n$-measurable random variable.
  The hypotheses that have been  modified  in the extended version are marked by the symbol $\ast$ below:

  \begin{theorem}[Extended Dvoretzky's theorem~\coqdef{QLearn.Dvoretzky}{Dvoretzky_DS_extended_simple_vec_theta}\label{thm:extended-dvor}] Assuming the following:
  \begin{description}
  \item [$\mathsf{H}_1:$] Let $(\Omega, \mathcal{F}, P)$ be a probability space

  \item[$\mathsf{H}_2:$] For $n = 1,2, \dots$
  
  \item[$\mathsf{H}_3:$] Let $\mathcal{F}_{n}$ be an increasing sequence of sub $\sigma$-fields
  of $\mathcal{F}$

  \item[$\mathsf{H}_4:$] Let $X_{n}$ be $\mathcal{F}_n$-measurable random
  variables taking values in $\R$.

  \item[$\ast \mathsf{H}_5:$] Let $T_{n}$ be $\mathcal{F}_{n}$-measurable $\R$-valued random variable

  \item[$\mathsf{H}_6:$] Let $W_{n}$ be $\mathcal{F}_{n+1}$-measurable $\R$-valued random variables
   such that:  \[ X_{n+1} = T(x_1, \dots, x_n) + W_{n}\]

  \item[$\mathsf{H}_7:$] $\mathbb{E}(W_n | \mathcal{F}_n) = 0$

  \item[$\mathsf{H}_8:$]  $ \sum_{n=1}^{\infty} \mathbb{E} W_n^2 < \infty$

  \item[$\ast \mathsf{H}_9:$] Let $\alpha_n, \beta_n, \gamma_n : \Omega \to \R$ be functions\footnote{Technically, Dvoretzky in his revisited paper \cite{DVORETZKY1986220} requires  $\alpha_n, \beta_n, \gamma_n$ to be $\mathcal{F}_{n}$-measurable, but we formalized the theorem without 
  having to assume this.} such that:

  \item[$\mathsf{H}_{10}:$] $\alpha_n \geq 0$

  \item[$\mathsf{H}_{11}:$] $\beta_n \geq 0$

  \item[$\mathsf{H}_{12}:$] $\gamma_n \geq 0$

  \item[$\ast \mathsf{H}_{13}:$] $\lim_{n\to\infty} \alpha_n  = 0$ with probability 1

  \item[$\ast \mathsf{H}_{14}:$] $\lim_{n \to \infty} \sum_{k=1}^{n} \beta_k < \infty$ with probability 1 

  \item[$\ast \mathsf{H}_{15}:$] $\lim_{n \to \infty} \sum_{k=1}^{n} \gamma_k = \infty$ with probability 1

  \item[$\mathsf{H}_{16}$:] Let $x_{*}$ be a point in $\R$ such that for all $n = 1,2,\dots$
  and for all $x_1, \dots, x_n \in \R$ we have:
      \begin{equation}\label{eq:Tn-bnd}
  | T_n(x_1, \dots, x_n) - x_{*} | \leq \max(\alpha_n, (1 + \beta_n)|x_n - x_{*}| - \gamma_n)
    \end{equation}
\end{description}
Then the sequence of random variables $X_1, X_2, \dots $ converges with probability 1
to $x_{*}$: 
  \begin{equation*}
    \quad \quad P\{\lim_{n\to\infty} X_{n} = x_{*}\} = 1
  \end{equation*}
\end{theorem}

We will now sketch the key pieces which go into the formalization of the Derman-Sacks proof of Theorem \ref{th:Dvoretzky}. 

\subsubsection{Overview of the proof}\label{sec:proof-overview}
The Derman-Sacks proof relies on a number of prerequisites 
in Probability Theory and Real Analysis. For example, the proof begins
by stating that we may replace the series $\sum_n \mathbb{E}W_n^2 < \infty$
by the series $\sum_n \frac{\mathbb{E}W_n^2}{\alpha_n^2} < \infty$ where $\alpha_n \to 0$.
This statement invokes a classical theorem of du Bois-Reymond \cite{bois1873neue}
which states:

\begin{theorem}\coqdef{utils.Sums}{no_worst_converge_iff}\label{thm:du-bois}
  Let $(a_n)$ be a sequence of nonnegative real numbers. The series 
  $\sum_n a_n$ converges if and only if there is another sequence 
  of positive real numbers $(b_n)$ such that $b_n \to \infty$ and 
  $\sum_n a_n b_n < \infty$. 
\end{theorem}
In other words, this theorem states that \textit{no worst convergent series 
exists} (see \cite{ash1997neither}). This elementary theorem did require 
some effort to formalize, in part because existing proofs such as 
the one in \cite{ash1997neither} require the sequence $(a_n)$ to consist 
only of positive terms, while our application (Dvoretzky's theorem) needed 
them to be non-negative. Additionally, we had to prove convergence 
of the product series without using the integral test (as used in \cite{ash1997neither}),
because it was unavailable in our library. Our final proof of \cref{thm:du-bois} involved a case analysis in 
which we case on whether the sequence $(a_n)$ was eventually positive or not \coqdef{utils.Sums}{eventually_pos_dec},
and we bypassed the need to use the integral test by using an exercise 
from Rudin's \textit{Principles of Mathematical Analysis} \cite{rudin1976principles}.

The main workhorse of the Derman-Sacks proof is the sequence 
$Z_n := W_n \ \mathrm{sgn} \ T_n$. First, they apply the following 
theorem\footnote{the proof of this theorem is a modification of Theorem 
6.2.1 in Ash's \textit{Probability and Measure Theory} \cite{ash2000probability}.}
to the sequence of random variables $(Z_n)$: 

\begin{theorem}[Lo\`eve \cite{zbMATH03068821}~\coqdef{QLearn.slln}{Ash_6_2_1_filter}]\label{thm:loeve}
 Let $X_1,X_2,\dots$ be a sequence of random variables adapted to a filtration 
 $(\mathcal{F}_n)_{n\in \mathbb{N}}$. 
 Assume that $\mathbb{E}[X_{n+1} \ | \ \mathcal{F}_n] = 0$ almost surely for all $n$ and also 
 that $\sum_{n=1}^\infty \mathbb{E}X_n^2$ converges. Then we have that $\sum_{n=1}^\infty X_n$ converges 
 almost surely. 
\end{theorem}
to conclude that the series $\sum_n Z_n$ converges almost surely. 
To apply this theorem we need to prove that $(Z_n)$ is adapted to the filtration $\mathcal{F}$, 
which critically uses the fact that $T_n : \mathcal{H}^n \to \mathcal{H}$ is a 
measurable function. (Here we take $\mathcal{H} = \R$.) The proof of the theorem
uses  $\mathbb{E}[X_{n+1} \ | \ \mathcal{F}_n] = 0$ to show that since the sequence is adapted,
we have  $\mathbb{E}[X_i X_j] = 0$ for all $i \neq j$. This depends on the ``factor out''
property of conditional expectation~\coqdef{ProbTheory.ConditionalExpectation}{is_conditional_expectation_factor_out}.

Next, it is shown that $|Z_n| \le \alpha_n$ almost surely for sufficiently
large $n$. This argument uses the Borel-Cantelli lemma \coqdef{ProbTheory.RandomVariableFinite}{Borel_Cantelli} 
and the Chebyshev inequality \coqdef{ProbTheory.Expectation}{Chebyshev_ineq_div_mean0}, 
both of which needed a significant amount of probability theory 
to be set up (see \cref{sec:form-prob-add}). Using this bound for $Z_n$ and the 
bound for $|T_n|$ in the hypothesis, an elementary argument shows that
\[ 
|X_{n+1}| \le \max(2 \alpha_n, |T_n| + Z_n) \le \max(2 \alpha_n, (1 + \beta_n)|X_n| + Z_n - \gamma_n)
\]
almost surely for sufficiently large $n$. 

Now, the conclusion $X_{n+1} \to 0$ almost surely follows by applying the following lemma:
\begin{lemma}\coqdef{QLearn.Dvoretzky}{DS_lemma1}
  Let $\{a_n\},\{b_n\},\{c_n\},\{\delta_n\}$ and $\{\xi_n\}$ be sequences of real
  numbers such that 
  \begin{enumerate}
    \item $\{a_n\},\{b_n\},\{c_n\},\{\xi_n\}$ are non-negative
    \item $\lim_{n \to \infty} a_n = 0, \ \sum_n b_n < \infty, \ \sum_n c_n = \infty, \ \sum_n \delta_n$ converges.
    \item For all $n$ larger than some $N_0$, $\xi_{n+1} \le \max(a_n, (1 + b_n)\xi_n + \delta_n - c_n)$ \label{it:lemma1_prop}
  \end{enumerate}
  then, $\lim_{n \to \infty}\xi_n = 0$.
\end{lemma}

The proof of the lemma is somewhat unusual since it involves running an iteration backwards: the property (\ref{it:lemma1_prop}) 
is applied repeatedly to derive an inequality between $\xi_{n+1}$ and $\xi_{N}$ for $n > N > N_0$ \coqdef{QLearn.Dvoretzky}{DS_1_helper}.
Besides using several properties of infinite products and list maximums, the final convergence result is an application of Abel's
descending convergence criterion \coqdef{utils.Sums}{Abel_descending_convergence} which says if the series $\sum_n b_n$ converges, and $a_n$ is a bounded descending sequence, then the series $\sum_n a_n b_n$ also converges.

We note that our formalization is firmly within the classical 
territory for a number of reasons: first of all,  while constructive measure theory and constructive analysis are both
actively researched topics (see \cite{coquand2008integrals,coquand2002metric,bishop1967foundations}) 
we are unaware if our main result (Dvoretzky's theorem) is constructively valid. 
Secondly, the theory of 
real numbers within the Coq standard library (which we use) uses non-computable axioms \cite{geuvers2000constructive}. 
Thirdly, as we remarked above, our proof of \cref{thm:du-bois} requires a case 
split on whether a particular sequence of real numbers is eventually zero or not, for which 
we use the axiom of constructive indefinite description.

\subsubsection{Variants of Dvoretzky's Theorem.}\label{sec:variants}

While Dvoretzky's theorem admits generalizations in many different ways, we chose 
to focus on formalizing the ones most suited for applications. 
\begin{enumerate}
  \item As already mentioned, we prove \cref{thm:extended-dvor} which is a generalization of 
  \cref{th:Dvoretzky} in which the sequences of numbers $\alpha_n, \ \beta_n, \ \gamma_n$
  are replaced by sequences of functions on the probability space. 
  This generalization is called the \textit{extended} Dvoretzky theorem \coqdef{QLearn.Dvoretzky}{Dvoretzky_DS_extended_alt_simple_vec_theta}.
  All conditions on the sequences $\alpha_n, \ \beta_n, \ \gamma_n$ now hold pointwise,
   almost everywhere. 
  \item   To apply \cref{thm:loeve} in the proof of \cref{th:Dvoretzky} we needed to
  prove that $(Z_n)$ is adapted to the filtration $\mathcal{F}$, which needed us to make
  assumptions on the functions $T_n$. These assumptions on $T_n$ can be modified and 
  generalized as:
  \begin{enumerate}
    \item  in the regular (non-extended) case, $T_n : \R^n \to \R$ are deterministic and measurable. \coqdef{QLearn.Dvoretzky}{Dvoretzky_DS_simple_vec}
    \item in the extended case, $T_n : \R^n \times \Omega \to \R$ are stochastic and $\mathcal{F}_n$-adapted. \coqdef{QLearn.Dvoretzky}{Dvoretzky_DS_extended_simple_vec}
  \end{enumerate} 
  Since Derman-Sacks do not explicitly state either assumption, we formalized 
  Dvoretzky's theorem under both assumptions. It should be noted that Dvoretzky's
  original paper \cite{dvoretzky1956} and his revisited paper \cite{DVORETZKY1986220} 
  treat both the above cases. 
  \item  We have also formalized a corollary of the extended Dvoretzky's theorem \coqdef{QLearn.Dvoretzky}{Dvoretzky_DS_extended_alt_simple_vec_theta} which proves that the 
  theorem holds in the context where the bound on $T$ in (\ref{eq:T-bound}) is assumed as follows with all other assumptions intact:
  \[
   | T_n(x_1, \dots, x_n) - x_{*} | \leq \max(\alpha_n, (1 + \beta_n - \gamma_n)|x_n - x_{*}|)  
  \]
  While this formulation is weaker compared to the original, it is convenient to have it for several
  applications of stochastic approximation theorems.  A proof of
  this corollary used a classical analysis result of Abel \cite{abel1828note}
  on the fact that the terms in a divergent sum-series could be multiplied by
  infinitesimally small series and the sum-series would still diverge
  \coqdef{utils.RealAdd}{no_best_diverge_iff}. This was addressed in 
  Dvoretzky's paper \cite[(5.1)]{dvoretzky1956}.
\end{enumerate}

\subsubsection{Formalized Probability Theory}\label{sec:form-prob-add}

As can be seen from its statement, Dvoretzky's theorem relies on a lot of Probability Theory. The statement alone makes references to events~\coqdef{ProbTheory.Event}{event}, $\sigma$-algebras~\coqdef{ProbTheory.Event}{SigmaAlgebra}, Borel $\sigma$-algebra~\coqdef{ProbTheory.BorelSigmaAlgebra}{borel_sa}, probability spaces~\coqdef{ProbTheory.ProbSpace}{ProbSpace}, the notion of an event holding almost everywhere~\coqdef{ProbTheory.Almost}{almost}, random variables~\coqdef{ProbTheory.RandomVariable}{RandomVariable}, expectation of a random variable~\coqdef{ProbTheory.Expectation}{Expectation}, conditional expectations~\coqdef{ProbTheory.ConditionalExpectation}{ConditionalExpectation},  filtrations~\coqdef{ProbTheory.SigmaAlgebras}{IsFiltration}, stochastic processes adapted to filtrations~\coqdef{ProbTheory.RandomVariable}{IsAdapted} all of which needed to be developed from scratch. These results ended up becoming a comprehensive library of probability theory which may itself be of independent interest. There are many other results proven in the library; here we highlight
two that are used in the Derman-Sacks proof: Chebyshev's inequality
and the Borel-Cantelli lemma.

Chebyshev's inequality~\coqdef{ProbTheory.Expectation}{Chebyshev_ineq_div_mean0} which states
that given a random variable $X$ and a positive constant $a$, the
probability of $\| X \| \ge a$ is less that or equal to the
expectation of $ X^2 / a^2$.
\begin{coq}
  Lemma Chebyshev_ineq_div_mean0
  (X : Ts -> R) (rv : RandomVariable dom borel_sa X) (a : posreal) :
  Rbar_le (ps_P (event_ge dom (rvabs X) a))
          (Rbar_div_pos 
          (NonnegExpectation (rvsqr X)) (mkposreal _ (rsqr_pos a))).
\end{coq}

Another is the Borel-Cantelli lemma~\coqdef{ProbTheory.RandomVariableFinite}{Borel-Cantelli} which states
that if the sum of probabilities of a sequence of events is finite,
then the probability of all but finitely many of them occuring is 0.

\begin{coq}
  Theorem Borel_Cantelli (E : nat -> event dom) :
 (forall (n:nat), sa_sigma (E n)) ->
  ex_series (fun n => ps_P (E n)) ->
  ps_P (inter_of_collection 
          (fun k => union_of_collection 
                      (fun n => E (n + k)))) = 0.
\end{coq}

In this theorem statement, \coqe{ex_series f}, defined in the Coquelicot library~\cite{Coquelicot},
asserts that the infinite series of partial sums
$\lim_{n\to\infty} \sum_{0\le i\le n} f(i)$ converges to a finite limit.


\chapter{Conclusions}
In this concluding chapter, we summarize our understanding and make some speculative comments on possible future work. 

\section{The Reinhardt Conjecture}
As we have seen, the Reinhardt Optimal Control problem has a remarkable amount of structure and has deep connections with hyperbolic geometry, Hamiltonian mechanics and the theory of chattering control. 
It is our belief that the Reinhardt conjecture has now been transformed from an impossible problem to a difficult, but approachable one. 

\subsection*{Stability of Log-Spirals}
One of the immediate extensions of our results would be to classify all the solutions that approach zero of the Fuller system that we have described in Section \ref{sec:fuller-system}. 
\[
z_3' = z_2 \quad z_2' = z_1 \quad z_1' = -i\frac{z_3}{|z_3|}
\]

We have found the inward log-spiral solution in Section \ref{sec:log-spiral-solutions} and numerical evidence suggests that this is an \textit{unstable} solution of this system, while the outward log-spirals are \textit{stable}.  

\subsection*{Ubiquity of Fuller's Phenomenon}
Fuller's problem was exhbited by Fuller~\cite{fuller1963study} as an oddity at first, but was later shown to be \emph{ubiquitous} in a very precise sense in a paper of Kupka~\cite{kupka2017ubiquity}: so long as the extended state space of our optimal problem is of sufficiently high dimension, one can find a Fuller trajectory as an extremal.

Recently, Zelikin, Lokutsievskii \& Hildebrand~\cite{zelikin2017typicality} show that for a linear-quadratic optimal problem with control variates in a two-dimensional simplex, the extremals perform infinite switchings in finite time, and their switches are chaotic in nature. Further, they prove that this behaviour is generic for piecewise smooth Hamiltonian systems near the junction where three hyper-surfaces meet in a codimension 2 manifold. The main innovation in Zelikin, Lokutsievskii \& Hildebrand~\cite{zelikin2017typicality} is the so-called \emph{descending system of Poisson brackets}, which is a clever change of coordinates of the generic system near the singularity made so that the results of the model problem are applicable. This method is illustrated in the very recent paper of Manita, Ronzhina \& Lokutsievskii~\cite{ronzhina2021neighborhood}. 

It remains to be seen whether these results can be derived for the Reinhardt problem, in particular because it is not quite clear what the \emph{order} of our singular locus is (the above results hold for second-order singular extremals\footnote{The definition of order of a singular extremal was itself was a contentious issue for a time, see Lewis~\cite{lewis1980definitions}.}) since the Hamiltonian maximization in the Reinhardt problem is a problem of fractional-linear optimization. See equation \eqref{eq:max-ham}.

However, it seems reasonable to expect that there is a ubiquity result which gives rise to the Fuller system which we derived by hand in Section \ref{sec:fuller-system}. We further hope for a Grobman-Hartman type result which would act as a bridge between the two types of dynamics: this mythical result would guarantee that the dynamics of the Reinhardt system is \emph{topologically conjugate} to the dynamics of the Fuller system in a neighbourhood of the singular locus. This result does not seem too far-fetched, since there are similar results in classical control theory, for example, the Brunovsky canonical form which proves local equivalence of affine control systems to linear control systems (see Theorem 3.9 of Elkin~\cite{elkin1999reduction}).

\subsection*{Four-cycle and Six-cycle}
In the 2017 article, Hales~\cite{hales2017reinhardt} describes the trajectories in the upper half-plane giving rise to the smoothed $(6k\pm2)$-gons. We have briefly described this in Section \ref{sec:6k+2-gons}. These trajectories (Pontryagin extremals) are given by bang-bang controls (where the control set is the 2-simplex $U_T$). The switching of the associated control vertices is in the order given by $e_3 \to e_2 \to e_1 \to e_3$ and $e_1 \to e_2 \to e_3 \to e_1$. 

The paper of Zelikin, Melnikov \& Hildebrand~\cite{ZelMelHil01} studies a Fuller-type problem in which the control set is the 2-simplex. In this paper, periodic solutions other than the ones listed above are also exhibited as being extremals. The trajectories exhibited have the control switchings as three 4-cycles: $e_i \to e_j \to e_i \to e_k  \to e_i$ with $j,k \ne i$ and two 6-cycles: $e_1 \to e_2 \to e_3 \to e_1 \to  e_2 \to e_3 \to e_1$ and in the other in the reverse order. 

While these trajectories are only extremals for the \emph{model problem} which this paper studies, we speculate that in light of the ubiquity results of Kupka~\cite{kupka2017ubiquity} and Zelikin, Lokutsievskii \& Hildebrand~\cite{zelikin2017typicality}, we should find similar periodic trajectories for the Reinhardt control problem near the singular locus. 

\subsection{Hypotrochoidal Control Problem}
The hypotrochoid result we have described in Appendix \ref{sec:hypotrochoids} might allow us a further speculation. We can compute the curvatures $\kappa_j(t)$ of the curves $\sigma_{2j}(t)$ defined in that section and compute their normalization and label them as being the \emph{controls}. This then shows that a hypotrochoid determines a control function in the control vector space $\{(u_0,u_1,u_2)~|~u_0+u_1+u_2=1\}$. 

We might then ask for an optimal control problem which has a particular hypotrochoid as a global optimizer and investigate how it might relate to the smoothed octagon. 



\section{Formal Verification of Optimal Control}

With the CertRL library and the formal proof of Dvoretzky's theorem, our work lays the foundations of formal Reinforcement Learning theory in the model-based and the model-free cases. 
While our results in the latter case are general, our intended application 
was formalizing machine learning theory, on which there is a  growing body of work \cite{tassarotti2021formal,DBLP:journals/corr/abs-2007-06776,markovInHOL,DBLP:conf/icml/SelsamLD17,DBLP:conf/aaai/Bagnall019,DBLP:journals/jar/BentkampBK19}. Our work is a step in this direction, providing future developers
of secure machine learning systems a library of basic reinforcement learning and also formalized stochastic approximation results.

In the former case we saw how the technique of metric coinduction and the Giry monad helped simplify the formalization of the value and policy iteration algorithms. It turns out that the power of the metric coinduction method goes beyond just simplifying proofs for Markov decision processes. See Kozen \& Ruozzi~\cite{KozenRuozzi} for other examples. 
Other applications of metric coinduction may be considered. Kozen's paper lists iterated function systems as a promising candidate.

What is missing in the latter case is an actual formal proof of convergence of a concrete model-free RL algorithm. The most promising candidate for this task is Q-learning, which is the prototypical model-free RL algorithm. The key ingredient in Q-learning convergence proofs is, as we remarked before, Stochastic Approximation theory, which Dvoretzky's theorem provides.  Indeed, convergence proofs of Q-Learning appeal to standard results of
stochastic approximation (see Watkins \& Dayan \cite{watkins1992q},
Jaakkola et al. \cite{jaakkola1994convergence}, Tsitsiklis
\cite{tsitsiklis1994asynchronous}). We plan to use our formalization of 
Dvoretzky's theorem to complete a convergence proof of the Q-learning algorithm. This would be a major milestone in the field of Safe RL.


\appendix     
\chapter{}

\section{The exceptional isomorphisms of \texorpdfstring{$\SL(\R)$}{SL2R} and  \texorpdfstring{$\sl(\R)$}{sl2R}}\label{sec:sl2-exceptional}
\begin{lemma}
We have the following 
\begin{itemize}
    \item isomorphisms of Lie algebras:
\begin{equation}\label{eq:sl2-isoms}
\sl(\R) \cong \sotwo \cong \su
\end{equation}
    \item isomorphism of Lie groups:
    \begin{equation}\label{eq:SL2-SU11}
        \SL(\R) \cong \SU
    \end{equation}
\end{itemize}
\end{lemma}
\begin{proof}
Consider the mapping:
\begin{align}
\phi : \sl(\R) &\to \sotwo  \nonumber \\
  \mattwo a b c {-a}  &\mapsto \left(
\begin{array}{ccc}
 0 & b-c & b+c \\
 c-b & 0 & -2 a \\
 b+c & -2 a & 0 \\
\end{array}
\right) \label{eq:sl2-so21}
\end{align}
One can check easily that $\phi$ has trivial kernel and preserves the Lie bracket. This map is a consequence of the fact that for every $g \in \SL(\R)$, $\Ad_g$ belongs to $\SOtwo$, since it acts on $\sl(\R)$ and leaves the signature $(2,1)$ trace inner product invariant. Taking differentials, the matrix of
the linear operator $\ad(g) : \sl(\R) \to \sl(\R)$ lies inside $\sotwo$. For a generic $g$, this matrix is given by \eqref{eq:sl2-so21}.

The next isomorphism is provided by the Cayley transform $C = \frac{1}{\sqrt{2}}\mattwo 1 i i 1$. We have:
\begin{align}
    C^{-1} \ \SL(\R) \ C &= \SU \\
    C^{-1} \ \sl(\R) \ C &= \su
\end{align} 

\end{proof}

We also collect the following properties of matrices in $\sl(\R)$: 
\begin{proposition}\label{prop:brack-brack}
For matrices $A,B \in \sl(\R)$ we have:
\[
[[B,A],A] = -2 \det(A) B - 2 A B A
\]
\end{proposition}

\begin{proposition}\label{prop:trace-quotient-sub}
For matrices $A,B,C,D \in \sl(\R)$, we have:
\[ 
\langle A,C\rangle\bracks{B}{D} - \langle B,C \rangle \bracks{A}{D} = -\frac{1}{2}\bracks{[A,B]}{[C,D]}
\]
\end{proposition}
\begin{proof}
Compute.
\end{proof}
\begin{proposition}\label{prop:brack-brack-inner-prod}
For matrices $A,B,C, D \in \sl$ such that $\bracks B C = 0$, we have:
\[
\bracks{[A,B]}{[C,D]} = 2 \bracks{A}{C} \bracks{B}{D}
\]
\end{proposition}
\begin{proof}
This is immediate from the previous proposition.
\end{proof}

\section{Lie-Poisson dynamics for \texorpdfstring{$X$}{X}.}\label{sec:X-lie-poisson}
Let $G$ be a Lie group with Lie algebra $\mathfrak{g}$.
Its vector space dual, $\mathfrak{g}^*$ can be equipped with a Poisson bracket called the $\pm$ \textit{Lie-Poisson bracket}: if $F,H$ are two smooth functions on $\mathfrak{g}^*$, then the bracket is given by:
\[
\{F,H\}(\nu) = \pm \left\langle \nu, \left[ \frac{\delta F}{\delta\nu}, \frac{\delta H}{\delta\nu} \right] \right\rangle_{*} \quad \nu \in \mathfrak{g}^*
\]

Here $\frac{\delta F}{\delta\nu} \in \mathfrak{g}^{**}\cong \mathfrak{g}$ is the functional derivative defined by:
\begin{equation}\label{eq:functional-derivative}
\lim_{\epsilon \to 0}\frac{1}{\epsilon}[F(\nu + \epsilon \delta\nu) - F(\nu)] = \left\langle \delta\nu, \frac{\delta F}{\delta\nu} \right\rangle_{*}
\end{equation}
for all $\delta\nu \in \mathfrak{g}^*$ and the pairing $\langle \cdot, \cdot\rangle_{*}$ is the natural pairing between a vector space and its dual. 

Hamilton's equations with respect to this bracket are called \textit{Lie-Poisson equations} and take the following form (Marsden \& Ratiu~\cite[Proposition~10.7.1]{marsden2013introduction}):
\begin{proposition}[Lie-Poisson equations]
Let $G$ be a Lie group. The equations of motion for a smooth Hamiltonian $H : \mathfrak{g}^* \to \R$ with respect to the $\pm$ Lie-Poisson brackets on $\mathfrak{g}^*$ are:
\begin{equation}\label{eq:lie-poisson-gen}
\frac{d\mu}{dt} = \mp \ad^*_{\delta H/ \delta \mu}\mu \quad \mu \in \mathfrak{g}^*
\end{equation}
\end{proposition}

Assume further that our Lie algebra $\mathfrak{g}$ is semisimple: it can be equipped with a nondegenerate inner product, which we denote by $\bracks \cdot \cdot$.
This inner product satisfies the following relation:
\begin{equation}\label{eq:nondegen-inner-prod}
\bracks {a} {[b,c]} = \bracks {[a,b]} c \quad a,b,c \in \mathfrak{g}
\end{equation}

Using this inner product, we can identify $\mathfrak{g}^*$ with $\mathfrak{g}$ as follows:
\begin{equation}\label{eq:semisimple-lg-ident}
y^*(\xi) = \bracks y \xi \qquad \xi, y\in \mathfrak{g} , \ y^* \in \mathfrak{g}^*
\end{equation}
in this way, $y^*$ is identified with $y$.

This identification maps the operator $\ad$ to $\ad^*$ and so equation \eqref{eq:lie-poisson-gen} becomes:
\[ 
\frac{d\mu}{dt} = \mp \ad_{\delta H/\delta \mu} \mu =  \mp \left[\frac{\delta H}{\delta \mu}, \mu \right]
\]

Armed with this background material, in this section we recast the dynamics for $X$ in our system, as given in Theorem $\ref{thm:X-dynamics}$, as the Lie-Poisson equation of a control-dependent Hamiltonian on the vector space $\sl(\R)^*$. 
To do this we shall need to exhibit a Hamiltonian function. Recall that we have defined $\bracks{A}{B} = \tr(AB)$ for matrices $A,B \in \sl(\R)$.

\begin{proposition}\label{lem:X-ham-frac-deriv}
If $h(X) = -\frac{\langle X,X \rangle}{2} \log \frac{\langle X,X \rangle}{\langle X, Z_0 \rangle}$ then \[X' = -\ad_{\delta h/\delta X}X = -\frac{\langle X,X \rangle}{2\langle X,Z_0 \rangle}\left[Z_0, X\right]\]
\end{proposition}
\begin{proof}
The function $h$ is well-defined since $\langle X, X \rangle = -2 $ and on the star-domain, by Corollary \ref{cor:sl2-star-condition}, we have that $\langle X,Z_0 \rangle < 0$. 

Now we have:
\begin{align*}
    \ad_{\delta h/\delta X}X &= \left[ \frac{\delta h}{\delta X}, X \right] \\
    &= \left[\frac{\delta}{\delta X}\left(-\frac{\langle X,X \rangle}{2} \log \frac{\langle X,X \rangle}{\langle X, Z_0 \rangle}\right) ,X\right]\\
    &= \left[ -X \log \frac{\langle X,X \rangle}{\langle X, Z_0 \rangle} - \frac{\bracks X X}{2} \left( \frac{\bracks X {Z_0}}{\bracks X X} \frac{2X \bracks X {Z_0} - \bracks X X Z_0}{{\bracks X {Z_0}}^2}\right), X\right] \\
    &= \frac{\bracks X X}{2 \bracks X {Z_0}}\left[ Z_0, X \right] = X'
\end{align*}
Thus, we see that the the dynamics for $X$ is Lie-Poisson with respect to the Hamiltonian $h(X) = -\frac{\langle X,X \rangle}{2} \log \frac{\langle X,X \rangle}{\langle X, Z_0 \rangle}$. 
\end{proof}
\begin{remark}\normalfont
Note that, with the parameterization of Section \ref{sec:lie-algebra-dynamics}, the Hamiltonian becomes $h(X) = \log(-\bracks{X}{Z_0})$.
\end{remark}

\section{Poisson Bracket on \texorpdfstring{$T^*(T\SL(\R))$}{T*TSL2R}}\label{sec:poisson-bracket}
The Lie-Poisson bracket helps us to write the Poisson bracket for the manifold of our control problem, viz., $T^*(T\SL(\R))$. 

Our control problem involves dynamics in the Lie group and the Lie algebra, and as we've described in Section \ref{sec:ROC}, our control problem is a left-invariant control problem at the group level. Also, we have seen in Proposition \ref{prop:left-invariant-hams} that the Hamiltonians of such problems are functions on the dual of the Lie algebra $\sl(\R)$ only.

As already mentioned, our extended state space is given by:

\[T^*(\SL(\R) \times \sl(\R)) \cong \left(\SL(\R) \times \sl(\R)\right)\times \left (\sl(\R) \times \sl(\R)\right)\]

where we have used Proposition \ref{prop:trivialization} making $T^*\SL(\R) \cong \sl(\R)^* \times \SL(\R)$ and the identification $\sl(\R)^* \cong \sl(\R)$ via the trace inner product as in Appendix \ref{sec:X-lie-poisson}. So, for $(g,X) \in \SL(\R) \times \sl(\R)$ the PMP costate variables denoted $(\Lambda_1,\Lambda_2) \in \sl(\R) \times \sl(\R)$. 

These remarks imply that the Hamiltonian is a function on $\sl(\R) \times \sl(\R) \times \sl(\R)$, since it is left-invariant. 

We now have the following theorem, which also appears in Jurdjevic \cite{jurdjevic2016optimal}:
\begin{theorem}\label{thm:extended-space-poisson-bracket}
If $F$ and $G$ are left-invariant smooth functions on $T^*(T\SL(\R))$, then their \emph{extended space Poisson bracket} is given by:
\[
\{F,G\}_{ex} := \left\langle \Lambda_1, \left[ \frac{\delta F}{\delta\Lambda_1}, \frac{\delta G}{\delta\Lambda_1} \right] \right\rangle
 + \bracks{\frac{\delta F}{\delta X}}{\frac{\delta G}{\delta \Lambda_2}} -  \bracks{\frac{\delta F}{\delta \Lambda_2}}{\frac{\delta G}{\delta X}}
\]
which is the sum of the Lie-Poisson bracket on $\sl(\R)^*$ and the canonical Poisson bracket on $T^*(\sl(\R))$. Here $\delta/\delta X$ denotes the functional derivative with respect to $X$.  
\end{theorem}
\begin{proof}
Let $\mathcal{H}$ be the Hamiltonian of the optimal control problem. We have seen (in Section \ref{sec:costate-variables}) that the dynamics for $\Lambda_1,\Lambda_2$ and $X$ are given by:
\begin{align*}
\Lambda_1' &= \ad^*_{\delta \H/\delta \Lambda_1}(\Lambda_1) \cong \left[\Lambda_1, \delta \H/\delta \Lambda_1 \right] \\
X' &= \frac{\delta \H}{\delta \Lambda_2} \\
\Lambda_2' &= -\frac{\delta \H}{\delta X}
\end{align*}

Now, if $F(\Lambda_1,\Lambda_2,X)$ is an arbitrary left-invariant smooth function, then its derivative along the lifted trajectory is computed by the chain rule as follows:
\begin{align*}
    \frac{dF}{dt} &= \mathbf{D}_1F(\Lambda_1,\Lambda_2,X)\cdot \frac{d \Lambda_1}{dt} + \mathbf{D}_2 F(\Lambda_1,\Lambda_2,X)\cdot \frac{d \Lambda_2}{dt} + \mathbf{D}_3 F(\Lambda_1,\Lambda_2,X)\cdot \frac{d X}{dt} \\ 
    &= \bracks{\frac{\delta F}{\delta \Lambda_1}}{\frac{d \Lambda_1}{dt}} + \bracks{\frac{\delta F}{\delta \Lambda_2}}{\frac{d \Lambda_2}{dt}} + 
     \bracks{\frac{\delta F}{\delta X}}{\frac{d X}{dt}} \\ 
     &= \left\langle \Lambda_1, \left[ \frac{\delta F}{\delta\Lambda_1}, \frac{\delta \H}{\delta\Lambda_1} \right] \right\rangle
 + \bracks{\frac{\delta F}{\delta X}}{\frac{\delta \H}{\delta \Lambda_2}} -  \bracks{\frac{\delta F}{\delta \Lambda_2}}{\frac{\delta \H}{\delta X}} \\ &= \{F,\H\}_{ex}
\end{align*}
which gives the required. Thus, the Poisson structure on $T^*(T\SL(\R)) \cong T^*(\SL(\R) \times \sl(\R))$ is the direct sum of the Poisson structures on $T^*\SL(\R)$ and $T^*\sl(\R)$. Since we restrict to left-invariant functions, the former is actually the Lie-Poisson structure.
\end{proof}

\section{Symplectic Structure of Coadjoint Orbits}\label{sec:kirillov}

On a Lie group $G$ with Lie algebra $\mathfrak{g}$, Kirillov~\cite{kirillov2004lectures} has defined a symplectic structure on the coadjoint orbit $\O_{\mu}:=\{\Ad^*_{g^{-1}}\mu~|~g \in G\}$ through $ \mu \in \mathfrak{g}^*$ (the linear dual of the Lie algebra $\mathfrak{g}$). This 2-form $\omega^K$ on $\O_\mu$ is given by:
\[
\omega_\mu^K(\ad_X^*\nu,\ad_Y^*\nu) := \bracks{\nu}{[X,Y]}, \quad \nu \in \O_\mu, \quad X,Y \in \mathfrak{g}
\]
where $\ad^*_X\nu,\ad^*_Y\nu \in T_{\nu}\O_{\mu}$. 
We specialize this general construction to our setting with $G = \SL(\R)$. 

Since the Lie algebra $\sl(\R)$ carries with it the nondegenerate trace inner product: $\bracks{X}{Y} = \mathrm{trace}(XY)$, this sets up a linear isomorphism $\sl(\R)^* \cong \sl(\R)$, which we use to transform the symplectic structure from coadjoint orbits to adjoint orbits.

In this section, we prove that the Kirillov symplectic structures on the adjoint orbit $\O_X \subset \sl(\R)$ and the symplectic structure on the Poincar\'{e} upper half-plane $\h$ are (anti)-equivalent. Recall that we have the following map:
\begin{align*}
     \Phi : \mathfrak{h} &\to \mathcal{O}_J \\
     z = x+iy &\mapsto \mattwo {x/y} {-(x^2 + y^2)/y} {1/y} {-x/y} =: X_z
\end{align*} 
transporting us from the upper half-plane to coadjoint orbit $\OX = \O_J$ of $\sl(\R)$.
\begin{lemma}
     The map $\Phi$ (defined in Lemma \ref{lem:def-phi}) is an anti-symplectomorphism.
\end{lemma}
\begin{proof}
Let $\omega$ be the symplectic form of the upper half-plane:
\[
\omega = \frac{dx \wedge dy}{y^2}     
\]
and let $\omega^K$ be the Kirillov 2-form on the coadjoint orbit $\OX$.  We have to show $\omega^K$ pulls back to the 2-form $\omega$ on the upper half-plane by $\Phi : \h \to \OX$. So, at a point $z = x+iy \in \h$ and tangent vectors $v,w \in T_z\h$:
\begin{align*}
     \Phi^*\omega^K_z(v,w) &= \omega^K_{\Phi(z)}((T_z\Phi(v),T_z\Phi(w)) \\
     &=\left\langle X_z,\left[ \mattwo {\frac{v_2}{2y}} {\frac{v_1y - v_2 x}{y})} {0} {\frac{v_2}{2y}},\mattwo {\frac{w_2}{2y}} {\frac{w_1y - w_2 x}{y})} {0} {\frac{w_2}{2y}} \right]\right\rangle \\
     &=-\frac{v_1 w_2 - v_2 w_1}{y^2} = -\omega_z(v,w)
\end{align*}
This proves that $\Phi^* \omega^K = -\omega$. 
\end{proof}

\section{Hypotrochoids.}\label{sec:hypotrochoids}

\begin{figure}[htbp]
    \centering
    \includegraphics[scale=0.5]{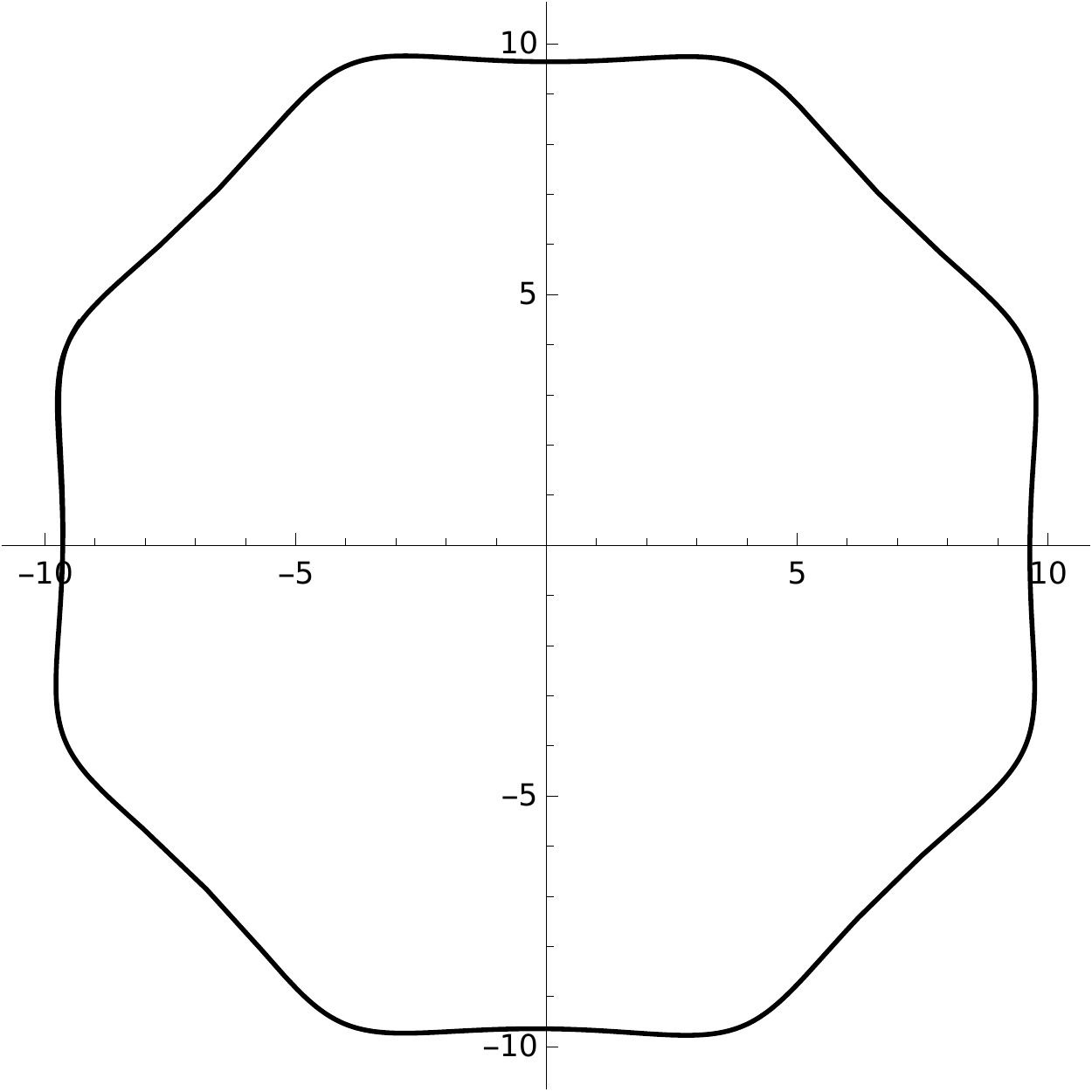}
    \caption{A hypotrochoid resembling the smoothed octagon.}
    \label{fig:smoothed-octagon-hypo}
\end{figure}

A hypotrochoid is a roulette curve which is traced by a point which is at a distance $d$ from the center of a circle of radius $r$ as it rolls without slipping on the inside of a circle of a fixed circle of radius $R$. 
The parametric equation of a hypotrochoid is given by:
\begin{align*}
    x(\theta) &= (R_1-r_1)\cos(\theta) +d_1 \cos\left(\frac{(R_1-r_1)}{r_1}\theta \right) \\
    y(\theta) &= (R_1-r_1)\sin(\theta) -d_1 \sin\left(\frac{(R_1-r_1)}{r_1}\theta \right)
\end{align*}

This section was motivated by the striking figure in Figure \ref{fig:smoothed-octagon-hypo}, which depicts a hypotrochoid with parameters $r_1 = 2.498, R_1= 2.855$ and $d_1=-10$. As we can see, this figure resembles the smoothed octagon in Figure \ref{fig:smoothed}. 

If $\zeta$ is a primitive cube root of unity, define \[\sigma_{2j}(\theta) := Re^{i\theta}\zeta^j +r e^{-i\rho \theta}\zeta^{-j}.\] which is a closed curve of period $2\pi/\rho$. We see that we recover a hypotrochoid from $\sigma_0$ by setting $\rho = (R_1-r_1)/r_1$, $R = (R_1-r_1)$ and $r = d_1$. 

In general, the smoothed octagon is given by a bang-bang control and hence cannot be a real analytic curve. But the following proposition shows that the hypotrochoid is realized by a curve in $\SL(\R)$. 

\begin{proposition}
Assuming for $|r| \ne |R|$, there exists a curve in $\SL(\R)$ which realizes the hypotrochoid.
\end{proposition}
\begin{proof}
We will prove the following of the curves $\sigma_{2j}(t)$:
\begin{align*}
\sigma_0(t) + \sigma_2(t) + \sigma_4(t) &= 0 \\
\Im(\bar{\sigma}_0(t)\sigma_2(t)) &= \mathrm{constant} \\
\sigma_{2j}(t + \frac{2\pi}{3\rho}) &= \sigma_{2j+2}(t)
\end{align*}
The first identity is due to the fact that $1 + \zeta + \zeta^2 = 0$. The second follows from:
\[
\bar{\sigma}_0(t)\sigma_2(t) = r^2 \zeta^2 + R^2 \zeta + 2 rR\Re(\zeta e^{it + i\rho t})
\]
The third follows from the definition of $\sigma_{2j}$. 
Identifying $\mathbb{C}$ with $\R^2$ we get that 
$\Im(\bar{\sigma}_0(t)\sigma_2(t) = \sigma_0(t) \wedge \sigma_2(t) = \mathrm{constant}$. 
This means that there is a constant $s > 0$ such that the curves $s\sigma_{2j}$ form a multicurve as in Definition \ref{def:multi-pt-multi-curve}. We can go through the same construction now as in Section \ref{sec:lie-group-dynamics}: define a curve $\tilde{g} \in \SL(\R)$ so that $s \sigma_{2j}(t) = \tilde{g}(t)e_{2j}^*$ to give the required. 
\end{proof}



    



\safebibliography{etdbib.bib}
\bibliographystyle{plain}

\end{document}